\newcommand{\imageformat}{bmp}
\renewcommand{\imageformat}{eps}
\renewcommand{\include}{\input}
\date{}

\title{Constructive Homological Algebra\\and Applications\\
       {\normalsize (2006 Genova Summer School)}}
\author{{\small \em Julio Rubio, Francis Sergeraert}}

\documentclass[a4paper,12pt]{article}
\usepackage{amsfonts}
\usepackage{amssymb}
\usepackage{amsmath}
\usepackage{amscd}
\usepackage{graphicx}
\usepackage[all]{xy}

\newtheorem{thr}{Theorem}
\newtheorem{mth}[thr]{Meta-Theorem}
\newtheorem{prp}[thr]{Proposition}
\newtheorem{crl}[thr]{Corollary}
\newtheorem{lmm}[thr]{Lemma}
\newtheorem{rmr}[thr]{Remark}
\newtheorem{dfn}[thr]{Definition}

\newtheorem{prb}[thr]{Problem}

\newtheorem{uos}[thr]{UOStated}
\newtheorem{clm}[thr]{Claim}
\newtheorem{ntt}[thr]{Notation}

\newcommand{\proof}{\noindent\textsc{Proof.}\ }
\newcommand{\QED}{\oneright{\rule{5pt}{5pt}}}

\newcommand{\ArHr}{\textrm{ArHr}}
\newcommand{\Cobar}{\textrm{Cobar}}
\newcommand{\Cone}{\textrm{Cone}}
\newcommand{\Ib}{\overline{I}}
\newcommand{\Ksz}{\textrm{Ksz}}
\newcommand{\ND}{^{ND}}
\newcommand{\Rsl}{\textrm{Rsl}}
\newcommand{\SSEH}{\mathcal{SS}_{EH}}
\newcommand{\Tor}{\textrm{Tor}}

\newcommand{\bC}{\mathbb{C}}
\newcommand{\bH}{\mathbb{H}}
\newcommand{\bN}{\mathbb{N}}
\newcommand{\bQ}{\mathbb{Q}}
\newcommand{\bR}{\mathbb{R}}
\newcommand{\bZ}{\mathbb{Z}}
\newcommand{\boldqm}{\mbox{\boldmath\(?\)}}
\newcommand{\boxtt}[1]{\mbox{\small\texttt{#1}}}
\newcommand{\cP}{\mathcal{P}}

\newcommand{\cU}{\mathcal{U}}
\newcommand{\coker}{\textrm{coker}\,}
\newcommand{\dinj}{\Delta^{\mathrm{inj}}}
\newcommand{\dsrj}{\Delta^{\mathrm{srj}}}
\newcommand{\ev}{\textrm{ev}}
\newcommand{\fG}{\mathfrak{G}}
\newcommand{\fR}{\mathfrak{R}}
\newcommand{\fRb}{\overline{\mathfrak{R}}}
\newcommand{\fh}{\mathfrak{h}}
\newcommand{\fk}{\mathfrak{k}}
\newcommand{\functorobj}{\hspace{2pt}<\hspace{0pt}?\hspace{0pt}>\hspace{2pt}}
\newcommand{\fm}{\mathfrak{m}}
\newcommand{\hC}{\widehat{C}}
\newcommand{\hd}{\widehat{d}}
\newcommand{\hdl}{\widehat{\delta}}
\newcommand{\id}[1]{\textrm{id}_{#1}}
\newcommand{\ideal}[1]{<\!\!#1\!\!>}
\newcommand{\im}{\textrm{im}\,}
\newcommand{\pr}{\textrm{pr}}
\newcommand{\oC}{\overline{C}}
\newcommand{\sKsz}{\textrm{\scriptsize Ksz}}
\newcommand{\senumerate}{\vspace{-10pt}\begin{enumerate}\sitemsep}
\newcommand{\sitemsep}{\setlength{\itemsep}{-3pt}}
\newcommand{\st}{{\ \ \textbf{st}\ \ }}
\newcommand{\thrtitle}[1]{\textbf{\emph{(#1)}}}
\newcommand{\und}[1]{\mbox{$\underline{\mathbf{#1}}$}}

\newcommand{\bast}{\mbox{\boldmath\(\ast\)}}
\newcommand{\bfPS}{\textbf{PS}}
\newcommand{\bfSP}{\textbf{SP}}
\newcommand{\cI}{\mathcal{I}}

\newcommand{\oneright}[1]{\hspace*{0pt}\hfill#1}

\newcommand{\lrdc}{\mbox{\,\(\Leftarrow\hspace{-9pt}\Leftarrow\hspace{-9pt}\Leftarrow\)\,}}
\newcommand{\rrdc}{\mbox{\,\(\Rightarrow\hspace{-9pt}\Rightarrow\hspace{-9pt}\Rightarrow\)\,}}
\newcommand{\eqvl}{\mbox{\(\lrdc\hspace{-10pt}\rrdc\)}}

\newcommand{\image}[3]%
  {\includegraphics[width=#1\pixel,height=#2\pixel]{#3.\imageformat}}
\newlength{\pixel}

\newcommand{\bmp}{\rule{0pt}{1pt}\\[-15pt]{\tiny.\dotfill.}\\*[-3pt]}
\newcommand{\bmpi}{\begingroup\footnotesize}
\newcommand{\empi}{\endgroup\rule{0pt}{1pt}\\*[-2pt]}
\newcommand{\empim}{\endgroup\hspace{5pt}\mbox{\small\(\maltese\)}\rule{0pt}{1pt}\\*[-2pt]}
\newcommand{\empimx}{\endgroup\hspace{5pt}\mbox{\small\(\maltese\)}\rule{0pt}{1pt}\\[-2pt]}
\newcommand{\empix}{\endgroup\rule{0pt}{1pt}\\[-2pt]}

\newcommand{\emp}{\rule{0pt}{1pt}\\*[-17pt]{\tiny.\dotfill.}\\*[-10pt]\rule{0pt}{1pt}}

\newlength{\ul}

\newcommand{\autoarrow}[1]
  {\begin{xy}<-#1\ul,0\ul>:<0\ul,-1\ul>::
   (0,-12) ; (0,12)
   **\crv{(40,-24) & (104,0) &(40,24)} ?>*\dir{>} ;
   (1,-12) ; (1,12)
   **\crv{(41,-24) & (105,0) &(41,24)} ?>*\dir{>} ;
   (0,-11) ; (0,13)
   **\crv{(40,-23) & (104,1) &(40,25)} ?>*\dir{>} ;
   (0,-24)*{} ; (70,0) *{} ; (0,24) *{}
   \end{xy}}

\begin{document}

\voffset=-2.5cm \hoffset=-0.9cm \sloppy

\maketitle

\section{Introduction.}

Standard  homological  algebra  is   not  \emph{constructive},  and  this  is
frequently the source of serious problems when \emph{algorithms} are looked
for. In particular the usual  exact  and  spectral  sequences  of homological
algebra frequently are in general not  sufficient  to obtain  some unknown
homology or homotopy group.   We will  explain it  is  not difficult to fill in
this gap, the main tools being on one hand, from a mathematical  point of view,
the so-called Homological Perturbation Lemma,  and on the  other hand,  from a
computational point  of view, Functional Programming.

We  will   illustrate  this   area  of  constructive   mathematics  by
applications in two domains:

\begin{itemize}
\item
Commutative  Algebra frequently  meets  homological objects,  in particular
when resolutions  are  involved (syzygies).   Constructive Homological Algebra
produces new  methods to process old problems such as homology of  Koszul
complexes and resolutions. The solutions so obtained are \emph{constructive}
and therefore more complete than the usual ones, an important point for their
concrete use.
\item
Algebraic  Topology  is  the historical  origin  of  Homological Algebra. The
usual exact and spectral sequences of Algebraic Topology can  be easily
transformed   into  new  effective  versions,  giving algorithms computing  for
example  unknown  homology  and  homotopy groups in wide standard contexts.  In
particular the effective version of  the Eilenberg-Moore spectral  sequence
gives  a very simple solution  for the  old Adams' problem: what algorithm
could  compute the homology groups of iterated loop spaces?
\end{itemize}

Thanks are due to Ana Romero who carefully proofread several sections of this
text.

\section{Standard Homological Algebra.}

We briefly recall in this section the minimal standard background of
homological algebra. We mainly concentrate on definitions and basic results.
Many good textbooks are available for the corresponding proofs, the main one
being maybe~\cite{MCLN2}. The only problem almost never considered in these
books is the relevant \emph{computability problem}. Besides giving the expected
background, our aim consists in making obvious why standard homological algebra
does not at all satisfy the modern \emph{constructiveness requirement}.

\subsection{Ingredients.}

Homological algebra is a general style of cooking where the main ingredients
are a \emph{ground ring} \(\fR\), \emph{chain-complexes}, \emph{chain groups},
\emph{boundary maps}, \emph{chains}, \emph{boundaries}, \emph{cycles},
\emph{homology classes}, \emph{homology groups}, \emph{exact sequences} and,
the last but not the least, \emph{spectral sequences}. In particular we do not
consider here the cohomological operations, where a good reference
is~\cite{MSTN}; this roughly defines the frontier between which is covered in
this text and which is not. Let us remark also that cohomological operations
would probably be filed by most algebraic topologists in Algebraic
\emph{Topology}, but it was explained in~\cite[Section 2]{RBSR7} why such a
discussion in fact does not make sense. In the same way, modern homological
algebra requires the notion of algebraic operad~\cite{MRSS}, a completely
different approach toward constructive algebraic topology, very interesting,
but which unfortunately did not yet produce significant concrete computer
programs. An operad is nothing but an algebra of generalized \emph{abstract}
cohomological operations.

Homological algebra was invented to systematically organize the algebraic
environment needed by the computation of the homology groups associated with
some topological objects. The first systematic presentation of Algebraic
Topology heavily based on homological algebra certainly is~\cite{ELST}, another
convenient reference for a detailed presentation and the relevant proofs of
most elementary facts. Now homological algebra is a fundamental tool in many
domains not directly connected to algebraic topology. Section~\ref{74280} here
devoted to the so-called \emph{Spencer cohomology}, where homological algebra
is applied to commutative algebra and local non-linear PDE systems, is a
typical example.

\subsection{Chain-complexes.}

\subsubsection{Definitions.}

The \emph{ground ring} \(\fR\) is an arbitrary unitary commutative ring; in the
topological case, an abelian \emph{group}, without any multiplicative structure
can also be considered, frequent when studying spectral sequences, because of
``coefficients'' that are other homology groups. In algebraic topology, \(\fR\)
is often \(\bZ\), the most general case because of the \emph{universal
coefficient theorem}~\mbox{\cite[V.11]{MCLN2}}: if you know the homology groups
with respect to the ground ring~\(\bZ\), you can easily deduce the same
homology groups with respect to any other ground ring of coefficients. But
because of the power of the \(\bZ\)-homology groups, they are of course the
most difficult to be computed. Other less ambitious possibilities are \(\fR =
\bQ\) or \(\bZ_p\) (\(p\) being a prime number); note that in algebraic
topology, \(\bZ_p\) does not denote the \(p\)-adic ring, it is simply \(\bZ_p
:= \bZ/p\bZ\); the rings \(\bQ\) and \(\bZ_p\) are in fact fields, making
easier certain calculations, and the last but not obvious step then consists in
reconstructing the \(\bZ\)-homology groups from the \(\bQ\) and \(\bZ_p\)
homology groups, the main tool being the Bockstein-Browder spectral
sequence~\cite[Chap.10]{MCCL}; a critical and interesting open problem consists
in obtaining a \emph{constructive} version of this spectral sequence.

\begin{uos}\footnote{Unless otherwise stated.}---
In these notes, unless otherwise stated, some underlying ring \(\fR\) is
assumed given. A module is therefore implicitly an \(\fR\)-module.
\end{uos}

In algebraic topology, the most useful ring is \(\bZ\) and you can assume this
convenient hypothesis. In commutative algebra, the ground ring will be most
often a field; a module is then a vector space, making some problems
significantly easier; but this apparent comfort is also misleading:
\emph{effective} homology is as useful in commutative algebra as in algebraic
topology.

\begin{dfn}---
\emph{A \emph{chain-complex} \(C_\ast\) is a pair of sequences \(C_\ast = (C_q,
d_q)_{q \in \bZ}\) where:
\begin{itemize}
\item
For every \(q \in \bZ\), the component \(C_q\) is an \(\fR\)-module, the
\emph{chain group} of degree~\(q\).
\item
For every \(q \in \bZ\), the component \(d_q\) is a module morphism \(d_q: C_q
\rightarrow C_{q-1}\), the differential map.
\item
For every \(q \in \bZ\), the composition \(d_q d_{q+1}\) is null: \(d_q d_{q+1}
= 0\).
\end{itemize}}
\end{dfn}
\[
\xymatrix@C30pt@1{\cdots & C_{q-2} \ar[l]_{d_{q-2}} & C_{q-1} \ar[l]_{d_{q-1}}
\ar@/^15pt/[ll]^0 & C_q \ar[l]_{d_{q}} \ar@/^15pt/[ll]^0 & C_{q+1}
\ar[l]_{d_{q+1}} \ar@/^15pt/[ll]^0 & C_{q+2} \ar[l]_{d_{q+2}} \ar@/^15pt/[ll]^0
& \cdots \ar[l]_{d_{q+3}} \ar@/^15pt/[ll]^0}
\]
\begin{dfn}---
\emph{If \(C_\ast = (C_q, d_q)_{q \in \bZ}\) is a chain-complex, the module
\(C_q\) is called the \emph{chain group} of degree \(q\) (in fact it is a
module, but the terminology \emph{chain group} is so traditional\ldots), or the
group of \emph{\(q\)-chains}. The image \(B_q = d_{q+1}(C_{q+1}) \subset C_q\)
is the (sub) group of \emph{\(q\)-boundaries}. The kernel \(Z_q = \ker(d_q)
\subset C_q\) is the group of \emph{\(q\)-cycles}. The relation \(d_q \circ
d_{q+1} = 0\) is equivalent to the inclusion relation \(B_q \subset Z_q\):
every boundary is a cycle but the converse in general is not true. The
``difference'' (quotient) \(H_q = Z_q/B_q\) is the \emph{homology
group}~\(H_q(C_\ast)\), again in fact a module.}
\end{dfn}

Another possible point of view consists in considering \(C_\ast = \oplus_q
C_q\) is a \emph{graded module} and the differential \(d: C_\ast \rightarrow
C_{\ast - 1}\) is a graded morphism of degree -1 satisfying \(d^2 = 0\).
According to the situation one point of view or other is more convenient, and
you must be able to immediately translate from one point of view to the other
one.

Most often, the chain groups in negative degree are null: \(q < 0 \Rightarrow
C_q = 0\), so that it becomes tempting to decide to index by \(q \in \bN\)
instead of \(q \in \bZ\), but experience shows it is not a good idea. The main
reason is that this requires specific definitions in degree 0, the cycle group
being then no longer defined, unless you decide to put an extra \(C_{-1} = 0\)
and the problem is transferred at -1\ldots\ In particular when we write down
corresponding programs, a choice \(q \in \bN\) would require specific code for
the particular case \(q = 0\), which quickly becomes painful, without any
advantage.

\begin{dfn}\label{80060}---
\emph{More generally, let \(C_\ast\) be a chain-complex and \(M\) a
\emph{coefficient group}, that is, an \(\fR\)-module. Then \(C_\ast\) and \(M\)
generate two other chain-complexes:
\begin{itemize}
\item
\(C_\ast \otimes_\fR M := (C_q \otimes_\fR M, d_q \otimes_\fR \id{M})\). The
corresponding cycles, boundaries and homology groups are then usually denoted
by \(Z_q(C_\ast ; M)\), \(B_q(C_\ast ; M)\) and \(H_q(C_\ast; M)\). We speak
then of homology groups ``with coefficients in \(M\)''.
\item
\(\textrm{Hom}(C_\ast ; M) := (\textrm{Hom}(C_q, M), d^q)\) with \(d^q\) the
morphism \(d^q: \textrm{Hom}(C_q, M) \rightarrow
\textrm{Hom}(C_{q\,\underline{+1}}, M)\) \emph{dual} to \(d_{q\,\underline{+
1}}\). The corresponding objects are then denoted with \(q\)-\emph{exponents}:
\(Z^q(C_\ast ; M)\), \(B^q(C_\ast ; M)\) and \(H^q(C_\ast ; M)\). In this case,
when the differential has degree \(+1\), it is common to call the complex a
\emph{co}chain-complex, to call the corresponding objects \emph{co}cycles,
\emph{co}boundaries, \emph{co}homology groups (not homology cogroups!).
\end{itemize}}
\end{dfn}

  Others
prefere to reverse the indices, deciding that \(C^q(C_\ast ; M) :=
\textrm{Hom}(C_{-q}, M)\); question of taste. the cohomological context will be
rarely considered in these notes.

\subsubsection{Simplicial complexes.}\label{70994}

\begin{dfn}---\label{66339}
\emph{A \emph{simplicial complex} \(K\) is a pair \(K = (V,S)\) where:
\begin{itemize}
\item
 The component \(V\) is a totally ordered\footnote{A more intrinsic definition
 does not require such an order, but the associated chain-complex is
 significantly bigger; it can always be reduced over a much smaller
 chain-complex, the definition of which requires a total order over the vertex
 set. See Sections 2 and 6 of Chapter VI of~\cite{ELST}.} set, the set of vertices of \(K\).
\item
 The component \(S\) is a set of \emph{non-empty finite} parts of \(V\), the
 \emph{simplices} of \(K\), satisfying the properties:
 \begin{itemize}
 \item
 For every \(v \in V\), the singleton \((v) \in S\).
 \item
 For every \(\sigma \in V\), then \(\emptyset \neq \sigma' \subset \sigma\)
 implies \(\sigma' \in V\).
 \end{itemize}
\end{itemize}
}
\end{dfn}

For example the small simplicial complex drawn here:
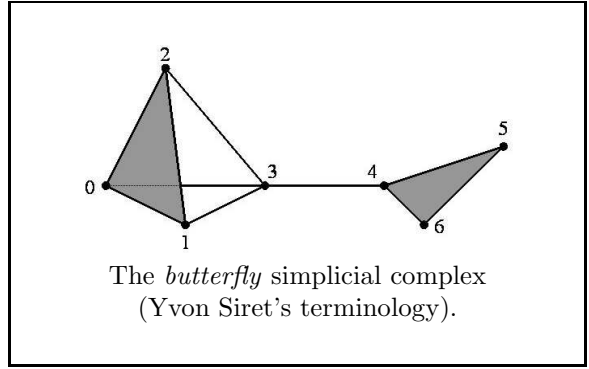
\begin{center}
\fbox{\parbox{0.5\textwidth}{
 \begin{center}
 \image{270}{130}{Butterfly}\\
 \footnotesize The \emph{butterfly} simplicial complex\\(Yvon Siret's terminology).
 \end{center}}}
\end{center}
is mathematically defined as the object \(B = (V,S)\) with:
\begin{center}
 \(V = (0,1,2,3,4,5,6)\)\\[10pt]
 \(S = \left\{
 \begin{array}{l}
 (0),(1),(2),(3),(4),(5),(6),\\
 (0,1),(0,2),(0,3),(1,2),(1,3),(2,3),(3,4),(4,5),(4,6),(5,6),\\
 (0,1,2),(4,5,6)
 \end{array}
 \right\}\)
\end{center}
In other words, the second component, the simplex list, gives the list of all
vertex combinations which are (abstractly) spanned by a simplex. The vertex set
\(V\) could be for example ordered as the integers are. Note also, because the
vertex set is ordered, the list of vertices of a simplex is also ordered, which
allows us to use a \emph{sequence} notation \((\cdots)\) and not a subset
notation \(\{\cdots\}\) for a simplex and also for the total vertex list \(V\).

A simplicial complex can be infinite. For example if \(V = \bN\) and \(S =
\{(n)\}_{n \in \bN} \cup \{(0,n)\}_{n \geq 1}\), the simplicial complex so
obtained could be understood as an infinite bunch of segments. Standard
algebraic topology proves that most ``sensible'' homotopy types can be modelled
as simplicial complexes, often infinite. We will see the notion of simplicial
set, roughly similar but more sophisticated, is also much more powerful to
reach this goal\footnote{There is here an amusing bug of terminology: the
notion of simplicial set, due to Sam Eilenberg, is more \emph{complex} than the
notion of simplicial\ldots\ complex.}.

\subsubsection{From simplicial complexes to chain-complexes.}\label{43994}

\begin{dfn}---
\emph{Let \(K = (V,S)\) be a simplicial complex. Then the set \(S_n(K)\) of
\(n\)-simplices of \(K\) is the set made of the simplices of cardinality
\(n+1\).}
\end{dfn}

For example the set of simplices \(S_0(K)\) is the set of singletons \(S_0(K) =
\{(v)\}_{v \in V}\). The set of 2-simplices of the butterfly \(B\) is
\(\{(0,1,2), (4,5,6)\}\); in the same case, the set of 1-simplices has ten
elements.

\begin{dfn}---
\emph{Let \(K=(V,S)\) be a simplicial complex. Then the \emph{chain-complex
\(C_\ast(K)\) canonically associated with \(K\)} is defined as follows. The
chain group \(C_n(K)\) is the free module generated by \(S_n(K)\). Let \((v_0,
\ldots, v_n)\) be an \(n\)-simplex, that is, a generator of \(S_n(K)\). The
boundary of this generator is then defined as:
\[
d_n((v_0, \ldots, v_n)) = (v_1, v_2, \ldots, v_n) - (v_0, v_2, v_3, \ldots,
v_n) + \cdots + (-1)^n (v_0, v_1, \ldots, v_{n-1})
\]
and this definition is linearly extended to \(C_n(K)\). }
\end{dfn}

A variant of this definition is important.

\begin{dfn}---
\emph{Let \(K = (V,S)\) be a simplicial complex. Let \(n \geq 1\) and \(0 \leq
i \leq n\) be two integers \(n\) and \(i\). Then the \emph{face operator}
\(\partial^n_i\) is the linear map \(\partial^n_i: C_n(K) \rightarrow
C_{n-1}(K)\) defined by:
\[
\partial^n_i((v_0, \ldots, v_n)) = (v_0, \ldots, v_{i-1}, v_{i+1}, \ldots,
v_n)\ :
\]
the \(i\)-th vertex of the simplex is removed, so that an \((n-1)\)-simplex is
obtained.}
\end{dfn}

\begin{rmr}---
The boundary operator \(d_n\) is the alternate sum:
\[
d_n := \sum_{i=0}^n (-1)^i \partial^n_i.
\]
\end{rmr}

This definition will be generalized in Section~\ref{18426} thanks to the notion
of simplicial \emph{set}.

Our butterfly example is then sufficient to understand the nature of the
notions of chain, cycle, boundary and homology class. An example of 1-chain is
\(c = (1,3) + (3,4) + (4,5) \in C_1(B)\); we here have chosen an example as
close as possible to the usual concrete notion of ``chain'', but \(c' = (0,2) +
(3,4) + (5,6)\) is a chain as well. The boundaries are \(d_1(c) = - (1) + (5)\)
and \(d_1(c') = + (2) - (0) + (4) - (3) + (6) - (5)\). The chain \(c_1 = (1,2)
+ (2,3) - (1,3)\) is a cycle, but is not a boundary, it is an ``interesting''
cycle, the homology class of which is non-null. On the contrary the cycle \(c_2
= (4,5) + (5,6) - (4,6)\) is trivial, it is the boundary of the 2-chain
\((4,5,6)\), and its homology class is null. If a cycle is homologous to 0, it
can be in general the boundary of several \emph{different} chains; for example,
in our butterfly, the 0-cycle \((3) - (1)\) is the boundary of the 1-chain
\((1,3)\), but also the boundary of \((1,2) + (2,3)\), a different 1-chain.

\subsubsection{Computing homology groups.}

Computing a homology group amounts to computing the relevant boundary matrices,
and to determine a kernel, an image and the quotient of the first one by the
second one. For example, if we want to compute the homology group \(H_1(B)\),
the 1-dimensional homology group of our butterfly, we have to describe the
kernel of \(d_1\):
\[
\begin{array}{rcl}
 \ker d_1 &=& \phantom{\ \oplus\ } \fR((0,1) + (1,2) - (0,2))\\
 && \ \oplus\  \fR((0,1) + (1,3) - (0,3))\\
 && \ \oplus\  \fR((0,2) + (2,3) - (0,3))\\
 && \ \oplus\  \fR((4,5) + (5,6) - (4,6))\\
\end{array}
\]
and the image of \(d_2\):
\[
\begin{array}{rcl}
 \im d_2 &=& \phantom{\ \oplus\ } \fR((0,1) + (1,2) - (0,2))\\
 && \ \oplus\  \fR((4,5) + (5,6) - (4,6)).\\
\end{array}
\]
Note in particular the tempting cycle \((1,2) + (2,3) - (1,3)\) is the
alternate sum of the first three ones in the description of \(\ker d_1\). So
that the homology group \(H_1(B)\) is \emph{isomorphic} to \(\fR^2\) with
\((0,1) + (1,3) - (0,3)\) and \((0,2) + (2,3) - (0,3)\) as possible
\emph{representants} of generators, but adding to such a representant an
arbitrary boundary gives another representant of the same homology class.

These computations quickly become complicated and it is then better -- or
necessary -- to be helped by a computer. Let us examine for example the case of
the real projective plane \(P^2\bR\). It can be proved the minimal
triangulation of \(P^2\bR\) as a simplicial complex is described by this
figure:
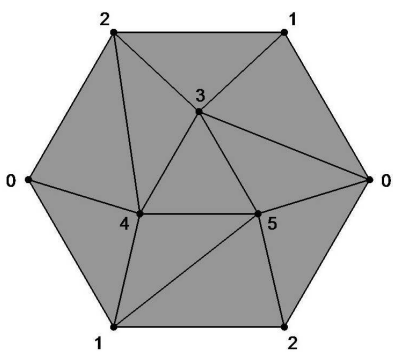
\begin{center}
\image{248}{216}{Projective-Plane}
\end{center}
The projective plane is the quotient of the 2-sphere by the antipodal relation.
Taking a hemisphere, that is, a disk, as a fundamental domain, we must then
identify two opposite points on the limit circle. Replacing the disk by the
homeomorphic hexagon, we obtain the figure above, the identification of
opposite points of the perimeter explaining the \emph{apparent} repetition of
the vertices 0, 1 and 2 and the corresponding edges.

This simplicial complex has six vertices, fifteen edges and ten triangles. The
1-skeleton is a complete graph with six vertices: any two vertices are
connected by an edge\footnote{A necessarily \emph{unique} edge in the context
of simplicial complexes.}. Computing by hand the homology groups of this
simplicial complex is a little lengthy. The Kenzo program obtains the result as
follows.

 \bmp
 \bmpi\verb|> (setf P2R|\empi
 \bmpi\verb|    (build-finite-ss|\empi
 \bmpi\verb|     '(v0 v1 v2 v3 v4 v5|\empi
 \bmpi\verb|       1 e01 (v1 v0) e02 (v2 v0) e03 (v3 v0) e04 (v4 v0) e05 (v5 v0)|\empi
 \bmpi\verb|         e12 (v2 v1) e13 (v3 v1) e14 (v4 v1) e15 (v5 v1) e23 (v3 v2)|\empi
 \bmpi\verb|         e24 (v4 v2) e25 (v5 v2) e34 (v4 v3) e35 (v5 v3) e45 (v5 v4)|\empix
 \bmpi\verb|       2 t013 (e13 e03 e01) t014 (e14 e04 e01) t024 (e24 e04 e02)|\empi
 \bmpi\verb|         t025 (e25 e05 e02) t035 (e35 e05 e03) t123 (e23 e13 e12)|\empi
 \bmpi\verb|         t125 (e25 e15 e12) t145 (e45 e15 e14) t234 (e34 e24 e23)|\empi
 \bmpi\verb|         t345 (e45 e35 e34))))|\empimx
 \bmpi\verb|Checking the 0-simplices...|\empi
 \bmpi\verb|Checking the 1-simplices...|\empi
 \bmpi\verb|Checking the 2-simplices...|\empi
 \bmpi\verb|[K1 Simplicial-Set]|\empi
 \emp

A Kenzo listing of this sort must be read as follows. The initial `\boxtt{>}'
is the Lisp prompt of this implementation. The user types out a Lisp statement,
here \boxtt{(setf...e35 e34))))} and the maltese cross \(\scriptstyle\maltese\)
(in fact not visible on the user screen) marks here the end of the Lisp
statement, the right number of closing parentheses is reached. The
corresponding Return key asks Lisp to evaluate the statement. Here a finite
simplicial \emph{set} is constructed according to the given description, it is
assigned to the symbol \boxtt{P2R}, and \emph{returned}, that is, displayed: it
is the Kenzo object \#1 (\boxtt{K1}) and it is a simplicial set; this is just a
small external display, the internal structure is not shown. Kenzo explains
beforehand it verifies the coherence of the definition of the simplicial set.

\label{24484}This construction of the projective plane is a little laborious.
In general the simplicial complexes are not used in ``good'' algebraic
topology; we have in fact used the more general notion of simplicial set. The
definition goes as follows: the \boxtt{build-finite-ss} Kenzo function is used,
which requires one argument, a list describing the finite
\underline{s}implicial \underline{s}et to be constructed; firstly the vertices
are given (six symbols \boxtt{v0} to \boxtt{v5}, then the edges (fifteen
symbols \boxtt{e01} to \boxtt{e45}) and for each of them their both ``faces''
(ends), and finally ten triangles and their faces (sides).

To arouse the interest for general simplicial sets, we immediately give the
minimal combinatorial definition of the projective plane as a simplicial
\emph{set}:

 \bmp
 \bmpi\verb|> (setf short-P2R|\empi
 \bmpi\verb|    (build-finite-ss|\empi
 \bmpi\verb|     '(v 1 e (v v) 2 t (e v e))))|\empim
 \bmpi\verb|Checking the 0-simplices...|\empi
 \bmpi\verb|Checking the 1-simplices...|\empi
 \bmpi\verb|Checking the 2-simplices...|\empi
 \bmpi\verb|[K6 Simplicial-Set]|\empi
 \emp

It is explained here only one vertex~`\boxtt{v}' is necessary, one
edge~`\boxtt{e}' and one triangle~`\boxtt{t}'. Both ends of the edge are the
unique vertex. The sides 0 and 2 of the triangle are the unique edge, and the
side 1 is collapsed on the vertex. It is clear \boxtt{P2R} and
\boxtt{short-P2R} are homeomorphic, and the second definition is much more
natural, but the underlying theory is not so easy.

The boundary matrices of \boxtt{P2R} are:

 \bmp
 \bmpi\verb|> (chcm-mat P2R 1)|\empim
 \bmpi\verb|========== MATRIX 6 lines + 15 columns =====|\empi
 \bmpi\verb|L1=[C1=-1][C2=-1][C3=-1][C4=-1][C5=-1]|\empi
 \bmpi\verb|L2=[C1=1][C6=-1][C7=-1][C8=-1][C9=-1]|\empi
 \bmpi\verb|L3=[C2=1][C6=1][C10=-1][C11=-1][C12=-1]|\empix
 \bmpi\verb|L4=[C3=1][C7=1][C10=1][C13=-1][C14=-1]|\empix
 \bmpi\verb|L5=[C4=1][C8=1][C11=1][C13=1][C15=-1]|\empi
 \bmpi\verb|L6=[C5=1][C9=1][C12=1][C14=1][C15=1]|\empi
 \bmpi\verb|========== END-MATRIX|\empi
 \emp

\noindent between degrees 1 and 0 and :

 \bmp
 \bmpi\verb|> (chcm-mat P2R 2)|\empim
 \bmpi\verb|========== MATRIX 15 lines + 10 columns =====|\empi
 \bmpi\verb|L1=[C1=1][C2=1]|\empi
 \bmpi\verb|L2=[C3=1][C4=1]|\empi
 \bmpi\verb|L3=[C1=-1][C5=1]|\empix
 \bmpi\verb|L4=[C2=-1][C3=-1]|\empix
 \bmpi\verb|L5=[C4=-1][C5=-1]|\empix
 \bmpi\verb|L6=[C6=1][C7=1]|\empix
 \bmpi\verb|L7=[C1=1][C6=-1]|\empix
 \bmpi\verb|L8=[C2=1][C8=1]|\empix
 \bmpi\verb|L9=[C7=-1][C8=-1]|\empix
 \bmpi\verb|L10=[C6=1][C9=1]|\empix
 \bmpi\verb|L11=[C3=1][C9=-1]|\empix
 \bmpi\verb|L12=[C4=1][C7=1]|\empix
 \bmpi\verb|L13=[C9=1][C10=1]|\empix
 \bmpi\verb|L14=[C5=1][C10=-1]|\empi
 \bmpi\verb|L15=[C8=1][C10=1]|\empi
 \bmpi\verb|========== END-MATRIX|\empi
 \emp

\noindent between degrees 2 and 1. Because large matrices can happen, a sparse
display is given; for example, for the last matrix, the row (line) 1 has only
two non null terms 1 in columns 1 and 2, the row 7 has a 1 in column 1 and a -1
in column 6, etc.

Computing the homology groups amounts to determining the kernel of the first
matrix, the image of the second one and the quotient of the kernel by the
image, a work a little painful. Significantly less painful for
\boxtt{short-P2R}:

 \bmp
 \bmpi\verb|> (chcm-mat short-P2R 1)|\empim
 \bmpi\verb|========== MATRIX 1 lines + 1 columns =====|\empi
 \bmpi\verb|L1=|\empi
 \bmpi\verb|========== END-MATRIX|\empix
 \bmpi\verb|> (chcm-mat short-P2R 2)|\empim
 \bmpi\verb|========== MATRIX 1 lines + 1 columns =====|\empi
 \bmpi\verb|L1=[C1=2]|\empi
 \bmpi\verb|========== END-MATRIX|\empi
 \emp

\noindent\label{49275} which means the chain-complex of \boxtt{short-P2R} is:
\[
0 \leftarrow \bZ \stackrel{0}{\leftarrow} \bZ \stackrel{\times 2}{\leftarrow}
\bZ \leftarrow 0
\]
if the ground ring is \(\bZ\) and it is then clear \(H_0 = \bZ\), \(H_1 =
\bZ_2\) and \(H_2 = 0\).

The homology groups of \boxtt{P2R} and \boxtt{short-P2R} can be computed by
Kenzo, for example the \(H_1\) groups.

 \bmp
 \bmpi\verb|> (homology P2R 1)|\empim
 \bmpi\verb|Homology in dimension 1 :|\empi
 \bmpi\verb|Component Z/2Z|\empi
 \bmpi\verb|---done---|\empix
 \bmpi\verb|> (homology short-P2R P2R 1)|\empim
 \bmpi\verb|Homology in dimension 1 :|\empi
 \bmpi\verb|Component Z/2Z|\empi
 \bmpi\verb|---done---|\empi
 \emp

The actual Kenzo display is more verbose and we keep here only the interesting
parts. You can do the same in degrees 0 and 2, both spaces have the same
homology groups, which ``confirms'' -- but does not prove -- both spaces are
homeomorphic.

\subsection{Chain-complex morphisms.}

\subsubsection{Definition.}

\begin{dfn}
--- \emph{Let \(A_\ast = \{A_q, d_q\}_q\) and \(B_\ast = \{B_q, d_q\}_q\) be two
chain-complexes\footnote{We do not hesitate to use the same symbol, \(d\) in
this case, for different\ldots\ differentials, the context being sufficient to
avoid any ambiguity.}. A \emph{chain-complex morphism} \(f: A_\ast \rightarrow
B_\ast\) is a collection of linear morphisms \(f = \{f_q: A_q \rightarrow
B_q\}_q\) satisfying the differential condition: for every \(q\), the relation
\(f_{q-1}d_q = d_q f_q\), or more simply \(df = fd\):
\[
\xymatrix{
 A_{q-1} \ar[d]_f & B_q \ar[l]_-d \ar[d]^f
 \\
 B_{q-1} & B_q \ar[l]^-d
 }
\]
is satisfied.}
\end{dfn}

More and more frequently, we will not indicate the indices of morphisms,
clearly implied by context. Also we use the same notation for a morphism and
some other morphisms directly deduced from the first one.

If \(f: A_\ast \rightarrow B_\ast\) is a chain-complex morphism, many other
maps are naturally induced; most often they are denoted by the same symbol,
\(f\) in this case. Because of the differential condition, the image of a cycle
is a cycle and we have induced maps \(f: Z_q(A_\ast) \rightarrow Z_q(B_\ast)\),
the same for the boundaries \(f: B_q(A_\ast) \rightarrow B_q(B_\ast)\), and for
homology classes and homology groups \(f: H_\ast(A_\ast) \rightarrow
H_\ast(B_\ast)\).

\subsubsection{Simplicial morphisms.}

\begin{dfn} ---
\emph{ Let \(K = (V, S)\)} and \(K' = (V', S')\) be two simplicial complexes. A
(simplicial) morphism \(f: K \rightarrow K'\) is a map \(f: V \rightarrow V'\)
satisfying the conditions:
\begin{itemize}
\item
The map \(f\) is compatible\footnote{Again a more general definition is
possible, without any order defined over \(V\) and \(V'\), see
Definition~\ref{70994}, but its use is then significantly more technical and
this matter is not directly our matter. Yet it is again a matter of
\emph{reductions}~\cite[Section 4]{ELMC1}!} with the orders defined over \(V\)
and \(V'\), see Definition~\ref{66339}. More precisely, if \(v \leq v'\) in
\(V\), then \(f(v) \leq f(v')\) in \(V'\)\footnote{In particular \(v < v'\) and
\(f(v) = f(v')\) is possible.}.
\item
If \(\sigma \in S\), then \(f(\sigma) \in S'\).
\end{itemize}
\end{dfn}

In other words, if \(v_0 < \cdots < v_k\) span a simplex of \(K\), then
\(f(v_0) \leq \cdots \leq f(v_k)\) span a simplex of \(V'\), but in the second
sequence, repetitions are allowed.

Now a simplicial morphism \(f: K \rightarrow K'\) induces a chain-complex
morphism again denoted by \(f: C_\ast(K) \rightarrow C_\ast(K')\). Only one
possible definition. If \(\sigma \in S_k(K)\) is a \(k\)-simplex of \(K\),
therefore a generator of \(C_k(K)\), two cases; if \(f(\sigma)\) again is a
\(k\)-simplex of \(K'\), that is, if there is no repetition in the images of
the vertices, then \(f(\sigma) := \ldots f(\sigma)\) where the left hand side
is understood in \(C_k(K')\) and the right hand one in \(S_k(K')\); if on the
contrary \(f(\sigma) \in S_{\mbox{\boldmath\(\ell\)}}(K')\) with \(\ell < k\),
then we decide \(f(\sigma) := 0\) in \(C_{\mbox{\boldmath\(k\)}}(K')\): the
image simplex is ``squeezed'' --- we will see later the appropriate terminology
is ``degenerate'', see Definition~\ref{62937} --- and this simplex do not
anymore contribute to homology. We advise the reader to verify the
chain-complex \emph{map} \(f: C_\ast(K) \rightarrow C_\ast(K')\) so defined is
compatible with the differentials, and therefore actually is a chain-complex
\emph{morphism}, the underlying sign game is instructive. Examples of
simplicial morphisms will be soon used in the proof of Theorem~\ref{38708}.

\subsection{Homotopy operators.}

\subsubsection{Definition and first properties.}

\begin{dfn}
--- \emph{Let \(A_\ast = \{A_q, d_q\}_q\) and \(B_\ast = \{B_q, d_q\}_q\) be two
chain-complexes. A homotopy operator \(h: A_\ast \rightarrow B_\ast\) is a
collection \(h = \{h_q: A_q \rightarrow B_{q+1}\}_q\) of linear maps. In other
words, it is a linear map \(h: A_\ast \rightarrow B_{\ast+1}\) of degree +1,
this degree being implicitly implied by the index `\(\ast+1\)' of
\(B_{\ast+1}\).}
\end{dfn}

In particular, no compatibility condition is required with the respective
differentials of \(A_\ast\) and \(B_\ast\). In the interesting cases, the
homotopy operator is rather ``seriously non-compatible'' with these
differentials.

\begin{dfn} ---
 \emph{Let \(f,g: A_\ast \rightarrow B_\ast\) be two chain-complex morphisms. A
homotopy operator \(h: A_\ast \rightarrow B_{\ast+1}\) is a \emph{homotopy
between \(f\) and \(g\)} if the relation \(g-f = dh + hd\) is satisfied.}
\end{dfn}

The next diagram shows there is a unique way to understand this relation when
you start from \(A_q\) and arrive at \(B_q\):
\[
\xymatrix{
 A_{q-1} \ar@<-2pt>[d]_f \ar@<2pt>[d]^g \ar[rd]|h
 & A_q \ar[l]_-d  \ar@<-2pt>[d]_f
 \ar@<2pt>[d]^g \ar[rd]|h & A_{q+1} \ar[l]_-d  \ar@<-2pt>[d]_f \ar@<2pt>[d]^g
 \\
 B_{q-1} & B_q \ar[l]^-d & B_{q+1} \ar[l]^-d
}
\]

\begin{prp} ---
 If two chain-complex morphisms \(f,g: A_\ast \rightarrow B_\ast\) are
 \emph{homotopic}, then the induced maps \(f,g: H_\ast(A_\ast) \rightarrow
 H_\ast(B_\ast)\) are \emph{equal}.
\end{prp}

\proof Let \(h\) be a homotopy between \(f\) and \(g\). If \(z\) is a
\(q\)-cycle representing the homology class \(\fh \in H_q(A_\ast)\), then the
relation \(gz - fz = dhz + hdz\) is satisfied; but \(z\) is a cycle and \(hdz =
0\), so that \(gz - fz = dhz\), which expresses the cycles \(fz\) and \(gz\)
representing the homology classes \(f\fh\) and \(g\fh\) are homologous, their
difference is a boundary; and therefore \(f\fh = g\fh\).\QED

\begin{dfn} ---
\emph{A \emph{homology equivalence} between two chain-complexes \(A_\ast\) and
\(B_\ast\) is a pair \((f,g)\) of chain-complex morphisms \(f: A_\ast
\rightarrow B_\ast\) and \(g: B_\ast \rightarrow A_\ast\) such that \(gf\) is
homotopic to \(\id{A_\ast}\) and \(fg\) is homotopic to \(\id{B_\ast}\).}
\end{dfn}

The terminology is not well stabilized, many authors use rather \emph{chain
equivalence}, or \emph{homotopy equivalence}. We feel more simple and clear our
terminology. We can also say that \(f: A_\ast \rightarrow B_\ast\) is a
homology equivalence if there exists a \emph{homological inverse} \(g: B_\ast
\rightarrow A_\ast\) such that the pair \((f,g)\) satisfies the above
definition.

\begin{prp} ---
If \(f: A_\ast \rightarrow B_\ast\) is a homology equivalence, then the induced
maps \(\{f_q: H_q(A_\ast) \rightarrow H_q(B_\ast)\}_q\) are isomorphisms.
\end{prp}

\proof The maps \(gf\) and \(fg\) are respectively \emph{homotopic} to
\(\id{A_\ast}\) and \(\id{B_\ast}\), so that the induced maps \(gf: H_q(A_\ast)
\rightarrow H_q(A_\ast)\) and \(fg: H_q(B_\ast) \rightarrow H_q(B_\ast)\) are
\emph{equal} to the corresponding identities. \QED

\subsubsection{Example.}

\begin{dfn} ---
\emph{If \(n \in \bN\) is a non-negative integer, we denote by
\(\underline{n}\) the initial segment of integers \(\underline{n} := (0, 1,
\ldots, n)\).}
\end{dfn}%

\begin{dfn} ---
\emph{The \emph{standard \(n\)-simplex} \(\Delta^n\) of dimension \(n\) is the
simplicial complex \((\underline{n}, \cP_\ast(\underline{n}))\) where
\(\cP_\ast(\underline{n})\) is the set of \emph{non-empty} subsets of
\(\underline{n}\).}
\end{dfn}

\begin{thr}\label{38708} ---
The homology groups of the standard simplex \(\Delta^n\) are null except
\(H_0(\Delta^n) = \fR\), the ground ring.
\end{thr}

\proof The result is obvious when \(n = 0\). Otherwise we can consider two
simplicial morphisms \(f: \Delta^0 \rightarrow \Delta^n\) and \(g: \Delta^n
\rightarrow \Delta^0\) where f(0) = 0 and \(g(i) = 0\) for every~\(i\). The
composition \(gf\) is the identity, the composition \(fg\) is not, but the
induced map \(fg: C_\ast(\Delta^n) \rightarrow C_\ast(\Delta^n)\) is homotopic
to the identity. The needed homotopy operator \(h: C_\ast(\Delta^n) \rightarrow
C_{\ast+1}(\Delta^n)\) is defined as follows; let \(\sigma = (i_0, \ldots,
i_k)\) a \(k\)-simplex generator of \(C_k(\Delta^n)\), that is, an ordered
sequence of \(k+1\) integers \(i_0 < \cdots < i_k\) of \(\underline{n}\). If
\(i_0 > 0\), we decide \(h(\sigma) = (0, i_0, \ldots, i_k)\); if on the
contrary \(i_0 = 0\), then we decide \(h(\sigma) = 0\). An interesting but
elementary computation then shows \(dh + hd = \id{C_\ast(\Delta^n)} - fg\). So
that the map \(fg: H_\ast(\Delta^n) \rightarrow H_\ast(\Delta^n)\) is simply
\emph{equal} to the identity and \(f: H_\ast(\Delta^0) \rightarrow
H_\ast(\Delta^n)\) is an isomorphism. \QED

\subsection{Exact sequences.}

\begin{dfn}
--- \emph{A chain-complex \(C_\ast = \{C_q, d_q\}_{q \in \bZ}\) is \emph{exact at degree
\(q\)} if \(\ker d_q = \im d_{q+1}\), in other words if \(H_q(C_\ast) = 0\), or
if \(Z_q(C_\ast) = B_q(C_\ast)\): every \(q\)-cycle is a \(q\)-boundary, no
``interesting'' cycle in degree \(q\). The chain-complex is \emph{exact} if it
is exact at every degree. In the same case, it is frequent also to state the
chain-complex is \emph{acyclic}; this does not mean there is no cycle, you must
understand there is no \emph{non-trivial} cycle, that is, a cycle which is not
a boundary. The expressions ``exact chain-complex'', ``acyclic chain-complex'',
``exact sequence'' are perfectly synonymous.}
\end{dfn}

\begin{prp}
--- Let \((C_\ast, d)\) be a chain-complex. If there exists a homotopy operator
\(h: C_\ast \rightarrow C_{\ast + 1}\) satisfying \emph{\(\id{} = dh + hd\)},
then the chain-complex \((C_\ast, d)\) is acyclic (or exact).
\end{prp}

\proof We can rewrite our relation \(\id{} - 0 = dh + hd\), that is, the
identity map is homotopic to the null map. The induced maps in homology
therefore are \emph{equal}. These induced maps are respectively the identity
maps and the null maps \(H_\ast(C_\ast) \rightarrow H_\ast(C_\ast)\). If \(M\)
is a module and if \(\id{M} = 0_M\), this is possible only if \(M = 0\). \QED

\begin{dfn}
--- {A \emph{short exact sequence} is a sequence of modules:
\[
0 \leftarrow C'' \stackrel{j}{\leftarrow} C \stackrel{i}{\leftarrow} C'
\leftarrow 0
\]
which is exact, that is in this case, the map \(i\) is injective, the map \(j\)
is surjective and \(\emph{\im} i = \ker j\).}
\end{dfn}

In particular the module \(C'\) is then canonically isomorphic to the kernel of
\(j\) and \(C''\) to the cokernel of \(i\). One says the central module \(C\)
is an \emph{extension} of \(C''\) \emph{by} \(C'\). In general there are
several possible extensions. For example if the ground ring is \(\bZ\), there
are two extensions of \(\bZ_6\) by \(\bZ_2\), namely \(\bZ_2 \oplus \bZ_6\) and
\(\bZ_{12}\) which are not isomorphic. The so-called \emph{extension problem},
how to determine in a particular case which is the right extension when the
left hand and right hand modules are known, is often a major problem in
homological algebra.

If a ``long'' sequence \(C_\ast\) is exact, there is no reason the short
sequence:
\[
0 \leftarrow C_{q-1} \stackrel{d_{q}}{\leftarrow} C_q
\stackrel{d_{q+1}}{\leftarrow} C_{q+1} \leftarrow 0
\]
is exact. To make it exact we must force \(d_{q+1}\) to be injective and
\(d_q\) to be surjective, and we obtain the short exact sequence:
\[
0 \leftarrow \im d_q \stackrel{d_{q}}{\leftarrow} C_q
\stackrel{d_{q+1}}{\leftarrow} C_{q+1}/\ker d_{q+1} \leftarrow 0
\]
but because of the exactness of the long sequence, we can write as well:
\[
0 \leftarrow \ker d_{q-1} \stackrel{d_{q}}{\leftarrow} C_q
\stackrel{d_{q+1}}{\leftarrow} \coker d_{q+2} \leftarrow 0
\]

This ``justifies'' the standard use of the long exact sequences: if a long
exact sequence \(C_\ast\) is given and if for every \(q\) the chain groups
\(C_{3q+1}\) and \(C_{3q+2}\) are known:
\[
\newcommand{\lc}[2]{\raisebox{-5pt}{\(\begin{array}{c}#1\\[-5pt]\textrm{\tiny #2}
\end{array}\)}}
\cdots \leftarrow \lc{C_{3q-2}}{known} \stackrel{d_{3q-1}}{\leftarrow}
\lc{C_{3q-1}}{known} \stackrel{d_{3q}}{\leftarrow} \lc{C_{3q}}{???}
\stackrel{d_{3q+1}}{\leftarrow} \lc{C_{3q+1}}{known}
\stackrel{d_{3q+2}}{\leftarrow} \lc{C_{3q+2}}{known} \leftarrow \cdots
\]
\emph{and also the morphisms \(d_{3q+2}\)}, then the chain group \(C_{3q}\) is
an extension of \(\ker d_{3q-1}\) by \(\coker d_{3q+2}\):
\[
0 \leftarrow \ker d_{3q-1} \leftarrow C_{3q} \leftarrow \coker d_{3q+2}
\leftarrow 0
\]
You understand you need to \emph{know} the maps \(d_{3q+2}\) for every \(q\) to
determine such kernels and cokernels, and when this is done, there remains an
extension problem.

In simple situations, this is easy. For example if every \(d_{3q+2}\) is known
to be an isomorphism, then kernels and cokernels are null and \(C_{3q} = 0\).
Another case is when every \(C_{3q-1}\) (resp. \(C_{3q+1}\)) is null, then
\(C_{3q} \cong C_{3q+1}\) (resp. \(C_{3q-1}\)).

But in general, the problem is highly non-trivial. Difference between
\emph{effective} homology and ordinary homology consists in particular in being
permanently \emph{vigilant} to be able to determine the maps \(d_{3q+2}\) and
to have sufficient data to solve the extension problem.

\subsection{The long exact sequence of a short exact sequence.}

It is a short exact sequence \emph{of chain-complexes} which produces a long
exact sequence.

\begin{thr}\label{15277}
\emph{\cite[II.4.1]{MCLN2}} --- Let:
\[
0 \leftarrow A_\ast \stackrel{j}{\leftarrow} B_\ast \stackrel{i}{\leftarrow}
C_\ast \leftarrow 0
\]
a short exact sequence \emph{of chain-complexes}. Then a canonical long exact
sequence of modules is obtained:
\[
\cdots \leftarrow H_{q-1}(C_\ast) \stackrel{\partial}{\leftarrow} H_{q}(A_\ast)
\stackrel{j}{\leftarrow} H_{q}(B_\ast) \stackrel{i}{\leftarrow} H_{q}(C_\ast)
\stackrel{\partial}{\leftarrow} H_{q+1}(A_\ast) {\leftarrow} \cdots
\]
\end{thr}

A short exact sequence of chain-complexes is a large diagram:
\[
\xymatrix{
   & \ar@{.>}[d] & \ar@{.>}[d] & \ar@{.>}[d] \\ 0 & A_{q+1} \ar[l] \ar[d] &
 B_{q+1} \ar[l]_j \ar[d] & C_{q+1} \ar[l]_i \ar[d] & 0 \ar[l] \\ 0 & A_{q}
 \ar[l] \ar[d] & B_{q} \ar[l]_j \ar[d] & C_{q} \ar[l]_i \ar[d] & 0 \ar[l] \\ 0
 & A_{q-1} \ar[l] \ar@{.>}[d] & B_{q-1} \ar[l]_j \ar@{.>}[d] & C_{q-1} \ar[l]_i
 \ar@{.>}[d] & 0 \ar[l] \\ &&&& }
\]
where all the horizontal short sequences are exact, and the three vertical
sequences are chain-complexes. It is understood \(i\) and \(j\) are
\emph{chain-complex morphisms}, that is, every square of our diagram is
commutative.

\proof See \cite[II.4.1]{MCLN2}. It is a matter of \emph{diagram chasing} in
our diagram. The \emph{connection morphism} for example \(\partial:
H_{q+1}(A_\ast) \rightarrow H_{q}(C_\ast)\) is of particular interest. It comes
from a diagram of objects:
\[
\xymatrix{
 {\fh_{q+1} \ni z_{q+1}} \ar@{|->}@{|->}[d] & c_{q+1} \ar@{|->}[l] \ar@{|->}[d] \\
 0 & b_q \ar@{|->}[l] \ar@{|->}[d] & {z_q \in \fh_q} \ar@{|->}[l] \ar@{|->}[d] \\
 & 0 & 0 \ar@{|->}[l]
}
\]
obtained as follows. Let \(\fh_{q+1} \in H_{q+1}(A_\ast)\) be a homology class
of \(A_\ast\) of degree \(q+1\). Let \(z_{q+1} \in A_{q+1}\) be a cycle
representing \(\fh_{q+1}\): the image in \(A_q\) by the vertical boundary map
is null. Because every \(j\) is surjective, we can find a chain \(c_{q+1} \in
B_{q+1}\) which is a \(j\)-preimage of \(z_{q+1}\). Then the vertical image
\(b_q\) of \(c_{q+1}\) must satisfy \(j(b_q) = 0\), for the left hand square is
commutative. Exactness of the horizontal row implies there is a unique
\(i\)-preimage \(z_q \in C_q\). The right hand square is also commutative. The
vertical image of \(b_q\) is null (\(dd = 0\)), so that, taking account of the
injectivity of \(i\), the vertical image of \(z_q\) is also null: \(z_q\) is a
cycle which defines a homology class \(\fh_q \in H_q(C_\ast)\). If \(c'_{q+1}\)
is another choice instead of \(c_{q+1}\) for a preimage of \(z_{q+1}\), then
this generates in the same way \(b'_q\), \(z'_q\) and \(\fh'_q\) but in fact
\(\fh_q = \fh'_q\), which results from the other diagram and analogous
arguments:
\[
\xymatrix{
 0 \ar@{|->}[d] & {c'_{q+1} - c_{q+1}} \ar@{|->}[l] \ar@{|->}[d] & {c''_{q+1}} \ar@{|->}[l]
 \ar@{|->}[d]
 \\
 0 & {b'_q - b_q} \ar@{|->}[l] & {z'_q - z_q} \ar@{|->}[l]
 }
\]
You must also prove the independance with respect to the choice of \(z_{q+1}
\in H_{q+1}(A_\ast)\), analogous exercise. The connexion map \(\partial:
\fh_{q+1} \mapsto \fh_q\) so defined is a module morphism --- exercise --- and
you must construct the other (induced) morphisms \(i\) and \(j\) of the long
exact sequence --- exercises. You must prove this long sequence actually
is\ldots\ exact. For example let us examine the exactness in
\(H_{q+1}(A_\ast)\). If ever \(\fh_{q+1}\) is the image of \(\fh'_{q+1} \in
H_{q+1}(B_\ast)\), we may choose \(c_{q+1} = z'_{q+1} \in \fh'_{q+1}\), it is a
cycle and \(b_q = 0\):
\[
\xymatrix{
 z_{q+1} & z'_{q+1} \ar@{|->}[l] \ar@{|->}[d]
 \\
 & b_q = 0 & z_q=0 \ar@{|->}[l]
}
\]
so that \(\fh_q = 0\). Conversely, if \(\fh_q = 0\), this means the final cycle
\(z_q\) is a boundary:
\[
\xymatrix{
 z_{q+1} & c_{q+1} \ar@{|->}[l] \ar@{|->}[d] & c'_{q+1} \ar@{|->}[d]
 \\
 & b_q & z_q \ar@{|->}[l]
}
\]
But this implies you have also this diagram:
\[
\xymatrix{
 0 & c''_{q+1} \ar@{|->}[l] \ar@{|->}[d] & c'_{q+1} \ar@{|->}[d]
 \ar@{|->}[l]
 \\
 & b_q & z_q \ar@{|->}[l]
}
\]
Now \(c_{q+1} - c''_{q+1}\) is another choice for \(c_{q+1}\), a choice which
is a (vertical) cycle, therefore defining a homology class \(\fh'_{q+1}\)
satisfying \(j(\fh'_{q+1}) = \fh_{q+1}\). If it is the first time you practice
this sport, you must carefully examine all the details of the other components
of the proof.\QED

We will see later, cf. Definition~\ref{87354}, that in \emph{effective
homology}, the analogous theorem needs a further hypothesis: the exactness
property of the short exact sequence of chain-complexes must be
\emph{effective}: an \emph{algorithm} must be present in the environment
\emph{returning} (producing) the various preimages which are required; it
happens it is always the case in the practical applications. And the
demonstration is then much easier and, very important, other \emph{algorithms}
are produced making \emph{effective} the exactness property of the resulting
long exact sequence.

\subsubsection{Examples.}\label{72185}

\begin{dfn} ---
\emph{A \emph{simplicial pair}} \((K, L)\) is a pair made of a simplicial
complex~\(K\) and a simplicial subcomplex \(L\) of \(K\).
\end{dfn}

The vertex set \(V_L\) of \(L\) is a subset \(V_L \subset V_K\) of the vertex
set of \(K\), the same for the simplices. For example let us define the
(simplicial) \((n-1)\)-sphere as the simplicial complex \(S^{n-1} =
(\underline{n}, \cP_\ast(\underline{n}) - \{\underline{n}\})\). A simplex is an
arbitrary subset of \(\underline{n}\) except the void subset \(\emptyset\) and
the full subset \(\underline{n}\). For example the 2-sphere is:
\[
 S^2 = (\underline{3}, \left\{\begin{array}{l}
 (0), (1), (2), (3),\\ (0,1), (0,2), (0,3),
(1,2), (1,3), (2,3),\\ (0,1,2), (0,1,3), (0,2,3), (1,2,3)
\end{array}\right\})
\]
It is the \emph{hollow} tetrahedron, while \(\Delta^3\) is the \emph{solid}
tetrahedron. In general \(S^{n-1}\) is a simplicial subcomplex of the standard
\(n\)-simplex \(\Delta^n\), and \((\Delta^n, S^{n-1})\) is a simplcial pair.

\begin{dfn} ---
\emph{Let \((K,L)\) be a \emph{simplicial pair}. The \emph{relative
chain-complex} \(C_\ast(K,L)\) is the quotient complex \(C_\ast(K,L) =
C_\ast(K) / C_\ast(L)\). The \emph{relative homology} \(H_\ast(K,L)\)
accordingly is \(H_\ast(K,L) := H_\ast(C_\ast(K)/C_\ast(L))\).}
\end{dfn}

The second component \(L\) is a simplicial subcomplex of the first one \(K\),
so that the corresponding chain-complex \(C_\ast(L)\) is a sub-chain-complex of
\(C_\ast(K)\), both differentials are compatible, which allows us to define the
quotient chain-complex \(C_\ast(K)/C_\ast(L)\) and the relative homology is the
homology of this quotient.

\begin{thr} ---
If \((K,L)\) is a simplicial pair, then a long exact sequence is obtained:
\[
\cdots \leftarrow H_{q-1}(L) \stackrel{\partial}{\leftarrow} H_q(K,L)
\stackrel{j}{\leftarrow} H_q(K) \stackrel{i}{\leftarrow} H_q(L)
\stackrel{\partial}{\leftarrow} H_{q+1}(K,L) \leftarrow \cdots
\]
\end{thr}

Note in particular the tempting result \(H_q(K,L) \cong H_q(K)/H_q(L)\) not
only in general is false, but it does not make sense: in general no inclusion
relation between \(H_q(L)\) and \(H_q(K)\). The inclusion relations \(C_q(L)
\subset C_q(K)\), \(Z_q(L) \subset Z_q(K)\) and \(B_q(L) \subset B_q(K)\) are
true, a canonical map \(H_q(L) \rightarrow H_q(L)\) therefore is defined, but
this map is not in general injective\footnote{\(3 \leq 6\) and \(2 \leq 6\) do
not imply \(3/2 \leq 6/6\).}.

\proof For every \(q\), we have a short exact sequence:
\[
0 \leftarrow C_q(K)/C_q(L) \stackrel{j}{\leftarrow} C_q(K)
\stackrel{i}{\leftarrow} C_q(L) \leftarrow 0
\]
But the maps \(i\) and \(j\) are compatible with the differentials of the
chain-complexes, so that we have in fact a short exact sequence \emph{of
chain-complexes}:
\[
0 \leftarrow C_\ast(K)/C_\ast(L) \stackrel{j}{\leftarrow} C_\ast(K)
\stackrel{i}{\leftarrow} C_\ast(L) \leftarrow 0
\]
and there remains to apply Theorem~\ref{15277}. \QED

\begin{prp}\label{35366} ---
Let \(S^{n-1}\) be the \((n-1)\)-sphere and \(\fR\) be the ground ring. Then,
if \(n \geq 2\), the homology groups \(H_q(S^{n-1})\) are null except
\(H_0(S^{n-1}) = H_{n-1}(S^{n-1}) =~\fR\).
\end{prp}

\proof Let us consider the pair \((\Delta^n, S^{n-1})\). Then all the chain
groups of \(C_\ast(\Delta^n)/C_\ast(S^{n-1})\) are null except
\(C_n(\Delta^n)/C_n(S^{n-1}) = \fR\): only the maximal simplex of \(\Delta^n\)
is not in \(S^{n-1}\). So that all the relative homology groups \(H_q(\Delta^n,
S^{n-1})\) are null except \(H_n(\Delta^n, S^{n-1}) = \fR\). Now in the long
exact sequence connecting \(H_\ast(\Delta^n)\) (known), \(H_\ast(\Delta^n,
S^{n-1})\) (known) and \(H_\ast(S^{n-1})\) (unknown), there are essentially two
interesting sections:
\[
\begin{array}{c}
 [H_0(\Delta^n, S^{n-1}) = 0] \leftarrow [H_0(\Delta^n) = \fR] \leftarrow
 [H_0(S^{n-1}) =\ ?] \leftarrow [H_1(\Delta^n, S^{n-1}) = 0]
 \\[5pt]
 \hspace{0pt}[H_{n-1}(\Delta^{n-1}) = 0] \leftarrow [H_{n-1}(S^{n-1}) =\ ?]
 \stackrel{\partial}{\leftarrow} [H_n(\Delta^n, S^{n-1}) = \fR] \leftarrow [H_n(\Delta^n) = 0]
\end{array}
\]
The extreme modules are null and, because of exactness, the central morphisms
are isomorphisms\footnote{What about the case \(n=1?\)}. \QED

Note also the connexion morphism \(\partial\) allows us to identify a canonical
representant (in fact unique up to sign, why?) for  a generator of
\(H_{n-1}(S^{n-1})\), namely the boundary of the maximal simplex \((0, \ldots,
n)\) of \(\Delta^n\); note this maximal simplex does not live in \(S^{n-1}\),
but its boundary does.

\subsection{About computability.}

All these didactical examples involve \emph{finite} simplicial complexes, so
that no theoretical computability problem here. However the benefit of the
various explained methods is already clear. For example let us take the
standard simplex \(\Delta^{10}\). If you want to compute \(H_5(\Delta^{10}) =
0\) by brute force, the boundary matrices to be considered are \(462 \times
462\) and \(462 \times 330\) so that proving kernel = image by ``simple''
computation is already a little serious. Moreover, when we will ask for
\emph{constructive} homology, see Section~\ref{76988}, we will have to be ready
to quickly return a boundary preimage for every cycle, for this cycle is
certainly homologous to 0. But Theorem~\ref{38708} gives immediately the
answer: this theorem in fact gives a \emph{reduction} (see
Definition~\ref{58661}) \(C_\ast(\Delta_n) \rrdc C_\ast(\Delta^0)\), so that
the homological problem for \(\Delta^n\) is equivalent to the same problem for
\(\Delta^0\), which problem is very simple.

This is a common situation. Even when the theoretical computability problem has
a trivial solution, an appropriate theoretical study of this computability
problem can produce dramatically better solutions. Another typical example is
the computation of the homology groups of the Eilenberg-MacLane spaces \(K(\pi,
n)\) for \(\pi\) an Abelian group of finite type. The general results quickly
sketched after Theorem~\ref{50468} prove the effective homology of these spaces
is computable. In the particular case where \(\pi\) is a \emph{finite} Abelian
group, the brute result is trivial, but the general method deduced from
Theorem~\ref{50468} remains essential for concrete computations. Let us
consider for example the group \(H_8(K(\bZ_2,4)) = \bZ_4\). A ``direct''
computation would require \(n_7 \times n_8\) and \(n_8 \times n_9\) matrices
with:
\[
\begin{array}{rcl}
n_7 &=& 34359509614
\\
n_8 &=& 1180591620442534312297
\\
n_9 &=& 85070591730234605240519066638188154620
\end{array}
\]
The method resulting from Theorem~\ref{50468} reduces the problem to a smaller
chain-complex with the analogous dimensions being \(n'_7 = 4\), \(n'_8 = 8\)
and \(n'_9 = 15\). The result is then obtained in less than 2 seconds with a
modest laptop, most computing time being devoted to \emph{compute} these small
matrices, which remains a non-trivial task.

But the most striking results which are obtained  in \emph{constructive}
homological algebra concern cases where the studied chain-complex defining
homology groups is \emph{not} of finite type. It is the general situation for
loop spaces leading to our simple solution for Adams' problem, see
Section~\ref{05508}. For Eilenberg-MacLane spaces, if you are interested by
\(H_8(K(\bZ, 4)) = \bZ_3 + \bZ\), then the corresponding numbers \(n_7\),
\(n_8\) and \(n_9\) are \emph{infinite}. Eilenberg and MacLane in their
wonderful papers~\cite{ELMC1,ELMC2} explained how to obtain an equivalent
chain-complex of finite type (in this case \(n'_7 = 1\) and \(n'_8 = n'_9 =
2\)) giving the right homology groups; it was the first historical case where
\emph{constructive} homological algebra was implicitly used, without any
constructive terminology\ldots\ The matter of these notes consists in a
systematic extension of these constructive methods, producing results with a
very general scope. The strong difference with the general style of
Eilenberg-MacLane's work is that we will have to \emph{keep in our environment}
the original \(K(\bZ, 4)\) itself, with a functional implementation, as a
\emph{locally effective} object, for example to be able to compute a spectral
sequence where this object is involved.

\section{Spectral sequences.}

\subsection{Introduction.}

The previous section explained how the long exact sequence of a short exact
sequence of chain-complexes can be used to determine some unknown homology
groups. The typical case being the last example: three chain-complexes are
present in the environment: \(C_\ast(\Delta^n)\), \(C_\ast(\Delta^n, S^{n-1})\)
and \(C_\ast(S^{n-1})\). We knew the homology groups \(H_\ast(\Delta^n)\) and
\(H_\ast(\Delta^n, S^{n-1})\) and the long exact sequence allowed us to obtain
the unknown groups \(H_\ast(S^{n-1})\).

This is the general process in the computation of homology groups, and the same
for homotopy groups in Algebraic Topology. Objects more and more complicated
are considered, and the invariants of the new objects are obtained from
invariants of simpler previous ones and of a careful study of the
``difference''.

But this description unfortunately is simplistic. For example if you know
\(H_\ast(K)\) and \(H_\ast(K,L)\), and you try to deduce \(H_\ast(L)\), the
long exact sequence:
\[
\cdots \leftarrow H_q(K,L) \stackrel{j}{\leftarrow} H_q(K)
\stackrel{i}{\leftarrow} H_q(L) \stackrel{\partial}{\leftarrow} H_{q+1}(K,L)
\stackrel{j}{\leftarrow} H_{q+1}(K) \leftarrow \cdots
\]
produces a short exact sequence:
\[
0 \leftarrow \ker j \stackrel{i}{\leftarrow} H_q(L)
\stackrel{\partial}{\leftarrow} \coker j \leftarrow 0
\]
and if you \emph{do not know} the exact nature of the map \(j\), you cannot
proceed; as soon as the situation becomes a little more complicated, it is the
\emph{most frequent} case. And even if you can determine the groups \(\ker j\)
and \(\coker j\), there remains a possible extension problem needing also other
informations to be solved.

\begin{clm} ---
Except in\ldots\ exceptional situations, the long exact sequence of homology
\underline{is not} an algorithm allowing one to compute an unknown group when
the four neighbouring groups are known.
\end{clm}

And most books about homological agebra do not give any explanations about this
lack in the theory; they even give frequently the unpleasant feeling that they
\emph{hide} this deficiency, but more probably the authors do not have a
sufficiently precise knowledge of the very nature of the constructive
requirement.

The present text is exactly devoted to provide the missing tools allowing one
to transform usual homological algebra into a modern \emph{constructive}
theory. Experience shows it is quite elementary, but two essential notions are
required. From an algorithmic point of view, \emph{higher-order functional
programming} is definitively necessary; fortunately, standard computer science
knows this matter from a long time, and the modern application tools are the so
called functional programming languages such as Lisp, ML, Haskell, with
powerful compilers. From an ``ordinary'' mathematical point of view, the
\emph{basic perturbation lemma} (Henri Cartan, Shih Weishu~\cite{SHIH}, Ronnie
Brown~\cite{BRWNR1}) is the key point.

Which probably explains the terrible delay of homological algebra with respect
to the modern constructive point of view is the fact that the elementary
results explained here to satisfy constructiveness cannot be seriously used
\emph{without machines and programs}. The analogy with commutative algebra
thirty years ago is striking. Noone can now hope to work a long time in
commutative algebra without using Groebner bases. Groebner bases are
elementary, but cannot be used without auxiliary machines and programs.
Groebner bases are quite elementary, the same for the homological perturbation
lemma. More precisely the basic theory of Groebner bases is elementary, but
looking for more and more efficient implementations, in particular for special
cases, remains an active research subject. And the situation is the same for
the homological perturbation lemma.

This section is devoted to a short presentation of the \emph{spectral sequence}
theory, and the situation for spectral sequences is the same as for the long
exact sequence. In \emph{exceptional} cases, a spectral sequence \emph{can be}
a process giving unknown homology groups when other homology groups are given,
but in the general situation, the constructive requirement is not satisfied: no
\emph{general algorithm} can be deduced from the spectral sequence theory. The
homological perturbation lemma will allow us to replace the usual spectral
sequences by \emph{effective} versions. Which is quite amazing in this case is
the fact that these effective versions are \emph{terribly simpler} to design
than the ordinary spectral sequences, but, think of the Groebner bases, these
effective spectral sequences cannot be used without the corresponding machines
and programs.

The spectral sequences are also used in commutative algebra, because of the
frequent presence of Koszul complexes playing an important role through their
homology groups. We will see the point of view presented here also gives very
interesting results in commutative algebra, mainly to compute the
\emph{effective} homology of Koszul complexes, richer than the ordinary
homology; for example this effective homology gives a direct method to compute
a resolution of the initial module.

\subsection{Notion of spectral sequence.}

One of the best references to attack the subject is~\cite[Chapter XI]{MCLN2}.
The \emph{didactic} quality of this text is the highest we know. In particular
MacLane begins to explain how to \emph{use} a spectral sequence before proving
its construction, a wise organization. We just give here a short presentation
of the general structure of spectral sequences, advising the reader to study
~\cite[Chapter XI]{MCLN2} for further details and results. The most complete
reference about spectral sequences of course is~\cite{MCCL}.

\begin{dfn} ---
\emph{A \emph{spectral sequence} is a collection \(\{E^r_{p,q}, d^r_{p,q}\}^{r
\geq r_0}_{p,q \in \bZ}\) satisfying the following properties:
\begin{itemize}
\item
Every \(E^r_{p,q}\) is an \(\fR\)-module (\(\fR\) is the underlying ground
ring).
\item
Every \(d^r_{p,q}\) is a morphism \(d^r_{p,q}: E^r_{p,q} \rightarrow E^r_{p-r,
q+r-1}\).
\item
Every composition \(d^r_{p,q} d^r_{p+r, q-r+1}\) is null, so that a homology
group \(H^r_{p,q} = \ker d^r_{p,q} / \im d^r_{p+r, q-r+1}\) is defined.
\item
For every \(r \geq r_0\), \(p,q \in \bZ\), an isomorphism \(H^r_{p,q} \cong
E^{r+1}_{p,q}\) is provided.
\end{itemize}}
\end{dfn}

A geometric representation of the notion of spectral sequence is very useful.
Look at this figure\footnote{Strongly inspired by the analogous scheme
of~\cite[Section XI.1]{MCLN2}, without any kind permission of
Springer-Verlag.}:

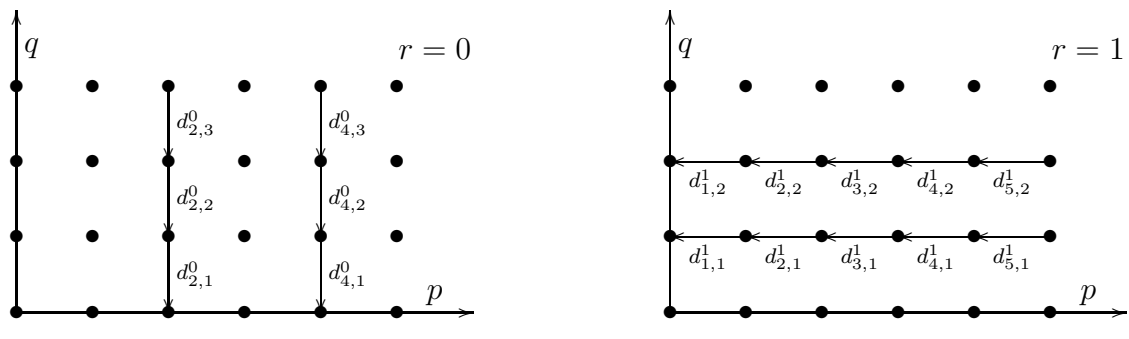
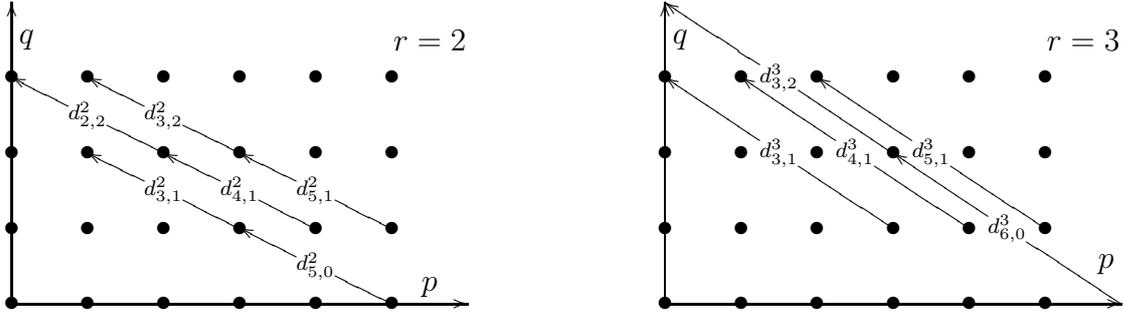
\begin{center}
 \mbox{\begin{xy}<1cm,0cm>:<0cm,1cm>::
 (0,0)*{\hspace{0pt}} ;
 (0,0)*{\bullet} ; (1,0)*{\bullet} ; (2,0)*{\bullet} ; (3,0)*{\bullet} ;
 (4,0)*{\bullet} ; (5,0)*{\bullet} ;
 (0,1)*{\bullet} ; (1,1)*{\bullet} ; (2,1)*{\bullet} ; (3,1)*{\bullet} ;
 (4,1)*{\bullet} ; (5,1)*{\bullet} ;
 (0,2)*{\bullet} ; (1,2)*{\bullet} ; (2,2)*{\bullet} ; (3,2)*{\bullet} ;
 (4,2)*{\bullet} ; (5,2)*{\bullet} ;
 (0,3)*{\bullet} ; (1,3)*{\bullet} ; (2,3)*{\bullet} ; (3,3)*{\bullet} ;
 (4,3)*{\bullet} ; (5,3)*{\bullet} ;
 (5.5,0.1) *!D{p} ; (0.1,3.5) *!L{q} ; (5.5,3.5)*{r=0} ;
 \ar (0,0);(6,0)
 \ar (0,0);(0,4)
 \ar^{d^0_{4,1}} (4,1) ; (4,0)
 \ar^{d^0_{4,2}} (4,2) ; (4,1)
 \ar^{d^0_{4,3}} (4,3) ; (4,2)
 \ar^{d^0_{2,1}} (2,1) ; (2,0)
 \ar^{d^0_{2,2}} (2,2) ; (2,1)
 \ar^{d^0_{2,3}} (2,3) ; (2,2)
 \end{xy}}
 \hfill
 \mbox{\begin{xy}<1cm,0cm>:<0cm,1cm>::
 (0,0)*{\hspace{0pt}} ;
 (0,0)*{\bullet} ; (1,0)*{\bullet} ; (2,0)*{\bullet} ; (3,0)*{\bullet} ;
 (4,0)*{\bullet} ; (5,0)*{\bullet} ;
 (0,1)*{\bullet} ; (1,1)*{\bullet} ; (2,1)*{\bullet} ; (3,1)*{\bullet} ;
 (4,1)*{\bullet} ; (5,1)*{\bullet} ;
 (0,2)*{\bullet} ; (1,2)*{\bullet} ; (2,2)*{\bullet} ; (3,2)*{\bullet} ;
 (4,2)*{\bullet} ; (5,2)*{\bullet} ;
 (0,3)*{\bullet} ; (1,3)*{\bullet} ; (2,3)*{\bullet} ; (3,3)*{\bullet} ;
 (4,3)*{\bullet} ; (5,3)*{\bullet} ;
 (5.5,0.1) *!D{p} ; (0.1,3.5) *!L{q} ; (5.5,3.5)*{r=1} ;
 \ar (0,0);(6,0)
 \ar (0,0);(0,4)
 \ar^{d^1_{5,1}} (5,1) ; (4,1)
 \ar^{d^1_{4,1}} (4,1) ; (3,1)
 \ar^{d^1_{3,1}} (3,1) ; (2,1)
 \ar^{d^1_{2,1}} (2,1) ; (1,1)
 \ar^{d^1_{1,1}} (1,1) ; (0,1)
 \ar^{d^1_{5,2}} (5,2) ; (4,2)
 \ar^{d^1_{4,2}} (4,2) ; (3,2)
 \ar^{d^1_{3,2}} (3,2) ; (2,2)
 \ar^{d^1_{2,2}} (2,2) ; (1,2)
 \ar^{d^1_{1,2}} (1,2) ; (0,2)
\end{xy}}
\\[1cm]
 \mbox{\begin{xy}<1cm,0cm>:<0cm,1cm>::
 (0,0)*{\hspace{0pt}} ;
 (0,0)*{\bullet} ; (1,0)*{\bullet} ; (2,0)*{\bullet} ; (3,0)*{\bullet} ;
 (4,0)*{\bullet} ; (5,0)*{\bullet} ;
 (0,1)*{\bullet} ; (1,1)*{\bullet} ; (2,1)*{\bullet} ; (3,1)*{\bullet} ;
 (4,1)*{\bullet} ; (5,1)*{\bullet} ;
 (0,2)*{\bullet} ; (1,2)*{\bullet} ; (2,2)*{\bullet} ; (3,2)*{\bullet} ;
 (4,2)*{\bullet} ; (5,2)*{\bullet} ;
 (0,3)*{\bullet} ; (1,3)*{\bullet} ; (2,3)*{\bullet} ; (3,3)*{\bullet} ;
 (4,3)*{\bullet} ; (5,3)*{\bullet} ;
 (5.5,0.1) *!D{p} ; (0.1,3.5) *!L{q} ; (5.5,3.5)*{r=2} ;
 \ar (0,0);(6,0)
 \ar (0,0);(0,4)
 \ar|{d^2_{5,1}} (5,1) ; (3,2)
 \ar|{d^2_{3,2}} (3,2) ; (1,3)
 \ar|{d^2_{4,1}} (4,1) ; (2,2)
 \ar|{d^2_{2,2}} (2,2) ; (0,3)
 \ar|{d^2_{5,0}} (5,0) ; (3,1)
 \ar|{d^2_{3,1}} (3,1) ; (1,2)
 \end{xy}}
 \hfill
 \mbox{\begin{xy}<1cm,0cm>:<0cm,1cm>::
 (0,0)*{\hspace{0pt}} ;
 (0,0)*{\bullet} ; (1,0)*{\bullet} ; (2,0)*{\bullet} ; (3,0)*{\bullet} ;
 (4,0)*{\bullet} ; (5,0)*{\bullet} ;
 (0,1)*{\bullet} ; (1,1)*{\bullet} ; (2,1)*{\bullet} ; (3,1)*{\bullet} ;
 (4,1)*{\bullet} ; (5,1)*{\bullet} ;
 (0,2)*{\bullet} ; (1,2)*{\bullet} ; (2,2)*{\bullet} ; (3,2)*{\bullet} ;
 (4,2)*{\bullet} ; (5,2)*{\bullet} ;
 (0,3)*{\bullet} ; (1,3)*{\bullet} ; (2,3)*{\bullet} ; (3,3)*{\bullet} ;
 (4,3)*{\bullet} ; (5,3)*{\bullet} ;
 (5.8,0.4) *!D{p} ; (0.1,3.5) *!L{q} ; (5.5,3.5)*{r=3} ;
 \ar (0,0);(6,0)
 \ar (0,0);(0,4)
 \ar|{d^3_{5,1}} (5,1) ; (2,3)
 \ar|{d^3_{4,1}} (4,1) ; (1,3)
 \ar|{d^3_{3,1}} (3,1) ; (0,3)
 \ar|{d^3_{6,0}} (6,0) ; (3,2)
 \ar|{d^3_{3,2}} (3,2) ; (0,4)
\end{xy}}
\end{center}

You must consider the integer parameter \(r\) as a discrete time, a spectral
sequence can be thought of as a dynamical system. The figures represent the
state of our system at times 0, 1, 2 and 3. Usually the initial time \(r_0\) is
0, 1 or 2 and we will mot mention it anymore. A convenient terminology consists
in considering a spectral sequence as a book where the page~\(r\) is visible at
time \(r\). The page~\(r\) is made of a collection of modules
\(\{E^r_{p,q}\}_{p, q \in \bZ}\); every morphism \(d^r_{p,q}\) starts from
\(E^r_{p,q}\) and goes to \(E^r_{p-r,q+r-1}\): the shift for the
\emph{horizontal} degree \(p\) is \(-r\), the page number, and the shift for
the \emph{total} degree \(p+q\) always is \(-1\), so that the shift for the
\emph{vertical} degree \(q\) necessarily is \(r-1\). On the above figures, only
a few differentials are displayed.

Because of the rule about the composition of two successive \(d^r_{p,q}\)'s,
every page is a collection of chain-complexes, where in the above
representation the (oriented) ``slope'' is \((q-1)/(-p)\). Therefore the page
\(r\) produces a collection of homology groups \(H^r_{p,q}\) and \(H^r_{p,q}\)
is isomorphic to \(E^{r+1}_{p,q}\), one usually says \(H^r_{p,q}\) ``is''
\(E^{r+1}_{p,q}\). In short, every page is a collection of chain-complexes and
the collection of corresponding homology groups ``is'' the next page, but it is
exactly at this point the constructiveness property in most situations fails,
point examined later.

Very frequently the spectral sequence is null outside some quadrant of the
\((p,q)\)-plane; for example, if \(p\ \textrm{or}\ q < 0 \Rightarrow E^r_{p,q}
= 0\), one says it is a \emph{first quadrant} spectral sequence; a \emph{second
quadrant} spectral sequence is null for \(p > 0\) or \(q < 0\).

\begin{dfn} ---
\emph{A spectral sequence \(\{E^r_{p,q}, d^r_{p,q}\}\) is \emph{convergent} if
for every \mbox{\(p, q \in \bZ\)} the relations \(d^r_{p,q} = 0 =
d^r_{p+r,q-r+1}\) holds for \(r \geq r_{p,q}\).}
\end{dfn}

If the convergence property is satisfied, then \(E^r_{p,q} = H^r_{p,q} \
\textrm{``}\!=\!\!\textrm{''}\ E^{r+1}_{p,q} = \cdots\) for \(r = r_{p,q}\).

\begin{dfn} ---
\emph{If a spectral sequence \(\{E^r_{p,q}, d^r_{p,q}\}\) is convergent,
\(E^\infty_{p,q} := \textrm{``}\lim\!\textrm{''}_{r \rightarrow \infty}
E^r_{p,q}\).}
\end{dfn}

As usual, only the isomorphism class of the limit is defined. For example a
first quadrant spectral sequence is necessarily convergent, because \(r > p
\Rightarrow d^r_{p,q} = 0\) and \(r
> q+1 \Rightarrow d^r_{p+r, q-r+1} = 0\). A second quadrant spectral sequence is
not necessarily convergent.

\begin{dfn} ---
\emph{Let \(\{H_n\}_{n \in \bZ}\) be a collection of modules, probably a
collection of interesting homology groups. The spectral sequence \(\{E^r_{p,q},
d^r_{p,q}\}\) \emph{converges towards} \(\{H_n\}_{n \in \bZ}\) if the spectral
sequence is convergent and if there exists a filtration \(\{H_{p,q}\}_{p+q=n}\)
of every \(H_n\) such that \(E^\infty_{p,q} \cong H_{p,q} / H_{p-1,q+1}\). The
collection \(\{H_n\}_{n \in \bZ}\) is then called the \emph{abutment} of the
spectral sequence.}
\end{dfn}

The filtration of \(H_n\) must be coherent, that is \(H_{p,q} \subset H_{p+1,
q-1}\), \mbox{\(\cap_{p+q = n} H_{p,q} = 0\)} and \(\cup_{p+q=n} H_{p,q} =
H_n\). It is an increasing filtration indexed on \(p\), but it is convenient to
recall the second index \(q\), which also implicitly implies the total degree
\(n = p+q\). For example for a first quadrant spectral sequence, the context
would imply \(0 = H_{-1,n+1} \subset H_{0,n} \subset H_{1,n-1} \subset \cdots
\subset H_{n,0} = H_n\).

There is a strange but convenient notation for such a convergence property:
\[
E^r_{p,q} \Rightarrow H_{p+q}
\]
The convergence is implicitly concerned by which happens when \(r \rightarrow
\infty\). The double arrow `\(\Rightarrow\)' instead of the simple one
`\(\rightarrow\)' recalls the convergence property is quite complex. The
ambiguous index of \(H_{p+q}\) means some filtration of \(H_n\) is involved
correlated to the double indexation of \(E^\infty_{p,q}\).

\subsection{The Serre spectral sequence.}

The Serre spectral sequence was invented in 1950 of course by Jean-Pierre
Serre, using anterior works of Jean Leray and Jean-Louis Koszul; this spectral
sequence allowed him to determine many homotopy groups, in particular sphere
homotopy groups. This spectral sequence concerns the \emph{fibrations}:
\[
F \hookrightarrow E \rightarrow B
\]
where \(F\) is the \emph{fibre} space, \(B\) the \emph{base} space and \(E\)
the \emph{total} space. These were initially topological spaces, but this
notion of fibration can be generalized to many other situations. The total
space \(E\) is to be considered as a \emph{twisted product} of the base space
\(B\) by the fibre space \(F\). The underlying twisting operator \(\tau\) is
defined by different means according to the context, but the idea is constant:
\(\tau\) explains how the twisted product \(E = F \times_\tau B\) is different
from the trivial product \(F \times B\), which product depends in turn on the
category we are working in. See for example~\cite[Section I.2]{STNR} for the
original case of the \emph{fibre bundles}; the twist then is a collection of
\emph{coordinate functions}.

\begin{thr} ---
 Let \(E = F \times_\tau B\) be a topological fibration with a base space \(B\)
 simply connected. Then a first quadrant spectral sequence \(\{E^r_{p,q},
 d^r_{p,q}\}_{r \geq 2}\) is defined with \(E^2_{p,q} = H_p(B ; H_q(F))\) and
 \(E^r_{p,q} \Rightarrow H_{p+q}(E)\).
\end{thr}

We are working in ``\emph{general}'' topology and there is a process called
\emph{singular homology} associating with every topological space \(X\), every
integer \(n\) and every abelian group \(\fG\) (here not necessarily a ring) a
homology group \(H_n(X ; \fG)\), the \mbox{\(n\)-th} homology group of \(X\)
with coefficients in \(\fG\). The process is strongly inspired by which had
been done in Section~\ref{43994}, but adapted to an arbitrary topological space
thanks to the notion of \emph{singular simplex}, see for example~\cite[Chapter
VII]{ELST}. We will not be concerned by the (interesting) definition of the
singular homology groups. It happens if the topological space \(X\) comes from
a simplicial complex, the simplicial homology groups and the singular homology
groups are canonically isomorphic.  The role of \emph{coefficients}, \(H_q(F)\)
here, simpler, was explained at Definition~\ref{80060}; note there is no
misprint: the coefficient group used to define \(H_p(B ; H_q(F))\) is in turn a
homology group \(H_q(F) := H_q(F ; \fR)\) if \(\fR\) is the underlying ground
ring.

The Serre spectral sequence establishes a rich set of relations between the
homology groups \(H_\ast(F)\), \(H_\ast(E)\) and \(H_\ast(B)\) of the fibre
space, total space and base space of a fibration, at least when the base space
is simply connected. It is frequently somewhat implicitly ``suggested'' this
spectral sequence is a process allowing one for example to \emph{compute} the
groups \(H_\ast(E)\) when the groups \(H_\ast(B)\) and \(H_\ast(F)\) are known.
But in general this is false. In general the differentials \(d^2_{p,q}\) are
unknown, and even if you know them, you will be able to compute the
\(E^3_{p,q}\)'s, but to continue the process, you need now the differentials
\(d^3_{p,q}\) and in general you do not have the necessary information to
compute them. And so on.

And if by any chance you reach the limit groups \(E^\infty_{p,q}\), you have
the group \(H_{0,n} = E^\infty_{0,n}\), but to determine the next component of
the filtration of \(H_n\), the exact sequence:
\[
0 \leftarrow E^\infty_{1,n-1} \leftarrow H_{1,n-1} \leftarrow H_{0,n}
\leftarrow 0
\]
shows \(H_{1,n-1}\) is the solution of an extension problem which can be very
difficult, we will show a typical example. And if you succeed, again an
extension problem for \(H_{2,n-2}\), and so on\ldots

\begin{clm} ---
Let \(F \hookrightarrow E \rightarrow B\) be a given fibration with \(B\)
simply connected. Except in\ldots\ exceptional situations, the Serre spectral
sequence \underline{is not} an algorithm allowing to compute \(H_\ast(E)\) when
\(H_\ast(B)\) and \(H_\ast(F)\) are known. More generally, except in
exceptional situations, the page \(r+1\) of a spectral sequence cannot be
deduced from the page \(r\) and the other available data.
\end{clm}

These negative appreciations of course must not reduce the interest of the
various known spectral sequences. The point of view used here is the following:
yes the spectral sequences are in many circumstances quite essential, yes they
allowed to obtain many very interesting results, but their general organisation
is not algorithmic; how this deficiency with respect to usual modern
mathematics could be corrected? In short, how a spectral sequence can be made
\emph{constructive}? It is our main concern.

\subsubsection{A positive example.}\label{11076}

When writing these notes, MacLane's excellent book~\cite{MCLN2} is not far and
instead of considering the loop spaces of spheres, the first example of this
book, we use the symmetrical example of \(B\bH_\ast = P^\infty \bH\), the
classifying space of the multiplicative group \(\bH_\ast\) of the quaternion
field \(\bH\), in other words the infinite quaternionic projective space. The
topological group \(\bH_\ast\) automatically generates \cite{MLNR3} a universal
principal fibration:
\[
\bH_\ast \hookrightarrow E\bH_\ast \rightarrow B\bH_\ast.
\]
This means our group \(\bH_\ast\) freely acts on the total space \(E\bH_\ast\),
the base space \(B\bH_\ast\) being the corresponding homogeneous space
\(B\bH_\ast = E\bH_\ast / \bH_\ast\). Saying the fibration is \emph{universal}
amounts to requiring the total space \(E\bH_\ast\) is \emph{contractible}, that
is, has the homotopy type of a point, which needs a few definitions to be
understood.

\begin{dfn} ---
\emph{Two continuous maps \(f_0, f_1: X \rightarrow Y\) are \emph{homotopic} if
there exists a continuous map \(F: X \times [0,1] \rightarrow Y\) such that
\(f_0(x) = F(x,0)\) and \(f_1(x) = F(x,1)\) for every \(x \in X\).}
\end{dfn}

In other words, two continuous maps \(f_0\) and \(f_1\) are homotopic if a
continuous \emph{deformation} \(F\) can be installed between them.

\begin{thr} \emph{\cite[Section VII.7]{ELST}} ---
If two continuous maps \(f, g: X \rightarrow Y\) are homotopic, then the
induced maps \(f_\ast, g_\ast: H_\ast(X;\fR) \rightarrow H_\ast(Y;\fR)\)
between singular homology groups, with respect to an arbitrary coefficient
group \(\fR\), are equal.
\end{thr}

\begin{dfn} ---
\emph{A continuous map \(f: X \rightarrow Y\) is a \emph{homotopy equivalence}
if there exists another continuous map \(g: Y \rightarrow X\) such that \(gf\)
is homotopic to \(\id{X}\) and \(fg\) is homotopic to~\(\id{Y}\).}
\end{dfn}

\begin{dfn} ---
\emph{Two topological spaces \(X\) and \(Y\) have the \emph{same homotopy type}
if there exists a homotopy equivalence \(f: X \rightarrow Y\).}
\end{dfn}

A homotopy equivalence \(f: X \rightarrow Y\) therefore induces isomorphisms
\(f: H_\ast(X) \stackrel{\cong}{\rightarrow} H_\ast(Y)\). The same homotopy
type requires isomorphic homology groups, but unfortunately the converse is
false: it is an \emph{open problem} to give \emph{computable} characteristic
conditions for homotopy equivalence. It is generally ``understood'' such a
condition is given by the so-called \emph{Postnikov-invariants} or
\emph{\(k\)-invariants}, but this is false~\cite{RBSR8}.

\begin{dfn} ---
\emph{A topological space \(X\) is \emph{contractible} if it has the homotopy
type of a point.}
\end{dfn}

\begin{dfn}---
\emph{If \(G\) is a topological group, a principal fibration:
\[
G \hookrightarrow EG \rightarrow BG
\]
is \emph{universal} if the total space \(EG\) is contractible~\cite{STNR,HSML}.
It is then proved the homotopy type of the so-called \emph{classifying space}
\(BG\) is well defined up to homotopy.}
\end{dfn}

A point \(\ast\)\footnote{Not to be confused with a generic index such as in
\(H_\ast(X)\).} is a (multiplicative) \emph{unit} in the topological world: the
product \(\ast \times X\) is canonically homeomorphic to \(X\). The total space
\(EG\) of a universal fibration is some twisted product \(EG = G \times_\tau
BG\), and because this product has the homotopy type of a point, the
classifying space \(BG\) can be understood as a ``twisted inverse'' of the
initial group \(G\), but \emph{up to homotopy}. Such a twisted inverse is
itself unique up to homotopy.

These classifying spaces \(BG\) are very important and the computation of their
homology groups as well. The dual notion of \emph{cohomology groups} of these
classifying spaces leads to the important notion of characteristic classes of
principal fibrations~\cite{MLST}. And it is essential to be able to compute the
homology groups of classifying spaces.

In the case of our quaternionic multiplicative group \(\bH_\ast\), a radial
homotopy easily allows one to prove the inclusion \(S^3 \hookrightarrow
\bH_\ast\) is a homotopy equivalence, so that the homology groups of
\(\bH_\ast\) and \(S^3\) are the same. Proposition~\ref{35366} proves the
\emph{simplicial} homology groups \(H_n(S^3; \bZ)\) are null except \(H_0(S^3)
= H_3(S^3)~=~\bZ\). And the isomorphism theorem between singular and simplicial
homology groups~\mbox{\cite[Section~VII.10]{ELST}} implies it is the same for
the singular homology groups, so that \(H_n(\bH_\ast) = 0\) except
\(H_0(\bH_\ast) = H_3(\bH_\ast) = \bZ\).

It is convenient to shorten \(\bH_\ast =: G\), \(E\bH_\ast =: E\) and
\(B\bH_\ast =: B\), so that the diagram:
\[
G \hookrightarrow E \rightarrow B
\]
denotes now our specific universal fibration around the topological group \(G =
\bH_\ast\).

Because the total space \(E\) is contractible, all its homology groups are null
except \(H_0(EG) = \bZ\). Knowing the groups \(H_\ast(G)\) and \(H_\ast(E)\),
the game now consists in guessing the groups \(H_\ast(BG)\).

The Serre spectral sequence of a fibration involving \(G\), \(E\) and \(B\)
describes \(E^2_{p,q} = H_p(B; H_q(G))\); in general the \emph{universal
coefficient theorem}~\cite[Section~V.11]{MCLN2} allows to deduce the groups
\(H_n(X; \fR)\), where \(\fR\) is an arbitrary abelian group, from the
\emph{integer} homology groups \(H_n(X; \bZ)\) most often denoted by \(H_n(X)\)
in short. Here the situation is simple: \(H_p(BG; H_q(G)) =~0\) except for \(q
= 0\) or 3 where \(H_p(BG; H_q(G)) = H_p(BG; \bZ)\). In particular \(H_0(BG,
H_q(G)) = 0\) except for \(q = 0\) or 3 where the value is \(\bZ\); this is
because \(BG\) is necessarily connected, which implies \(H_0(BG; \bZ) = \bZ\).

The initial state of our study is the known state of the page 2 of our spectral
sequence:

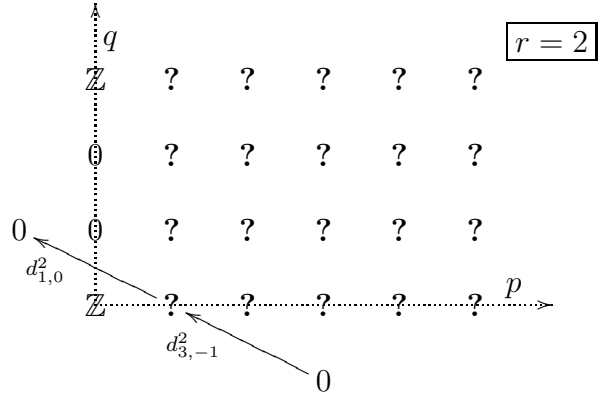
\[\begin{xy}<1cm,0cm>:<0cm,1cm>::
 (0,0)*{\hspace{0pt}} ;
 (-1,1)*{0} ; (3,-1)*{0} ;
 (0,0)*{\bZ} ; (1,0)*{\boldqm} ; (2,0)*{\boldqm} ; (3,0)*{\boldqm} ;
 (4,0)*{\boldqm} ; (5,0)*{\boldqm} ;
 (0,1)*{0} ; (1,1)*{\boldqm} ; (2,1)*{\boldqm} ; (3,1)*{\boldqm} ;
 (4,1)*{\boldqm} ; (5,1)*{\boldqm} ;
 (0,2)*{0} ; (1,2)*{\boldqm} ; (2,2)*{\boldqm} ; (3,2)*{\boldqm} ;
 (4,2)*{\boldqm} ; (5,2)*{\boldqm} ;
 (0,3)*{\bZ} ; (1,3)*{\boldqm} ; (2,3)*{\boldqm} ; (3,3)*{\boldqm} ;
 (4,3)*{\boldqm} ; (5,3)*{\boldqm} ;
 (5.5,0.1) *!D{p} ; (0.1,3.5) *!L{q} ; (6,3.5)*{\fbox{\(r=2\)}} ;
 \ar@{.>} (0,0);(6,0)
 \ar@{.>} (0,0);(0,4)
 \ar^<(0.8){d^2_{1,0}} (0.8,0.1) ; (-0.8,0.9)
 \ar^<(0.8){d^2_{3,-1}} (2.8,-0.9) ; (1.2,-0.1)
 \end{xy}\]

\hspace*{\parindent}It is a first quadrant spectral sequence, so that the
\(d^2\)-arrows arriving and starting from \(E^2_{1,0}\) necessarily are null.
This entails \(E^3_{1,0} = \ker d^2_{1,0} / \im d^2_{3,-1} = E^2_{1,0} / 0 =
E^2_{1,0}\). The same for the next \(r\)'s, and \(E^2_{1,0} = E^3_{1,0} =
\cdots = E^\infty_{1,0}\). At the abutment of the spectral sequence, we know
all the \(H_n(EG)\) are null for \(n>0\), so that certainly all the
corresponding \(E^\infty_{p,q} = H_{p,q}(EG) / H_{p-1,q+1}(EG)\) also are null.
This implies that when \(E^r_{p,q}\) becomes fixed, that is when \(r
> \max(p,q+1)\), the relation \(E^r_{p,q} = 0\) is satisfied: for every \((p,q)\) with \(p\)
or \(q > 0\), the spectral group \(E^r_{p,q}\) must ``die''.

But for \(E^2_{1,0}\), it must be already died at time \(r=2\), otherwise
\(E^2_{1,0} = E^\infty_{1,0} \neq 0\). We have proved \(E^2_{1,0} = H_1(BG) =
0\). Now \(E^2_{1,q} = H_1(BG; H_q(G)) = 0\), because of the universal
coefficient theorem. So that we obtain this partial description for the page~3
of our spectral sequence.

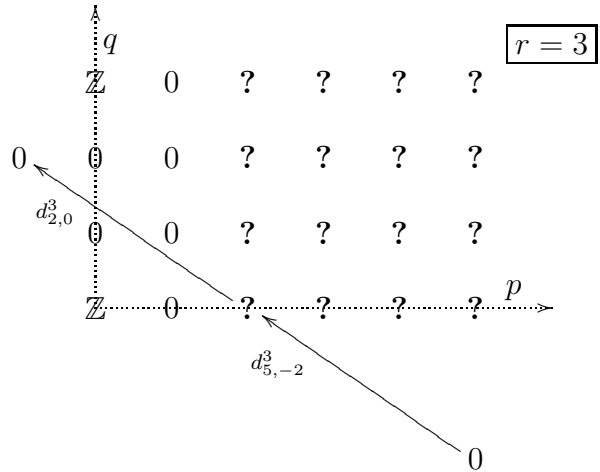
\[\begin{xy}<1cm,0cm>:<0cm,1cm>::
 (0,0)*{\hspace{0pt}} ;
 (-1,2)*{0} ; (5,-2)*{0} ;
 (0,0)*{\bZ} ; (1,0)*{0} ; (2,0)*{\boldqm} ; (3,0)*{\boldqm} ;
 (4,0)*{\boldqm} ; (5,0)*{\boldqm} ;
 (0,1)*{0} ; (1,1)*{0} ; (2,1)*{\boldqm} ; (3,1)*{\boldqm} ;
 (4,1)*{\boldqm} ; (5,1)*{\boldqm} ;
 (0,2)*{0} ; (1,2)*{0} ; (2,2)*{\boldqm} ; (3,2)*{\boldqm} ;
 (4,2)*{\boldqm} ; (5,2)*{\boldqm} ;
 (0,3)*{\bZ} ; (1,3)*{0} ; (2,3)*{\boldqm} ; (3,3)*{\boldqm} ;
 (4,3)*{\boldqm} ; (5,3)*{\boldqm} ;
 (5.5,0.1) *!D{p} ; (0.1,3.5) *!L{q} ; (6,3.5)*{\fbox{\(r=3\)}} ;
 \ar@{.>} (0,0);(6,0)
 \ar@{.>} (0,0);(0,4)
 \ar^<(0.8){d^3_{2,0}} (1.8,0.1) ; (-0.8,1.9)
 \ar^<(0.8){d^3_{5,-2}} (4.8,-1.9) ; (2.2,-0.1)
 \end{xy}\]
 This argument can be repeated for the column 2, starting this time from
 \(E^2_{2,0}\), and also for the column 3, starting from \(E^2_{3,0}\) and in this
 case, taking account of \(E^2_{1,1} = 0\). We obtain \(H_2(BG) = H_3(BG) =
 0\). But for \(E^r_{4,0}\), there is something new when \(r=4\):
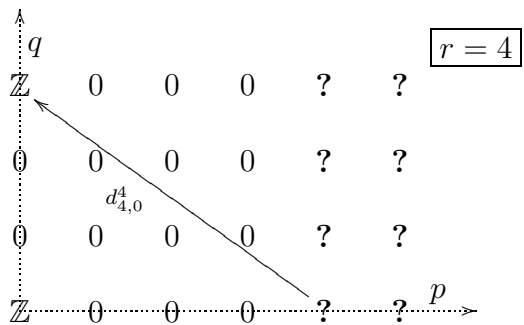
\[\begin{xy}<1cm,0cm>:<0cm,1cm>::
 (0,0)*{\hspace{0pt}} ;
 (0,0)*{\bZ} ; (1,0)*{0} ; (2,0)*{0} ; (3,0)*{0} ;
 (4,0)*{\boldqm} ; (5,0)*{\boldqm} ;
 (0,1)*{0} ; (1,1)*{0} ; (2,1)*{0} ; (3,1)*{0} ;
 (4,1)*{\boldqm} ; (5,1)*{\boldqm} ;
 (0,2)*{0} ; (1,2)*{0} ; (2,2)*{0} ; (3,2)*{0} ;
 (4,2)*{\boldqm} ; (5,2)*{\boldqm} ;
 (0,3)*{\bZ} ; (1,3)*{0} ; (2,3)*{0} ; (3,3)*{0} ;
 (4,3)*{\boldqm} ; (5,3)*{\boldqm} ;
 (5.5,0.1) *!D{p} ; (0.1,3.5) *!L{q} ; (6,3.5)*{\fbox{\(r=4\)}} ;
 \ar@{.>} (0,0);(6,0)
 \ar@{.>} (0,0);(0,4)
 \ar^<(0.6){d^4_{4,0}} (3.8,0.2) ; (0.2,2.8)
 \end{xy}\]
 The group in position \((0,3)\) starts non-null at time 2: \(E^2_{0,3} = \bZ\)
 and it must also die. Because the columns 2 or 3 are null, this group can be
 killed only at time 4, which implies the arrow \(d^4_{4,0}: E^4_{4,0}
 \rightarrow E^4_{0,3}\) necessarily is an isomorphism. So that \(E^4_{4,0} =
 \bZ\) and going back to time 2, \(H_4(BG) := H_4(BG, \bZ) = H_4(BG, H_0(G)) =
 \bZ\). The conclusion is the following: The column 4 for \(r = 2,3,4\) is made of \(E^r_{4,q} =
 \bZ\) for \(q = 0\) or 3, \(E^r_{4,q} = 0\) for \(q \neq 0\) and 3. But 4 is also the last time
 where it is possible to kill \(E^4_{4,3} = \bZ\), which implies by the same
 argument \(H_8(BG) = \bZ\). Finally we have proved:
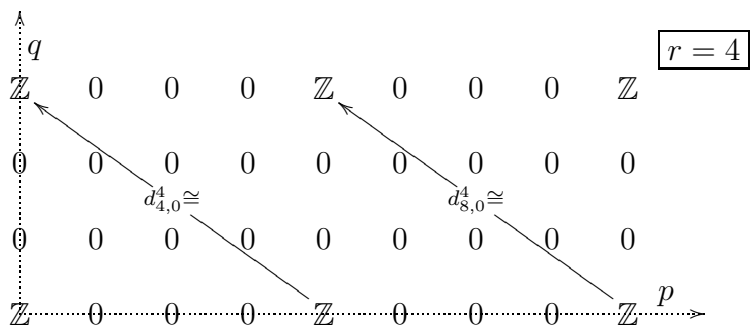
\[\begin{xy}<1cm,0cm>:<0cm,1cm>::
 (0,0)*{\hspace{0pt}} ;
 (0,0)*{\bZ} ; (1,0)*{0} ; (2,0)*{0} ; (3,0)*{0} ;
 (4,0)*{\bZ} ; (5,0)*{0} ; (6,0)*{0} ; (7,0)*{0} ; (8,0)*{\bZ} ;
 (0,1)*{0} ; (1,1)*{0} ; (2,1)*{0} ; (3,1)*{0} ;
 (4,1)*{0} ; (5,1)*{0} ; (6,1)*{0} ; (7,1)*{0} ; (8,1)*{0} ;
 (0,2)*{0} ; (1,2)*{0} ; (2,2)*{0} ; (3,2)*{0} ;
 (4,2)*{0} ; (5,2)*{0} ; (6,2)*{0} ; (7,2)*{0} ; (8,2)*{0} ;
 (0,3)*{\bZ} ; (1,3)*{0} ; (2,3)*{0} ; (3,3)*{0} ;
 (4,3)*{\bZ} ; (5,3)*{0} ; (6,3)*{0} ; (7,3)*{0} ; (8,3)*{\bZ} ;
 (8.5,0.1) *!D{p} ; (0.1,3.5) *!L{q} ; (9,3.5)*{\fbox{\(r=4\)}} ;
 \ar@{.>} (0,0);(9,0)
 \ar@{.>} (0,0);(0,4)
 \ar|{d^4_{4,0} \cong} (3.8,0.2) ; (0.2,2.8)
 \ar|{d^4_{8,0} \cong} (7.8,0.2) ; (4.2,2.8)
 \end{xy}\]

 \begin{thr}---
 The homology groups \(H_{4n+k}(B\bH_\ast)\) are null for \(k = 1, 2\) or \(3\) and the
 groups \(H_{4n}(B\bH_\ast)\) are all equal to \(\bZ\).
 \end{thr}

Because of the very specific situation, this is a (rare) case where the
spectral sequence can be entirely described: \(E^r_{p,q} = \bZ\) if [\(p = q =
0\) and \(2 \leq r \leq \infty\)] or [\(p = 4n\) and \(q = (0\) or \(3\)) and
\(2 \leq r \leq 4\)]. Otherwise every \(E^r_{p,q} = 0\). The only non null
\(d^r_{p,q}\)'s occur for \(r=4\), \(p=4n\) and \(q=0\) and they are
isomorphisms \(d^4_{4n,0}: E^4_{4n,0} = \bZ \stackrel{\cong}{\rightarrow}
E^4_{4n-4,3} = \bZ\).

\subsubsection{A negative example.}\label{62859}

Jean-Pierre Serre got one of the 1954 Fields Medal, mainly for his computations
of sphere homotopy groups, where the principal tool was his famous spectral
sequence. To illustrate the non-constructive nature of this spectral sequence,
we describe the beginning of his computations, up to the first point where the
method failed.

If \((X, \ast)\) is a \emph{based} topological space, that is, some \emph{base
point} \(\ast \in X\) is given, the set of homotopy classes of continuous maps
\(\pi_n(X) := [(S^n, \ast), (X, \ast)]\) has a natural commutative group
structure for \(n \geq 2\) and it is a popular sport in algebraic topology to
find out the groups \(\pi_n(S^k)\). It is not hard to prove \(\pi_n(S^k) = 0\)
for \(n < k\), \(\pi_n(S^n) = \bZ\), and the first event in the story was the
amazing discovery by Hopf in 1935 that \(\pi_3(S^2) = \bZ\). In 1937,
Freudenthal proved \(\pi_4(S^2) = \bZ_2\) (in algebraic topology, it is common
to denote \(\bZ_2\) the quotient group \(\bZ/\bZ_2\), not the \(p\)-adic
ring!), and Serre at the beginning of the fifties computed many sphere homotopy
groups; in particular he proved \(\pi_6(S^3)\) has 12 elements, but could not
choose between \(\bZ_{12}\) and \(\bZ_2 + \bZ_6\). We want to describe the
point where the spectral sequence method fails.

To compute \(\pi_4(S^3)\), we can proceed as follows. We consider a fibration:
\[F_2 \hookrightarrow X_4 \rightarrow S^3\] where \(F_2 := K(\bZ, 2)\) is an
\emph{Eilenberg-MacLane space}, in this case a connected space where every
homotopy group is null except \(\pi_2(F_2) = \bZ\); such a space is well
defined up to homotopy, it happens we can take \(F_2 = P^\infty \bC\). The
beginning of the spectral sequence uses the homology of \(S^3\), null except
\(H_0(S^3) = H_3(S^3) = \bZ\) and the homology of \(F_2\), null except
\(H_{2q}(F_2) = \bZ\) for every \(q \geq 0\). The critical page of the spectral
sequence is the page \(r=3\):
\[\begin{xy}<1cm,0cm>:<0cm,1cm>::
 (0,0)*{\hspace{0pt}} ;
 (0,0)*{\bZ} ; (1,0)*{0} ; (2,0)*{0} ; (3,0)*{\bZ} ;
 (0,1)*{0} ; (1,1)*{0} ; (2,1)*{0} ; (3,1)*{0} ;
 (0,2)*{\bZ} ; (1,2)*{0} ; (2,2)*{0} ; (3,2)*{\bZ} ;
 (0,3)*{0} ; (1,3)*{0} ; (2,3)*{0} ; (3,3)*{0} ;
 (0,4)*{\bZ} ; (1,4)*{0} ; (2,4)*{0} ; (3,4)*{\bZ} ;
 (3.5,0.1) *!D{p} ; (0.1,4.5) *!L{q} ; (4,4.5)*{\fbox{\(r=3\)}} ;
 \ar@{.>} (0,0);(4,0)
 \ar@{.>} (0,0);(0,5)
 \ar|<(0.75){d^3_{3,0} = \times 1} (2.8,0.2) ; (0.2,1.8)
 \ar|<(0.75){d^3_{3,2} = \times 2} (2.8,2.2) ; (0.2,3.8)
 \end{xy}\]

 Our fibration is not completely defined, we have not explained how the
 twisting operator \(\tau\) of \(X_4 = F_2 \times_\tau S^3\) is defined. We do
 not want to give the details, but the twisting operator \(\tau\) is entirely
 defined\footnote{Up to sign.} by the fact the arrow \(d^3_{3,0}\) is an
 isomorphism. It is then necessary to know the arrows \(d^3_{3,2q}\); in this
 \emph{particular case}, a \emph{specific tool} gives the solution; examining the
 multiplicative structure of the analogous spectral sequence in cohomology, it
 can be proved the arrow \(d^3_{3,2q}: \bZ \rightarrow \bZ\) is the
 multiplication by \(q+1\). This implies the \(E^3_{3,2q}\) die and
 \(E^r_{0,2q} = \bZ_{q}\) for \(4 \leq r \leq \infty\) and \(q \geq 2\). So
 that the Serre spectral sequence entirely gives the homology groups \(H_0(X_4)
 = \bZ\), \(H_{2n}(X_4) = \bZ_n\) for \(n \geq 2\) and the other \(H_n(X_4)\)
 are null. In particular, please believe the Hurewicz theorem~\cite[Section~IV.7]{WHTH1}
 and the long exact sequence of homotopy~\cite[Section~IV.8]{WHTH1} imply \(\pi_4(S^3)
 = \pi_4(X_4) = H_4(X_4) = \bZ_2\), a result known by Freudenthal.

Conclusion: the computation of \(H_\ast(X_4)\) needs more information than
which is given by the spectral sequence itself, information coming from the
multiplicative structure of the \(X_4\)-cohomology.

To compute \(\pi_5(S^3)\), we must consider a new fibration:
\[
F_3 \hookrightarrow X_5 \rightarrow X_4
\]
where \(F_3 = K(\bZ_2, 3)\) again is an Eilenberg-MacLane space, with every
homotopy group null except \(\pi_3(F_3) = \bZ_2\), chosen because \(\pi_4(X_4)
= \bZ_2\). We cannot give the details allowing us to use the spectral sequence
in this case, but the next figure gives an idea of the complexity of the
situation\footnote{The details of this spectral sequence which are shown here
have been obtained thanks to Ana Romero's program~\cite{RMRS}, a good
illustration of its possibilities.}.
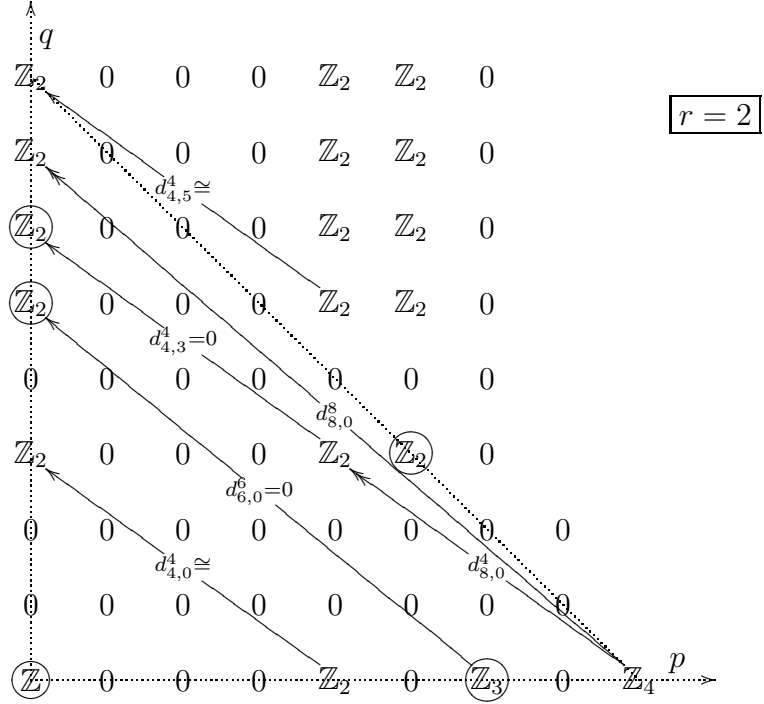
\[\begin{xy}<1cm,0cm>:<0cm,1cm>::
 (0,0)*{\hspace{0pt}} ;
 (0,0)*+{\bZ}*\frm{o} ; (1,0)*{0}     ; (2,0)*{0} ; (3,0)*{0} ; (4,0)*{\bZ_2} ;
 (5,0)*{0}   ; (6,0)*+{\bZ_3}*\frm{o} ; (7,0)*{0} ; (8,0)*{\bZ_4} ;
 (0,1)*{0} ; (1,1)*{0} ; (2,1)*{0} ; (3,1)*{0} ;
 (4,1)*{0} ; (5,1)*{0} ; (6,1)*{0} ; (7,1)*{0} ;
 (0,2)*{0} ; (1,2)*{0} ; (2,2)*{0} ; (3,2)*{0} ;
 (4,2)*{0} ; (5,2)*{0} ; (6,2)*{0} ; (7,2)*{0} ;
 (0,3)*{\bZ_2} ; (1,3)*{0} ; (2,3)*{0} ; (3,3)*{0} ;
 (4,3)*{\bZ_2} ; (5,3)*+{\bZ_2}*\frm{o} ; (6,3)*{0} ;
 (0,4)*{0} ; (1,4)*{0} ; (2,4)*{0} ; (3,4)*{0} ;
 (4,4)*{0} ; (5,4)*{0} ; (6,4)*{0} ;
 (0,5)*+{\bZ_2}*\frm{o} ; (1,5)*{0} ; (2,5)*{0} ; (3,5)*{0} ;
 (4,5)*{\bZ_2} ; (5,5)*{\bZ_2} ; (6,5)*{0} ;
 (0,6)*+{\bZ_2}*\frm{o} ; (1,6)*{0} ; (2,6)*{0} ; (3,6)*{0} ;
 (4,6)*{\bZ_2} ; (5,6)*{\bZ_2} ; (6,6)*{0} ;
 (0,7)*{\bZ_2} ; (1,7)*{0} ; (2,7)*{0} ; (3,7)*{0} ;
 (4,7)*{\bZ_2} ; (5,7)*{\bZ_2} ; (6,7)*{0} ;
 (0,8)*{\bZ_2} ; (1,8)*{0} ; (2,8)*{0} ; (3,8)*{0} ;
 (4,8)*{\bZ_2} ; (5,8)*{\bZ_2} ; (6,8)*{0} ;
 (8.5,0.1) *!D{p} ; (0.1,8.5) *!L{q} ; (9,7.5)*{\fbox{\(r=2\)}} ;
 \ar@{.>} (0,0);(9,0)
 \ar@{.>} (0,0);(0,9)
 \ar@{.} (8,0);(0,8)
 \ar|{d^4_{4,0} \cong} (3.8,0.2) ; (0.2,2.8)
 \ar|{d^6_{6,0} = 0} (5.8,0.2) ; (0.2,4.8)
 \ar|{d^4_{8,0}}@{->>} (7.8,0.2) ; (4.2,2.8)
 \ar|{d^4_{4,3}=0} (3.8,3.2) ; (0.2,5.8)
 \ar|{d^8_{8,0}}@{->>} (7.8,0.2) ; (0.2,6.8)
 \ar|{d^4_{4,5} \cong} (3.8,5.2) ; (0.2,7.8)
 \end{xy}\]

We show the page \(r = 2\) and all the arrows which are necessary to determine
the \(E^\infty_{p,q}\) for \(p + q \leq 8\). Up to \(p + q \leq 8\), the
\(E^2_{p,q}\) which remain definitively alive are circled, the others die, and
in particular \(E^2_{8,0}\) will die in two steps at times \(r = 4\) and 8.

The twisting operator of the fibration is the unique one giving \(d^4_{4,0} =
\id{\bZ_2}\) and \(H_3(X_5) = H_4(X_5) = 0\). No choice for \(d^6_{6,0}\), it
is necessarily the null map, so that \(E^7_{0,5} = E^\infty_{0,5} =
H_{0,5}(X_5) = H_5(X_5) = \bZ_2\). Again the Hurewitz theorem and the long
homotopy exact sequence imply \(H_5(X_5) = \pi_5(X_5) = \pi_5(X_4) = \pi_5(S^3)
= \bZ_2\); it was the first important result obtained by Serre.

It happens the arrow \(d^4_{8,0}\) is the only non-null arrow from \(\bZ_4\) to
\(\bZ_2\); this implies the next arrow \(d^4_{4,3}\) is null. It was the last
possible event for \(E^r_{0,6}\), so that \(\bZ_2 = E^7_{0,6} = E^\infty_{0,6}
= H_{0,6}\). Another ingredient for \(H_6(X_5)\) is \(E^2_{6,0} =
E^\infty_{6,0} = \bZ_3\). Therefore two stages in the filtration of
\(H_6(X_5)\) at the abutment, which gives the short exact sequence:
 \[
 0 \leftarrow \bZ_3 \leftarrow H_6(X_5) \leftarrow \bZ_2 \leftarrow 0
 \]
 The group \(H_6(X_5)\) is an extension of \(\bZ_3\) by \(\bZ_2\), and
 fortunately there exists a unique extension \(H_6(X_5) = \bZ_6\).

 Please believe that \(d^8_{8,0}\) kills \(E^8_{0,7}\) and \(E^8_{8,0} = \ker
d^4_{8,0} = \bZ_2\); in particular \(H_7(X_5) = 0\). Also \(d^4_{4,5}\) in
particular kills \(E^4_{0,8}\) which implies \(H_8(X_5) = E^\infty_{5,3} =
E^2_{5,3} = \bZ_2\).

We have obtained the sequence \((\bZ, 0, 0, 0, 0, \bZ_2, \bZ_6, 0, \bZ_2)\) for
the first homology groups of \(X_5\).

Jean-Pierre Serre was able to obtain all the necessary ingredients for the
various \(d^r_{p,q}\) which play an essential role in the beginning of this
spectral sequence. The main ingredients are the multiplicative structure in
cohomology and more generally the module structure with respect to the Steenrod
algebra \(\mathcal{A}_2\), a subject not studied here.

Let us study now the next fibration:
\[
F_4 \hookrightarrow X_6 \rightarrow X_5
\]
with \(F_4 = K(\bZ_2, 4)\) and the twisting operator is chosen to have
\(\pi_5(X_6) = 0\). The part of the spectral sequence interesting for us is
simple, we need only the part \(p + q \leq 6\):
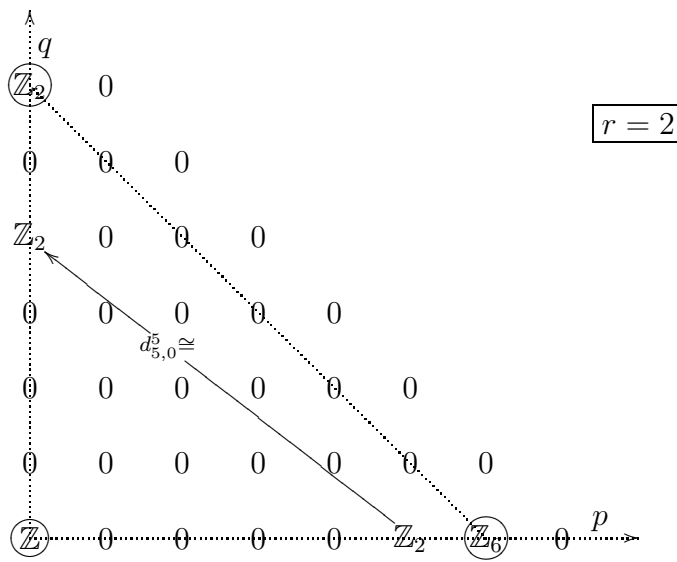
\[\begin{xy}<1cm,0cm>:<0cm,1cm>::
 (0,0)*{\hspace{0pt}} ;
 (0,0)*+{\bZ}*\frm{o} ; (1,0)*{0}     ; (2,0)*{0} ; (3,0)*{0} ; (4,0)*{0} ;
 (5,0)*{\bZ_2}   ; (6,0)*+{\bZ_6}*\frm{o} ; (7,0)*{0} ;
 (0,1)*{0} ; (1,1)*{0} ; (2,1)*{0} ; (3,1)*{0} ;
 (4,1)*{0} ; (5,1)*{0} ; (6,1)*{0} ;
 (0,2)*{0} ; (1,2)*{0} ; (2,2)*{0} ; (3,2)*{0} ;
 (4,2)*{0} ; (5,2)*{0} ;
 (0,3)*{0} ; (1,3)*{0} ; (2,3)*{0} ; (3,3)*{0} ;
 (4,3)*{0} ;
 (0,4)*{\bZ_2} ; (1,4)*{0} ; (2,4)*{0} ; (3,4)*{0} ;
 (0,5)*{0} ; (1,5)*{0} ; (2,5)*{0} ;
 (0,6)*+{\bZ_2}*\frm{o} ; (1,6)*{0} ;
 (7.5,0.1) *!D{p} ; (0.1,6.5) *!L{q} ; (8,5.5)*{\fbox{\(r=2\)}} ;
 \ar@{.>} (0,0);(8,0)
 \ar@{.>} (0,0);(0,7)
 \ar@{.} (6,0);(0,6)
 \ar|<(0.65){d^5_{5,0} \cong} (4.8,0.2) ; (0.2,3.8)
 \end{xy}\]

 The same argument as before produces a short exact sequence:
 \[
 0 \leftarrow \bZ_6 \leftarrow H_6(X_6) \leftarrow \bZ_2 \leftarrow 0
 \]
but this time two possible extensions, the trivial one \(\bZ_2 + \bZ_6\) and
the twisted one~\(\bZ_{12}\). And the Serre spectral sequence does not give any
information, given the available data, which allows us to choose the right
extension. The conclusion of Serre was only: ``The group \(\pi_6(S^3)
=\pi_6(X_6) = H_6(X_6)\) has 12 elements''. Two years later, Barratt and
Paechter, using a \emph{quite specific method}, proved the group \(\pi_6(S^3)\)
in fact contains an element of order~4, so that finally \(\pi_6(S^3) =
\bZ_{12}\), it is the non-trivial extension which is the right one.
See~\cite{BRPC} and also~\cite[pp.105-110]{SERR3}.

The modern process to determine homotopy groups consists in using the Adams
spectral sequence and the numerous other related spectral sequences. Some exact
sequences, in particular the chromatic exact sequence, are also very useful.
The basic reference about these methods is the marvelous book~\cite{RVNL}. It
is a marvelous book, numerous important and spectacular results are obtained,
but no spectral sequence in this book is made \emph{constructive}.

\section{Effective homology.}

\subsection{Notion of \emph{constructive} mathematics.}

Standard mathematics are based on Zermelo-Fraenkel (ZF) axiomatics. When
\emph{existence} results are involved, another axiomatics, the
\emph{constructive} logic, allows the user to express the results in a more
precise way; in this constructive context, one carefully distinguishes the
situation where some existence result should\ldots exist~(\(!\)) from the other
situation where a \emph{constructive} process is exhibited producing a copy of
the object the existence of which is stated.

The most common example allowing a novice to understand the difference is the
following. Question: does there exist two irrational real numbers \(\alpha\)
and \(\beta\) such that \(\alpha^\beta\) is rational? Let us inspect \(\gamma =
\sqrt{2}^{\sqrt{2}}\). If \(\gamma\) is rational, then \(\alpha = \beta =
\sqrt{2}\) is a solution. Otherwise \(\gamma\) is irrational, but then \(\alpha
= \gamma\) and \(\beta = \sqrt{2}\) is a solution, for
\((\sqrt{2}^{\sqrt{2}})^{\sqrt{2}} = 2\) is rational. This solution is correct
in ZF, but is not in constructive logic. The point is the following; in the
existence statement:
\[
(1)  \hspace{1cm} (\exists \alpha \in \bR-\bQ) (\exists \beta \in
\bR-\bQ)(\alpha^\beta \in \bQ)
\]
you \emph{did not} give a process allowing the user to \emph{construct} such a
pair \((\alpha, \beta)\). You have only produced two candidate solutions
\((\sqrt{2}, \sqrt{2})\) and \((\sqrt{2}^{\sqrt{2}}, \sqrt{2})\) and an
argument explaining that one of both candidate solutions must satisfy the
required property; but which one, this remains \emph{unknown}: you are not able
to produce \emph{one} genuine solution.

In constructive logic, if \(P\) is a predicate, \(\neg\neg P\)\footnote{\(\neg
=\) not.} is not equivalent to \(P\). In the above example, we have only
proved:
\[
(2) \hspace{1cm} \neg \neg (\exists \alpha \in \bR-\bQ) (\exists \beta \in
\bR-\bQ)(\alpha^\beta \in \bQ)
\]
Let us detail this point. The precise interpretation of \(\neg P\) is \(P
\Rightarrow \bot\) to be read: \(P\) implies a contradiction. Typically, \(\neg
(\sqrt{2} \in \bQ)\), because if some rational \(p/q\) is a square root of
\(2\), Euclid's analysis of the prime decompositions of \(p\) and \(q\)
generates a contradiction. Proving the double negation (2) consists in proving
the statement: \[ (3) \hspace{1cm} (\exists \alpha \in \bR-\bQ) (\exists \beta
\in \bR-\bQ)(\alpha^\beta \in \bQ) \Rightarrow \bot
\]
implies a contradiction, that is:
\[
(4) \hspace{1cm} ((\exists \alpha \in \bR-\bQ) (\exists \beta \in
\bR-\bQ)(\alpha^\beta \in \bQ) \Rightarrow \bot) \Rightarrow  \bot
\]

Let us assume this statement \((3)\). We \emph{then} prove firstly
\(\sqrt{2}^{\sqrt{2}} \in \bR - \bQ\). In fact applying \((3)\) to \(\alpha =
\beta = \sqrt{2}\) known irrational (Euclid), the hypothesis
\(\sqrt{2}^{\sqrt{2}} \in \bQ\) generates a contradiction, which is the
\emph{very definition} of \(\sqrt{2}^{\sqrt{2}} \in \bR - \bQ\). We can again
apply \((3)\) this time to \(\alpha = \sqrt{2}^{\sqrt{2}}\), now known \(\in
\bR - \bQ\), and \(\beta = \sqrt{2}\); the computation \(\alpha^\beta = 2\)
proves \(\alpha^\beta \in \bQ\), so that we have proved \((3)\) implies a
contradiction; in other words we have proved \((4)\), that is, (2).

On the contrary, our discussion is not a \emph{constructive} proof of \((1)\),
so that (1) and~(2) are not equivalent. Mathematicians usually think ``not-not
= yes'', but if the existence is constructively understood, you see \(\neg\neg
P\) is not necessarily equivalent to~\(P\). You see also the constructive
interpretation of \((2)\) gives a better interpretation of the ZF statement
(1): constructive mathematics is more precise and richer than ZF mathematics,
and mainly closer to the actual world. Consider these statements:

\begin{enumerate}
\item[(1)]
It is false there is no book about constructive analysis in this library.
\item[(2)]
The upper shelf to the left of the east window at the second floor of the
library has a book about constructive analysis\footnote{Maybe the famous book
by Bishop and Bridges, Springer-Verlag, 1985.}.
\end{enumerate}

Are these statements equivalent?

A constructive interpretation of existence quantifiers is an elegant way to
implicitly require as far as possible algorithms producing the objects whose
existence is stated. Sometimes it is possible, sometimes not; sometimes the
problem is open.

To be complete about our example around the \(\sqrt{2}\)'s, we must mention
that in fact a famous theorem of Gelfond and Schneider proves \(a^b\) is
transcendant as soon as \(a\) and \(b\) are algebraic, \(a \neq 0,1\) and \(b
\in \bR - \bQ\). So that \(\sqrt{2}^{\sqrt{2}}\) is transcendant and \((\alpha,
\beta) = (\sqrt{2}^{\sqrt{2}}, \sqrt{2})\) is this time, thanks to Gelfond and
Schneider, a constructive solution of our problem. But the proof is a long
story!

Another \emph{constructive} solution\footnote{Communicated by Thierry Coquand.}
is quite elementary; it can be obtained as follows: take \(\alpha = \sqrt{2}\)
and \(\beta = 2 \log_2 3\); Euclid knew \(\alpha \notin \bQ\) and if he had
known the definition of \(\log_2\), it would have been able to prove \(\beta
\notin \bQ\) as well. And \(\alpha^\beta = 3\).

\subsection{Existential quantifiers and homological algebra.}

Let \(C_\ast\) be a chain-complex and \(n\) be some integer. Let us study the
statement \(H_n(C_\ast) = 0\). By definition \(H_n(C_\ast) = Z_n(C_\ast) /
B_n(C_\ast)\), and \(H_n(C_\ast) = 0\) means any \(n\)-cycle is an
\(n\)-boundary. In a still more detailed way:
\[
(\forall c \in C_n)((dc = 0) \Rightarrow ((\exists c' \in C_{n+1}) (dc' = c)))
\]
And the critical question is the following: what about the exact status of the
existential quantifier?

In ordinary homological algebra, no constructiveness property is required for
this quantifier and, because constructing this preimage \(c'\) is most often a
little difficult, standard homological algebra is in a sense a catalog of
methods allowing you to prove some homology group is null without exhibiting an
algorithm constructing a boundary preimage for a cycle. For example if you can
insert your homology group in an exact sequence where both close groups are
null, then you know your group is null too, and in ordinary homological
algebra, this is enough.

But this habit has a severe drawback. For example we have explained how
Jean-Pierre Serre was unable to choose between \(\bZ_2 + \bZ_6\) and
\(\bZ_{12}\) when computing the group \(\pi_6{S^3}\). The homology groups
\(E^2_{6,0}\) and \(E^2_{0,6}\) of his spectral sequence did not give any
information about the nature, trivial or not, of the extension of \(\bZ_6\) by
\(\bZ_2\). We will give later an analysis of this difficulty: it comes from a
lack of representants of homology classes. When it is claimed \(H_2(C_\ast) =
\bZ_6\), it is in fact an unfortunate shorthand for: there \emph{exists} an
isomorphism\ \raisebox{3pt}{\xymatrix@1{H_2(C_\ast) \ar@<-1.5pt>[r]_-\phi &
\bZ_6 \ar@<-1.5pt>[l]_-\psi}}; but this claimed existence most often is not
constructive. To make it constructive, you must be able to
\emph{construct}~\(\psi\), in other words you must be able to \emph{construct}
the homology classes in front of the elements of \(\bZ_6\), for example by
exhibiting cycles \((z_i)_{0 \leq i \leq 5}\) representing them. It is not
finished, you must next construct~\(\phi\); let \(\fh \in H_2(C_\ast)\); most
often the homology class \(\fh\) is given through a cycle \(z\), and because
\(\psi\) is assumed available, defining \(\phi(\fh)\) amounts to identify which
\(z_i\) is homologous to \(z\). Let us assume in a particular case it is
\(z_5\); this means \((\exists c \in C_3)(dc = z - z_5)\), again an existential
quantifier.

We will explain how it is possible, and elementary, to systematically organize
homological algebra in a constructive style. It is not hard and very useful.
The fuzzy classical tools such as exact and spectral sequences will easily so
become \emph{algorithms} allowing you to compute wished homology and homotopy
groups. Of course you must remain lucid about the complexity of the algorithms
so obtained, but there is an interesting intermediate work level where these
algorithms will produce results otherwise unreachable.

\subsection{The homological problem for a chain-complex.}\label{76988}

We translate the constructiveness requirement roughly described in the previous
section into a definition. This definition, a little heavy but unavoidable, is
essentially \emph{temporary}. It will be soon replaced by the notion of
\emph{reduction}.

\begin{dfn}\label{73027} ---
\emph{Let \(\fR\) be a ground ring and \(C_\ast\) a chain-complex of
\(\fR\)-modules. A \emph{solution \(S\) of the homological problem for
\(C_\ast\)} is a set \(S = (\sigma_i)_{1 \leq i \leq 5}\) of five
\emph{algorithms}:
\begin{enumerate}
\item
\(\sigma_1: C_\ast \rightarrow \{\bot,\top\}\) (\(\bot\) = false, \(\top\) =
true) is a predicate deciding for every \(n \in \bZ\) and every \(n\)-chain \(c
\in C_n\) whether \(c\) is an \(n\)-cycle or not, in other words whether \(dc =
0\) or \(dc \neq 0\), whether \(c \in Z_n(C_\ast)\) or not.
\item
\(\sigma_2: \bZ \rightarrow \{\fR\textrm{-modules}\}\) associates to every
integer \(n\) some \(\fR\)-module \(\sigma_2(n)\) in principle isomorphic to
\(H_n(C_\ast)\). The image \(\sigma_2(n)\) will \emph{model} the isomorphism
class of \(H_n(C_\ast)\) in an effective way to be defined.
\item
The algorithm \(\sigma_3\) is indexed by \(n \in \bZ\); for every \(n \in
\bZ\), the algorithm \(\sigma_{3,n}: \sigma_2(n) \rightarrow Z_n(C_\ast)\)
associates to every \(n\)-homology class \(\fh\) coded as an element \(\fh \in
\sigma_2(n)\) a cycle \(\sigma_{3,n}(\fh) \in Z_n(C_\ast)\) \emph{representing}
this homology class.
\item
The algorithm \(\sigma_4\) is indexed by \(n \in \bZ\); for every \(n \in
\bZ\), the algorithm \(\sigma_{4,n}: C_n \supset Z_n(C_\ast) \rightarrow
\sigma_2(n)\) associates to every \(n\)-cycle \(z \in Z_n(C_\ast)\) the
homology class of \(z\) coded as an element of \(\sigma_2(n)\).
\item
The algorithm \(\sigma_5\) is indexed by \(n \in \bZ\); for every \(n \in
\bZ\), the algorithm \(\sigma_{5,n}: ZZ_n(C_\ast) \rightarrow C_{n+1}\)
associates to every \(n\)-cycle \(z \in Z_n(C_\ast)\) \emph{known as a
boundary} by the previous algorithm, a boundary preimage \mbox{\(c \in
C_{n+1}\)}: \mbox{\(dc = z\)}. In particular \(ZZ_n(C_\ast) := \ker
\sigma_{4,n}\)
\end{enumerate}}
\end{dfn}

Several complements are necessary to clarify this definition.

The computational context needs some method to \emph{code} on our theoretical
or concrete machine the chain-complex~\(C_\ast\) and the homology groups
\(H_n(C_\ast)\); and also their elements. We will see a \emph{locally
effective} representation of \(C_\ast\) will be enough; this subtle notion,
very important, in fact most often ordinarily underlying, is detailed in the
next section.

In most important cases, the set of interesting isomorphism classes of
\(\fR\)-modules is countable, and some simple process defines a relevant
isomorphism class as a finite machine object. If \(\fR\) is a principal ring,
\(\bZ\) for example, an \(\fR\)-module \emph{of finite type} \(H\) may be
described as a sequence \(H = (d_1, \ldots, d_r) \in \fR^r\) for some \(r\),
the pseudo-rank, the sequence~\(H\) satisfying the \emph{divisor condition}:
\(d_1\) divides \(d_2\), which divides \(d_3\) and so on up to~\(d_r\). For
example the \(\bZ\)-module \(\bZ^2 + \bZ_6 + \bZ_{15}\) would be represented as
the sequence \((3, 30, 0, 0)\). This representation is perfect: the
correspondance between isomorphism classes and representations is
bijective\footnote{In a different context, the presentation of a group by
finite sets of generators and relators is not perfect: no \emph{effective}
canonical presentation, because of the G\"odel-Novikov-Rabin theorem.}. An
element of such an \(\fR\)-module~\(H\) is then coded as a simple machine
object using the standard structured types.

As usual, an isomorphism class is defined through a \emph{representant} of this
class, but to make complete such a representation, an isomorphism must also
\emph{effectively} be given between the original group and the representant of
the isomorphism class: this is the role of \(\sigma_3\) and \(\sigma_4\). In
our context, \(\sigma_{3,n}\) describes the isomorphism from the
\emph{representant} \(\sigma_2(n)\) of the homology group to the \emph{genuine}
homology group \(H_n(C_\ast)\), an element of the last group being in turn
represented by a cycle. The algorithm \(\sigma_{4,n}\) is the reciprocal. Note
the map \(\sigma_{3,n}\) cannot be in general a module morphism. In the
chain-complex \(0 \leftarrow \bZ \stackrel{\times 2}{\longleftarrow} \bZ
\leftarrow 0\) null outside degrees 0 and~1, \(Z_0(C_\ast) = \bZ\) and
\(H_0(C_\ast) = \bZ_2\). The map \(\sigma_{4,0}\) is surjective and a morphism,
but the map \(\sigma_{3,0}\), a section of the previous one, cannot be a module
morphism. This unpleasant possible behavior will soon be avoided thanks to the
notion of \emph{reduction}.

Observe a homology class is represented in two different ways, and it is
important to understand the subtle difference. An ``actual'' \(n\)-homology
class is represented by a cycle \(z \in Z_n(C_\ast)\), while its image
\(\sigma_{4,n}(z)\) represents the same element in the model \(\sigma_2(n)\) of
the isomorphism class of the homology group.

The algorithm \(\sigma_5\) is in particular a \emph{certificate} for the
claimed properties of \(\sigma_3\) and \(\sigma_4\), but its role is not at all
limited to this authentication. We will see it is the main ingredient allowing
us to make constructive the usual exact and spectral sequences.

\subsection{Notion of locally effective object.}

When you use a simple pocket computer, this computer is able to compute for
example the sum of two integers \(a\) and \(b\) for a large set of integers
\(\bZ' \subset \bZ\). This situation is quite common, but not precisely enough
analyzed. We will describe this situation by a convenient terminology; we will
say the computer contains a \emph{locally effective} version of the standard
ring \(\bZ\).

The mathematical ring \(\bZ\) is a large set provided with a few operators. On
your computer, you can ask for \(2 + 3\) and the answer is 5. You \emph{enter}
(input) two particular elements of \(\bZ\) and another one has been computed,
the right terminology being: `5' has been \emph{returned} (output). Any
analogous computation can be done, at least if it is possible to enter the
arguments, when they are not too large. Note no \emph{global} description of
\(\bZ\) is given by your computer. But for \emph{arbitrary} integers \(a\) and
\(b\), the computer can effectively compute \(a + b\). We will use in such a
situation the following expression: the addition on \(\bZ\) is \emph{locally
effective}; this expression is a little inappropriate, no topology here to
justify the adverb ``locally'', but experience shows it is very convenient. In
a detailed way, we mean there is no \emph{global} implementation of the
addition; the possible \emph{global} properties of the addition, for example
associativity, commutativity, are unreachable by your computer, but this does
not prevent you from using it fruitfully. It is not frequent to need a global
property of the addition, most often we use only ``local'', more precisely
\emph{elementwise}, properties. For example for the specific elements 2 and 3,
the sum is~5.

For two \emph{arbitrary} elements \(a, b \in \bZ\), the computer can compute
\(a + b\); really arbitrary? Not exactly. Not many computers could accept for
example \(a = 10^{10 ^{10}}\). The user of such a concrete locally effective
implementation of \(\bZ\) usually knows he must be sensible about input size.
The specific problem met here most often is a problem of memory size, or
technical bounds. In computational algebra systems allowing you to handle the
so-called \emph{extended} integers, with a claimed arbitrary number of digits,
you are yet limited by the memory size of your machine. For a specific
computation you could after all buy more memory to succeed\footnote{But memory
extensions, except for Turing machines, have their own technical bounds!}. On
pocket computers, technical limits are most often given, maybe you are limited
to integers with less than ten decimal digits. From a theoretical mathematical
point of view, these constraints are most often neglected, without any serious
drawback, at least in a first step. The underlying implicit statement is: when
you will use \emph{concrete} implementations of locally effective objects, be
careful, you can meet memory limitations, otherwise the results will be
correct. And this is enough in a first step. Of course, time and space
complexity is an important subject, theoretically as much as practically, but
we decide it is another subject, which of course will be quickly present in
concrete calculations.

Another point is to be considered. If you try to enter the ``arbitrary
integer'' \boxtt{234hello567}, we hope your computer or computational algebra
system complains! Another formalism is here necessary. The universe \(\cU\) is
the set of all the \emph{objects} that can be handled\footnote{Without taking
account of size limitations! We will not make this precision anymore.} by a
machine. The set of ``legal'' integers is a small subset of it; the computer
scientists use the notion of \emph{type} to formalize this point. Specific
machine objects, more precisely specific machine \emph{predicates}, can be used
to verify whether an object is an integer or not. Which allows the machine or
the program, when it is safely organized, to detect an incoherent input.
Situations are quite different according to concrete implementations. The
simplest pocket computers do not have alphabetic keys. More sophisticated ones
have and almost always detect our incorrect integer. If you use an intermediary
programming language, according to the language, \boxtt{234hello567} is an
object or not: in Lisp yes, in C not. In Lisp this object is accepted but it is
a symbol, which cannot be a legal argument for addition, a type error is in
principle detected. In C this character string does not denote any machine
object and the compiler or more rarely the interpreter will detect an
incoherent input, most often being unable to guess what your intention could
be.

These technical but unavoidable considerations will be formalized here by
characteristic functions. A locally effective object will contain a
\emph{membership} predicate, that is an algorithm \(\chi: \cU \rightarrow
\{\bot, \top\}\) allowing the user or the program, if necessary, to verify the
object it must process has the right type\footnote{Such a characteristic
function is a universal predicate, and an interesting question is to construct
the type of the universal predicates, in other words the type of types! If you
study a little more this matter, you will quickly rediscover G\"odel's
incompleteness theorem.}, that is, actually is a member of the underlying set.

Another more subtle predicate must also be used. On most simple pocket
computers, instead of keying \boxtt{2}, you could enter as well
\boxtt{0.002E3}, two \emph{different} notations are possible, because \(2 =
0.002 \times 10^3\). And this is a permanent problem when implementing
mathematical objects: \emph{different} machine objects can \emph{code} the
\emph{same} mathematical object. Sometimes it is an extremely technical point:
the integer object \boxtt{2} somewhere in the machine is or is not ``equal'' to
another object again \boxtt{2} but somewhere else in the machine\footnote{Only
Common Lisp correctly handles this matter, see the Lisp functions \texttt{eq}
and \texttt{eql}.}. Sometimes, such a decision depends on the technical choice
of the user: if you have to implement \(\bZ_5 := \bZ/5\bZ\), you can decide to
implement an object as an integer \boxtt{0} or \boxtt{1} or \boxtt{2} or
\boxtt{3} or \boxtt{4}, why not; but sometimes it is much better to decide to
represent an element of \(\bZ_5\) by an arbitrary machine integer, taking care
that in fact \boxtt{12} and \boxtt{17} represent the same element of \(\bZ_5\).

We will not give more details about this notion of locally effective object.
The numerous examples studied in this text are sufficient illustrations.

\subsection{Notion of effective object.}

On the contrary, we must sometimes be able to ``know everything'' about an
object, including the \emph{global} properties. For example if you intend to
compute some homology group \(H_n(C_\ast)\) of the chain-complex \(C_\ast\),
you must know the \emph{global} nature of \(C_k\) for \(k = n-1, n, n+1\), and
you must know also the differentials \(d_{k}\) and \(d_{k+1}\) in such a way
you can compute \(\ker d_k\), \(\im d_{k+1}\) and finally the looked-for
homology group.

If the chain-complex is only \emph{locally effective}, these calculations in
general are not possible, you must have more information about your
chain-complex. We will say a chain-complex is \emph{effective} when every chain
group \(C_n\) is of finite type. Then the a priori locally effective
implementation of a boundary operator \(d_n\) becomes effective and the
homology group can be computed. Instead of painful abstract definitions, we
prefer to illustrate this point by a typical Kenzo example.

Let us assume we are interested by \(H_7(K(\bZ, 3))\). The Eilenberg-MacLane
space \(K(\bZ,3)\) has the following characteristic property: its homotopy
groups are null except \(\pi_3 K(\bZ, 3) = \bZ\). The Kenzo program can
construct it:

 \bmp
 \bmpi\verb|> (setf KZ3 (k-z 3))|\empim
 \bmpi\verb|[K11 Abelian-Simplicial-Group]|\empi
 \emp

The simplicial set \(K(\bZ, 3)\) is locally effective, and in principle it is
not possible to deduce from its implementation its homology groups. But the
Kenzo program is intelligent enough to use the \emph{definition} of \(K(\bZ,
3)\) to undertake sophisticated computations giving the result. Look at the
(Kenzo) \emph{definition} of this object.

 \bmp
 \bmpi\verb|> (dfnt KZ3)|\empim
 \bmpi\verb|(CLASSIFYING-SPACE [K6 Abelian-Simplicial-Group])|\empi
 \emp

It is the \emph{classifying space} of another simplicial group. Using this
definition and others, Kenzo can compute the homology group.

 \bmp
 \bmpi\verb|> (homology KZ3 7)|\empim
 \bmpi\verb|Homology in dimension 7 :|\empi
 \bmpi\verb|Component Z/3Z|\empi
 \bmpi\verb|---done---|\empi
 \emp

But let us play now to \emph{hide} the definition. We reinitialize --
\boxtt{cat-init} -- the environment, otherwise it would not be sufficient.

 \bmp
 \bmpi\verb|> (cat-init)|\empim
 \bmpi\verb|---done---|\empix
 \bmpi\verb|> (setf KZ3 (k-z 3))|\empim
 \bmpi\verb|[K11 Abelian-Simplicial-Group]|\empix
 \bmpi\verb|> (setf (slot-value KZ3 'dfnt) '(hidden-definition))|\empim
 \bmpi\verb|(HIDDEN-DEFINITION)|\empix
 \bmpi\verb|> (homology KZ3 7)|\empim
 \bmpi\verb|Error: I don't know how to determine the effective homology of: [K11|\empi
 \bmpi\verb|Abelian-Simplicial-Group] (Origin: (HIDDEN-DEFINITION)).|\empi
 \emp

This is due to the fact that the chain-complex associated with our \(K(\bZ,
3)\) is only \emph{locally effective}: no \emph{global} information is
reachable:

 \bmp
 \bmpi\verb|> (basis KZ3 7)|\empim
 \bmpi\verb|Error: The object [K11 Abelian-Simplicial-Group] is locally-effective.|\empi
 \emp

\noindent and in fact the basis is infinite. Let us reinstall the right
definition:

 \bmp
 \bmpi\verb|> (setf (slot-value KZ3 'dfnt) `(classifying-space ,(k 6)))|\empim
 \bmpi\verb|(CLASSIFYING-SPACE [K6 Abelian-Simplicial-Group])|\empi
 \emp

\noindent The basis of the chain-complex is still unreachable:

 \bmp
 \bmpi\verb|> (basis KZ3 7)|\empim
 \bmpi\verb|Error: The object [K11 Abelian-Simplicial-Group] is locally-effective.|\empi
 \emp

\noindent but the homology group is computable:

 \bmp
 \bmpi\verb|> (homology KZ3 7)|\empim
 \bmpi\verb|Homology in dimension 7 :|\empi
 \bmpi\verb|Component Z/3Z|\empi
 \bmpi\verb|---done---|\empi
 \emp

How this is possible? It is here the \emph{heart} of our subject. Because of
the correct definition, Kenzo is able to construct the \emph{effective
homology} of \(K(\bZ, 3)\). Taking account of \boxtt{efhm} =
\underline{Ef}fective \underline{H}o\underline{m}ology:

 \bmp
 \bmpi\verb|> (efhm KZ3)|\empim
 \bmpi\verb|[K265 Equivalence K11 <= K255 => K251]|\empi
 \emp

This homology equivalence is the key point, it is an equivalence between the
\emph{locally effective} chain-complex \boxtt{K11} \(= C_\ast(K(\bZ, 3))\) and
the \emph{effective} chain-complex \boxtt{K251} which cannot be detailed at
this point.

 \bmp
 \bmpi\verb|> (basis (K 11) 7)|\empim
 \bmpi\verb|Error: The object [K11 Abelian-Simplicial-Group] is locally-effective.|\empix
 \bmpi\verb|> (basis (K 251) 7)|\empim
 \bmpi\verb|(<<Abar[7 <<Abar[2 S1][2 S1][2 S1]>>]>>)|\empi
 \emp

In fact there is only one generator in \(C_7(\boxtt{K251})\), which does not
prevent the chain-complex \(\boxtt{K251}\) from being homology equivalent to
\(\boxtt{K11}\), the \(C_7\) of which being on the contrary not at all of
finite type. And Kenzo, knowing this equivalence, computes in fact the homology
group of \(\boxtt{K251}\) when \(H_7(K(\bZ, 3))\) is asked for.

The \emph{effective homology theory} is essentially a systematic method
combining locally effective chain-complexes with effective chain-complexes
through homology equivalences. A locally effective chain-complex is too
``vague'' to allow us to compute its homology groups, but it is so possible to
implement infinite objects such as our Eilenberg-MacLane space \(K(\bZ, 3)\).
The effective chain-complexes are objects where homology groups can be
elementary computed, but only simple objects of finite type can be so
implemented. Homology equivalences will allow us to settle \emph{bridges}
between both notions, making homological algebra \emph{effective}.

\subsection{Reductions.}

Definition \ref{73027} is relatively complex and the notion of \emph{reduction}
is an interesting intermediate organization allowing the topologist to work on
the contrary in a convenient environment, from a traditional mathematical point
of view and also when computer implementations are planned.

\begin{dfn}\label{58661}
--- \emph{A \emph{reduction} \(\rho: \hC_\ast \rrdc C_\ast\) is a diagram:
\[
 \rho = \raisebox{4pt}{\framebox{\xymatrix@1{ \scriptstyle h\ \textstyle
 \autoarrow{0.8}\ \hC_\ast \ar@<-10\ul>[r]_-f & C_\ast \ar@<-10\ul>[l]_-g }}}
\]
where:
\begin{enumerate}
\item
\(\hC_\ast\) and \(C_\ast\) are chain-complexes.
\item
\(f\) and \(g\) are chain-complex morphisms.
\item
\(h\) is a homotopy operator (degree +1).
\item These relations are satisfied:
\begin{enumerate}
\item \(fg = \id{C_\ast}\).
\item \(gf + dh + hd = \id{\hC_\ast}\).
\item \(fh = hg = hh = 0\).
\end{enumerate}
\end{enumerate}
}
\end{dfn}

A reduction is a particular homology equivalence between a \emph{big}
chain-complex \(\hC_\ast\) and a \emph{small} one \(C_\ast\). This point is
detailed in the next proposition.

\begin{prp}
--- Let \(\rho: \hC_\ast \rrdc C_\ast\) be a reduction. This reduction is
equivalent to a decomposition: \(\hC_\ast = A_\ast \oplus B_\ast \oplus
C'_\ast\):
\begin{enumerate}
\item
\(\hC_\ast \supset C'_\ast = \emph{\im} g\) is a subcomplex of \(\hC_\ast\).
\item
\(A_\ast \oplus B_\ast = \ker f\) is a subcomplex of \(\hC_\ast\).
\item
\(\hC_\ast \supset A_\ast = \ker f \cap \ker h\) is not in general a subcomplex
of \(\hC_\ast\).
\item
\(\hC_\ast \supset B_\ast = \ker f \cap \ker d\) is a subcomplex of
\(\hC_\ast\) with null differentials.
\item
The chain-complex morphisms \(f\) and \(g\) are inverse isomorphisms between
\(C'_\ast\) and \(C_\ast\).
\item
The arrows \(d\) and \(h\) are module isomorphisms of respective degrees -1
and~+1 between \(A_\ast\) and \(B_\ast\).
\end{enumerate}
\end{prp}

In other words a reduction is a compact and convenient form of the following
diagram.
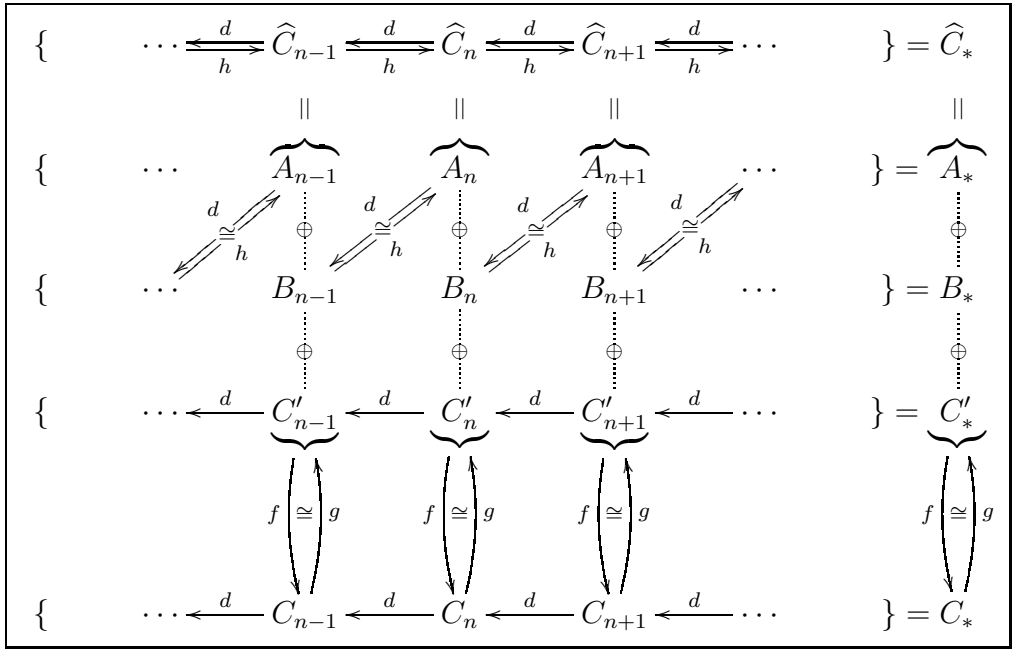
\[
 \fbox{\newcommand{\diagiso}{\ar@<-5\ul>[ld]
 |-{\phantom{\cong}}_-{\raisebox{10\ul}{\(\scriptstyle d\)}}
 \ar@{<-}@<5\ul>[ld]|-{\phantom{\cong}}^-h
 \ar@{}[ld]|-{\cong}}
 \xymatrix@C100\ul{ \{ & \cdots
 \ar@<-5\ul>[r]_-h \ar@{<-}@<5\ul>[r]^-d & \hC_{n-1}
 \ar@{}[d]|{||} \ar@<-5\ul>[r]_-h \ar@{<-}@<5\ul>[r]^-d &
 \hC_n \ar@{}[d]|{||} \ar@<-5\ul>[r]_-h
 \ar@{<-}@<5\ul>[r]^-d & \hC_{n+1} \ar@{}[d]|{||}
 \ar@<-5\ul>[r]_-h \ar@{<-}@<5\ul>[r]^-d & \cdots & \} =
 \hC_\ast \ar@{}@<40\ul>[d]|{||}
 \\
 \{ & \cdots & {\overbrace{A_{n-1}}}
 \ar@{.}[d]|{\oplus} \diagiso & {\overbrace{A_n}}
 \ar@{.}[d]|{\oplus} \diagiso & {\overbrace{A_{n+1}}}
 \ar@{.}[d]|{\oplus} \diagiso & \cdots \diagiso & \} =
 \overbrace{A_\ast} \ar@{.}@<40\ul>[d]|{\oplus}
 \\
 \{ & \cdots & B_{n-1} \ar@{.}[d]|{\oplus} &
 B_n \ar@{.}[d]|{\oplus} & B_{n+1}
 \ar@{.}[d]|{\oplus} & \cdots & \} = B_\ast
 \ar@{.}@<40\ul>[d]|{\oplus}
 \\
 \{ & \cdots & {\underbrace{C'_{n-1}}} \ar@/_/[dd]_f \ar[l]_-d
 \ar@{}[dd]|{\cong} & {\underbrace{C'_n}} \ar[l]_-d \ar@/_/[dd]_f
 \ar@{}[dd]|{\cong} & {\underbrace{C'_{n+1}}} \ar[l]_-d \ar@/_/[dd]_f
 \ar@{}[dd]|{\cong} & \cdots \ar[l]_-d & \} =
 \underbrace{C'_\ast} \ar@/_/@<40\ul>[dd]_f \ar@{}@<40\ul>[dd]|{\cong}
 \\\\
 \{ & \cdots & C_{n-1} \ar[l]_-d \ar@/_/[uu]_g &
 C_n \ar[l]_-d \ar@/_/[uu]_g & C_{n+1} \ar[l]_-d
 \ar@/_/[uu]_g & \cdots \ar[l]_-d & \} = C_\ast
 \ar@<-40\ul>@/_/[uu]_g}}
\]

It is a simple exercise of elementary linear algebra to prove the equivalence
between the above diagram and the initial reduction. Every chain group
\(\hC_n\) is then decomposed into three components, \(A_n\) made of chains in
canonical bijection with \(B_{n-1}\) thanks to \(d\) and \(h\). We can consider
\(A_n\) is a collection of \(n\)-chains ready to explain the elements of
\(B_{n-1}\) are not only cycles, but also boundaries. \(B_n\) is a collection
of cycles known as boundaries, because of the bijection between \(A_{n+1}\) and
\(B_n\) again through \(d\) and \(h\). Finally the component \(C'_n\) is a copy
of \(C_n\) and their homological natures therefore are the same.

A reduction \(\rho: \hC_\ast \rrdc C_\ast\) is a decomposition \(\hC_\ast =
\ker f \oplus\,C'_\ast\) in two components; no specific information about the
second one other than \(C'_\ast \cong C_\ast\); but the first one \(\ker f\) is
\emph{acyclic}, for the restriction of the relation \(\id{\hC_\ast} = gf + dh +
hd\) to \(\ker f\) is simply \(\id{\ker f} = dh + hd\); note in particular
\(\ker f\) is a subcomplex of \(\hC_\ast\): \(f\) is a chain-complex morphism,
that is, \(df = fd\), which implies \(d(\ker f) \subset \ker f\). Note also
\(h(\hC_\ast) \subset \ker f\), a consequence of \(fh = 0.\) The component
\(\ker f\), known as acyclic, is in turn decomposed in two components, \(\ker f
= A_\ast + B_\ast\) with \(A_\ast = \ker f \cap \ker h\) and \(B_\ast = \ker f
\cap \ker d\). This can be considered as a \emph{Hodge decomposition} of
\(\hC_\ast\), describing in a detailed way why the homology groups of
\(\hC_\ast\) and \(C_\ast\) are canonically isomorphic.

\begin{thr}---
Let \(\rho = (f,g,h): \hC_\ast \rrdc C_\ast\) be a reduction where the
chain-complexes \(\hC_\ast\) and \(C_\ast\) are locally effective. If the
homological problem is solved in the small chain-complex \(C_\ast\), then the
reduction \(\rho\) induces a solution of the homological problem for the big
chain-complex \(\hC_\ast\).
\end{thr}

\proof Let us examine the criteria of Definition~\ref{73027}.

1. Let \(c \in \hC_\ast\); the chain-complex \(\hC_\ast\) is locally effective
and the ``local'' calculation \(dc\) can be achieved, which allows you to
determine whether the chain \(c\) satisfies \(dc = 0\) or not, whether \(c\) is
a cycle or not.

2. The known relations \(\id{C_\ast} = fg\) and \(\id{\hC_\ast} = gf + dh +
hd\) imply \(f\) and \(g\) are inverse homology equivalences. The homology
groups \(H_n(\hC_\ast)\) and \(H_n(C_\ast)\) are \emph{canonically} isomorphic.
Let \(\sigma_\ast\) be the algorithms provided by the solution of the
homological problem for \(C_\ast\) and let us call \(\widehat{\sigma}_\ast\)
the algorithms to be constructed for~\(\hC_\ast\). We can choose in particular
\(\widehat{\sigma}_{2,n} = \sigma_{2,n}\), the last \emph{equality} being a
\emph{genuine} one.

3. The chain morphism \(f\) induces an isomorphism between \(H_n(\hC_\ast)\)
and \(H_n(C_\ast)\). This allows us to choose \(\widehat{\sigma}_{3,n}(z) :=
\sigma_{3,n}(f(z))\).

4. In the same way, choose \(\widehat{\sigma}_{4,n}(\fh) :=
g(\sigma_{4,n}(\fh))\).

5. Finally, if \(z \in \hC_n\) is a cycle known homologous to zero, a boundary
preimage is \(\widehat{\sigma}_{5,n}(z) := h(z) + g(\sigma_{5,n}(f(z)))\). In
fact: \(d(hz + g(\sigma_{5,n}(f(z)))) = dhz + gd\sigma_{5,n}(f(z)) = dhz + gfz
= z - hdz = z\), for \(g\) is a chain-complex morphism, \(\sigma_{5,n}\) finds
boundary preimages, and \(z\) is a cycle.\QED

\begin{crl}---
If \(\rho = (f,g,h): \hC_\ast \rightarrow C_\ast\) is a reduction where
\(\hC_\ast\) is \emph{locally effective} and \(C_\ast\) is \emph{effective},
then this reduction produces a solution of the homological problem for
\(\hC_\ast\).
\end{crl}

\proof The small chain-complex \(C_\ast\) is effective and a solution of the
homological problem for \(C_\ast\) therefore is elementary.\QED

\begin{prp}---
Let \(\rho = (f,g,h): \hC_\ast \rightarrow C_\ast\) be a reduction, where the
homological problem is solved for \(\hC_\ast\). Then the homological problem is
also solved for the small chain complex \(C_\ast\).
\end{prp}

\proof The small chain complex being a sub-chain-complex of the big one, the
situation is more comfortable. The only point deserving a little attention is
the search of a boundary preimage for a \(C_\ast\)-cycle known being a
boundary: exercise. \QED

\begin{crl}---
If \(\varepsilon: C_\ast \eqvl C'_\ast\) is an equivalence between two
chain-complexes, a solution of the homological problem for \(C'_\ast\) gives a
solution of the same problem for \(C_\ast\). In particular, if \(C'_\ast\) is
effective, the homological problem is solved for \(C_\ast\).
\end{crl}

The reader probably wonders why, in presence of such a reduction \(\rho:
\hC_\ast \rrdc C_\ast\), the user continues to give some interest to
\(\hC_\ast\). The big chain-complex \(\hC_\ast\) is the direct sum of the small
one \(C_\ast\) and \(\ker f\), the last component not playing any role from a
homological point of view. The point is the following: frequently we have to
work with chain-complexes which carry more structure than a chain-complex
structure. For example if the chain-complex comes from a simplicial set or
complex, there is another structure, the simplicial structure which is present,
and the chain-complex structure in this case is \emph{underlying}; and it is
frequent the chain-complex structure can be \emph{reduced} but the simplicial
structure \emph{not}. So that you must continue to play with the big
chain-complex \(\hC_\ast\) and its further simplicial structure, but when the
subject is homology, you can transfer the work to the small chain-complex
\(C_\ast\). And the planned work is always of this sort: playing simultaneously
with big objects provided with sophisticated structures, most often not
significantly \emph{reducible}, and their small homological reductions.

\subsection{Kenzo example.}

We want to \emph{concretely} illustrate how reductions between locally
effective and effective chain-complexes allow a user to obtain and use the
corresponding solution of a homological problem.

The mathematical underlying theory will be explained later in
Section~\ref{64720} and we use here Example~\ref{25180} of this section. We
consider the polynomial ring \(\fR = \bQ[t,x,y,z]_0\) and in this ring the
ideal:
\[
I = \ \ideal{t^5 - x, t^3 y - x^2, t^2 y^2 - x z, t^3 z - y^2, t^2 x - y, t x^2
- z, x^3 - t y^2, y^3 - x^2 z, x y - t z}.
\]
It happens the homology of the \emph{Koszul complex} \(\Ksz(\fR/I)\) reflects
deep properties of the ideal \(I\). The Koszul complex is a \(\bQ\)-vector
space of infinite dimension, but yet an algorithm can compute its
\emph{effective} homology. Kenzo constructs the ideal as a list of generators,
each generator being a combination (\boxtt{cmbn}) of monomials, each monomial
being a list of exponents, for example \boxtt{(3 0 1 0)} codes \(t^3y\).

 \bmp
 \bmpi\verb|> (setf ideal|\empi
 \bmpi\verb|    (list|\empi
 \bmpi\verb|     (cmbn 0 1 '(5 0 0 0) -1 '(0 1 0 0))|\empi
 \bmpi\verb|     (cmbn 0 1 '(3 0 1 0) -1 '(0 2 0 0))|\empi
 \bmpi\verb|[... 6 lines deleted ...]|\empi
 \bmpi\verb|     (cmbn 0 1 '(0 1 1 0) -1 '(1 0 0 1))))|\empimx
 \bmpi\verb|(|\empi
 \bmpi\verb|----------------------------------------------------------------------{cmbn 0}|\empi
 \bmpi\verb|<1 * (5 0 0 0)> <-1 * (0 1 0 0)>|\empi
 \bmpi\verb|------------------------------------------------------------------------------|\empi
 \bmpi\verb|[... other lines deleted ...]|\empi
 \bmpi\verb|)|\empi
 \emp

The display is simply the list of generators, only the first one is given here.
The Koszul complex \(\Ksz(\fR/I)\) is then constructed.

 \bmp
 \bmpi\verb|> (setf ksz (k-complex/gi 4 ideal))|\empim
 \bmpi\verb|[K5 Chain-Complex]|\empi
 \emp

Kenzo returns \boxtt{K5}, the Kenzo object \boxtt{\#5}, a chain-complex. The
ideal in fact is as well generated by the toric generators \(x-t^5\),
\(y-t^7\), \(z-t^{11}\); we will see how the \emph{effective} homology of the
Koszul complex can \emph{discover} this fact. Three generators and four
variables, the quotient is certainly of infinite \(\bQ\)-dimension. If we ask
for the \(\bQ\)-basis of the Koszul complex in degree 2 for example, an error
is returned.

 \bmp
 \bmpi\verb|> (basis ksz 2)|\empim
 \bmpi\verb|Error: The object [K5 Chain-Complex] is locally-effective.|\empi
 \emp

Several procedures in Kenzo can compute the effective homology of \boxtt{K5}.
In particular the procedure \boxtt{koszul-min-rdct} computes the \emph{minimal}
effective homology as a reduction.

 \bmp
 \bmpi\verb|> (setf mrdct (koszul-min-rdct ideal "H"))|\empim
 \bmpi\verb|[K778 Reduction K5 => K763]|\empi
 \emp

The reduction is assigned to the symbol \boxtt{mrdct}, a reduction of the
chain-complex~\boxtt{K5} over the chain-complex \boxtt{K763}. You observe
several hundreds of Kenzo objects, chain-complexes, morphisms, reductions,
equivalences, \ldots, have been necessary to obtain the result, but this work
of automatic writing of programs is very fast, less than half a second for our
modest laptop. The small chain-complex \boxtt{K763} is effective. The Lisp
statement \boxtt{(mapcar ...)} gives the list of \(\bQ\)-dimensions from 0 to
4.

 \bmp
 \bmpi\verb|> (mapcar|\empi
 \bmpi\verb|     #'(lambda (i) (length (basis (k 763) i)))|\empi
 \bmpi\verb|     '(0 1 2 3 4))|\empim
 \bmpi\verb|(1 3 3 1 0)|\empi
 \emp

Let us look for the first generator in degree 2 and compute its differential.

 \bmp
 \bmpi\verb|> (first (basis (k 763) 2))|\empim
 \bmpi\verb|H-2-1|\empix
 \bmpi\verb|> (? (k 763) 2 *)|\empim
 \bmpi\verb|----------------------------------------------------------------------{CMBN 1}|\empi
 \bmpi\verb|------------------------------------------------------------------------------|\empi
 \emp

The generator is the symbol \boxtt{H-2-1} and its differential is null. The
esoteric Lisp statement ``\boxtt{(? (k 763) 2 *)}'' is to be understood as
follows: as already observed, ``\boxtt{(k 763)}'' returns the Kenzo object
\boxtt{K763}, a chain complex. The functional operator `\boxtt{?}' makes the
differential of this chain complex work in this case on a generator of degree
2, namely `\boxtt{*}', that is, the last object returned by the Lisp
interpreter, the symbol \boxtt{H-2-1}.

In fact the same behaviour can be observed for the eight basis elements: the
differential is the null-morphism of degree -1. This property is characteristic
of the \emph{minimal} effective homology of our Koszul complex. So that the
elements of the list \boxtt{(1 3 3 1 0)} are the \emph{Betti} numbers of the
Koszul complex. The first \boxtt{3} informs us for example the minimal number
of generators for our ideal is 3, while the ideal was defined with 9
generators.

The chain-complex \boxtt{K763} is nothing but a model for ``the''
\emph{homology} of our Koszul complex \boxtt{K5}. The homology class
\boxtt{h-2-1} is represented by the cycle \(g(\boxtt{h-2-1})\) if \(g\) is the
\(g\)-component of the reduction \boxtt{K778} = \boxtt{mrdct} = \((f,g,h)\).

 \bmp
 \bmpi\verb|> (g mrdct 2 'h-2-1)|\empimx
 \bmpi\verb|----------------------------------------------------------------------{CMBN 2}|\empi
 \bmpi\verb|<-1 * ((0 2 0 0) (1 1 0 0))>|\empi
 \bmpi\verb|<1 * ((4 0 0 0) (1 0 0 1))>|\empi
 \bmpi\verb|<-1 * ((0 0 0 0) (0 1 0 1))>|\empi
 \bmpi\verb|------------------------------------------------------------------------------|\empi
 \emp

\noindent which cycle would be denoted by \(- x^2 \, dt.dx + t^4 \, dt.dz -
dx.dz\) in the standard notation explained Section~\ref{66039}. You see not
only the homology groups are computed, but representants of homology classes
can be exhibited.

Let us play now with cycles and boundary preimages. If we take a random element
of the Koszul complex, in general it is not a cycle.

 \bmp
 \bmpi\verb|> (? ksz 2 '((2 0 0 0) (1 1 0 0)))|\empi
 \bmpi\verb|----------------------------------------------------------------------{CMBN 1}|\empi
 \bmpi\verb|<-1 * ((0 0 1 0) (1 0 0 0))>|\empi
 \bmpi\verb|<1 * ((3 0 0 0) (0 1 0 0))>|\empi
 \bmpi\verb|------------------------------------------------------------------------------|\empi
 \emp

\noindent The differential of \(t^2 \, dt.dx\) is not null; this object is not
a cycle. Now the demonstrator goes for a moment into the wings of his theater
and comes back with the object z1. Is it a cycle?

 \bmp
 \bmpi\verb|> (setf z1|\empi
 \bmpi\verb|    (cmbn 2|\empi
 \bmpi\verb|          1 '((1 0 1 9) (1 1 0 0))|\empi
 \bmpi\verb|         -1 '((0 2 0 0) (1 1 0 0))|\empi
 \bmpi\verb|         -1 '((1 1 0 9) (1 0 1 0))|\empi
 \bmpi\verb|          1 '((4 0 0 0) (1 0 0 1))|\empix
 \bmpi\verb|          1 '((2 0 0 9) (0 1 1 0))|\empi
 \bmpi\verb|         -2 '((1 1 0 0) (0 1 1 0))|\empi
 \bmpi\verb|          2 '((2 0 0 0) (0 1 0 1))|\empi
 \bmpi\verb|         -1 '((0 0 0 0) (0 1 0 1))|\empi
 \bmpi\verb|         -2 '((0 0 0 0) (0 0 1 1))))|\empimx
 \bmpi\verb|----------------------------------------------------------------------{CMBN 2}|\empi
 \bmpi\verb|<1 * ((1 0 1 9) (1 1 0 0))>|\empi
 \bmpi\verb|[... Lines deleted  ...]|\empi
 \bmpi\verb|------------------------------------------------------------------------------|\empix
 \bmpi\verb|> (? ksz z1)|\empim
 \bmpi\verb|----------------------------------------------------------------------{CMBN 1}|\empi
 \bmpi\verb|------------------------------------------------------------------------------|\empi
 \emp

The combination \((tyz^9 - x^2) \, dt.dx - txz^9 \, dt.dy + t^4 \, dt.dz +
(t^2z^9 - 2 tx) \, dx.dy + (2 t^2 - 1) \, dx.dz - 2 \, dy.dz\) is a cycle of
degree 2. What about its homology class?

 \bmp
 \bmpi\verb|> (f mrdct z1)|\empim
 \bmpi\verb|----------------------------------------------------------------------{CMBN 2}|\empi
 \bmpi\verb|<1 * H-2-1>|\empi
 \bmpi\verb|<-2 * H-2-3>|\empi
 \bmpi\verb|------------------------------------------------------------------------------|\empi
 \emp

We obtain the homology class by applying the \(f\)-component of the reduction
to the cycle; the homology class is \(\boxtt{h-2-1} - 2\ \boxtt{h-2-3}\). The
demonstrator again goes into the wings and comes back with another cycle
\boxtt{z2}.

 \bmp
 \bmpi\verb|> (setf z2|\empi
 \bmpi\verb|    (cmbn 2|\empi
 \bmpi\verb|          1 '((1 0 1 9) (1 1 0 0))|\empi
 \bmpi\verb|         -1 '((1 1 0 9) (1 0 1 0))|\empi
 \bmpi\verb|          1 '((2 0 0 9) (0 1 1 0))))|\empim
 \bmpi\verb|----------------------------------------------------------------------{CMBN 2}|\empi
 \bmpi\verb|<1 * ((1 0 1 9) (1 1 0 0))>|\empi
 \bmpi\verb|[... 2 lines deleted ...]|\empi
 \bmpi\verb|------------------------------------------------------------------------------|\empix
 \bmpi\verb|> (? ksz z2)|\empim
 \bmpi\verb||\empi
 \bmpi\verb|----------------------------------------------------------------------{CMBN 1}|\empi
 \bmpi\verb|------------------------------------------------------------------------------|\empix
 \bmpi\verb|> (f mrdct z2)|\empim
 \bmpi\verb|----------------------------------------------------------------------{CMBN 2}|\empi
 \bmpi\verb|------------------------------------------------------------------------------|\empi
 \emp

This time the cycle is \(tyz^9 \, dt.dx - txz^9 \, dt.dy + t^2 z^9 \, dx.dy\),
but its homology class is null. To obtain a boundary preimage, because the
homology is minimal, it is sufficient to apply the \(h\)-component of the
reduction.

 \bmp
 \bmpi\verb|> (h mrdct z2)|\empi
 \bmpi\verb|----------------------------------------------------------------------{CMBN 3}|\empi
 \bmpi\verb|<1 * ((1 0 1 8) (1 1 0 1))>|\empi
 \bmpi\verb|<-1 * ((1 1 0 8) (1 0 1 1))>|\empi
 \bmpi\verb|<1 * ((2 0 0 8) (0 1 1 1))>|\empi
 \bmpi\verb|------------------------------------------------------------------------------|\empi
 \emp

The claimed preimage is \(tyz^8 \, dt.dx.dz - txz^8 dt.dy.dz + t^2z^8 \,
dx.dy.dz\). To verify this claim, we compute the difference between the
original \boxtt{z2} and the boundary of the preimage.

 \bmp
 \bmpi\verb|> (2cmbn-sbtr (cmpr ksz) z2 (? ksz *))|\empim
 \bmpi\verb|----------------------------------------------------------------------{CMBN 2}|\empi
 \bmpi\verb|------------------------------------------------------------------------------|\empi
 \emp

 A comparison operator between generators is necessary to compute such a
 difference, it is the reason why the first argument is the comparison operator
 (\boxtt{cmpr}) of the Koszul complex (\boxtt{ksz}). The result is null, OK!

 These small computations illustrate how any homological question in the Koszul
 complex is \emph{effectively} solved, thanks to the reduction \boxtt{mrdct}.
 Even if the chain-complex is not of finite \(\bQ\)-type.
 There remains to understand how it is possible to construct the critical
 reduction, more generally the necessary equivalence.

\subsection{Homological Perturbation theory.}

\subsubsection{Presentation.}

The most important tool allowing us to efficiently work with reductions is the
so-called \emph{basic perturbation lemma}, a ``lemma'' which would be better
called the \emph{fundamental theorem of homological algebra}. We intend to
construct and study objects that are in a sense \emph{recursively} constructed,
that is, constructed from previous objects already studied. And we need tools
to study the new objects using the informations that are known for the previous
ones.

Typically, many topological spaces can be described as the total space of a
fibration. This total space \(E\) is then presented as a twisted product of two
other spaces: \(E := F \times_\tau B\); the space \(B\) (resp. \(F\)) is the
base space (resp. fibre space) and instead of the ordinary product \(F \times
B\), some important modification in the construction of the product, following
the instructions given by the \emph{twisting function} \(\tau\), allows one to
construct a different space, for some reason or other. For example in
Section~\ref{62859} we have constructed \(X_4\) and \(X_5\) as twisted products
\(X_4 = K(\bZ, 2) \times_\tau S^3\) and \(X_5 = K(\bZ_2, 3) \times_{\tau'}
X_4\) where \(\tau\) and \(\tau'\) were chosen to ``kill'' the first non-null
homotopy group of \(S^3\) and \(X_4\).

So that the game rule is the following. Given: the homological nature of \(F\)
and \(B\). Problem: How to determine the \emph{same} information for \(E = F
\times_\tau B\)? In this case, the Eilenber-Zilber theorem gives the homology
of the \emph{non-twisted} product \(E' = F \times B\); and if an appropriate
hypothesis is satisfied for \(B\) (simple connectivity), then the basic
perturbation lemma allows to consider the twisted product \(E\) as a
\emph{perturbation} of the non-twisted product \(E'\) and to obtain the
looked-for homological information for \(E\). This will be our \emph{effective
version} of the Serre spectral sequence.

\begin{dfn}---
\emph{Let \((C_\ast, d)\) be a chain-complex. A collection of module morphisms
\(\delta = (\delta_n: C_n \rightarrow C_{n-1})_{n \in \bZ}\) is called a
\emph{perturbation} of the differential \(d\) if the sum \(d + \delta\) is also
a differential, that is, if \((d + \delta)^2 = 0\).}
\end{dfn}

Such a perturbation produces a \emph{new} chain-complex \((C_\ast, d +
\delta)\) and in general the homological nature of the chain-complex is so
deeply\dots\ perturbed. Two theorems are available in this area. The first one,
called the \emph{easy} perturbation lemma, is trivial but useful. The second
one, called the \emph{basic} perturbation lemma (BPL) is not trivial at all: in
a sense it gives more information than some spectral sequences, typically the
Serre and Eilenberg-Moore spectral sequences. The BPL was discovered by Shih
Weishu~\cite{SHIH} to overcome some gaps in the Serre spectral sequence, and
Ronnie Brown gave the abstract modern form~\cite{BRWNR1}.

\subsubsection{Easy perturbation lemma.}\label{33336}

\begin{prp}\label{73716}---
Let \(\rho = (f,g,h): (\hC_\ast, \hd) \rrdc (C_\ast, d)\) be a reduction and
let \(\delta: C_\ast \rightarrow C_{\ast - 1}\) be a perturbation of the
differential \(d\) of the small chain-complex. Then a ``new'' reduction \(\rho
= (f, g, h): (\hC_\ast, \hd + \hdl) \rrdc (C_\ast, d + \delta)\) can be
constructed above the perturbed the chain-complex.
\end{prp}

\proof The differential of the small chain-complex is perturbed, so that a
priori the components \(f\) and \(g\) of the reduction \(\rho\) are no more
compatible with the differentials \(\hd\) and \(d + \delta\). But the reduction
\(\rho\) induces a decomposition \(\hC_\ast = \ker f \oplus C'_\ast\) where
\(C'_\ast = \im g\) is a copy of the small chain-complex \(C_\ast\); so that it
is enough to copy also the perturbation, that is, to introduce the perturbation
\(\hdl = g \delta f\) of \(\hd\). The nature of \(\ker f\) is not modified and
the previous components \(f\), \(g\) and \(h\) of the reduction \(\rho\) can be
let unchanged. This is the reason why the new reduction is not so ``new'', it
is the \emph{same} reduction between \emph{different} chain-complexes! \QED

\subsubsection{Basic perturbation lemma.}

The situation is now dramatically harder: we intend to perturb the differential
of the \emph{big} chain-complex of the reduction. In general it is not possible
to coherently perturb the differential of the small chain-complex, even by
modifying the reduction itself. For example, let \(\hC_\ast\) be the ``big''
chain-complex where \(\hC_n = 0\) except \(\hC_0 = \hC_1 = \bZ\) and \(d_1 =
\id{\bZ}\). This chain-complex is acyclic, which implies there is a reduction
\(\rho = (0, 0, h): \hC_\ast \rrdc 0\) over the null chain-complex. If you
introduce the perturbation \(\hdl_1 = -\id{\bZ}\), then the differential
becomes null, the chain-complex is no more acyclic and it is not possible to
perturb coherently the differential of the null chain-complex, which
differential in fact cannot be actually ``perturbed''. This simple example
shows some further hypothesis is necessary to make possible a coherent
perturbation for the small chain-complex and for the reduction.

\begin{thr}\label{07404}{\bf (Basic Perturbation Lemma)} ---
Let \(\rho = (f,g,h): (\hC_\ast, \hd) \rrdc (C_\ast, d)\) be a reduction and
let \(\hdl\) be a perturbation of the differential~\(\hd\) of the big
chain-complex. We assume the \emph{nilpotency hypothesis} is satisfied: for
every \(c \in \hC_n\), there exists \(\nu \in \bN\) satisfying \((h\hdl)^\nu(c)
= 0\). Then a perturbation \(\delta\) can be defined for the differential \(d\)
and a new reduction \(\rho' = (f', g, h'): (\hC_\ast, \hd + \hdl) \rrdc
(C_\ast, d + \delta)\) can be constructed.
\end{thr}

The nilpotency hypothesis states the composition \(h\hdl\) is pointwise
nilpotent. Note the differential of the small chain-complex is modified but
also the components \((f,g,h)\) of the reduction which become something else
\((f',g',h')\): we will have to perturb these components as well.

Which is magic in the BPL is the fact that a sometimes complicated perturbation
of the ``big'' differential can be accordingly reproduced in the ``small''
differential; in general it is not possible, unless the nilpotency hypothesis
is satisfied.

\proof Because of the nilpotency  condition, the following series have, for
each element which they work on, only a  finite number of non-null terms and
their sums are defined:

\label{85606} $$ \phi = \sum_{i=0}^\infty (-1)^i (h\widehat{\delta})^i ;
\hspace{2cm} \psi = \sum_{i=0}^\infty (-1)^i (\widehat{\delta} h)^i. $$

The operators $\phi$ and $\psi$ have degree 0  and trivially satisfy a few
relations; \emph{these relations are  the only ones that are from now on
utilized}:
\begin{quotation}
$\phi h = h \psi$~;

$\widehat{\delta} \phi = \psi \widehat{\delta}$~;

$\phi = 1-h \widehat{\delta} \phi = 1 - \phi h \widehat{\delta}
      = 1 - h \psi \widehat{\delta}$~;

$\psi = 1 - \widehat{\delta} h \psi = 1 - \psi \widehat{\delta} h
      = 1 - \widehat{\delta} \phi h$.
\end{quotation}

The reduction \(\rho' = (f', g', h'): (\hC_\ast, \hd') \rrdc (C_\ast, d')\) to
be constructed is then simply defined by:
\begin{quotation}
$\widehat{d}'  = \hd + \hdl$ is the new differential of
$\hC_\ast$~;\label{35472}

$d' =  d  +  \delta$ is the new differential of  $C_\ast$ where
\label{80476}$\delta  = f\widehat{\delta} \phi g = f \psi \widehat{\delta} g$~;

$f' = f \psi$~;

$g' = \phi g$~;

$h' = \phi h = h \psi$~.
\end{quotation}

\begin{lmm}\label{81951}
--- Let $(C_\ast, d)$ be a  chain-complex and let $h$
be an operator on $C_\ast$ of degree $+1$, satisfying the relations:
\begin{quotation}
hh = 0~;

hdh = h.
\end{quotation}
Then $D=dh+hd$ is a projector which splits the  chain-complex $C_\ast$ into the
direct sum  of chain-complexes $\ker D  \oplus  \emph{\im} D$ where the second
one is acyclic. More precisely, if \(\gamma\) is the canonical inclusion \(\ker
D \rightarrow C_\ast\), then \((\textrm{\emph{id}} - D, \gamma, h): C_\ast
\rrdc \ker D\) is a reduction.
\end{lmm}

\proof The  operator $D$ is a  projector, because  of the computation: $D^2 =
\mbox{\((d h + h d)^2\)} = d h d h + h d h d = d h + h  d = D$ (because $h h =
0$ and $d d = 0$). The operator $D$  and  therefore also $\id{}-D$ are
chain-complex morphisms~: $d(dh+hd) =  dhd  =  (dh+hd)  d$ (because $dd=0$).
The operator \(h\) also commutes with \(D\) and therefore preserves
\(\ker(\id{}-D)\); it is null on \(\ker D\), for \((dh + hd) = 0\) implies
\(h(dh + hd) = h = 0\).\QED

\textsc{Proof of Theorem continued.} In   the theorem, the operator $h$   does
satisfy these relations with respect to $\widehat{d}$,  because $hh=0$ is
explicitly required among the reduction properties   and $h \widehat{d}  h =
(1-\widehat{d} h - gf)h = h$ (because $hh=0$ and  $fh=0$). The projection $D =
\widehat{d} h+h\widehat{d}$  is also the  difference $1-gf$, and therefore the
complementary projection $1-D$ is the composition $gf$.

The new homotopy operator $h'$ has  been defined by $h'=\phi h=h\psi$. Firstly,
we  naturally  obtain   from the  definition  of   $h'$  the definitions of
$f'$, $g'$ and $\delta$.

The new  operator  $h'$  satisfies  also the  relations  $h'h'=0$  and
$h'\widehat{d}'h'=h'$.    In fact    $h'h'=\phi   h    h \psi=0$   and
$h'\widehat{d}'h'=  \phi  h(\widehat{d}+\widehat{\delta})h\psi  = \phi
h\widehat{d}  h\psi+\phi h\widehat{\delta}  h\psi   =\phi h  \psi+\phi
h(1-\psi)=\phi h = h'$ (because $\widehat{\delta} h \psi = 1-\psi$).

We        then  obtain     from   the      lemma      the  fact   that
$D'=\widehat{d}'h'+h'\widehat{d}'$  is a  projector;  let us denote by $\pi=gf$
the complementary   projector   of $D$  and $\pi'=1-D'$  the complementary
projector of $D'$.

We already know  the relations $hh=h'h'=0$. Furthermore $hh'=hh\psi=0$ and
$h'h=\phi  h h=0$. In fact  any composition of  an operator of type $h$  with
an  operator of   type $\pi$  is    null. Firstly $\pi h   = (1-\widehat{d} h-h
\widehat{d})h=h-h\widehat{d}      h=h-h=0$  and $h\pi=h(1-\widehat{d}
h-h\widehat{d})=h-h\widehat{d}   h=h-h=0$.  Next $\pi  h'=\pi h\psi=0$ and
$h'\pi=\phi h\pi=0$. Then $\pi'h'=h'\pi'=0$ is      proved      like $\pi
h=h\pi=0$.           Finally $\pi'h=(1-\widehat{d}'h'-h'\widehat{d}')h$~; but
$h'h=0$   and $\widehat{d}'=\widehat{d}+\widehat{\delta}$,   therefore
$\pi'h=h-\phi h(\widehat{d}+\widehat{\delta})h    =h-\phi      h\widehat{d}
h-\phi h\widehat{\delta} h  =  h-\phi h-(1-\phi)h=0$  (because  $h\widehat{d}
h=h$ and  $\phi  h\widehat{\delta}=1-\phi$).   In   the   same   way
$h\pi'=h(1-\widehat{d}'h'-h'\widehat{d}')                            =
h-h(\widehat{d}+\widehat{\delta})h\psi=h-h\widehat{d} h\psi-h\widehat{\delta}
h\psi=h-h\psi-h(1-\psi)=0$.

Let us now  consider the compositions $\pi\pi'\pi$  and $\pi'\pi\pi'$. Firstly
$\pi\pi'\pi  =   \pi  (1-\widehat{d}'h'-h'\widehat{d}')  \pi = \pi^2=\pi$,
because    $\pi  h'  =   h'\pi=0$.    In  the   same way
$\pi'\pi\pi'=\pi'(1-\widehat{d}   h-h\widehat{d})\pi'   =\pi'^2=\pi'$.
Therefore the operators $\pi$ and $\pi'$ are inverse morphisms between the
images of  $\pi'$ and $\pi$~;  they are only homomorphisms {\em of graded
modules},  in    general  non  compatible  with  the  natural differentials of
the  respective images. But the  image of  $\pi$ has a bijective mapping
towards the small graded module  $C_\ast$ through $f$ and $g$, so that a
composition provides an isomorphism of graded  modules between $C_\ast$ and the
image of  $\pi'$ which allows  us to install a new differential  on  $C_\ast$
deduced  from  the differential  of $\im{\pi'}$, restriction of
$\widehat{d}'=\widehat{d}+\widehat{\delta}$.

Firstly let us  note that $h'g=\phi   hg=0$, and that  $fh'=fh\psi=0$. Taking
account of what was explained in the  previous paragraph, it is natural  to
define $g'=\pi'g=(1-\widehat{d}'h'-h'\widehat{d}')g  = g-\phi  h\widehat{d}
g-\phi  h\widehat{\delta}   g =-\phi hgd+(1-\phi h\widehat{\delta})g=\phi g$.
Then the ``projection''  $f'$ will be the composition of the   actual
projection $\pi'$ with    the composition $f\pi$. But  $f\pi=f(1-\widehat{d}
h-h\widehat{d})=f-f\widehat{d} h-fh\widehat{d}= f-dfh-fh\widehat{d}=f$  and  we
obtain $f'=f\pi\pi'=f\pi'= f(1-\widehat{d}'h'-h'\widehat{d}')=f-f\widehat{d}
h\psi-f\widehat{\delta} h\psi=-\widehat{d} fh\psi+f(1-\widehat{\delta} h\psi)
=f\psi$.   We  have obtained   the  announced formulas  for the desired
reduction components $f'$ and $g'$.

The  new differential to be  installed on the graded module underlying $C$
remains to    be  determined. We  naturally  compute:   $d+\delta=
f\pi(\widehat{d}+\widehat{\delta})\pi'g                              =
f(\widehat{d}+\widehat{\delta})\pi    g      =       f\widehat{d}\pi'g
+f\widehat{\delta}\phi                                               g
=f\widehat{d}(1-\widehat{d}'h'-h'\widehat{d}')g+f\widehat{\delta}\phi g =
f\widehat{d}                        g
-f\widehat{d}\widehat{d}'h'g-dfh'\widehat{d}'g+f\widehat{\delta}\phi
g=f\widehat{d}  g + f\widehat{\delta}\phi g =d+f\widehat{\delta}\phi g =
d+f\psi\widehat{\delta}  g$~;    we    must   therefore   choose
$\delta=f\widehat{\delta}\phi g=f\psi\widehat{\delta} g$.

The basic perturbation lemma is proved. \QED

We will frequently -- not always -- use the basic perturbation lemma in the
following context.

\subsection{Objects with effective homology.}

An object \emph{with effective homology} is a complex object made of a
\emph{locally effective} object -- the object under study, an \emph{effective}
object -- namely an effective chain-complex describing the homological nature
of the object under study, both objects being connected by an appropriate
homology equivalence. Because of the latter, the homological problem for the
underlying object is solved.

\begin{dfn}---
\emph{A \emph{strong homology equivalence}, in short an \emph{equivalence}
\(\varepsilon: C_\ast \eqvl D_\ast\) between two chain-complexes is a pair of
reductions connecting \(C_\ast\) and \(D_\ast\) through a third chain-complex
\(\hC_\ast\):
\[
\varepsilon = \fbox{\(C_\ast \stackrel{\rho_\ell}{\lrdc} \hC_\ast
\stackrel{\rho_r}{\rrdc} D_\ast\)}
\]}
\end{dfn}

Because of the fundamental importance of this sort of equivalence, this will be
simply called in this text an \emph{equivalence}. If the homological problem is
solved for~\(D_\ast\), it is also solved for \(C_\ast\).

\begin{dfn}\label{93815}---
\emph{An \emph{object with effective homology} \(X\) is a quadruple \(X = (X,
C_\ast X, EC_\ast, \varepsilon)\) where: \begin{itemize}
\item
\(X\) is a \emph{locally effective} object, the homological nature of which is
under study.
\item
\(C_\ast X\) is the (locally effective) chain-complex canonically associated
with \(X\) when the homological nature of \(X\) is studied.
\item
\(EC_\ast\) is an \emph{effective} chain-complex.
\item
Finally \(\varepsilon\) is an equivalence \(\varepsilon: C_\ast X \eqvl
EC_\ast\).
\end{itemize} }
\end{dfn}

Typically the object \emph{under study} \(X\) could be an infinite simplicial
complex; if it is infinite, we must content ourselves with a \emph{locally
effective} implementation. Then \(C_\ast X\) is the chain-complex canonically
associated with it (Section~\ref{70994}); it is not of finite type and it is
also implemented as a \emph{locally effective} chain-complex. In many
situations, the homology groups of this chain-complex yet are of finite type:
so that some \emph{effective} chain-complex can have the right homology groups.
The last but not the least, an \emph{equivalence} between the genuine
chain-complex associated with our object and our effective chain-complex will
play an essential role in the next constructions. In most cases, the
\emph{basic perturbation lemma} will be the main tool constructing new
equivalences from others already constructed.

A good didactic simple example of object with effective homology, didactic but
very useful, is the Koszul complex \(\Ksz_\fR(\fR)\), see Section~\ref{91102}.
If \(\fk\) is the underlying ground field, then the Koszul complex has an
infinite \(\fk\)-dimension; but it is a resolution of \(\fk\) and its homology
is only \(\fk\) in dimension 0, nothing else. The reduction \(\Ksz_\fR(\fR)
\rrdc \fk\) which will be constructed is the equivalence
component~\(\varepsilon_\fR\). So that the quadruple \((\Ksz_\fR(\fR),
\Ksz_\fR(\fR), \fk, \varepsilon_\fR)\) is a version \emph{with effective
homology} of the Koszul complex. In this case, and this is not seldom, the
object under study is the chain-complex itself.

The main result of Effective Homology Theory is the following ``meta-theorem''.

\begin{mth}---
Let \(X_1, \ldots, X_k\) be a collection of objects and \(\phi\) some
``reasonable'' \emph{constructor} \(\phi: (X_1, \ldots, X_k) \mapsto X\). Then
a \emph{version with effective homology} \(\phi_{EH}\) can be obtained,
constructing a \emph{version \(X_{EH}\) with effective homology} of the result
\(X\) of the construction when \emph{versions with effective homology} of the
\(X_i\)'s are given:
\[
\phi_{EH}: ((X_1, C_\ast X_1, EC{^X_1}_\ast, \varepsilon_1), \ldots, (X_k,
C_\ast X_k, EC{^X_k}_\ast, \varepsilon_k)) \mapsto (X, C_\ast X, EC^X_\ast,
\varepsilon).
\]
\end{mth}

The nature of \emph{constructive} homological algebra is now simply defined:
please transform the standard theorems of homological algebra into instances of
this meta-theorem. The version with effective homology \(\phi_{EH}\) of the
constructor \(\phi\) is a collection sometimes sizeable of algorithms
constructing algorithm components of the result \(X_{EH} = (X, C_\ast X,
EC^X_\ast, \varepsilon)\) from the algorithm components of the data \(X_{i,EH}
= (X_i, C_\ast X_i, EC{^X_i}_\ast, \varepsilon_i)\). An \emph{algorithm}
constructing \emph{algorithms} from other \emph{algorithms} requires
\emph{functional programming}; this wonderful tool is theoretically known since
Church's work in logic~\cite{CHRC}, a theoretical work leading to the currently
most complete programming language, Common Lisp.

\section{Constructive Homology and Commutative Algebra.}\label{74280}

\subsection{Presentation.}

The homological framework was not available at Hilbert's time, but among his
famous results in Commutative Algebra, typically the theorems about
\emph{syzygies}, many of them in fact have a \emph{homological} nature. Henri
Cartan and Sam Eilenberg~\cite{CREI} understood the algebraic tools of
Algebraic Topology can be organized to be fruitfully used in other domains, for
example in Commutative Algebra: it was the birth of the subject
\emph{Homological Algebra}.

We explain in this section how the point of view of \emph{constructive}
homological algebra gives new insights about some homological domains of
commutative algebra. The following theme is particularly convenient. A
classical theorem, the bicomplex spectral sequence, allows one to prove the
equivalence of both definitions of torsion groups:
\[
H_\ast(\Rsl_\fR(M) \otimes_\fR N) =: \Tor^\fR_\ast(M, N) := H_\ast(M
\otimes_\fR \Rsl_\fR(N)).
\]
\(\fR\) is a commutative unitary ground ring, \(M\) and \(N\) are two
\(\fR\)-modules. The torsion groups of \(M\) and \(N\) are defined by taking
for example an \(\fR\)-resolution \(\Rsl(M)\) of \(M\) and computing the tensor
product \(\Rsl(M) \otimes_\fR N\); the last chain-complex in general is no
longer exact, and its homology groups are the torsion groups. If you do the
symmetric work with a resolution of \(N\), the result is the same; the result
is also independent of the chosen resolutions, so that these torsion groups
express deep abstract relations between the modules \(M\) and \(N\).

We intend to illustrate that a systematic \emph{constructive} point of view in
these homological notions produces new methods and also allows their users to
have a more global understanding of the various studied properties. The last
but not the least, most often the proofs are more elementary! We will so obtain
the striking result: there is a perfect direct equivalence between the
\emph{effective} homology of \(\Ksz(M)\), the Koszul complex of an
\(\fR\)-module \(M\), and a resolution \(\Rsl_\fR(M)\) of the same module with
respect to the ground ring \(\fR\).

\subsection{Koszul complex.}\label{66039}

\begin{uos}---
In this section about Commutative Algebra, the ground ring \(\fR\) is \(\fR =
\fk[x_1, \ldots, x_m]_0\), the usual polynomial ring with \(m\) variables,
\emph{localized} at \mbox{\(0 \in \fk^m\)}. We denote by \(V\) the ``abstract''
vector space \(V = \fk^m\) provided with the basis \((dx_1, \ldots, dx_m)\).
\end{uos}

The ground \emph{field} \(\fk\) is an arbitrary commutative field, in
particular the case of a finite characteristic is covered without any extra
work. An element of \(\fR\) is a ``quotient'' \(P/Q\) of two polynomials, the
second one being non-null at 0. It happens the denominators, because of the
context, will not play any role, but the general correct framework is the case
of \(\fR\) a \emph{regular local ring}\footnote{The \(dx_i\)'s in the
forthcoming definition of the Koszul complex are essentially a dual
\(\fk\)-basis of \(\fm_0 /\fm_0^2\) for \(\fm_0\) the maximal ideal of
\(\fk[x_1, \ldots, x_m]\) at 0. The Koszul complex so defined analyzes the
\emph{local} properties of an \(\fR\)-module at 0. Using the ordinary
polynomial ring \(\fRb = \fk[x_1, \ldots, x_m]\) does not work. Consider for
example the ideal \(I = \ideal{x}\) of \(\fk[x]\) and the \(\fk[x]\)-module \(M
= \fk[x]/I\). The homology of the Koszul complex is trivial, whereas the module
is not; using instead the local ring \(\fk[x]_0\) as ground ring, then \(M_0 =
\fk[x]_0/I\) is trivial: the Koszul complex analyses the initial module at 0
\emph{only}.}. The basic reference about local rings is~\cite{SERR5}. To make
significantly more readable the exposition, we prefer to consider only the case
of \(\fR = \fk[x_1, \ldots, x_m]_0\).

\begin{dfn}---
\emph{The \emph{Koszul complex} \(\Ksz(M)\) of the \(\fR\)-module \(M\) is a
chain-complex of \(\fR\)-modules constructed as follows. The chain group in
degree \(n \geq 0\) is \(\Ksz_n(M) = \mbox{\(M \otimes_\fk \wedge^n V\)}\) and
the differential \(d: \Ksz_n(M) \rightarrow \Ksz_{n-1}(M)\) is defined by the
formula:
\[
\begin{array}{rcl}
d(\alpha \, dx_{i_1} \ldots dx_{i_n}) &=& \alpha x_{i_1} \, dx_{i_2} \ldots
dx_{i_n}
\\
&& - \alpha x_{i_2} \, dx_{i_1} . dx_{i_3} \ldots dx_{i_n}
\\
&& + \cdots
\\
&& + (-1)^{n-1} \alpha x_{i_n} \, dx_{i_1} \ldots dx_{i_{n-1}}.
\end{array}
\]
}
\end{dfn}

Observe we write simply \(\alpha dx_2.dx_4.dx_5\) instead of \(\alpha \otimes
(dx_2 \wedge dx_4 \wedge dx_5)\) if \(\alpha \in M\). The definition can be
generalized to an arbitrary collection of elements \((\alpha_1, \ldots,
\alpha_p)\) of \(\fR\) instead of the ``variables'' \((x_1, \ldots, x_m)\); the
differential of \(dx_i\) (\(1 \leq i \leq p\)) is then \(\alpha_i\).

The usual sign game shows the Koszul complex actually is a chain-complex.
Furthermore this will be also a \emph{consequence} of a recursive construction
given soon.

\subsection{Geometrical interpretation.}\label{69368}

The construction of a Koszul complex is a little strange, but becomes more
natural if we give a geometrical interpretation, in fact historically at the
origin of this notion~\cite{KSZL}. This interpretation is never used later in
this text.

In our environment, you must think of the ring \(\fR\) as a topological group,
used as a \emph{structural group} to construct fibrations. The exterior algebra
\(\wedge V\) is also a \emph{coalgebra} for the shuffle coproduct:
\[
\Delta(v_1 \wedge \cdots \wedge v_k) = \sum (-1)^{\sigma} \, (v_{\sigma_1}
\wedge \cdots \wedge v_{\sigma_\ell}) \otimes (v_{\sigma_{\ell + 1}} \wedge
\cdots \wedge v_{\sigma_k})
\]
where the sum is taken with respect to all the shuffles \(((\sigma_1 < \cdots <
\sigma_\ell), (\sigma_{\ell + 1} < \cdots < \sigma_k)\) for \(0 \leq l \leq
k\). The coalgebra structure of \(\wedge V\) gives it a flavor of
\emph{topological space}, think of the Alexander-Whitney coproduct over the
singular chain-complex of a topological space.

In the particular case \(M = \fR\), the Koszul complex \(\Ksz(\fR)\) can be
viewed as a principal fibration, the ``base space'' being \(\wedge V\) and the
``structural group'' \(\fR\). This is made more explicit in the notation
\(\Ksz(\fR) := \fR \otimes_t \wedge V\) to be understood as follows: the Koszul
complex is a twisted (index \(t\) of \(\otimes_t\)) product of the base space
\(\wedge V\) by the structural group \(\fR\), the twist \(t\) being defined by
a \emph{twisting cochain} \(t \in H^1(\wedge V ; \fR)\); in the particular case
of the Koszul complex, this twisting cochain is null outside the degree 1
component \(\wedge^1 V\) of \(\wedge V\) and \(t(dx_i) := x_i\); see for
example~\cite[\S~30]{MAY} for the general definition of the notion of twisting
cochain. Such a twisting cochain is the translation in the algebraic framework
of the \emph{coordinate functions}, more precisely of the coordinate changes
defining a fibre bundle~\cite[Section~I.2]{STNR}.

Finally, if \(M\) is an arbitrary \(\fR\)-module, it can be understood as a
topological space provided with an action \(M \otimes_\fk \fR \rightarrow M\),
which allows us to interpret the Koszul complex \(\Ksz(M) = M \otimes_t \wedge
V = M \otimes_\fR (\fR \otimes_t \wedge V)\) as the fibre bundle canonically
associated with the principal bundle \(\fR \otimes_t \wedge V\).

The chain-complex \(\Ksz(\fR)\) is acyclic, and we will see the \emph{homotopy
operator} proving this fact will play a quite essential role in our study. So
that \(\Ksz(\fR)\) has the ``homotopy type'' of a point; in other words the
fibration \(\fR \otimes_t \wedge V\) is the \emph{universal \(\fR\)-fibration};
in algebraic language \(\Ksz(\fR)\) is a \emph{resolution} of the ground field
\(\fk\).

\subsection{Tensor products of chain-complexes.}\label{05496}

We give here a few general technical results about tensor products of
chain-complexes and reductions. The ring \(\fR\) \emph{in this subsection} is
again an arbitrary commutative unitary ring.

\begin{dfn} \emph{\textbf{(Koszul convention)} ---
Let \(C_\ast\) and \(D_\ast\) be two graded modules \(C_\ast = \oplus_n C_n\)
and \(D_\ast = \oplus_n D_n\). A natural graduation is induced over \(T_\ast =
C_\ast \otimes D_\ast = \oplus_n (\oplus_{p+q=n} C_p \otimes D_q)\). If \(f:
C_\ast \rightarrow C'_{\ast + k}\) and \(g: D_\ast \rightarrow D'_{\ast +
\ell}\) are graded morphisms of respective degrees \(k\) and \(\ell\), then the
tensor product \(f \otimes g: (C \otimes D)_\ast \rightarrow (C' \otimes
D')_{\ast + k + l}\) is defined by \((f \otimes g)(a \otimes b) :=
(-1)^{\ell|a|} f(a) \otimes g(b)\) if \(a\) is homogeneous of degree \(|a|\).}
\end{dfn}

We think the necessary permutation of \(g\) (degree \(\ell\)) and \(a\) (degree
\(|a|\)) generates a signature \((-1)^{\ell |a|}\).

\begin{dfn} ---
\emph{Let \((C_\ast, d)\) and \((C'_\ast, d')\) be two chain-complexes of
\(\fR\)-modules. The tensor product \((C_\ast, d) \otimes (C'_\ast, d')\) is a
chain-complex defined as the module \((C_\ast \otimes C'_\ast)\) provided with
the differential \(d_{C_\ast \otimes C'_\ast} := d \otimes \id{C'_\ast} +
\id{C_\ast} \otimes d'\) where the Koszul convention must be applied. The
identity being of degree 0 and a differential of degree -1, this implies \(d(a
\otimes b) = da \otimes b + (-1)^{|a|} a \otimes db\).}
\end{dfn}

The tensor product operator is an important functor and we must be able to
define the tensor product of two reductions. It is better to start with the
composition of reductions.

\begin{prp}\label{28461} ---
Let \(\rho = (f,g,h): C_\ast \rrdc C'_\ast\) and \(\rho' = (f',g',h'): C'_\ast
\rrdc C''_\ast\) be two reductions. These reductions can be \emph{composed},
producing the reduction \(\rho'' = (f'', g'', h''): C_\ast \rrdc C''_\ast\)
with:
\[
\begin{array}{rcl}
f'' &=& f' f;
\\
g'' &=& g g';
\\
h'' &=& h + g h' f.
\end{array}
\]
\end{prp}

\proof Exercise. \QED

\begin{prp}\label{49390} ---
Let \(\rho = (f, g, h): C_\ast \rrdc D_\ast\) and \(\rho' = (f',g',h'): C'_\ast
\rrdc D'_\ast\) be two reductions. Then a tensor product:
\[
\rho'' = (f'', g'', h'') : C_\ast \otimes C'_\ast \rrdc D_\ast \otimes D'_\ast
\]
can be defined, with:
\[
\begin{array}{rcl}
f'' &=& f \otimes f' ;
\\
g'' &=& g \otimes g' ;
\\
h'' &=& h \otimes \id{C'_\ast} + gf \otimes h'
\end{array}
\]
\end{prp}

\proof Compose the reductions
\[
\begin{array}{rcl} \rho \otimes \id{C'_\ast} &:&
C_\ast \otimes C'_\ast \rrdc D_\ast \otimes C'_\ast
\\
\id{D_\ast} \otimes \rho' &:& D_\ast \otimes C'_\ast \rrdc D_\ast \otimes
D'_\ast.
\end{array}
\]\QED

Note the lack of symmetry in the result; you could replace the intermediate
complex by \(C_\ast \otimes D'_\ast\) and \(h'' = h \otimes \id{C'_\ast} + gf
\otimes h'\) by \(h'' = \id{C_\ast} \otimes h' + h \otimes g'f'\).

\subsection{Cones of chain-complexes.}

The \emph{cone constructor} is important in homological algebra, and we study
here the most elementary properties. We will meet the first application of the
BPL.

\begin{dfn}---
\emph{Let \(C _\ast\) and \(D_\ast\) be two chain-complexes and \(\phi: C_\ast
\leftarrow D_\ast\) be a chain-complex morphism. Then the cone of \(\phi\)
denoted by \(\Cone(\phi)\) is the chain-complex \(\Cone(\phi) = A_\ast\)
defined as follows. First \(A_n := C_n \oplus D_{n-1}\); then the boundary
operator is given by the matrix:
\[
d_{A_\ast} := \left[\begin{array}{cc}
 d_{C_\ast} & \phi \\ 0 & -d_{D_\ast}
 \end{array}\right]
\]
}
\end{dfn}

We prefer to turn to the left the arrow from \(D_\ast\) to \(C_\ast\), because
a cone is in fact a particular case of a bicomplex and experience shows it is
convenient to keep one's organisation as homogeneous as possible. The diagram
clearly explaining the nature of a cone is the following.
\[
\xymatrix{ \cdots & D_{n-2} \ar[l]_-{-d_D}\ar[dl]|{\phi} \ar@{.}[d]|{\oplus} &
D_{n-1} \ar[l]_-{-d_D}\ar[dl]|{\phi} \ar@{.}[d]|{\oplus} & D_n
\ar[l]_-{-d_D}\ar[dl]|{\phi} \ar@{.}[d]|{\oplus}& \cdots
\ar[l]_-{-d_D}\ar[dl]|{\phi}
\\
\cdots & C_{n-1} \ar[l]^-{d_C} & C_n \ar[l]^-{d_C} & C_{n+1} \ar[l]^-{d_C} &
\cdots \ar[l]^-{d_C}}
\]

You see the morphism \(\phi\) contributes to the differential of the cone. If
you do not change the sign of \(d_{D_\ast}\) in the cone the rule \(d \circ d=
0\) would not be satisfied. With our sign choice:
\[
\left[\begin{array}{cc} d & \phi \\ 0 & -d \end{array}\right] \
\left[\begin{array}{cc} d & \phi \\ 0 & -d \end{array}\right] =
\left[\begin{array}{cc} d^2 & d \phi - \phi d \\ 0 & d^2 \end{array}\right] =
0.
\]
for the initial differentials satisfy \(d^2 = 0\) and \(\phi\) is a
chain-complex morphism satisfying \(d_C \phi = \phi d_D\). In fact the Koszul
convention has been applied: the suspension operator \(\sigma\) which increases
the degree by 1 is implicitly applied to the elements of~\(D_\ast\) and this
suspension operator has degree +1. So that Koszul teaches us that \(d(\sigma c)
= - \sigma(dc)\) is the good choice: the morphisms \(\sigma\) (suspension,
degree +1) and \(d\) (differential, degree -1) have been permuted.

Studying carefully the next simple application of the BPL (basic perturbation
lemma) gives an excellent understanding of this wonderful Theorem strangely
called ``lemma''. This application is not difficult; all the applications have
the same style and this one is the simplest one. Consider this particular case
as the ideal didactic situation to \emph{learn} how to use the BPL; the other
applications are not more difficult, even in more or less terrible
environments.

\begin{thr}\label{47942} \textbf{\emph{(Cone Reduction Theorem)}} ---
Let \(\rho = (f,g,h) : C_\ast \rrdc D_\ast\) and \(\rho' = (f',g',h'): C'_\ast
\rrdc D'_\ast\) be two reductions and \(\phi: C_\ast \leftarrow C'_\ast\) a
chain-complex morphism. Then these data define a canonical reduction:
\[
\rho'' = (f'', g'', h'') : \Cone(\phi) \rrdc \Cone(f \phi g').
\]
\end{thr}

\proof This would be trivial if \(\phi = 0\): in such a case, we have also \(f
\phi g' = 0\) and the cones are simple direct sums (with a suspension applied
over \(C'_\ast\) and \(D'_\ast\)) and defining a direct sum of reductions is
trivial. Now look carefully at this diagram:
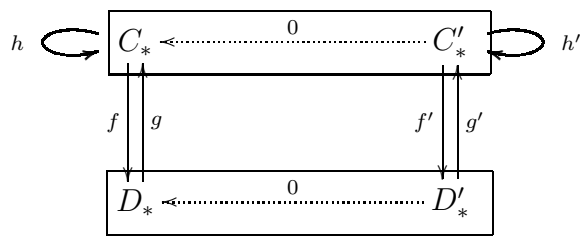
\[{\xymatrix @R150\ul @C100\ul{
  {\scriptstyle h} \ar@{}[r]|-{\autoarrow{1}}
  &
  *+{C_\ast} \ar@<-10\ul>[d]_{f}
  &&&
  *+{C'_\ast} \ar@{.>}[lll]_{{0}} \ar@<-10\ul>[d]_{f'}
  \ar@{*{}*{}*{}}@/_50\ul/[lll]_{\phantom{h \phi h'}}
  &
  {\scriptstyle h'} \ar@{}[l]|-{\autoarrow{-1}}
  \\
  &
  *+{D_\ast} \ar@<-10\ul>[u]_{g}
  &&&
  *+{D'_\ast} \ar@{.>}[lll]_{{0}} \ar@<-10\ul>[u]_{g'}
  \save "1,2"."1,5"*+[F]\frm{} \restore
  \save "2,2"."2,5"*+[F]\frm{} \restore
  }}
  \]
The rectangular boxes intend to visualize the cone constructions, simple direct
sums when the chain-complex morphisms are null. The suspensions applied to the
right-hand chain-complexes are not shown. Each chain-complex of these (trivial)
cones is a direct sum, so that the morphisms of our initial reduction are
represented by \(2 \times 2\) matrices:
\[\begin{array}{ccccc}
  \left[\begin{array}{cc} d_C & {0} \\ 0 & -d'_{C'}
  \end{array}\right]
  &
  \left[\begin{array}{cc} d_D & {0} \\ 0 & -d'_{D'} \end{array}\right]
  &
  \left[\begin{array}{cc} f & 0 \\ 0 & f' \end{array}\right]
  &
  \left[\begin{array}{cc} g & 0 \\ 0 & g' \end{array}\right]
  &
  \left[\begin{array}{cc} h & 0 \\ 0 & -h' \end{array}\right]
  \\
  d_{\textrm{\scriptsize top}} & d_{\textrm{\scriptsize bottom}} &
  f \oplus f' & g \oplus g' & h \oplus -h'
  \end{array}
\]
A homotopy operator has degree +1 and the Koszul convention must also be
applied between suspension and homotopy. Now if we install the right
morphism~\(\phi\) on the top cone, the reduction is nomore valid, the top
differential is modified and there is no reason the pairs \(f \oplus f'\) and
\(g \oplus g'\) are compatible with the new top differential. It is exactly in
such a situation the BPL is to be used. In this case the \emph{perturbation} to
be applied to \(d_{\textrm{\scriptsize top}}\) is:
\[
  \left[\begin{array}{cc} d_C & {0} \\ 0 & -d'_{C'}
  \end{array}\right] + \left[\begin{array}{cc} 0 & {\phi} \\ 0 & 0
  \end{array}\right] \mapsto \left[\begin{array}{cc} d_C & {\phi} \\ 0 & -d'_{C'}
  \end{array}\right]
\]

This is frequent in applications of the BPL, the perturbations are extra arrows
installed in the diagram after the starting situation, here only one arrow
\(\phi\). The BPL can be used only if the nilpotency condition is satisfied.
The composition \(h\hd\) of Theorem~\ref{07404} is here:
\[
\left[\begin{array}{cc} h & {0} \\ 0 & -h'_{C'}
  \end{array}\right] \circ
  \left[\begin{array}{cc} 0 & {\phi} \\ 0 & 0
  \end{array}\right]
  =
  \left[\begin{array}{cc} 0 & {h \phi} \\ 0 & 0
  \end{array}\right]
\]
which is clearly nilpotent. Instead of formal computations, verifying the
nilpotency condition is most often the following game: follow a perturbation
arrow, then a homotopy arrow, then a perturbation arrow, and so on. You must
show this treasure hunt, in general several possible choices at each step,
terminates after a finite number of steps, whatever your choices are. Here the
longest path is \(h \phi h'\) and it is not possible to extend this path of
length 3, the nilpotency condition is therefore satisfied.

We remind you of the magic Shih's formulas in the general framework of
Theorem~\ref{07404}, in particular with the notations of Theorem~\ref{07404}:
\[\phi = \sum_{i=0}^\infty (-1)^i
(h\widehat{\delta})^i ; \hspace{2cm} \psi = \sum_{i=0}^\infty (-1)^i
(\widehat{\delta} h)^i.
\]
\[
\delta  = f\widehat{\delta} \phi g = f \psi \widehat{\delta} g ;\hspace{0.5cm}
f' = f \psi ;\hspace{0.5cm} g' = \phi g ;\hspace{0.5cm} h' = \phi h = h \psi.
\]

Applying these formulas to our particular situation gives:
\[
\phi_{\textrm{\scriptsize Shih}} =   \left[\begin{array}{cc} 1 & {-h \phi} \\ 0
& 1 \end{array}\right];\hspace{1cm} \psi_{\textrm{\scriptsize Shih}} =
\left[\begin{array}{cc} 1 & {-\phi h'} \\ 0 & 1 \end{array}\right];
\]
and then, with our current notations, except \(\delta\) being the perturbation
to apply to the bottom cone:
\[
\delta = \left[\begin{array}{cc} 0 & {f \phi g'} \\ 0 & 0 \end{array}\right];
f'' = \left[\begin{array}{cc} f & {-f \phi h'} \\ 0 & f' \end{array}\right];
g'' = \left[\begin{array}{cc} g & {-h \phi g'} \\ 0 & g' \end{array}\right];
h'' = \left[\begin{array}{cc} h & {h \phi h'} \\ 0 & -h' \end{array}\right];
\]

In other words we have successfully constructed the right new reduction between
\(\Cone(\phi)\) and \(\Cone(f \phi g')\):
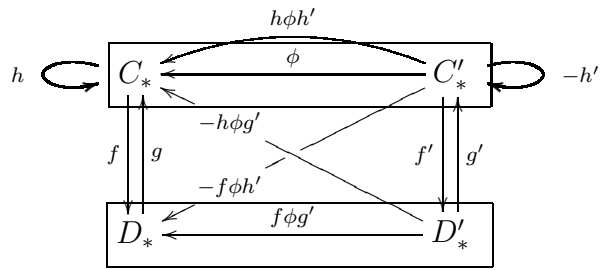
\[
{\label{37174}{\xymatrix @R150\ul @C100\ul{
  {\scriptstyle h} \ar@{}[r]|-{\autoarrow{1}}
  &
  *+{C_\ast} \ar@<-10\ul>[d]_{f}
  &&&
  *+{C'_\ast} \ar[lll]_{{\phi}} \ar@<-10\ul>[d]_{f'}
  \ar[dlll]|{\hole}|(0.7)*+{\scriptstyle {-f \phi h'}}
  \ar@/_50\ul/[lll]_{{h \phi h'}}
  &
  {\scriptstyle -h'} \ar@{}[l]|-{\autoarrow{-1}}
  \\
  &
  *+{D_\ast} \ar@<-10\ul>[u]_{g}
  &&&
  *+{D'_\ast} \ar[lll]_{{f \phi g'}} \ar@<-10\ul>[u]_{g'}
  \ar[ulll]|(0.7)*+{\scriptstyle {-h \phi g'}}
  \save "1,2"."1,5"*+[F]\frm{} \restore
  \save "2,2"."2,5"*+[F]\frm{} \restore
  }}}
\]\QED

Of course, most often this theorem is proved \emph{without} using the BPL, but
experience shows it is not so easy to guess the right compositions and the
right signs. Once the BPL is understood, it is easier to use it to prove the
cone theorem.

\subsection{Resolutions.}\label{45205}

In this section which has a general scope, the ring \(\fR\) is an arbitrary
unitary commutative ring.

\begin{dfn} ---
\emph{Let \(M\) be an \(\fR\)-module. A \emph{free \(\fR\)-resolution} of
\(M\), in short a \emph{resolution} of \(M\), is a chain-complex \(\Rsl(M)\)
null in negative degrees, made of \emph{free} \(\fR\)-modules, every
differential is an \(\fR\)-morphism, every homology group \(H_n(\Rsl(M))\) is
null except \(H_0(\Rsl(M))\): an \(\fR\)-isomorphism \(\varepsilon:
H_0(\Rsl(M)) \cong M\) is \emph{given}.}
\end{dfn}

Note the isomorphism is a component of the data defining the resolution;
strictly speaking the resolution is the \emph{pair} \((\Rsl(M), \varepsilon)\).
You can also consider the isomorphism \(\varepsilon\) as coming from a morphism
called \emph{augmentation} \(\overline{\varepsilon}: \Rsl_0(M) \rightarrow M\).
If you ``add'' \(\Rsl_{-1}(M) := M\) and this augmentation, you obtain the
exact sequence:
\[
0 \leftarrow M \stackrel{\overline{\varepsilon}}{\leftarrow} \Rsl_0(M)
\leftarrow \Rsl_1(M) \leftarrow \cdots
\]
but the good point of view is not to include \(M\) which must be isomorphic to
the \(H_0\)-group of the resolution, with a \emph{given} isomorphism. The
functional notation \(\Rsl(M)\) is justified by the fact such a resolution is
unique up to homotopy, a point not very important here.

Another detail about notations in this context must be given. Sometimes the
module \(M\) is better considered as a chain-complex concentrated in
dimension~0, with null differentials, in particular when we will soon consider
the notion of \emph{effective} resolution; to emphasize this point of view, we
sometimes use the \(\ast\)-notation \(M_\ast := [\cdots
\stackrel{0}{\leftarrow}  0 \stackrel{0}{\leftarrow} M \stackrel{0}{\leftarrow}
0 \stackrel{0}{\leftarrow} \cdots]\). Both points of view have their own
interest and it is not always possible to keep a constant notation.

We want to make \emph{effective} (or constructive) the definition of a
resolution. We want to make explicit a contracting homotopy \emph{proving} the
very nature of our resolution. In general there is no hope these homotopy
operators are \(\fR\)-morphisms. To keep some linear behaviour, \emph{we assume
now} our ground ring \(\fR\) is a \(\fk\)-algebra with respect to a commutative
field \(\fk\). It is the case for the usual rings of commutative algebra, for
example for the ring \(\fk[x_1, \ldots, x_m]_0\).

\begin{dfn} ---
\emph{Let \(M\) be an \(\fR\)-module. An \emph{effective} resolution
\(\Rsl(M)\) is a resolution with a \emph{\((\fR, \fk, \fk)\)-reduction} \(\rho
= (f,g,h): \Rsl(M) \rrdc M_\ast\) where the small chain-complex \(M_\ast\) is
made from \(M\) concentrated in degree 0.}
\end{dfn}

The prefix \((\fR, \fk, \fk)\) for our reduction means we require \(f\) is an
\(\fR\)-morphism, but \(g\) and \(h\) in general are only \(\fk\)-morphisms.

\noindent\textsc{Example}. Let us consider \(\fR = \fk[x]_0\) (one variable).
The evaluation \(\ev_0(P) = P(0)\) gives a structure of \(\fR\)-module to
\(\fk\): \((P, k) \mapsto P(0)k\). What about a resolution of \(\fk\)?

The Koszul complex \(\Ksz(\fR)\) is in this case very simple; it is a
chain-complex concentrated in degrees 0 and 1:
\[
0 \leftarrow \fR \stackrel{d_1}{\leftarrow} \fR.dx \leftarrow 0
\]
with \(d_1(P.dx) = Px\) (\(P \times x\), not P(x)). This \(d_1\) is injective,
and \(H_1 = 0\). The image is the maximal ideal \(\fm\) and therefore \(H_0 =
\fR/\fm = \fk\). The Koszul complex is a \emph{resolution of \(\fk\)}.

But we are not happy with this result, we prefer \emph{effective} resolutions.
Can this resolution be made effective? It is not hard. The projection \(f:
\Ksz(\fR) \rightarrow \fk_\ast\) is given by the composition \(\Ksz_0(\fR)
\rightarrow \Ksz_0(\fR)/d_1(\Ksz_1(\fR)) = \fR/\fm = \fk\). The inclusion \(g:
\fk \rightarrow \Ksz_0(\fR) = \fR\) is the canonical inclusion \(\fk
\hookrightarrow \fR\) \emph{which is not} an \(\fR\)-morphism. Finally the
homotopy operator must be defined in degree 0, the most natural choice being
\(h_0(P) = (P-P(0))/x\), which is not an \(\fR\)-morphism either: \(h_0(1) =
0\) and \(h_0(x) = 1 \neq x h_0(1) = 0\). But \(h_0\) is \(\fk\)-linear.

In a sense we want to extend the elementary study of this example to the
general case. We want to prove that, if \(\fR = \fk[x_1, \ldots, x_m]_0\), then
\(\Ksz(\fR)\) is an \emph{effective free \(\fR\)-resolution} of the ground
field \(\fk\). The proof is inductive, easy if the polynomial ring is not
localized~\cite[Proposition~VII.2.1]{MCLN2}, a little harder but also a little
more interesting in the localized case\footnote{We could also use the
\emph{flatness} property of \(\fR\) as \(\fRb\)-module, but an \emph{effective}
flatness is required; see~\cite[III.5]{MNRR} for the right definition.}. We
must precisely connect our various rings for different numbers of variables.

\begin{ntt}---
The number \(m\) of variables we are interested in is fixed. If \(0 \leq q \leq
m\), we denote by \(I_{q}\) the ideal of \(\fR = \fk[x_1, \ldots, x_m]_0\)
generated by the variables \(x_{q+1}, \ldots, x_m\): \(I_q =\ \ideal{x_{q+1},
\ldots, x_m}\). The quotient ring \(\fR / I_q\) is denoted by \(\fR_q\). We
denote by \(V_q\) the \(\fk\)-vector space of dimension \(m-q\) generated by
the distinguised basis \((dx_{q+1}, \ldots, dx_m)\).
\end{ntt}

The ring \(\fR_q\) is the same as \(\fR\) except any occurence of a variable
\(x_r\) with \(r > q\) is cancelled. So that \(\fR_q\) is the analogous local
ring but with \(q\) variables only. In particular \(\fR = \fR_m\) and \(\fk =
\fR_0\). If \(q \leq r\), canonical morphisms \(f_{q,r}: \fR_r \rightarrow
\fR_q\) and \(g_{q,r}: \fR_q \rightarrow \fR_r\) are defined. The first one is
a projection, it is also an evaluation process consisting in replacing the
variables \(x_{q+1}, \ldots, x_r\) by 0; it is an \(\fR_i\)-morphism for every
\(i\), in particular for \(i = m\). The second one is a canonical inclusion, it
is an \(\fR_i\)-morphism only for \(i \leq q\).

\begin{dfn}---
\emph{The definition of the Koszul complex is extended as follows. We denote by
\(\Ksz^q(M)\) the sub-chain-complex \(\Ksz^q_k(M) = M \otimes_\fk \wedge^k
V_q\) of \(\Ksz(M)\).}
\end{dfn}

The only difference between \(\Ksz^q(M)\) and \(\Ksz(M)\) is that in the first
case a \(dx_i\) with \(i \leq q\) is excluded.

\begin{thr}\label{37681}---
\(\emph{\Ksz}(\fR)\) is an effective free \(\fR\)-resolution of the
\(\fR\)-module \(\fk\).
\end{thr}

It is the particular case \(q = 0\) of the next theorem to be proved by
decreasing induction.

\begin{thr}---
\(\emph{\Ksz}^q(\fR)\) is an effective free \(\fR\)-resolution of the
\(\fR\)-module \(\fR_q\).
\end{thr}

Note strictly speaking such a statement is improper. When we claim some object
is \emph{effective}, we mean some collection of algorithms, more or less
difficult to be constructed, will allow us to justify the qualifier.

\proof The theorem is obvious for \(q = m\): the chain-complex \(0 \leftarrow
\fR \leftarrow 0\) concentrated in degree 0 is a resolution of \(\fR\).

Let us assume the theorem is proved for \(q\) and let us prove it for \(q-1\).
A reduction \(\rho_{q} = (f_{q}, g_{q}, h_{q}): \Ksz^{q}(\fR) \rrdc \fR_{q}\)
is available.

Our simple example above is easily adapted to prove:
\begin{lmm}---
The chain-complex
\[
0 \leftarrow \fR_{q} \stackrel{\times x_q}{\longleftarrow} \fR_{q} \leftarrow 0
\]
is an effective free resolution of \(\fR_{q-1}\).
\end{lmm}

It is a sophisticated and precise way to express the map \(\times x_q\) is
injective and its cokernel is \(\fR_{q-1}\). The relevant reduction is made of
the projection \(f_{q-1,q}\) which is an \(\fR\)-morphism, the injection
\(g_{q-1,q}\) which is an \(\fR_{q-1}\)-morphism only, and the homotopy
operator \(h_0(\alpha) = (\alpha - \alpha(x_q=0)/x_q)\) which is an
\(\fR_{q-1}\)-morphism. \QED

\textsc{Proof of Theorem continued}. Thanks to the reduction \(\rho_{q}\), the
object \(\Ksz^{q}(\fR)\) is ``above'' \(\fR_{q}\). The morphism \(\times x_q\)
is trivially lifted into a chain-complex morphism: \(\times x_q: \Ksz^{q}(\fR)
\leftarrow \Ksz^{q}(\fR)\); the source and the target of this morphism are
reduced through \(\rho_{q}\) over \(\fR_ {q}\) and we can apply the Cone
Reduction Theorem~\ref{47942}. Combining with the other reduction already
available, we obtain:
\[
\Cone(\Ksz^{q}(\fR) \stackrel{\times x_q}{\longleftarrow} \Ksz^{q}(\fR)) \rrdc
\Cone(\fR_{q,\ast} \stackrel{\times x_q}{\longleftarrow} \fR_{q,\ast}) \rrdc
\fR_{q-1}
\]
where the \(\fR_{q,\ast}\) terms are understood as chain-complexes concentrated
in degree 0. Composing both reductions (Proposition~\ref{28461}) gives the
result if we can identify the first cone with \(\Ksz^{q-1}(\fR)\). This cone is
made of two copies of \(\Ksz^{q}(\fR)\); to distinguish them, let us recall the
right hand one \(dx_q.\Ksz^{q}(\fR)\), that is, for every term of this
\(\Ksz^{q}(\fR)\), let us put a symbol \(dx_q\) between the coefficient in
\(\fR\) and the exterior part in \(\wedge V_{q}\). This increases the Koszul
degree in the chain-complex by +1, but by chance the right hand term in a cone
is suspended. When you compute the differential of \(\alpha \, dx_q. \ldots\)
in a Koszul complex, the contribution of \(dx_q\) corresponds here to our
\(\times x_q\) morphism, the other terms come from the differential of
\(\Ksz^{q}(\fR)\). In fact with another sign, but the sign of the differential
in the right hand component of a cone is also changed. Conclusion: there is a
natural canonical isomorphism of chain-complex \(\Cone(\Ksz^{q}(\fR)
\stackrel{\times x_q}{\longleftarrow} \Ksz^{q}(\fR)) \cong \Ksz^{q-1}(\fR)\).
\QED

A novice can be troubled by the following observation: more \(\fR_q\) is small,
more \(\Ksz^q(\fR)\) is big? The point is that if \(\fR_q\) is smaller, then
the ``difference'' between the ground ring \(\fR\) and \(\fR_q\) is bigger, so
that the resolution is logically more complicated. The proof start from \(\fR\)
and goes up to \(\fk\) through the various~\(\fR_q\).

From a computational point of view, it is important to make explicit the
homotopy component \(h\) of the reduction \(\Ksz(\fR) \rrdc \fk\). Using the
detailed formula given when proving the Cone Reduction Theorem~\ref{47942}, it
is easy to prove our homotopy operator is given by the formula:
\[\label{79290}
h(\alpha.\lambda) = \sum_{q=1}^m ((\alpha(x_1, \ldots, x_q, 0, \ldots, 0) -
\alpha(x_1, \ldots, x_{q-1}, 0, \ldots, 0)/x_q) \ dx_q . \lambda
\]
if \(\alpha \in \fR\) and \(\lambda \in \wedge V\) with the common
interpretations inside the exterior algebra \(\wedge V\): if ever \(dx_q\) is
present in \(\lambda\), then \(dx_q . \lambda = 0\); and if \(dx_q\) is at a
wrong place, putting it at the right place can need a sign change. It is
amusing to study the particular case where \(\alpha\) is a monomial \(\alpha =
x_{i_1}^{j_1} \cdots x_{i_k}^{j_k}\) with \(i_1 < \cdots < i_k\) and \(j_k >
0\). If \(dx_{i_k}\) is present in \(\lambda\), the result is 0; otherwise you
replace \(j_k\) by \(j_k-1\) in the monomial and you insert a \(dx_{i_k}\) in
\(\lambda\) at the right place with the right sign. In concrete programming,
this can be run very efficiently. We will see this algorithm is by far the most
used in the resulting programs. Because of the proverb: the difference between
effective homology and ordinary homology consists in using the \emph{explicit}
homotopy operators.

\subsection{Koszul complex with effective homology}\label{91102}

The reduction constructed in the previous section can be understood as
describing the \emph{effective homology} of our Koszul complex.

\begin{thr}\label{64499}---
If \(\fR\) is the ring \(\fR = \fk[x_1, \ldots, x_m]_0\), the Koszul complex
\(\emph{\Ksz}(\fR)\) ``is'' an object with effective homology.
\end{thr}

\proof As usual, strictly speaking, the statement is improper: the statement
part ``is an object \ldots'' is a shorthand; in fact it is claimed some process
allows us to complete the object under study, the Koszul complex, as a
quadruple satisfying the required rules of Definition~\ref{93815}. This
quadruple is \((\Ksz(\fR), \Ksz(\fR), \fk_\ast, \rho)\) where~\(\rho\) is the
reduction of the previous section considered as the equivalence:
\[
\Ksz(\fR) \stackrel{=}{\lrdc} \Ksz(\fR) \stackrel{\rho}{\rrdc} \fk_\ast.
\]\QED

\subsection{Torsion groups.}

\begin{dfn} ---
\emph{Let \(\fR = \fk[x_1, \ldots, x_m]_0\) be our traditional ring and let
\(M\) and \(N\) be two \(\fR\)-modules. The torsion groups \(\Tor^\fR_i(M, N)\)
are defined as follows. Let \(\Rsl(M)\) and \(\Rsl(N)\) be two (free)
\(\fR\)-resolutions of \(M\) and \(N\). Then:
\[
H_\ast(\Rsl_\fR(M) \otimes_\fR N) =: \Tor^\fR_\ast(M, N) := H_\ast(M
\otimes_\fR \Rsl_\fR(N)).
\]
}
\end{dfn}

It is not obvious the definition is \emph{coherent}, that is, the result does
not depend on the choice of using a resolution of \(M\) or \(N\), and does not
depend either on the choice of the resolution itself. The usual argument uses
the bicomplex spectral sequence, and it is a good opportunity to introduce the
effective version of this spectral sequence.

\begin{dfn}---
\emph{A first quadrant \emph{bicomplex} is a diagram of modules:
\[
\xymatrix{
 &  \ar[d] &  \ar[d] &
 \\
 \cdots & C_{p-1,q} \ar[l]\ar[d]_{d''} & C_{p,q} \ar[l]_-{d'}\ar[d]_{d''} &
 \cdots \ar[l]
 \\
 \cdots & C_{p-1,q-1} \ar[l]\ar[d] & C_{p,q-1} \ar[l]_-{d'}\ar[d] & \cdots \ar[l]
 \\
 & & &
 }
\]
with \(C_{p,q} = 0\) if \(p\) or \(q < 0\). Furthermore every horizontal is a
chain-complex (\(d'd' = 0\)), every vertical is a chain-complex (\(d''d'' =
0\)), and every square is \emph{anti}-commutative: \(d'd'' + d''d' = 0\). The
totalization of this bicomplex is a \emph{simple} chain-complex \((T_n, d_n)\)
where \(T_n = \oplus_{p+q=n} C_{p,q}\) and the differential \(dc\) of a chain
\(c \in C_{p,q} \subset T_{p+q}\) is \(dc = d'c \oplus d''c \in C_{p-1,q}
\oplus C_{p,q-1} \subset T_{p+q-1}\).}
\end{dfn}

The relations required for \(d'\) and \(d''\) are exactly the necessary
relations which do make the totalization a chain-complex. The bicomplex
spectral sequence gives a relation between the homology of every column (for
example) and the homology of the totalization. Other similar definitions can be
given for other quadrants or for the whole \((p,q)\)-plane.

\begin{thr}\textbf{\emph{(Bicomplex Spectral Sequence)}}---
If \((C_{p,q}, d'_{p,q}, d''_{p,q})\) is a first quadrant bicomplex, a spectral
sequence \((E^r_{p,q}, d^r_{p,q})\) can be defined with \(E^0_{p,q} =
C_{p,q}\), \(d^0_{p,q} = d''_{p,q}\), and \(E^1_{p,q} = H''_{p,q}\) is the
``vertical'' homology group of the \(p\)-column at index \(q\). Furthermore
this spectral sequence converges to the homology of the totalization:
\[
E^r_{p,q} \Rightarrow H_{p+q}(T_\ast).
\]
\end{thr}

Of course you can exchange the role of rows and columns and obtain
\emph{another} spectral sequence where \(E^1_{p,q} = H'_{p,q}\) is this time
the homology of the \(q\)-row at the index \(p\), converging exactly toward the
same homology groups \(H_\ast(T_\ast)\). For the proof, see for
example~\cite[Section XI.6]{MCLN2} where the \(E^2_{p,q}\) are also computed,
quite elementary. Which is a little harder is \(E^r_{p,q}\) for \(r > 2\), a
problem which gets \emph{constructive} answers with our constructive methods.

\begin{dfn}\label{12080}---
\emph{A first quadrant multicomplex \((C_{p,q}, d^r_{p,q})\)} is a collection
of \((p,q)\)-modules as for a bicomplex, but a large collection of
\((r,p,q)\)-arrows defined for \(r \geq 0\) and every \((p,q)\); the
\(r\)-index describes the horizontal shift: \(d^r_{p,q}: C_{p,q} \rightarrow
C_{p-r,q+r-1}\). Furthermore, the natural totalization of these data must be a
chain-complex.
\end{dfn}

There are striking analogies with spectral sequences but also important
differences. The most important difference is the following: the parameter
\(r\) is nomore a ``time'' parameter, that is, all the arrows \(d^r_{p,q}\)
\emph{coexist at the same time}, and anyway no time in this definition! The
modules \(C_{p,q}\) do not depend on \(r\) like the \(E^r_{p,q}\) of a spectral
sequence.

To define the totalization of a multicomplex, the process is as follows. As for
a bicomplex, \(T_n = \oplus_{p+q=n} C_{p,q}\). A component of the differential
starts from every \(C_{p,q}\) and goes to every \(C_{p-r, q+r-1}\) for \(r \geq
0\). A formula can be written:
\[
(d: T_n \rightarrow T_{n-1}) = \oplus_{p+q=n} (\oplus_{r \geq 0} d^r_{p,q})
\]
where the first \(\oplus\) takes account of the expression of the \emph{source}
as a direct sum and the second \(\oplus\) is analogous for the \emph{target}.
The relevant explicative diagram maybe is this one.
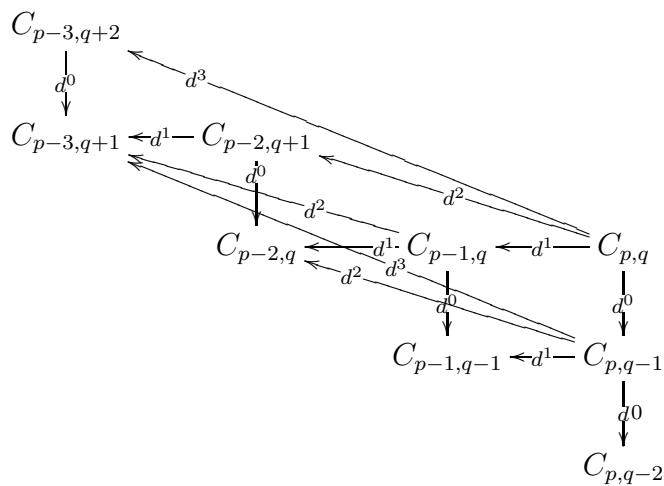
\[
\xymatrix{
 C_{p-3,q+2} \ar[d]|-{d^0} & & &
 \\
 C_{p-3,q+1} & C_{p-2,q+1} \ar[l]|-{d^1} \ar[d]|<(0.15){d^0}
 \\
 & C_{p-2,q} & C_{p-1,q} \ar[d]|-{d^0} \ar[l]|<(0.1){d^1} \ar[ull]|<(0.25){d^2} & C_{p,q}
 \ar[d]|-{d^0} \ar[l]|-{d^1} \ar[ull]|-{d^2} \ar[uulll]|<(0.75){d^3}
 \\
 & & C_{p-1,q-1} & C_{p,q-1} \ar[d]|-d^0 \ar[l]|-{d^1} \ar[ull]|<(0.7){d^2} \ar[uulll]|<(0.35){d^3}
 \\
 & & & C_{p,q-2}
}
\]
This diagram intends to study what happens when starting from \(C_{p,q}\).
Observe there is a unique way in our treasure hunt diagram starting from
\(C_{p,q}\) and arriving at \(C_{p,q-2}\). This implies necessarily \(d^0 d^0 =
0\) and therefore \emph{the columns are chain-complexes}. There are two ways
reaching \(C_{p-1,q-1}\) and we find again the anticommutativity property of
the squares. But now there are three ways leading to \(C_{p-2,q}\) and the
totalisation will actually be a differential only if \(d^0d^2 + d^1d^1 + d^2d^0
= 0\). And so on. In general \(\sum_{i = 0}^k d^{k-i}_{p-i,q+i-1} d^{i}_{p,q} =
0\) is required for every \((p,q)\) and every \(k \geq 0\).

\begin{thr}\label{61564} \textbf{\emph{(Bicomplex Reduction Theorem)}} ---
Let \((C_{p,q}, d'_{p,q}, d''_{p,q})\) be a bicomplex and \(T_\ast\) be its
totalization. Let \(\rho_p = (f_p, g_p, h_p): C_{p,\ast} \rrdc D_{p,\ast}\) be
a reduction of the (p)-column \emph{given} for every \(p\). Then a
multi-complex \((D_{p,q}, d^r_{p,q})\) can be defined with the following
property: let \(U_\ast\) be the totalization of this multicomplex; the
reductions \(\rho_p\) defines a ``total reduction'' \(\rho: T_\ast \rrdc
U_\ast\).
\end{thr}

Note the Cone Reduction Theorem~\ref{47942} is in fact a particular case: if
the bicomplex is null for \(p \geq 2\), which remains is simply the cone of the
columns~0 and~1; more precisely, in column~1, you must consider the
chain-complex with an opposite (vertical) differential; the morphism defining
the cone is given by the \(d'_{1,q}\) arrows.

We explain after the proof in which circumstance this theorem is mainly used.

\proof The proof is also a simple extension of the proof for the Cone Reduction
Theorem. You consider firstly the same bicomplex but with all the horizontal
differentials cancelled: \(d'_{\ast, \ast} = 0\). Then the \emph{different}
totalization, let us call it \((T'_\ast, d_{T'})\), is nothing but the direct
sum of the columns. The given reductions of the columns produce a reduction
\(\rho' = \oplus_p \rho_p: (T'_\ast, d_{T'}) \rrdc (U'_\ast, d_{U'})\) with
\(U'_\ast = \oplus_p D_{p,\ast}\). This being observed, let us reinstall now
the right horizontal arrows over \(T'_\ast\) to obtain again \(T_\ast\); this
can be viewed as a global \emph{perturbation} of the differential \(d_{T'}\) to
obtain the differential \(d_T\). Can we apply the BPL?

We must verify the nilpotency hypothesis. We must prove the composition \(h\hdl
=\) homotopy-perturbation is pointwise nilpotent. Let us start from
\(C_{p,q}\). The perturbation \(\hdl = d'_{p,q}\) in this case leads to
\(C_{p-1,q}\), the homotopy operator to \(C_{p-1,q+1}\). If we repeat, we go to
\(C_{p-2,q+2}\), and so on, and after \(p\) steps, we reach \(C_{0,p+q}\), but
here the perturbation is null and the nilpotency hypothesis is satisfied. The
role of the \emph{first quadrant} property then is clear: the snake path must
lead to some 0 module, whatever the starting point is.
\[
\xymatrix{
 \textrm{[nilpotency]} & C_{0,p+q} \ar[l]_-{\hdl = 0}
 \\
 & C_{0,p+q-1} \ar[u]_h & C_{1,p+q-1} \ar[l]_\hdl
 \\
 & & C_{1,p+q-2} \ar[u]_h & C_{2,p+q-2} \ar[l]_\hdl
 \\
 & & & \ar@{.>}[u]_h
}
\]

Examining the nilpotency hypothesis gives also a good idea about the nature of
the BPL series:
\[
\phi = \sum_{i=0}^\infty (-1)^i (h\widehat{\delta})^i ; \hspace{2cm} \psi =
\sum_{i=0}^\infty (-1)^i (\widehat{\delta} h)^i.
\]
They are the sums of all the terms obtained following our snake paths, starting
in the horizontal direction for \(\phi\), in the vertical direction for
\(\psi\). The next diagram shows the terms corresponding to \(i = 2\), called
\(\phi_2\), in dashed arrows, and \(\psi_2\), solid arrows.
\[
\xymatrix{
 C_{p-2,q+2} & C_{p-1,q+2} \ar[l]_-{\widehat{\delta}}
 \\
 C_{p-2,q+1} \ar@{-->}[u]^-h & C_{p-1,q+1} \ar[u]^-h \ar@{-->}[l]_-{\widehat{\delta}} &
 C_{p,q+1} \ar[l]_-{\widehat{\delta}}
 \\
 & C_{p-1,q} \ar@{-->}[u]_-h & C_{p,q} \ar[u]_-h
 \ar@{-->}[l]^-{\widehat{\delta}}
 \ar@{-->}@/^20pt/[uull]|{\phi_2} \ar@/_20pt/[uull]|<(0.4){\psi_2} }
\]
Then these series \(\phi\) and \(\psi\) have to be combined with the original
\(f\), \(g\), \(h\), \(d_{T'}\) and \(d_{U'}\) to produce the looked-for
reduction between \(T_\ast\), given, and \(U_\ast\), to be constructed. This is
the role of BPL. In particular, taking account of the formula \(\delta =
f\widehat{\delta} \phi g\) for the resulting perturbation on the small complex
\(U'_\ast\) transforming it into the perturbed small one \(U_\ast\), you see
the path to be followed to obtain a component of this differential \(d_U\). Let
\(D_{p,q}\) be a starting point. First you go up from \(D_{p,q}\) following
\(g_{p,q}\) arriving at \(C_{p,q}\). Then follow for example the \(\phi_2\)
path of the above figure going from \(C_{p,q}\) to \(C_{p-2,q+2}\). Then again
an arrow \(d'_{p-2,q+2}: C_{p-2,q+2} \rightarrow C_{p-3,q+2}\). And finally get
back in \(U_\ast\) by \(f_{p-3,q+2}: C_{p-3,q+2} \rightarrow D_{p-3,q+2}\).
This composition \(f_{p-3,q+2} d'_{p-2,q+2} \phi_2 g_{p,q}\) is the arrow
\(d^3_{p,q}: D_{p,q} \rightarrow D_{p-3,q+2}\) of the multicomplex \((D_{p,q},
d^r_{p,q})\).
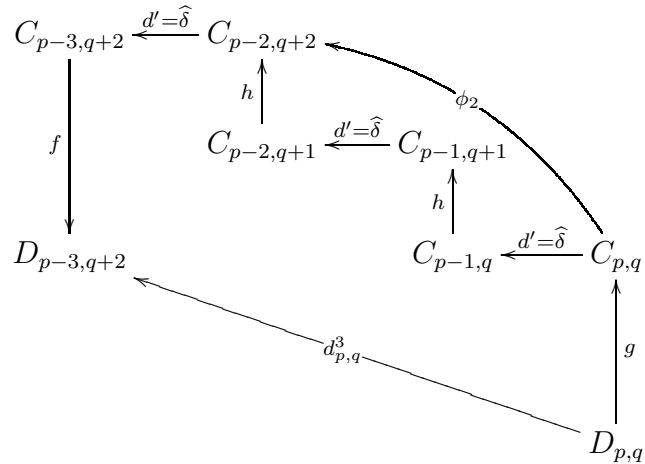
\[
\xymatrix{
 C_{p-3,q+2} \ar[dd]_f & C_{p-2,q+2} \ar[l]_-{d' = \widehat{\delta}}
 \\
 & C_{p-2,q+1} \ar[u]^-h & C_{p-1,q+1} \ar[l]_-{d' = \widehat{\delta}}
 \\
 D_{p-3,q+2} & & C_{p-1,q} \ar[u]^-h & C_{p,q} \ar[l]_-{d' = \widehat{\delta}}
 \ar@/_20pt/[uull]|{\phi_2}
 \\
 \\
 & & & D_{p,q} \ar[uu]_g \ar[uulll]|{d^3_{p,q}}
}
\]
Do the same for every shift \(r\) and you obtain the multicomplex \((U_\ast,
d^r_{p,q})\). The components of the final reduction \(T_\ast \rrdc U_\ast\) are
quite analogous: every component starting from a \(C_{pq}\) or a \(D_{p,q}\)
and going to another one is made of a snake path plus a few simple components
added at the departure and/or the arrival. \QED

In what context this theorem can be used? We will see several different
contexts where this theorem is a key tool. The simplest one is the following.
If ever the~\(D_{p,q}\) are \emph{effective}, the homology groups of \(U_\ast\)
are elementarily computable. We will meet many cases where the ``main''
bicomplex is only locally effective, so that the homology groups of the
totalization in general are not reachable. But frequently we can obtain for
example the \emph{effective homology} of every column. Our theorem will then
give us an equivalence between the initial totalisation and another one coming
from a multicomplex where the components are on the contrary \emph{effective}.
Then it will be possible to compute the homology groups.

Another classical use of the bicomplex spectral sequence theorem concerns the
case where \emph{every} column and row is ``almost'' exact. We will do the same
in our constructive framework, obtaining a \emph{constructive} result. Without
any spectral sequence.

\begin{thr}---
Let \((C_{p,q}, d'_{p,q}, d''_{p,q})\) be a first quadrant bicomplex satisfying
the following properties.
\begin{itemize}
\item
Every \((p)\)-column is exact except at \((p,0)\) producing a homology group
\(H''_{p,0}\).
\item
Every \((q)\)-row is exact except at \((0,q)\) producing a homology group
\(H'_{0,q}\).
\item We assume \emph{reductions} are available:
\[
(C_{p,\ast}, d'') \rrdc H''_{p,0} ; \hspace{1cm} (C_{\ast,q}, d') \rrdc
H'_{0,q}.
\]
\end{itemize}
Then an equivalence can be \emph{constructed}:
\[
(H''_{p,0}, d')_p \eqvl (H'_{0,q}, d'')_q.
\]
\end{thr}

The statement of the theorem needs a few explanations. In the column direction
for example, every colum is a chain-complex and requiring its exactness makes
sense. The vertical exactness is required in any position \((p,q)\) with \(q
> 0\). In position \((p,0)\), the arrow \(d''_{p,1}: C_{p,1} \rightarrow
C_{p,0}\) is not necessarily surjective, which defines the homology group
\(H''_{p,0} = C_{p,0}/ d''_{p,1}(C_{p,1})\). Now the horizontal arrow
\(d'_{p,0}\) induces a map \(d'_{p,0}: H''_{p,0} \rightarrow H''_{p-1,0}\) and
this produces a chain-complex \((H''_{p,0}, d')_p\). The same for the rows. The
classical result obtained in this case is that the homology groups of both
complexes \((H''_{p,0}, d')_p\) and \((H'_{0,q}, d'')_q\) are isomorphic. Here,
using the reductions of the statement, we \emph{construct} an equivalence
between these complexes, which of course implies the isomorphism between
homology groups.

\proof Let \(T_\ast\) be the totalisation of our bicomplex. Applying the
Bicomplex Reduction Theorem produces a reduction:  \(T_\ast \rrdc (H''_{p,0},
d')_p\). Doing the same with the rows finally gives:
\[
(H''_{p,0}, d')_p \lrdc T_\ast \rrdc (H'_{0,q}, d'')_q
\]
\QED

It is the first example where a natural \emph{equivalence} is obtained, instead
of a \emph{reduction}. This is frequent.

We are finally ready to prove the \emph{constructive} coherence of the Torsion
groups.

\begin{thr}---
Let \(\fR = \fk[x_1, \ldots, x_m]\) be the localized polynomial ring and  \(M\)
and~\(N\) two \(\fR\)-modules. Let \emph{\(\Rsl(M)\)} and \emph{\(\Rsl(N)\)} be
some \emph{effective} free resolutions. An explicit equivalence can be
installed between \emph{\(\Rsl_\ast(M) \otimes_\fR N\)} and \emph{\(M
\otimes_\fR \Rsl_\ast(N)\)}.
\end{thr}

\proof Reductions \(\Rsl(M) \rrdc M_\ast\) and \(\Rsl(N) \rrdc N_\ast\) are
available. We consider the bicomplex \(\Rsl(M) \otimes_\fR \Rsl(N)\). The
Koszul convention implies the totalization actually is a chain-complex. If we
examine the \((p)\)-column, the left factor \(\Rsl_p(M)\) of the tensor product
\(C_{p,q} = \Rsl_p(M) \otimes_\fR \Rsl_q(N)\) is independent of \(q\). The
\(\fR\)-module \(\Rsl_p(M)\) is free of rank \(r_p\) so that the tensor product
\(\Rsl_p(M) \otimes_\fR \Rsl_\ast(N)\) is nothing but the direct sum of \(r_p\)
copies of \(\Rsl_p(N)\). In particular the reduction \(\Rsl_\ast(N) \rrdc N\)
becomes a reduction \(\Rsl_p(M) \otimes \Rsl_\ast(N) \rrdc \Rsl_p(M) \otimes
N\): the homology of every \((p)\)-column is concentrated at \((p,0)\). We are
exactly in the situation of the previous theorem, obtaining an explicit
equivalence:
\[
\Rsl_\ast(M) \otimes_\fR N \lrdc T_\ast \rrdc M \otimes_\fR \Rsl_\ast(N).
\]
\QED

Note the equivalence depends for example on the chosen isomorphisms \(\Rsl_p(M)
\cong \fR^{r_p}\). The homotopy operators of \(\Rsl_\ast(N)\) in general
\emph{are not} \(\fR\)-morphisms and the expressions of the induced homotopy
operator over \(\Rsl_p(M) \otimes_\fR \Rsl_\ast(N)\) can so be modified: the
splitting into \(r_p\) components \emph{is not} intrinsic.

\section{Effective homology of Koszul complexes.}\label{64720}

\subsection{Presentation.}

The Koszul complexes play an important role when studying the \emph{formal
integrability problem} of PDE systems. The data in this case is a number of
(independent) variables \(m\), a ground field \(\fk = \bR\) or \(\bC\), the
ring \(\fR = \fk[x_1, \ldots, x_m]_0\) and an \(\fR\)-module of finite type
\(M\) coming from the PDE system. The nature of the PDE system then strongly
depends on the torsion groups \(\Tor_\ast(M, \fk)\)~\cite{GLSS}.

We explain in this Section how constructive homological algebra gives
completely new methods to study this problem. Usually the torsion groups are
computed as follows. First construct a finite free \(\fR\)-resolution of \(M\),
it is the classical Hilbert's syzygy problem. Efficient theoretical and
concrete methods are available, but they are rather technical; the Groebner
basis techniques are necessary. Then the tensor product \(\Rsl(M) \otimes_\fR
\fk\) is a finite chain-complex of finite dimensional \(\fk\)-vector spaces,
the homology groups of which can be elementarily computed.

But it happens the theoretical result at the origin of this computation comes
from the symmetric definition \(\Tor_\ast(M, \fk) = H_\ast(M \otimes_\fR
\Ksz(\fk)) = H_\ast(\Ksz(M))\). If these torsion groups are sufficiently null,
then the module \(M\) is \emph{involutive}, which expresses that ``good''
coordinate systems can be used to study the algebraic nature of \(M\). As
usual, the homological condition allows the user to claim there \emph{exists}
good coordinate systems. Making \emph{constructive} such a statement is a
natural goal; if such constructive results are obtained, we can reasonably hope
to be able to \emph{concretely} use the nice results of~\cite{GLSS}.

To conveniently explain how our constructive methods can be used, we choose a
framework a little simpler; the translation in the general framework is very
easy. This framework is also chosen to allow us to give simple machine
demonstrations with the current available Kenzo programs.

\begin{uos}---
In this section, the ground field \(\fk\) is an arbitrary commutative field; in
the Kenzo demonstrations, \(\fk = \bQ\). The ring \(\fR\) is as before
\(\fk[x_1, \ldots, x_m]_0\). Instead of an \(\fR\)-module \(M\), we consider an
ideal \(I = \ideal{g_1, \ldots, g_n} \subset \fR\) and the corresponding module
\(M = \fR/I\). We intend to construct a version \emph{with effective homology}
of \(\Ksz(M) = \Ksz(\fR/I)\).
\end{uos}

The Groebner methods will play also an essential role, but with a completely
different organization, significantly simpler and more conceptual from a
theoretical point of view, at least when the general style of
\emph{constructive} homological algebra is understood.

\subsection{Constructive homological algebra and short exact sequences of
chain-complexes.}

Theorem~\ref{15277} explains how a short exact sequence of chain-complexes
produces a long exact sequence of the corresponding homology groups. This exact
sequence is implicitly assumed solving computational problems when you know the
homology groups of two chain-complexes and you want to obtain the homology
groups of the third one. Section~\ref{72185} was devoted to elementary positive
examples, but we saw later, Section~\ref{62859}, that in general exact
sequences lead to \emph{extension problems} which can be really difficult.

This section will replace the long exact sequence of a short exact sequence of
chain-complexes by simple \emph{constructive} results, which systematically
avoid this difficulty. We will see why constructive homological algebra is also
in particular a \emph{general solving method} for extension problems.

We work in this subsection in a \emph{quite general framework}. The ground ring
\(\fR\) is an arbitrary unitary commutative ring, and in fact its
multiplicative structure is never used, it could be simply an Abelian group. No
ground field is concerned, except if \(\fR\) is itself a field\ldots

We begin with an easy extension of the Cone Reduction Theorem~\ref{47942}.

\begin{thr} \textbf{\emph{(Cone Equivalence Theorem)}} ---
Let \(\phi: C_{\ast, EH} \leftarrow C'_{\ast, EH}\) be a chain-complex morphism
between two chain-complexes \emph{with effective homology}. Then a general
algorithm computes a version \emph{with effective homology}
\(\emph{\Cone}(\phi)_{EH}\) of the cone.
\end{thr}

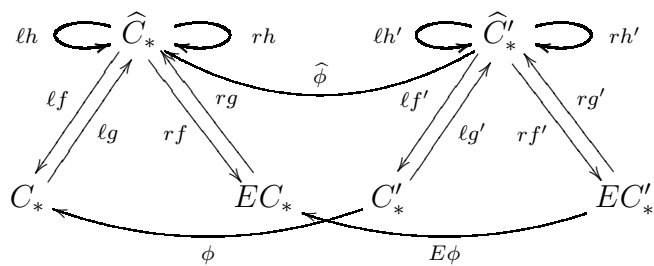
\[
\xymatrix @R150\ul @C80\ul {\scriptstyle \ell h \ar@{}[r]|{\autoarrow{1}} &
 {\widehat{C}_\ast} \ar@<-10\ul>[ld]_{\ell f} \ar@<-10\ul>[rd]_{rf} &
 \scriptstyle r h \ar@{}[l]|{\autoarrow{-1}} & \scriptstyle \ell h'
 \ar@{}[r]|{\autoarrow{1}} & {\widehat{C}'_\ast} \ar@<-10\ul>[ld]_{\ell f'}
 \ar@<-10\ul>[rd]_{rf'} \ar@/^80\ul/[lll]_{\widehat{\phi}} & \scriptstyle rh'
 \ar@{}[l]|{\autoarrow{-1}}
 \\
 {C_\ast} \ar@<-10\ul>[ru]_{\ell g} & & {EC_\ast} \ar@<-10\ul>[lu]_{r g} &
 {C'_\ast} \ar@<-10\ul>[ru]_{\ell g'} \ar@/^50\ul/[lll]^\phi
 & & {EC'_\ast} \ar@<-10\ul>[lu]_{r g'} \ar@/^50\ul/[lll]^{E\phi}}
\]

\proof We start with two equivalences \(C_\ast \lrdc \widehat{C}_\ast \rrdc
EC_\ast\) and \(C'_\ast \lrdc \widehat{C}'_\ast \rrdc EC'_\ast\) and the
morphism \(\phi: C_\ast \leftarrow C'_\ast\). In the figure above, \(\ell\) =
left and \(r\) = right, this for each given equivalence. The morphism \(\phi\)
naturally induces ``parallel'' morphisms \(\widehat{\phi} := (\ell g) \, \phi
\, (\ell f'): \widehat{C}_\ast \leftarrow \widehat{C}'_\ast\) and then \(E\phi
:= (rf) \, (\ell g) \, \phi \, (\ell f') \, (rg'): EC_\ast \leftarrow
EC'_\ast\).

As usual we can consider \(\phi\) is a perturbation of the differential of
\(\Cone(C_\ast \stackrel{0}{\leftarrow} C'_\ast)\). Using the Easy Perturbation
Lemma~\ref{33336} produces a reduction \(\Cone(\phi) \lrdc
\Cone(\widehat{\phi})\) where in fact the morphism \(\widehat{\phi}\) is
\emph{produced} by the lemma. Applying in the same way the Basic Perturbation
Lemma between \(\Cone(\widehat{C}_\ast \stackrel{0}{\leftarrow}
\widehat{C}'_\ast)\) and \(\Cone(\widehat{C}_\ast
\stackrel{\widehat{\phi}}{\leftarrow} \widehat{C}'_\ast)\) produces in turn a
new reduction \(\Cone(\widehat{\phi}) \rrdc \Cone(E\phi)\) where again, the
morphism \(E\phi\) is in fact \emph{produced} by the BPL. Combining these
reductions gives the looked-for equivalence:
\[
{%
\Cone(C_\ast \stackrel{\phi}{\leftarrow} C'_\ast) \lrdc \Cone(\widehat{C}_\ast
\stackrel{\widehat{\phi}}{\leftarrow} \widehat{C}'_\ast) \rrdc \Cone(EC_\ast
\stackrel{E\phi}{\leftarrow} EC'_\ast)}
\]\QED

You see how our perturbation lemmas are used. Some process is applied to the
left hand term of an equivalence, here the cone construction. This process
induces something analogous over the central chain complex of the given
equivalence, thanks to the easy perturbation lemma. The left hand reduction is
not here modified, this reduction is only used to \emph{copy} the perturbation
into the central chain complex. Then the actual Basic Perturbation Lemma is
applied to take account of the perturbation in the central chain complex to
replace the right hand reduction by a new appropriate reduction; in general the
differential of the right hand chain complex is modified.

\begin{dfn}\label{87354}---
\emph{An \emph{effective short exact sequence} of chain-complexes is a diagram:
\[
 \xymatrix @C100\ul @1{ *+++{0} & *+++{A_\ast} \ar[l]_-0 \ar@<10\ul>[r]^-\sigma &
 *+++{B_\ast} \ar@<10\ul>[r]^-\rho \ar@<10\ul>[l]^-j & *+++{C_\ast}
 \ar@<10\ul>[l]^-i & *+++{0} \ar[l] }
\]
 where \(i\) and \(j\) are \emph{chain-complex} morphisms, \(\rho\)
(retraction) and \(\sigma\) (section) are \emph{graded module} morphisms
satisfying:
\begin{itemize}
\item \(\rho i = \id{C_\ast}\);
\item \(i \rho + \sigma j = \id{B_\ast}\);
\item \(j \sigma = \id{A_\ast}\).
\end{itemize}
 }
\end{dfn}

It is an exact sequence in both directions, but to the left it is an exact
sequence of \emph{chain-complexes}, the exact sequence we are mainly interested
in, and to the right it is only an exact sequence of \emph{graded modules}, no
compatibility in general with the differentials. The components \(\rho\) and
\(\sigma\) are nothing but a homotopy operator describing a \emph{reduction} to
0 of our ``total'' chain-complex: you can think of this exact sequence as a
bicomplex with only three columns non-null. As usual for the homotopy
operators, weak properties are only required, and here for example it is not
required \(\rho\) and \(\sigma\) are compatible with the differentials.
Otherwise the chain-complex \(B_\ast\), \emph{differential included}, would be
the direct sum of \(A_\ast\) and \(C_\ast\), a trivial situation without any
interest. The exactness expresses \(i\) is injective, \(j\) is surjective, and
\(\rho\) and \(\sigma\) define a sum decomposition \(B_\ast = \im i \oplus \ker
\rho = \im \sigma \oplus \ker j\), but this decomposition \emph{is not} in
general a subcomplex decomposition, making the hoped-for results non-trivial.

\begin{thr}\label{52775} \textbf{\emph{(SES Theorems)}} ---
Let
\[
 \xymatrix @C100\ul @1{ *+++{0} & *+++{A_\ast} \ar[l]_-0 \ar@<10\ul>[r]^-\sigma &
 *+++{B_\ast} \ar@<10\ul>[r]^-\rho \ar@<10\ul>[l]^-j & *+++{C_\ast}
 \ar@<10\ul>[l]^-i & *+++{0} \ar[l] }
\]
be an \emph{effective} short exact sequence of chain-complexes. Then three
general algorithms are available:
\[
\begin{array}{c}
\textrm{\emph{SES}}_1: (B_{\ast,EH}, C_{\ast,EH}) \mapsto A_{\ast,EH}
\\
\textrm{\emph{SES}}_2: (A_{\ast,EH}, C_{\ast,EH}) \mapsto B_{\ast,EH}
\\
\textrm{\emph{SES}}_3: (A_{\ast,EH}, B_{\ast,EH}) \mapsto C_{\ast,EH}
\end{array}
\]
producing a version \emph{with effective homology} of one chain-complex when
versions \emph{with effective homology} of both others are given.
\end{thr}

SES = Short Exact Sequence. Observe the process is perfectly \emph{stable}: the
\emph{type} of the result is exactly the same as for the given objects. The
obtained object can then be used later in another exact or spectral sequence,
and so on.

\proof

Let us begin with the SES\(_1\) case.

\begin{lmm}\label{47522} ---
The effective exact sequence produces a reduction: \(\Cone(i) \rrdc A_\ast\).
\end{lmm}
\proof It is again a simple application of BPL. We mentioned, when describing
the notion of \emph{effective} short exact sequence, that \(\rho\)  and
\(\sigma\) are ``weak'' morphisms. This negative property is no longer an
obstacle if we cancel the differentials of our three chain-complexes. Let us
call \(A^0_\ast, B^0_\ast, C^0_\ast\) these chain-complexes with null
differentials. It is easy to obtain the looked-for reduction in this simple
case. It is:
\[
\rho^0 = (f^0, g^0, h^0) : \Cone(i: B^0_\ast \leftarrow C^0_\ast) \rrdc
A^0_\ast
\]
The morphism \(f^0: \Cone(i) \rightarrow A^0_\ast\) is the projection defined
by \(j: A^0_\ast \leftarrow B^0_\ast\), null on the \(C^0_\ast\) component of
the cone. The morphism \(g^0\) is defined by the section \(\sigma\) with values
in the \(B^0_\ast\) component of the cone. Finally the homotopy operator
\(h^0\) is the retraction \(\rho: B^0_\ast \rightarrow C^0_\ast\) inside the
cone. The reduction properties are direct consequence of the relations
satisfied by \(i, j, \rho\) and \(\sigma\). Note the \emph{components} of our
new cone have null differentials, but the cone itself has the component \(i\)
non null except if \(C_\ast = 0\).

Now we reinstall the right differentials \emph{over the cone}. Two components
for the perturbation \(\hdl\), a differential in general non trivial over
\(B_\ast\) and another one over~\(C_\ast\). Combined with the initial homotopy
operator of our reduction, we see \((h^0 \hdl)^2\) is null. The nilpotency
condition is satisfied.

Using Shih's formula for the new reduction, we obtain the reduction:
\[
\rho = (f, g, h) : \Cone(i: B_\ast \leftarrow C_\ast) \rrdc A_\ast
\]
with \(f = f^0 = j\) and \(h = h^0 = \rho\) not modified, but with \(g = \sigma
- \rho d_{B_\ast} \sigma\). Furthermore, the new differential to install on the
small chain-complex is by chance the initial differential \(d_{A_\ast}\) of
\(A_\ast\). \QED

\noindent\textsc{Proof of Theorem continued.} Consider the sequence:
\[
A_\ast \lrdc \Cone(i) \lrdc \Cone(\widehat{i}) \rrdc \Cone(Ei).
\]
The central and the right hand reductions are produced by the Cone Equivalence
Theorem, using the available equivalences describing the chain-complexes
\(B_\ast\) and \(C_\ast\) as chain-complexes with effective homology. The left
hand reduction is produced by the lemma just proved. Composing the central and
the left hand reductions gives another reduction and an equivalence between
\(A_\ast\) and \(\Cone(Ei)\) is obtained, describing also \(A_\ast\) as a
chain-complex with effective homology.

The case SES\(_3\) is symmetric and left to the reader.

Let us finally consider the case SES\(_2\), different.

\begin{lmm}---
The \emph{effective} short exact sequence generates a {connection}
\emph{chain-complex} morphism \(\chi: A_\ast \rightarrow C_\ast^{[1]}\).
\end{lmm}

The ``exponent '' \([1]\) explains the \emph{suspension functor} is applied to
the chain-complex \(C_\ast\): the degree of an element is increased by 1 and
the differential is replaced by the opposite.

\proof The connection morphism is defined as the composition \(\chi = \rho d
\sigma\) where the differential cannot be anything else than \(d = d_B\); this
differential has degree~\mbox{-1} and is the cause of the suspension. We must
verify the compatibility of this claimed \emph{chain-complex} morphism with the
differentials of \(A_\ast\) and \(C_\ast^{[1]}\).

Let us consider an element \(a \in A_n\), then its lifting \(\sigma a\) in
\(B_\ast\), and let us try to use \(d_B d_B = 0\) and also \(\sigma j + i \rho
= \id{}\). First:
\[
\begin{array}{lcr@{\hspace*{1cm}}l}
 d \sigma a &=& \sigma j d \sigma a + i \rho d \sigma a & (\sigma j + i \rho =
 \id{})
 \\
 &=& \sigma d a + i \rho d \sigma a & (j d = d j\ \textrm{and} j \sigma = \id{})
\end{array}
\]
Let us apply again \(d_B\):
\[
\begin{array}{rcl@{\hspace*{1cm}}l}
 0 &=& d \sigma d a + d i \rho d \sigma a
 \\
 &=& \sigma j d \sigma d a + i \rho d \sigma d a + \sigma j d i \rho d \sigma a
 + i \rho d i \rho d \sigma a & (\sigma j + i \rho = \id{})
 \\
 &=& 0 + i \rho d \sigma d a + 0 + i d \rho d \sigma a
\end{array}
\]
for \(jd = dj\), \(j \sigma = \id{}\), \(d d = 0\), \(j d = dj\) again
 and \(j i = 0\). The morphism \(i\) is injective, which implies:
 \[
 d (\rho d \sigma) a = - (\rho d \sigma) (d a).
 \]
 \QED

This looks a little magic, but in fact, as in ordinary magic, there is an
explanation. The central \(B_\ast\) is, \emph{as graded module}, the direct sum
of \(A_\ast\) and \(C_\ast\). If you think of an element of \(B_\ast\) as
having two components, one in \(A_\ast\) and the other one in \(C_\ast\), then
you obtain an expression of the differential of \(d_{B_\ast}\) as working in
\(A_\ast \oplus C_\ast\); the differential is a \(2 \times 2\) matrix of maps,
the component \(C_\ast \rightarrow A_{\ast - 1}\) being null because \(i\) and
\(j\) are compatible with the differentials and \(ji = 0\); the component
\(A_\ast \rightarrow C_\ast\) is our connection map. We so obtain a cone
diagram:
\[
 \xymatrix{ \cdots & A_{n-2} \ar[l]_-d \ar[dl]_-\chi & A_{n-1} \ar[l]_-d
 \ar[dl]_-\chi & A_n \ar[l]_-d \ar[dl]_-\chi & \cdots \ar[l]_-d \ar[dl]_-\chi
 \\
 \cdots & C_{n-2} \ar[l]_-d & C_{n-1} \ar[l]_-d & C_n \ar[l]_-d & \cdots
 \ar[l]_-d
 }
\]
The total differential of this diagram is null if and only if every
parallelogram is anticommutative.

\textsc{Proof of Theorem continued.} The previous diagram also explains in
fact~\(B_\ast\) is canonically isomorphic to \(\Cone(\chi)\). Using the Cone
Equivalence Theorem and the equivalences describing the effective homology of
\(A_\ast\) and \(C_\ast\), we obtain the looked-for equivalence between
\(B_\ast\) and an effective chain-complex. \QED

\subsection{Solution for monomial ideals.}\label{46887}

We come back to our goal in commutative algebra: computing the \emph{effective}
homology of \(\Ksz(\fR/I)\) for \(I\) an ideal of \(\fR = \fk[x_1, \ldots,
x_m]_0\). Our ring is Noetherian and \(I\) is described by a finite set of
generators \(I =\ \ideal{g_1, \ldots, g_n}\). Our work is decomposed in three
steps:
\begin{enumerate}
\item
Using a Groebner basis, we replace \(I\) by \(I'\) a \emph{monomial} ideal to
be considered as a \emph{good simple} approximation of \(I\);
\item
A recursive process over \(I'\) using a number of times the BPL gives a simple
solution for \(I'\);
\item
Applying again the BPL between \(I\) and \(I'\) will give the solution for the
ideal~\(I\).
\end{enumerate}

Step 1 is standard. You choose a coherent monomial order over \(\fR\), then a
reduced Groebner basis is canonically defined for our ideal \(I\). We assume
our expression \(I =\ \ideal{g_1, \ldots, g_n}\) just uses this Groebner basis.

The ideal \(I'\) is obtained by replacing every generator \(g_i\) by its
\emph{leading term} \(g_i'\): \(I' :=\ \ideal{g'_1, \ldots, g'_n}\). This
process is interesting for two reasons:
\begin{itemize}
\item
The monomial ideal \(I'\), because it is monomial, is more comfortable.
\item
Both ideals \(I\) and \(I'\) are ``close'' to each other: the \emph{graded
modules} \(\fR/I\) and \(\fR/I'\) are \emph{canonically isomorphic}.
\end{itemize}
Of course the multiplicative structures of \(\fR/I\) and \(\fR/I'\) are
different, but the isomorphism between the underlying \emph{graded modules}
will be enough when applying the BPL to process this difference.

\emph{In the rest of this section, we assume our ideal \(I\) is monomial: every
generator~\(g_i\) of \(I =\ \ideal{g_1, \ldots, g_n}\) is a monomial of \(\fR =
\fk[x_1, \ldots, x_m]_0\)}.

The recursive process then consists in obtaining the result for the simpler
ideal \(J =\ \ideal{g_2, \ldots, g_n}\), the generator \(g_1\) being removed.
What about the exact nature of the relation between \(I\) and \(J\)? We must
use the notion of \emph{quotient} of two ideals; the quotient \(I_1:I_2\) of
two ideals \(I_1\) and \(I_2\) is \((I_1 : I_2) := \{a \in \fR \st a I_2
\subset I_1\}\).

\begin{prp}\label{37176}---
An ideal \(I =\ \ideal{g_1, \ldots, g_n}\ \subset \fR\) produces an
\emph{effective} short exact sequence of \(\fR\)-modules:
\[
0 \leftarrow \frac{\fR}{\ideal{g_1, \ldots, g_n}} \leftarrow
\frac{\fR}{\ideal{g_2, \ldots, g_n}} \leftarrow \frac{\fR}{\ideal{g_2, \ldots,
g_n} : \ideal{g_1}} \leftarrow 0.
\]
\end{prp}


\proof Exercise. \QED

In this exercise, please observe the initial monomorphism \({\fR}/{\ideal{g_2,
\ldots, g_n}} \leftarrow {\fR}/({\ideal{g_2, \ldots, g_n} : \ideal{g_1}})\) is
defined by the \emph{multiplication by \(g_1\)}, while the terminal epimorphism
\({\fR}/{\ideal{g_1, \ldots, g_n}} \leftarrow {\fR}/{\ideal{g_2, \ldots,
g_n}}\) is the canonical projection; this remark will be important later. These
monomorphism and epimorphism are \(\fR\)-module morphisms. To make effective
the exact sequence, a section \(\sigma\) and a retraction \(\rho\) are needed,
see Definition~\ref{87354}. The section is the ``brute'' lifting which sends a
non-null monomial -- in fact its class modulo the ideal -- to the same. The
retraction examines whether a monomial is divisible by \(g_1\); if yes the
retraction gives the quotient by \(g_1\), otherwise the result is null. These
section and retraction are \(\fk\)-linear but not at all \(\fR\)-morphisms.

\begin{prp}---
If the given generators of \(\ideal{g_1, \ldots, g_n}\) are \emph{monomials},
then \(\ideal{g_2, \ldots, g_n} : \ideal{g_1}\ =\ \ideal{g'_2, \ldots, g'_n}\)
with \(g'_i =\ \textrm{\emph{lcm}}(g_1, g_i)/g_1\) for \(i \geq 2\).
\end{prp}

\proof Exercise. \QED

Which implies if, thanks to the short exact sequence, a recursive process
reduces some work for \(\fR/\!\!\ideal{g_1, \ldots, g_n}\) to the analogous
work for \(\fR/\!\!\ideal{g_2, \ldots, g_n}\) and \(\fR/\!\!\ideal{g'_2,
\ldots, g'_n}\), there remains to \emph{start} the recursive process.

\begin{crl}\label{17184}---
A general algorithm computes:
\[
\left[\emph{\Ksz}\left(\frac{\fR}{\ideal{g_2, \dots, g_n}}\right)_{\!\!\!EH},
\emph{\Ksz}\left(\frac{\fR}{\ideal{g'_2, \dots, g'_n}}\right)_{\!\!\!EH}\right]
\mapsto \emph{\Ksz}\left(\frac{\fR}{\ideal{g_1, g_2, \dots,
g_n}}\right)_{\!\!\!EH}
\]
when the generators \(g_1, g_2, \ldots, g_n\) are monomials, when \(g'_i =\
\textrm{\emph{lcm}}(g_1, g_i)/g_1\), where \(\emph{\Ksz}(\cdots)_{EH}\) is a
version \emph{with effective homology} of the Koszul complex
\(\emph{\Ksz}(\cdots)\).
\end{crl}

\proof The constructor \(M \mapsto \Ksz(M)\) is a functor from \(\fR\)-modules
to chain-complexes. An \(\fR\)-module morphism \(f: M \rightarrow N\) generates
a chain-complex morphism \(f := \Ksz(f): \Ksz(M) \rightarrow \Ksz(N)\).
Applying this functor, the \emph{effective} short exact sequence of
\(\fR\)-\emph{modules} of Proposition~\ref{37176} becomes an \emph{effective}
short exact sequence of \emph{chain-complexes}:
\[
0 \leftarrow \Ksz\left(\frac{\fR}{\ideal{g_1, \dots}}\right)_{\!\!\!EH}
\leftarrow \Ksz\left(\frac{\fR}{\ideal{g_2, \dots}}\right)_{\!\!\!EH}
\leftarrow \Ksz\left(\frac{\fR}{\ideal{g'_2, \dots}}\right)_{\!\!\!EH}
\leftarrow 0
\]
Applying the SES\(_1\) case of Theorem~\ref{52775} gives the result. \QED

We noted in Proposition~\ref{37176} the section \(\sigma\) for example is only
\(k\)-linear; then \(\Ksz(\sigma)\) is defined but is not compatible with
differentials; it is only a graded-module morphism.

The recursive process is now installed: computing the effective homology of a
Koszul complex \(\Ksz(\fR/\!\!\ideal{g_1, \ldots, g_n})\) is reduced to two
analogous problems with one generator less. What about the starting point of
this induction? The minimal case is 0 generator; we must determine
\(\Ksz(\fR/\!\!<>)_{EH} = \Ksz(\fR)_{EH}\). This was done at
Theorem~\ref{64499}, which theorem was a translation of Theorem~\ref{37681}.
Combining this remark with the above corollary gives the main result of this
section.

\begin{thr}---
A general algorithm computes:
\[
\ideal{g_1, \ldots, g_n}\ \mapsto \emph{\Ksz}\left(\frac{\fR}{\ideal{g_1,
\ldots, g_n}}\right)_{\!\!\!EH}
\]
where \(g_1, \ldots, g_n\) are \emph{monomial} generators in our localized
polynomial ring \(\fR\). \QED
\end{thr}

The \emph{homological problem} for the chain-complex \(\Ksz(\fR/\!\!\ideal{g_1,
\ldots, g_n})\) is solved in the \emph{monomial case}. How to obtain the same
result in the general case?

\subsection{Installing a general multigrading.}

It was explained at the beginning of the previous section we intend to apply
again the BPL to process the difference between an arbitrary ideal \(I\) and
its monomial approximation \(I'\). The required nilpotency hypothesis needs a
careful use of monomial orders. Two ingredients are necessary.

  On one hand we must \emph{delocalize} the problem, replacing the localized ring
\(\fR = \fk[x_1, \ldots, x_m]_0\) by the ordinary polynomial ring \(\fRb =
\fk[x_1, \ldots, x_m]\). We will prove later that if \(I\) is an ideal of
\(\fR\), then \(H_\ast(\Ksz_\fR(\fR / I)) \cong H_\ast(\Ksz_{\fRb}(\fRb/\Ib))\)
if \(\Ib = I \cap \fRb\), so that instead of studying the problem of \(I\)
inside \(\fR\), we can study the case of \(\Ib\) and \(\fRb\); furthermore this
isomorphism between different homology groups will be \emph{constructive}, and
a solution for the homological problem of \(\Ksz_\fR(\fR / I)\) is equivalent
to a solution for \(\Ksz_{\fRb}(\fRb/\Ib)\): we can get rid of the
denominators.

\begin{dfn}---
\emph{An \(\fRb\)-ideal \(\Ib \subset \fRb\) is \emph{localized at} \(0 \in
\fk^m\) if \((\Ib \fR) \cap \fRb = \Ib\).}
\end{dfn}

In this definition \(\Ib \fR\) is the \(\fR\)-ideal generated by \(\Ib \subset
\fRb \subset \fR\). If \(I\) is an ideal of the local ring \(\fR\), then \(\Ib
= I \cap \fRb\) is an \(\fRb\)-ideal localized at 0, and all those ideals are
obtained in this way. The inclusion \(\Ib \subset (\Ib \fR) \cap \fRb\) is
always satisfied, but the ideal \(\Ib =\ \ideal{1-x}\ \subset \fk[x]\) is not
localized at 0, for \((\Ib \fR) \cap \fRb = \fRb \neq \Ib\).

On the other hand, in order to be able to use the Groebner techniques in our
context, we must define and handle carefully multigradings and monomial orders.
Once for all, we choose a \emph{Groebner monomial order}. If \(x_1^{\alpha_1}
\cdots x_m^{\alpha_m}\) is a monomial, its \emph{multigrading} is an
\(m\)-tuple \(\mu(x_1^{\alpha_1} \cdots x_m^{\alpha_m}) := [\alpha_1, \ldots,
\alpha_m]\). We consider an ideal \(\Ib \subset \fRb\) defined by a reduced
Groebner basis \(\ideal{g_1, \ldots, g_n}\) for the chosen monomial order. The
leading term of \(g_i\) is \(g'_i\), a monomial, and a canonical \(\fk\)-vector
space isomorphism is defined between \(\fRb/\Ib\) and \(\fRb/\Ib'\) if \(\Ib'
=\ \mbox{\(\ideal{g'_1, \ldots, g'_n}\)} \).

The ideal \(\Ib'\) is monomial and \(\fRb / \Ib'\) is multigraded. The Koszul
complex \(\Ksz(\fRb / \Ib')\) is also multigraded if we decide:
\[
\mu(x_1^{\alpha_1} \cdots x_m^{\alpha_m} \, dx_1^{\beta_1} \cdots
dx_m^{\beta_m}) := [\alpha_1 + \beta_1, \cdots, \alpha_m + \beta_m].
\]
where \(\alpha_i \in \bN\) and \(\beta_i \in \{0,1\}\). In particular the
differential of the Koszul complex is multigraded: a differential is made of
terms where a \(dx_i\) is replaced by a \(x_i\), which does not change the
multigrading.

The work of the previous section for Koszul complexes of monomial ideals can be
repeated without any change for the case of \(\fRb\) and \(\Ib'\) in the
present section instead of \(\fR\) and \(I\) in the previous section. In
particular we must use the initial reduction \((f,g,h): \Ksz(\fRb) \rrdc
\fk_\ast\) defined exactly in the same way. Note the three components of the
reduction are also multigraded: the components \(f\) and~\(g\) are trivial
except for elements of null multigrading; and the homotopy operator~\(h\), see
its detailed construction at page~\pageref{79290}, does the contrary of the
differential: every term is obtained by replacing some \(x_i\) by the
corresponding \(dx_i\). Using an obvious terminology, we can state:

\begin{prp}---
The reduction \(\emph{\Ksz}(\fRb) \rrdc \fk_\ast\) constructed as in
Section~\ref{45205} is multigraded.
\end{prp}

Note in particular taking or not the denominators \emph{does not}
``significantly'' change the effective homology of the Koszul complex of the
ground ring.

For a monomial ideal \(\Ib'\), we must apply a few times the Cone Reduction
Theorem~\ref{47942} to compute the effective homology of the corresponding
module, as explained in the previous section. We need multigraded versions of
the Basic Perturbation lemma and its applications, in particular for the Cone
Reduction Theorem.

\begin{thr} ---
If the data \(\rho = (f,g,h): (\hC_\ast, \hd) \rrdc (C_\ast, d)\) and \(\hdl:
\hC_\ast \rightarrow \hC_{\ast-1}\) of the Basic Perturbation Lemma~\ref{07404}
are multigraded, the resulting reduction \(\rho' = (f', g', h'): (\hC_\ast, \hd
+ \hdl) \rrdc (C_\ast, d + \delta)\) is also multigraded.
\end{thr}

\proof In this statement, the underlying chain-complexes are multigraded, and
the various \emph{given} operators respect the multigrading. The theorem
asserts the same for the new reduction. Given the explicit formulas for the
components of the new reduction, the proof is obvious. \QED

In the same way, if a morphism \(\phi: C_\ast \leftarrow C'_\ast\) is a
multigraded morphism between multigraded chain-complexes, and if the effective
homology of both complexes is given and is multigraded too, then the effective
homology of \(\Cone(\phi)\) computed by Theorem~\ref{47942} is also
multigraded: the three components of both reductions describing the effective
homology of the cone are multigraded, their source and target as well.

Remember the main step when computing the effective homology of \(\Ksz(\fRb /
\Ib')\) consists in using the \emph{effective} short exact sequence of
\(\fRb\)-modules:
\[
0 \leftarrow \frac{\fRb}{\ideal{g'_1, \ldots, g'_n}}
\stackrel{\textrm{\scriptsize{pr}}}{\longleftarrow}\frac{\fRb}{\ideal{g'_2,
\ldots, g'_n}} \stackrel{\times g'_1}{\longleftarrow} \frac{\fRb}{\ideal{g'_2,
\ldots, g'_n} : \ideal{g'_1}} \leftarrow 0.
\]
The generators are monomials, so that the epimorphism `pr', the canonical
projection, is multigraded. The monomorphism `\(\times g'_1\)' is the
multiplication by a monomial, it is also multigraded if you \emph{shift the
multigrading} of the initial module \({\fRb} / ({\ideal{g'_2, \ldots, g'_n} :
\ideal{g'_1}})\) by \(\mu(g'_1)\).

Starting from the \emph{multigraded} effective homology \(\Ksz(\fRb) \rrdc
\fk\), applying repetitively this process produces a version with
\emph{multigraded} effective homology of the Koszul complex of our monomial
module:
\[
\Ksz(\fRb / \Ib') \lrdc \hC_\ast \rrdc EC_\ast
\]
where the three modules are multigraded, and the six reduction components as
well. The following theorem is proved.

\begin{thr}\label{18292}---
An algorithm computes:
\[
\Ib' \mapsto [\Ksz(\fRb/\Ib') \lrdc \hC_\ast \rrdc EC_\ast]
\]
where \(\Ib'\) is a monomial ideal of \(\fRb = \fk[x_1, \ldots, x_m]\) and the
result is a multigraded equivalence between the corresponding Koszul complex
and an effective multigraded chain-complex of finite-dimensional \(\fk\)-vector
spaces.\QED
\end{thr}

\textsc{Example.} We consider the toy example \(\Ib' =\ \ideal{x^2, y^3}\
\subset \fRb = \bQ[x,y]\); effective homology must in particular compute
homology groups for the Koszul complex with representants for homology classes.
The recursive process will lead to consider \(\Ib'_0 =\ <>\), then \(\Ib'_1 =\
\ideal{y^3}\) and finally \(\Ib'\).

The result for \(\Ib'_0\) was obtained at Theorem~\ref{64499}. The Koszul
complex \(\Ksz_\ast(\fRb)\) is an \(\fRb\)-resolution of \(\fk = \bQ\) and its
effective homology is a diagram:
\[
\Ksz_\ast(\bQ[x,y]) \stackrel{=}{\lrdc} \Ksz_\ast(\bQ[x,y])
\stackrel{\rho}{\rrdc} \bQ_\ast.
\]
with \(\bQ_\ast\) the chain-complex with only \(\bQ\) in degree 0. Only one
homology group, in degree 0, isomorphic to \(\bQ\); the representant of a
generator is obtained by taking the image of the generator of \(\bQ_\ast\) in
\(\Ksz_\ast(\bQ[x,y])\), it is the ``base point'' of this Koszul complex,
namely \(1 \in \Ksz_0(\bQ[x,y])\).

The next figure extracts the important parts of the effective homology of
\(\Ksz_\ast(\bQ[x,y]/\!\!\ideal{y^3})\) when looking for a representant of the
generator of the homology in degree 1.
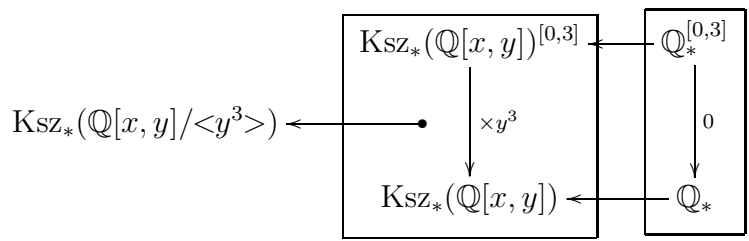
\[
\xymatrix@R10pt{
 & \Ksz_\ast(\bQ[x,y])^{[0,3]} \ar[dd]^-{\times y^3} & \bQ_\ast^{[0,3]} \ar[l]
 \ar[dd]^-0
 \\
 \Ksz_\ast(\bQ[x,y]/\!\!\ideal{y^3}) & {\phantom{\rule{30pt}{10pt}}} \ar@{{*}->}[l] &
 \\
 & \Ksz_\ast(\bQ[x,y]) & \bQ_\ast \ar[l]
 \save "1,2"."3,2"*+[F]\frm{} \restore
 \save "1,3"."3,3"*+[F]\frm{} \restore
}
\]
The boxes are cone chain-complexes produced by Corollary~\ref{17184}. The null
morphism between both copies of \(\bQ_\ast\) is the image of `\(\times y^3\)'
between both copies of \(\Ksz_\ast(\bQ[x,y])\). The right hand box is the
\emph{effective} chain-complex describing the homology of
\(\Ksz_\ast(\bQ[x,y]/\!\!\ideal{y^3})\). The central box settles the necessary
connection between the right hand box and
\(\Ksz_\ast(\bQ[x,y]/\!\!\ideal{y^3})\). The exponents \([0,3]\) show the
multigrading shift when necessary; in this way the `\(\times y^3\)' map is
multigraded. The next diagram displays at the right place elements of the nodes
of the previous diagram.
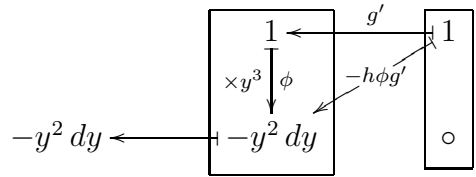
\[
\xymatrix@C40pt{
  & 1 \ar@{|->}[d]_-{\times y^3}^-{\phi} & 1 \ar@{|->}[l]_<(0.35){g'}
 \ar@{|->}[ld]|<(0.35){-h \phi g'}
  \\
  -y^2\,dy & -y^2\,dy  \ar@{|->}[l] & \circ
 \save "1,2"."2,2"*+[F]\frm{} \restore
 \save "1,3"."2,3"*+[F]\frm{} \restore
}
\]
The right hand 1 is the generator of the ``abstract'' homology in degree 1. Its
image in the intermediate box is obtained as explained in the cone reduction
diagram of page~\pageref{37174}: we have indicated in the present diagram the
relevant arrows labelled \(\phi\), \(g'\) and \(-h \phi g'\) of the generic
cone diagram. Note in particular the role of the contraction \(h\) when
obtaining the component \(- y^2\,dy\). The conclusion is : a generator of the
homology in degree 1 is the cycle \(-y^2\,dy \in
\mbox{\(Z_1(\Ksz_\ast(\bQ[x,y]/\!\!\ideal{y^3}))\)}\).

Now we must use the next short exact sequence to take care of the generator
\(x^2\) of the ideal \(\Ib'\):
\[
0 \leftarrow \frac{\fRb}{\ideal{x^2,y^3}}
\stackrel{\textrm{\scriptsize{pr}}}{\longleftarrow}\frac{\fRb}{\ideal{y^3}}
\stackrel{\times x^2}{\longleftarrow} \frac{\fRb}{\ideal{y^3}} \leftarrow 0.
\]
For in this case, \(\ideal{y^3}:\ideal{x^2}\ =\ \ideal{y^3}\). The available
work above, combined with the Cone Reduction Theorem gives the effective
homology of \(\Cone(\times x^2)\) when applied to the corresponding Koszul
complexes. Which cone can be reduced over
\(\Ksz_\ast(\fRb/\!\!\ideal{x^2,y^3})\). The main components of the result are
in the diagram:
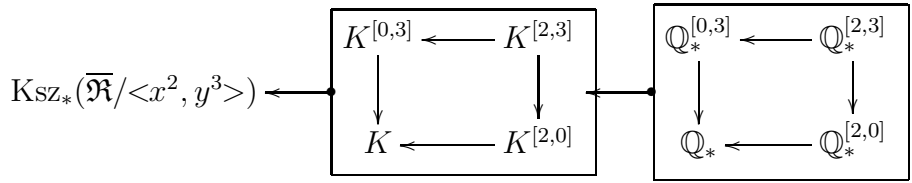
\[
\xymatrix@R0pt{
 & K^{[0,3]} \ar[dd] & K^{[2,3]} \ar[l] \ar[dd] & \bQ_\ast^{[0,3]} \ar[dd] &
 \bQ_\ast^{[2,3]} \ar[l] \ar[dd]
 \\
 \Ksz_\ast(\fRb /\!\! \ideal{x^2, y^3}) & {\phantom{\rule{30pt}{10pt}}} \ar@{{*}->}[l] &
 {\phantom{\rule{30pt}{10pt}}} & {\phantom{\rule{30pt}{10pt}}} \ar@{{*}->}[l]
 \\
 & K & K^{[2,0]} \ar[l] & \bQ_\ast & \bQ_\ast^{[2,0]} \ar[l]
 \save "1,2"."1,3"."3,2"."3,3"*+[F]\frm{} \restore
 \save "1,4"."1,5"."3,4"."3,5"*+[F]\frm{} \restore
}
\]
where \(K\) is a shorthand for \(\Ksz_\ast(\fRb)\). The cones we have to work
with are now \emph{cones of cones}, which explains why the boxes representing
these cones are now square boxes; because of the recursive organization, the
tower of cones can have in the general case an arbitrary number of floors.
Playing here the same game as before for the homology generator in degree 2
leads to the diagram:
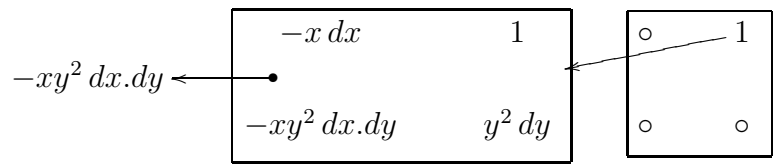
\[
\xymatrix@R0pt{
 & {-x\,dx} & 1 & \circ & 1 \ar[lld]
 \\
 {-xy^2\,dx.dy} & {\phantom{\rule{30pt}{10pt}}} \ar@{{*}->}[l] &
 {\phantom{\rule{30pt}{10pt}}}
 \\
 & {-xy^2\,dx.dy} & {y^2\,dy} & \circ & \circ
 \save "1,2"."1,3"."3,2"."3,3"*+[F]\frm{} \restore
 \save "1,4"."1,5"."3,4"."3,5"*+[F]\frm{} \restore
}
\]
Please try to do it by hand; it is not so hard but in particular when you have
to mix which has been done at the  previous level for \(\ideal{y^3}\) with the
new equivalence to be constructed, things become quickly relatively complex.
And if you have an ideal with many generators, of course a machine program is
necessary.

The Kenzo program can process these calculations. Let us make Kenzo construct
the previous diagram. First we define the ideal. Every generator \(x^\alpha
y^\beta\) is coded as the integer list \((\alpha\ \beta)\) and the ideal as a
list of generators.

 \bmp
 \bmpi\verb|> (setf I '((2 0) (0 3)))|\empim
 \bmpi\verb|((2 0) (0 3))|\empi
 \emp

Constructing the corresponding Koszul complex.

 \bmp
 \bmpi\verb|> (setf K (k-complex/i 2 I))|\empim
 \bmpi\verb|[K3 Chain-Complex]|\empi
 \emp

Constructing the \underline{ef}fective \underline{h}o\underline{m}ology of the
Koszul complex, assigned to the symbol \boxtt{EH}.

 \bmp
 \bmpi\verb|> (setf EH (efhm K))|\empim
 \bmpi\verb|[K98 Equivalence K3 <= K71 => K74]|\empi
 \emp

Kenzo automatically organizes the recursion process and returns an
\emph{equivalence} between the Koszul complex \(\boxtt{K3} =
\Ksz_\ast(\fRb/\!\!\ideal{x^2,y^3})\) and the \emph{effective} chain-complex
\boxtt{K74} via another chain-complex \boxtt{K71}, only \emph{locally
effective}, namely the left hand square of the above diagram.

The chain-complex \boxtt{K74} is effective and we can ask for the basis for
example in degree 2.

 \bmp
 \bmpi\verb|> (basis (k 74) 2)|\empim
 \bmpi\verb|(<Con1 <Con1 Z-GNRT>>)|\empi
 \emp

The rank is 1 and the unique generator is in a cone of cones. Because this
complex is effective, the homology groups are elementarily computed, for
example in degree 2. The function \boxtt{homology-gen} returns a list of
generators for this homology.

 \bmp
 \bmpi\verb|> (homology-gen (k 74) 2)|\empi
 \bmpi\verb|(|\empi
 \bmpi\verb|----------------------------------------------------------------------{CMBN 2}|\empi
 \bmpi\verb|<1 * <Con1 <Con1 Z-GNRT>>>|\empi
 \bmpi\verb|------------------------------------------------------------------------------|\empi
 \bmpi\verb|)|\empi
 \emp

Only one generator, we find again the same generator presented as a combination
of degree 2 (\boxtt{\{CMBN 2\}}) with one term, coefficient 1, and generator
\boxtt{<Con1 <Con1 Z-GNRT>>}.

Now, important intermediate step, we want to lift this generator of ``abstract
homology'' into \boxtt{K71}. We first extract this generator from the one
element list, then apply the \boxtt{rg} component (\boxtt{rg} = right hand
\(g\)) of our equivalence. The Lisp symbol `\boxtt{*}' points to the last
result returned.

 \bmp
 \bmpi\verb|> (first *)|\empim
 \bmpi\verb|----------------------------------------------------------------------{CMBN 2}|\empi
 \bmpi\verb|<1 * <Con1 <Con1 Z-GNRT>>>|\empi
 \bmpi\verb|------------------------------------------------------------------------------|\empix
 \bmpi\verb|> (rg EH *)|\empim
 \bmpi\verb|----------------------------------------------------------------------{CMBN 2}|\empi
 \bmpi\verb|<-1 * <Con0 <Con0 ((1 2) (1 1))>>>|\empi
 \bmpi\verb|<1 * <Con0 <Con1 ((0 2) (0 1))>>>|\empi
 \bmpi\verb|<-1 * <Con1 <Con0 ((1 0) (1 0))>>>|\empi
 \bmpi\verb|<1 * <Con1 <Con1 ((0 0) (0 0))>>>|\empi
 \bmpi\verb|------------------------------------------------------------------------------|\empi
 \emp

You easily recognize the element which was displayed in the left hand box of
the last diagram, taking account for example of the translation \(\boxtt{((1 2)
(1 1))} = x y^2 \, dx.dy\). There remains to go to our Koszul complex, applying
this time the \boxtt{lf} component of the equivalence.

 \bmp
 \bmpi\verb|> (lf EH *)|\empim
 \bmpi\verb|----------------------------------------------------------------------{CMBN 2}|\empi
 \bmpi\verb|<-1 * ((1 2) (1 1))>|\empi
 \bmpi\verb|------------------------------------------------------------------------------|\empi
 \emp

The representant cyle is \(- x y^2 \, dx.dy\). Computing the \emph{ordinary}
homology of such a simple Koszul complex is elementary; computing the
\emph{effective} homology is already not so easy, think of the six morphisms
defining the equivalence \boxtt{K98} between our Koszul complex \boxtt{K3} and
the effective chain-complex \boxtt{K74}; think also of the right differential
to be installed in the cones of cones. For an ideal with more variables and
more generators, this cannot reasonably be obtained without a machine. Let us
for example consider the ideal \(\Ib' = \mbox{\(< \!\! v^3 w^3 x^3
z^2\!,\)}\hspace*{5pt} \hfill v^2 w^3 x y z^3,\hfill v w^3 x^3 y z^2,\hfill v
w^3 x^2 y^2 z^2,\hfill v^3 w^3 x y^2 z,\hfill v^3 w^3 x^3 y,\hfill w^2 x^3 y^3
z^2,\hfill v^2 w^3 x^2 y^3, \\ v^2 w x y^3 z^3, v^2 x^3 y^2 z^3, v^2 w^2 x^2
y^2 z^3, v^3 w^2 x^3 y z^3 \!\!>\) of \(\bQ[v,w,x,y,z]\). Working exactly as
before, we use Kenzo to construct the ideal, the corresponding Koszul complex
and its effective homology. The 3-homology has rank 9 and we extract the
abstract generator number 7, for which a representant cycle is computed in the
Koszul complex.  You may observe the numerical notation of monomials by number
lists is quickly more readable than the usual one.

 \bmp
 \bmpi\verb|> (setf I '((3 3 3 0 2) (2 3 1 1 3) (1 3 3 1 2) (1 3 2 2 2)|\empi
 \bmpi\verb|            (3 3 1 2 1) (3 3 3 1 0) (0 2 3 3 2) (2 3 2 3 0)|\empi
 \bmpi\verb|            (2 1 1 3 3) (2 0 3 2 3) (2 2 2 2 3) (3 2 3 1 3)))|\empim
 \bmpi\verb|[ligne deleted]|\empix
 \bmpi\verb|> (setf K (k-complex/i 5 I))|\empim
 \bmpi\verb|[K3 Chain-Complex]|\empix
 \bmpi\verb|> (setf eh (efhm K))|\empim
 \bmpi\verb|[K1235 Equivalence K3 <= K1208 => K1211]|\empix
 \bmpi\verb|> (seventh (homology-gen (k 1211) 3))|\empim
 \bmpi\verb|----------------------------------------------------------------------{CMBN 3}|\empi
 \bmpi\verb|<1 * <Con0 <Con0 <Con1 <Con0 <Con1 <Con1 Z-GNRT>>>>>>>|\empi
 \bmpi\verb|------------------------------------------------------------------------------|\empix
 \bmpi\verb|> (rg EH *)|\empim
 \bmpi\verb|----------------------------------------------------------------------{CMBN 3}|\empi
 \bmpi\verb|<-1 * <Con0 <Con0 <Con0 <Con0 <Con0 <Con0 <Con0 <Con0 <Con0 ...|\empi
 \bmpi\verb|                  ... <Con0 <Con0 <Con0 ((2 1 2 3 2) (0 1 1 0 1))>>>>>>>>>>>>>|\empi
 \bmpi\verb|[... 8 lines deleted ...]|\empi
 \bmpi\verb|<1 * <Con0 <Con0 <Con1 <Con0 <Con1 <Con1 ((0 0 0 0 0) (0 0 0 0 0))>>>>>>>|\empi
 \bmpi\verb|------------------------------------------------------------------------------|\empix
 \bmpi\verb|> (lf EH *)|\empim
 \bmpi\verb|----------------------------------------------------------------------{CMBN 3}|\empi
 \bmpi\verb|<-1 * ((2 1 2 3 2) (0 1 1 0 1))>|\empi
 \bmpi\verb|<1 * ((2 1 3 2 2) (0 1 0 1 1))>|\empi
 \bmpi\verb|<-1 * ((2 2 2 2 2) (0 0 1 1 1))>|\empi
 \bmpi\verb|------------------------------------------------------------------------------|\empi
 \emp

Therefore a cycle representing our seventh 3-homology class is \(- v^2 w x^2
y^3 z^2 \, dw.dx.dz + v^2 w x^3 y^2 z^2 \, dw.dy.dz - v^2 w^2 x^2 y^2 z^2 \,
dx.dy.dz\). Note the terrible cone towers which are involved. But all these
calculations are almost instantaneous, and having tools doing them conveniently
will soon become mandatory in modern homological algebra.

It is not obvious the last element actually is a homology generator, but we can
at least verify it is a cycle!

 \bmp
 \bmpi\verb|> (? K *)|\empim
 \bmpi\verb|----------------------------------------------------------------------{CMBN 2}|\empi
 \bmpi\verb|------------------------------------------------------------------------------|\empi
 \emp

When the first argument is a chain-complex, here \boxtt{K}, a symbol pointing
to our Koszul complex \(\Ksz_\ast(\fRb/\Ib')\), and the second argument is a
combination, here `\boxtt{*}' a symbol pointing to the last result returned,
this combination must be an element of the chain-complex and the operator
`\boxtt{?}' computes the boundary of this combination in this chain-complex.
Absence of terms between both horizontal dash lines means the result is null.

\subsection{Case of a non-monomial ideal.}

The work around the \emph{monomial} ideal \(\Ib'\) in the previous section was
undertaken because we hope to be able to apply the BPL to obtain a version with
effective homology \(\Ksz(\fRb/\Ib)_{EH}\) in the general case from
\(\Ksz(\fRb/\Ib')_{EH}\) now available. Thanks to a Grobner basis of \(\Ib\),
the \(\fk\)-vector spaces \(\fRb/\Ib\) and \(\fRb/\Ib'\) are canonically
isomorphic, which implies the corresponding Koszul complexes are also
isomorphic as graded \(\fk\)-vector spaces. Only the differentials are
different, we are in a situation where the BPL is applicable.

What about the nilpotency condition? A simple example explains better what
happens than a generic description. Let us take \(\Ib =\ \ideal{x - t^3, y -
t^5}\ \subset \bQ[x,y,t]\); the DegRevLex reduced Groebner basis for the order
\(x > y > t\) is \(\Ib = \mbox{\(< x t^2 - y,\)}\ \mbox{\(t^3 - x,\)}\ x^2 - y
t\!\!>\) and the associated monomial ideal is \(\Ib' =\ \mbox{\(\ideal{x t^2,
t^3, x^2}\)}\). Both quotients \(\fRb/\Ib\) and \(\fRb/\Ib'\) are isomorphic
\(\fk\)-vector spaces with basis \(\cup_{\alpha \in \bN}\{y^\alpha, x y^\alpha,
t y^\alpha, x t y^\alpha, t^2 y^\alpha\}\). A generator of a Koszul complex is
a product of such a basis element and a combination of \(dx, dy\) and \(dt\)
without any repetition. The differential is obtained by successively replacing
the various~\(d\textit{?}\) by \(\textit{?}\) with the right signs. If ever the
resulting coefficient is not in our basis, two cases:
\begin{enumerate}
\item
In the monomial case, the reduction modulo the ideal cancels the corresponding
term.
\item
In the initial non-monomial case, a reduction modulo the ideal in general
generates other monomials.
\end{enumerate}

For example \(d_{\textrm{\scriptsize Ksz}}(x t \, dt) = x t^2\) which is not in
our basis, hence to be reduced; modulo \(\Ib'\), the result is null; modulo
\(\Ib\), because of the generator \(x t^2 - y\), the result is non null, it is
\(y\). \emph{The main point is here}: because of the structure of the Groebner
basis, the multigrading of the result is certainly \emph{strictly less} than
the multigrading of the initial monomial. In our small example, the
multigrading of \(x t \, dt\) is \([1,0,2]\) while the multigrading of \(y\) is
\([0,1,0] < [1,0,2]\) for DegRevLex in \(\bQ[x,y,t]\).

\begin{prp}\label{10322}---
Let \(\Ib \subset \fRb\) an ideal. Some Groebner monomial order is given for
the multigrading. Cancelling the trailing terms of the corresponding reduced
Groebner basis defined an ``approximate'' monomial ideal \(\Ib'\), allowing us
to identify \emph{as multigraded \(\fk\)-vector spaces} the Koszul complexes
\(\emph{\Ksz}_\ast(\fRb / \Ib)\) and \(\emph{\Ksz}_\ast(\fRb / \Ib')\). Then
the perturbation difference between both differentials \emph{strictly decreases
the multigrading}. \QED
\end{prp}

Theorem~\ref{18292} constructs an equivalence:
\[
\Ksz_\ast(\fRb/\Ib') \lrdc \hC'_\ast \rrdc EC'_\ast
\]
with an effective chain-complex \(EC'_\ast\). We would like to construct:
\[
\Ksz_\ast(\fRb/\Ib) \lrdc \hC_\ast \rrdc EC_\ast
\]

As usual, applying the Easy Perturbation Lemma~\ref{73716} between
\(\Ksz_\ast(\fRb/\Ib')\) and \(\Ksz_\ast(\fRb/\Ib)\) will produce the wished
chain-complex \(\hC_\ast\), the same graded vector space as \(\hC'_\ast\) but
with another differential. Then applying the serious Basic Perturbation
Lemma~\ref{07404} produces a new effective chain-complex \(EC_\ast\), the same
graded vector space as \(EC'_\ast\) with another differential.

The only critical point is the nilpotency hypothesis. The initial equivalence
produced by Theorem~\ref{18292} is entirely made of objects, differentials,
morphisms, homotopy operators that are, thanks to the \emph{multigrading shift}
process, \emph{multigraded}. The initial perturbation between Koszul complexes
on the contrary \emph{strictly} decreases the multigrading. The easy
perturbation lemma copies this perturbation into \(\hC'_\ast\) using
multigraded morphisms; therefore the differential perturbation to be applied to
\(\hC'_\ast\) to obtain \(\hC_\ast\) also \emph{strictly decreases the
multigrading}. Now the composition homotopy-perturbation \(h\hdl\) which must
be proved locally nilpotent is made of a multigraded map and another map which
strictly decreases the multigrading; the composition also strictly decreases
the multigrading.

A monomial order defines a well-founded order in the multigrading set; every
strictly decreasing sequence goes to the minimal element, the multigrading of 0
often decided to be \(-\infty\) and our composition \(h \hdl\) is nilpotent for
any argument.

\begin{thr}\label{91916}---
An algorithm computes:
\[
\Ib \mapsto [\emph{\Ksz}(\fRb/\Ib) \lrdc \hC_\ast \rrdc EC_\ast]
\]
where \(\Ib\) is an  ideal of \(\fRb = \fk[x_1, \ldots, x_m]\), and the result
is an equivalence between the corresponding Koszul complex and an effective
chain-complex of finite-dimensional \(\fk\)-vector spaces.\QED
\end{thr}

\subsection{Coming back to the local ring.}

There remains to come back to our local ring \(\fR = \fk[x_1, \ldots, x_m]_0\).
The only difference between the elements of \(\fRb\) and \(\fR\) is that
denominators are allowed for~\(\fR\), on condition such a denominator is non
null at 0. A canonical inclusion is defined \(\fRb \subset \fR\). The ring
\(\fRb\) is factorial and an element of \(\fR\) can be written in a unique
irreducible form \(p/(1-m)\) with \(p \in \fRb\) and \(m \in \fm_0\), the
maximal ideal of \(\fRb\) at 0.

\begin{thr}---
Let \(I\) be an ideal of \(\fR\) and \(\Ib = I \cap \fRb\). The injection
\(\lambda: \fRb \hookrightarrow \fR\) induces an injection \(\lambda:
\emph{\Ksz}_{\fRb}(\fRb / \Ib) \hookrightarrow \emph{\Ksz}_{\fR}(\fR / I)\)
which in turn induces an \emph{isomorphism}:
\[
\lambda: H_\ast(\emph{\Ksz}_{\fRb}(\fRb / \Ib))
\stackrel{\cong}{\longrightarrow} H_\ast(\emph{\Ksz}_{\fR}(\fR / I)).
\]
\end{thr}

In short, the denominators do not play any role in the homological nature of
these Koszul complexes.

\proof Let us qualify as \emph{polynomial} a chain element of the chain-complex
\(\Ksz_{\fRb}(\fRb / \Ib)\): all the coefficients are (equivalence classes of)
\emph{polynomials}. Every polynomial has a (total) degree and also an
\emph{order}, the smallest degree of a non-null monomial component, which
definitions are extended to polynomial chains, without taking account of the
``differential'' terms in \(\wedge V\). The \(\fk\)-vector space
\(H_\ast(\Ksz_{\fRb}(\fRb / \Ib))\) has a finite dimension; choosing cycles
representing some \emph{generators} of this homology, the degree of every
generator is \(< k\) for some \(k \in \bN\). We will carefully examine which
happens when objects are reduced modulo \(\fm_0^k\).

The space of all cycles, \emph{after reduction} modulo \(\fm_0^k\), is also a
finite dimensional vector space where the classes modulo \(\fm_0^k\) of
boundaries are a supplementary of the space generated by the chosen cycles
representing the generators of homology: \(Z' = H \oplus B'\) if \(H\) is the
vector space generated by our (exact) representants, if \(Z'\) (resp. \(B'\))
is the set of all cycles (resp. boundaries) \emph{truncated} at degree \(k\).
In particular any cycle \(z\) of order \(\geq k\) certainly is a boundary; in
fact the homology class of \(z\) is obtained as follows: you truncate the cycle
\(z\) at degree \(k\), obtaining an element \(z' \in Z'\) and the homology
class ``is'' the \(H\)-component \(h\) of \(z' = h + b'\). But if the cycle has
an order \(\geq k\), then \(z' = 0\).

Let us take now a ``local'' cycle \(z \in Z_\ast(\Ksz_{\fR}(\fR / I))\).
Reducing to the same denominator the various components of \(z\), this cycle
can be written \(z = \overline{z}/(1-m)\) with \(\overline{z} \in
Z_\ast(\Ksz_{\fRb}(\fRb / \Ib))\) and \(m \in \fm_0\). For \(\overline{z} =
(1-m)z\) again is a cycle: the differential is a \emph{module} morphism. Now
\(\overline{z}/\mbox{\((1-m)\)} = (\overline{z} + m \overline{z} + \cdots +
m^{k-1} \overline{z}) + m^k \overline{z}/(1-m)\). Because of its order, the
\emph{numerator} \(m^k \overline{z}\) is a boundary in the polynomial Koszul
complex, which allows to express also the fraction \mbox{\(m^k
\overline{z}/(1-m)\)} as a boundary in the localized Koszul complex, again
because the boundary operator is a \emph{module} morphism. The sum of the other
terms \(\overline{z} + m \overline{z} + \cdots + m^{k-1} \overline{z}\) is
polynomial, it is again a cycle and its homology class \(\fh\) in the
polynomial Koszul complex is defined; the previous study shows the homology
class of \(z\) in the localized Koszul complex is \(\lambda(\fh)\) and
\(\lambda\) at the homological level is surjective.

Let us take now \(\overline{z} \in Z_\ast(\Ksz_{\fRb}(\fRb / \Ib))\) and assume
\(\overline{z}\) is a boundary in \(\Ksz_{\fR}(\fR / I)\), that is with the
same calculation as before: \(\overline{z} = d(c/(1-m)) = d(c + mc + \cdots +
m^{k-1}c) + d(m^kc/(1-m)))\) with \(c\) a \emph{polynomial} chain in
\(\Ksz_{\fRb}(\fRb / \Ib)\). So that in fact the last term \(d(m^kc/(1-m))\) is
a difference between polynomial chains and it is also polynomial; furthermore
the computation of \(d(d(m^kc/(1-m)))\) can be done as well in the localized
Koszul complex; the result is null (\(dd = 0\)) and our pseudo-fraction
\(d(m^kc/(1-m))\) is also a cycle in the polynomial Koszul complex. Because of
the order, this polynomial cycle is a boundary \emph{in the polynomial Koszul
complex} and finally the cycle \(\overline{z}\) is a boundary in the polynomial
Koszul complex. In other words the map \(\lambda\) at the homological level is
injective. \QED

It is not very hard to transform this proof into a \((\fRb, \fk, \fk)\)-linear
reduction \(\Ksz_\fR(\fR / I) \rrdc \Ksz_{\fRb}(\fRb / \Ib)\). But the most
appropriate conclusion is the following: the homological problems for
\(\Ksz_{\fRb}(\fRb / \Ib))\) and \(\Ksz_\fR(\fR / I)\) are
\emph{constructively} equivalent. We have seen in the previous section how the
homological problem for \(\Ksz_{\fRb}(\fRb / \Ib)\) is solved thanks to two
essential ingredients: Groebner basis and BPL.

\begin{thr}---
The homological problem of \(\Ksz_\fR(\fR / I)\) is solved. \QED
\end{thr}

\subsection{Effective homology \(\Leftrightarrow\) Effective resolution.}

Let \(I\) be an ideal of our local ring \(\fR = \fk[x_1, \ldots, x_m]_0\). We
know how to compute the effective homology of \(\Ksz(\fR/I)\). We intend now to
use this information to obtain an \emph{effective} \(\fR\)-resolution of
\(\fR/I\). Conversely, an effective resolution naturally gives the effective
homology of the corresponding Koszul complex.

As before, it is better to work with \(\Ib = I \cap \fRb\). Elementary
arguments show an \(\fRb\)-resolution \(\Rsl_{\fRb}(\fRb/\Ib)\) induces an
\(\fR\)-resolution \(\Rsl_\fR(\fR/I) := \Rsl_{\fRb}(\fRb/\Ib) \otimes_{\fRb}
\fR\). In particular the \(\fRb\)-module \(\fR\) is \emph{flat}.

The connection between effective homology of \(\Ksz(\fRb/\Ib)\) and effective
resolution of \(\fRb/\Ib\) is the Aramova-Herzog bicomplex~\cite{ARHR}.

\begin{dfn}---
\emph{Let \(\Ib\) be an ideal of \(\fRb\). The \emph{Aramova-Herzog bicomplex}
\(\ArHr(\fRb/\Ib)\) of \(\fRb/\Ib\) is \(\ArHr(\fRb/\Ib) := \fRb/\Ib
\otimes_\fk \wedge V \otimes_\fk \fRb\) provided with both differentials coming
from both Koszul complexes present in its definition.}
\end{dfn}

We recall \(V\) is the \(\fk\)-vector space \(\fm_0 / \fm_0^2\) provided with
the canonical basis \((dx_1, \ldots, dx_m)\), the ideal \(\fm_0\) being the
maximal ideal at 0 of \(\fRb\). We can see \(\ArHr(\fRb/\Ib) := \fRb/\Ib
\otimes \wedge V \otimes \fRb = \Ksz(\fRb/\Ib) \otimes \fRb\) and the
\emph{vertical} differential \(\partial''\) of our bicomplex is \(\partial'' :=
d_{\sKsz(\fRb/\Ib)} \otimes \id{\fRb}\). In the same way, appropriately
swapping the factors \(\wedge V\) and \(\fRb\), we can interpret
\(\ArHr(\fRb/\Ib) := \fRb/\Ib \otimes \wedge V \otimes \fRb = \fRb/\Ib \otimes
\Ksz(\fRb)\) and the \emph{horizontal} differential is \(\partial' :=
\id{\fRb/\Ib} \otimes d_{\sKsz(\fRb)}\). See the following diagram where the
ground ring \(\fRb\) is split into its homogeneous components \(\fRb_p\), this
index \(p\) defining the horizontal grading, which implies the horizontal
differential has degree \((0, +1)\). In the same way, the central factor
\(\wedge V\) is split into homogeneous components \(\wedge^v V\), the vertical
degree is \(q = v+p\) but this time the vertical differential has degree \((-1,
0)\). The total degree therefore is \(v\). The bicomplex is null outside the
strip \(0 \leq v \leq m\), that is, \(q \in [p\,..\, p+m]\).

\[
 \xymatrix{
 \ar@{.>}[d] & \ar@{.>}[d] & \ar@{.>}[d] & \ar@{.>}[d]
 \\
 \fRb/\Ib \otimes \wedge^3 \otimes \fRb_0 \ar[r]^{\partial'} \ar[d]_{\partial''} &
 \fRb/\Ib \otimes \wedge^2 \otimes \fRb_1 \ar[r]^{\partial'} \ar[d]_{\partial''} &
 \fRb/\Ib \otimes \wedge^1 \otimes \fRb_2 \ar[r]^{\partial'} \ar[d]_{\partial''} &
 \fRb/\Ib \otimes \wedge^0 \otimes \fRb_3 \ar[r] \ar[d] & 0
 \\
 \fRb/\Ib \otimes \wedge^2 \otimes \fRb_0 \ar[r]^{\partial'} \ar[d]_{\partial''} &
 \fRb/\Ib \otimes \wedge^1 \otimes \fRb_1 \ar[r]^{\partial'} \ar[d]_{\partial''} &
 \fRb/\Ib \otimes \wedge^0 \otimes \fRb_2 \ar[r] \ar[d] & 0
 \\
 \fRb/\Ib \otimes \wedge^1 \otimes \fRb_0 \ar[r]^{\partial'} \ar[d]_{\partial''} &
 \fRb/\Ib \otimes \wedge^0 \otimes \fRb_1 \ar[r] \ar[d] & 0
 \\
 \fRb/\Ib \otimes \wedge^0 \otimes \fRb_0 \ar[r] \ar[d] & 0 \\
 0
 }
 \]

If we see \(\ArHr(\fRb/\Ib) = \fRb/\Ib \otimes \Ksz_{\fRb}(\fRb)\), using the
fact the Koszul complex of the ground ring \(\Ksz_{\fRb}(\fRb)\) is acyclic
(Theorem~\ref{37681}, or more precisely the variant for \(\fRb\), easier), we
will construct a reduction \(\ArHr(\fRb/\Ib) \rrdc \fRb/\Ib\). Considering now
the symmetric factorization \(\ArHr(\fRb/\Ib) = \Ksz(\fRb/\Ib) \otimes \fRb\),
using the effective homology \(\Ksz(\fRb/\Ib) \eqvl H\), that is, an
equivalence between the Koszul complex and some \emph{effective} chain-complex
\(H\), we will construct an equivalence \(\ArHr(\fRb/\Ib) \eqvl H \otimes
\fRb\) with an appropriate differential for \(H \otimes \fRb\) coming again
from the BPL. Combining this reduction and this equivalence will produce an
equivalence \(\fRb/\Ib \eqvl H \otimes \fRb\) which is the looked-for
resolution. And the whole process can be reversed, starting from
\emph{effective} resolution\emph{s}, going to \emph{effective}
homolog\emph{ies} of the Koszul complex.

Let us recall the possible geometrical interpretation of the Koszul complex
given Section~\ref{69368}. A natural analogous interpretation can be given
here. The Koszul complex \(\Ksz(M) = M \otimes_t \wedge V\) is the total space
of a fibration \(M \hookrightarrow M \otimes_t \wedge V \rightarrow \wedge V\).
We can also consider the symmetric Koszul complex \(\Ksz'(\fR) = \wedge V\,
_t\!\! \otimes \fR\) with an analogous fibration. Combining both fibrations
gives the diagram:
\[
\xymatrix{
 & & *++{M} \ar@{^{(}->}[d]
 \\
 & M \otimes_t \wedge V \, _t\!\!\otimes \fR \ar[r] \ar[d] & M \otimes_t \wedge
 V \ar[d]
 \\
 *++{\fR} \ar@{^{(}->}[r] & \wedge V \, _t\!\!\otimes \fR \ar[r] & \wedge V
 }
\]
where the Aramova-Herzog bicomplex \(M \otimes_t \wedge V \, _t\!\!\otimes
\fR\) is the \emph{pullback} of the vertical fibration by the horizontal map
\(\xymatrix@1{\wedge V \, _t\!\!\otimes \fR \ar[r] & \wedge V}\); but the
original space of this map is contractible (Theorem~\ref{37681}), so that it
has the homotopy type of a point, and the pullback, up to homotopy, is nothing
but the base fiber of the vertical fibration, that is, the module \(M\). It is
this homotopy equivalence which is systematically exploited by Aramova and
Herzog. Note the vertical arrow between \(M \otimes_t \wedge V\) and \(\wedge
V\) is not actually defined and there is only some ``analogy'' with the
projection of a topological fibration.

\vspace{10pt}

\noindent \textsc{First reduction.}
\\*
The chain-complex \(\Ksz(\fRb)\) is acyclic. More precisely, every
``horizontal'' \emph{sub}complex \((\oplus_{v+p = q} \wedge^v V \otimes \fRb_p,
\partial')\) at ordinate \(q\) is acyclic, except for \(q = 0\) where
\(\wedge^0 V \otimes \fRb_0 = \fk\). Appling the functor \(\fRb/\Ib\ \otimes
\hspace{-5pt} <?\hspace{-8pt}>\) gives a reduction \((\ArHr(\fRb/\Ib),
\partial') \rrdc \fRb/\Ib\), where, quite important, the vertical differential
\(\partial''\) has temporarily been cancelled. Reinstalling the vertical
differential is BPL's job. Verifying the nilpotency hypothesis is the following
game: you start from \(\ArHr_{p,q}(\fRb/\Ib)\), a homotopy operator expressing
\(\partial'\) is contractible leads you to \(\ArHr_{p-1,q}(\fRb/\Ib)\), the
perturbation \(\partial''\) goes to \(\ArHr_{p-1,q-1}(\fRb/\Ib)\), the next
homotopy operator goes to \(\ArHr_{p-2,q-1}(\fRb/\Ib)\) and so on. Finally you
get out from the diagram at \(\ArHr_{0,q-p}(\fRb/\Ib)\) after having run \(p\)
steps of a stairs leftdownward.\vspace{10pt}

\noindent \textsc{Second equivalence}.
\\*
The analogous work for the second interpretation of the Aramova-Herzog
bicomplex works as follows. We now consider \(\ArHr(\fRb/\Ib) = \Ksz(\fRb/\Ib)
\otimes \fRb\). Theorem~\ref{91916} constructs an equivalence \(\Ksz(\fRb/\Ib)
\eqvl H\) where \(H\) is a chain-complex of finite type, called \(H\) for it
describes the ``abstract'' homology of the Koszul complex. This equivalence can
be applied to every vertical of the Aramova-Herzog bicomplex, which produces an
equivalence \((\ArHr(\fRb/\Ib), \partial'') \eqvl H \otimes \fRb\) with \(d_{H
\otimes \fRb} = d_H \otimes \id{\fRb}\).

There remains to reinstall the horizontal differential \(\partial'\), again
under the responsability of BPL. The nilpotency check runs the same stairs as
before, but in the \emph{reverse direction}, and this time we \emph{do not
reach} any void part of the bicomplex. But the vertical homotopy operator comes
from Theorem~\ref{91916} and the details of the proof show this homotopy
operator \emph{does not} increase the Groebner multidegree: in the monomial
case, this operator is multigraded and respects the multigrading; in the
general case, Shih's magic formula \(h' = h \psi = h \sum_{i=0}^\infty (-1)^i
(\widehat{\delta} h)^i\), see page~\pageref{85606}, gives the result because of
Proposition~\ref{10322}.

We were speaking here of the multigrading of \(\ArHr(\fRb/\Ib) = \Ksz(\fRb/\Ib)
\otimes \fRb\) deduced from the left hand factor \(\Ksz(\fRb/\Ib)\), neglecting
the right hand factor \(\fRb\). If we consider the relevant perturbation
\(\partial'\), every term of \(\partial'(\kappa \otimes v \otimes \rho)\) is
obtained by replacing some \(dx_i\) in \(v \in \wedge V\) by the corresponding
\(x_i\) to be installed as a multiplier in \(\rho \in \fRb\). This
\emph{strictly} decreases the ``left hand'' multigrading of \(\ArHr(\fRb/\Ib) =
\Ksz(\fRb/\Ib) \otimes \fRb\). The nilpotency condition is satisfied.

Applying the BPL is allowed, which gives an equivalence:
\[
(\ArHr(\fRb/\Ib), \partial' \oplus \partial'') \eqvl (H \otimes \fRb, d')
\]
with a \emph{new} differential \(d' \neq d_h \otimes \id{\fRb}\) except in
trivial cases. Combining the first reduction and the second equivalence gives
the next theorem.

\begin{thr}---
The Aramova-Herzog bicomplex \(\ArHr(\fRb/\Ib)\) produces an equivalence:
\[
\fRb/\Ib \eqvl (H \otimes \fRb, d')
\]\QED
\end{thr}

The left hand term of this equivalence is without any differential, and more
exactly is a chain-complex concentrated in differential degree 0. Our
equivalence is nothing but a resolution \((H \otimes \fRb, d')\) for
\(\fRb/\Ib\). The component \(H\) is a free (\(!\)) \(\fk\)-vector space and
the tensor product \(H \otimes \fRb\) therefore is a \emph{free}
\(\fRb\)-module. The possibly sophisticated differential \(d'\), sophisticated
but \emph{automatically} produced by BPL, describes the main part of the
resolution.

By the way, why the differential \(d'\) is an \(\fRb\)-morphism? The BPL
constructs this differential as a combination of compositions whose ingredients
can be:
\begin{itemize}
\item
\(d \otimes \id{\fRb}\) for \(d\) the \(H\)-differential; the second factor
\(\otimes \id{\fRb}\) ensures the \(\fRb\)-linearity.
\item
\(g \otimes \id{\fRb}\) for \(g: H \rightarrow \fRb/\Ib\) the second component
of the effective homology of \(\fRb/\Ib\); same argument.
\item
\(h \otimes \id{\fRb}\) for \(h: \fRb/\Ib \rightarrow \fRb/\Ib\) the third
component of the effective homology of \(\fRb/\Ib\); same argument.
\item
\(\id{\fRb/\Ib} \otimes \partial'\), but \(\partial'\) is \(\fRb\)-linear.
\item
\(f \otimes \id{\fRb}\) for \(f: \fRb/\Ib \rightarrow H\) the first component
of the effective homology of \(\fRb/\Ib\); same argument as above.
\end{itemize}%
On the contrary, when constructing the homotopy \emph{effectively} describing
the acyclicity property of \((H \otimes \fRb, d')\), the contracting homotopy
of \(\wedge V \otimes \fRb\) is used, which homotopy \emph{is not}
\(\fRb\)-linear.

\subsection{Examples.}

\subsubsection{The minimal non-trivial example.}

Let \(\fR = \fk[x]\) (one variable) and \(M = \fR/\!<\!x^2\!>\). And let us
assume we do not know~(!) the minimal resolution. Here the ideal is monomial
and the steps 1 and~3 of our algorithm are void. The effective homology of the
Koszul complex:
\[
\textrm{Ksz}(M) = [\cdots\leftarrow  0 \leftarrow M \leftarrow M.dx \leftarrow
0 \leftarrow \cdots]
\]
is made of the chain-complex:
\[
H = [\cdots \leftarrow 0 \leftarrow \fk_0 \stackrel{0}{\leftarrow} \fk_1
\leftarrow 0 \leftarrow \cdots]
\]
(where \(\fk_0\) and \(\fk_1\) are copies of the ground field \(\fk\) with
respective homological degrees 0 and 1) \emph{and} of the maps \(\rho =
(f,g,h)\) with:
\begin{enumerate}
\item
\(f:M \rightarrow \fk_0\) is defined by \(f(1) = 1_0, f(x) = 0\).
\item
\(f: M.dx \rightarrow \fk_1\) is defined by \(f(1.dx) = 0, f(x.dx) = 1_1\).
\item
\(g: \fk_0 \rightarrow M\) is defined by \(g(1_0) = 1\).
\item
\(g: \fk_1 \rightarrow M.dx\) is defined by \(g(1_1) = x.dx\).
\item
\(h: M \rightarrow M.dx\) is defined by \(h(1) = 0, h(x) = 1.dx\).
\end{enumerate}

We must guess the right differential on \(\Rsl(M) = (H \otimes_\fk \fR, d =
?)\). The only non-trivial differential \(d_{\textrm{\scriptsize Rsl}(M)}(1_1
\otimes 1_\fR)\) comes from a unique non-null term in the series \((\Sigma)\),
following the path:
\[
1_1 \otimes 1_\fR \stackrel{g \otimes \textrm{\scriptsize id}_\fR}{\longmapsto}
x \otimes dx \otimes 1_\fR \stackrel{\partial'}{\mapsto} x \otimes 1 \otimes x
\stackrel{-h \otimes \textrm{\scriptsize id}_\fR}{\longmapsto} -1 \otimes dx
\otimes x \stackrel{\partial'}{\longmapsto} -1 \otimes 1 \otimes x^2
\stackrel{f \otimes \textrm{\scriptsize id}_\fR}{\longmapsto} -1_0 \otimes x^2
\]
and, surprise, we find the resolution \(1_1 \otimes 1_\fR \mapsto -1_0 \otimes
x^2 \). You find it is a little complicated for a so trivial particular case?
The point is the following: this example in a sense is \emph{complete}, the
most general case is not harder, you have here all the ingredients of the
general solution, nothing more is necessary.

\subsubsection{First Aramova-Herzog example.}

In the paper~\cite{ARHR}, Aramova and Herzog consider the toy example of the
ideal \(I = \mbox{\(<\!x_1 x_3, x_1 x_4, x_2 x_3, x_2 x_4\!>\)}\) in \(\fR =
\fk[x_1,x_2,x_3,x_4]\). The ideal is monomial and again, steps 1 and 3 of our
algorithm are void. The Betti numbers of \(\Ksz(\fR/I)\) are \((1,4,4,1)\) and
the effective homology of \(\Ksz(\fR/I)\) is a diagram:
\[\rho = \raisebox{4pt}{\framebox{\xymatrix@1{
 \scriptstyle h\ \textstyle \autoarrow{0.8}\ \textrm{Ksz}(\fR/I) \ar@<-10\ul>[r]_-f & H
 \ar@<-10\ul>[l]_-g }}}
 \]
 where \(H\) is the chain-complex with null differentials:
 \[
 \cdots \stackrel{}{\longleftarrow} \fk \stackrel{0}{\longleftarrow} \fk^4
 \stackrel{0}{\longleftarrow} \fk^4 \stackrel{0}{\longleftarrow} \fk \longleftarrow \cdots
 \]
 The arrows \(f\) and \(g\) are chain-complex morphisms satisfying \(fg =
\id{H}\), the self-arrow \(h\) is a homotopy between \(gf\) and
\(\id{\textrm{\scriptsize Ksz}(\fR/I)}\), that is, \(\id{\textrm{\scriptsize
Ksz}(\fR/I)} = gf + dh + hd\), and finally, the composite maps \(fh\), \(hg\)
and \(h^2\) are null. These maps smartly express the big chain-complex
\(\Ksz(\fR/I)\) as the direct sum of the small one \(H\), in this case with
trivial differentials, and an acyclic one (\(\ker f\)) \emph{with} an
\emph{explicit} contraction~\(h\). Our Kenzo program~\cite{DRSS} computes this
effective homology in a negligible time with respect to input-output. In
particular the map \(g\) defines representants for the alleged homology
classes, the map \(f\) is a projection which in particular sends cycles to
their homology classes, and \(h\) is the main component of a
\emph{constructive} proof of these claims.

The minimal resolution of \(\fR/I\) is \(\Rsl(\fR/I) = H \otimes \fR\) where a
non-trivial differential must be installed. Let us apply our formula to the
unique generator \(\fh_{3,1} \otimes 1_\fR\) of \(H_3 \otimes \fR\). Kenzo
chooses \(g(\fh_{3,1}) = x_2\,dx_1.dx_3.dx_4 - x_1\,dx_2.dx_3.dx_4\) and:
\[
\partial' (g\otimes 1_\fR) (\fh_{3,1} \otimes 1_\fR) =
\begin{array}[t]{l}
x_2\,dx_3.dx_4 \otimes x_1
\\
- x_1\,dx_3.dx_4 \otimes x_2
\\
+ (-x_2\,dx_1.dx_4 + x_1\,dx_2.dx_4) \otimes x_3
\\
+ ( x_2\,dx_1.dx_3 - x_1\,dx_2.dx_3) \otimes x_4
\end{array}
\]
Kenzo is a little luckier than Aramova and Herzog, for he had chosen:
\[
\begin{array}{rcl}
g(\fh_{2,1}) &=& - x_2 \, dx_1.dx_3 + x_1 \, dx_2.dx_3
\\
g(\fh_{2,2}) &=& - x_1 \, dx_3.dx_4
\\
g(\fh_{2,3}) &=& - x_2 \, dx_1.dx_4 + x_1 \, dx_2.dx_4
\\
g(\fh_{2,4}) &=& - x_2 \, dx_3.dx_4
\end{array}
\]
which is enough to imply:
\[
d(\fh_{3,1}) = - \fh_{2,1} \otimes x_4 + \fh_{2,2} \otimes x_2 + \fh_{2,3}
\otimes x_3 - \fh_{2,4} \otimes x_1
\]
that is, except for legal minor differences, directly the same result as
Aramova and Herzog.

Let us now \emph{force} Kenzo to choose Aramova and Herzog's representants for
the homology classes of \(H_2\). This amounts to replacing the component \(g\)
in degree 2 by another one \(g' = g + d\alpha\) for \(\alpha\) a map \(\alpha:
H_2 \rightarrow \Ksz_3(\fR/I)\) chosen to give the new representants. The cycle
\(- x_2 \, dx_1.dx_i + x_1 \, dx_2.dx_i\) (\(i = 3\) or 4) is homologous to the
cycle \(- x_i \, dx_1.dx_2\) (sign error in \cite{ARHR}) thanks to the boundary
preimage \(dx_1.dx_2.dx_i\). So that we transform Kenzo's choices to Aramova
and Herzog's choices by taking \(\alpha(\fh_{2,1}) = - dx_1.dx_2.dx_3\),
\(\alpha(\fh_{2,3}) = - dx_1.dx_2.dx_4\) and \(\alpha(\fh_{2,i}) = 0\) for \(i
= 2\) or 4.

The component \(f\) of the reduction does not change, but the homotopy \(h_2\)
must be replaced by \(h'_2 = h_2 (\id{} - d \alpha f_2)\). Repeating the same
computation, taking account of \(g_3 = g'_3\), now the homotopy term \((h'_2
\otimes \id{\fR}) \partial' (g_3 \otimes \id{\fR})(\fh_{3,1}) = dx_1.dx_2.dx_4
\otimes x_3 - dx_1.dx_2.dx_3 \otimes dx_4\) \emph{is not null}, so that we must
continue the expansion of the series \((\Sigma)\). We find:
\[
\begin{array}{rcl}
- \partial' (h'_2 \otimes \id{\fR}) \partial' (g \otimes \id{\fR}) (\fh_{3,1})
&=& - dx_2.dx_4 \otimes x_1 x_3 + dx_1.dx_4 \otimes x_2 x_3
\\
&&   + dx_2.dx_3 \otimes x_1 x_4 - dx_1.dx_3 \otimes x_2 x_4
\end{array}
\]
but applying \(f\) or \(h'\) to the left hand factors of the tensor products
this time gives~0 and the final result is the same: Aramova-Herzog's conclusion
is so justified; the \emph{possible} pure nature of the looked-for resolution,
known in advance after examining the Koszul cycles, may also be used to cancel
the examination of the critical homotopy operator, but we will see our method
can be applied in much more general situations, even in a non-homogeneous
situation. In more complicated situations, the result could have been
different: ``the'' minimal resolution is unique only up to chain-complex
isomorphism and this set of isomorphisms is very large. In this particular
case, many triangular perturbations can for example be applied to the simple
expression found for \(d(\fh_{3,1})\) without changing its intrinsic nature,
and in parallel the same for ``the'' effective homology of the Koszul complex.

Another comment is also necessary. After all, any (correct) choice for the
representants \(g(\fh_{2,i})\) is possible, so that why it would not be
possible to prefer Kenzo's choices to the initial unfortunate choices by
Aramova and Herzog? The point is the following: a resolution is not only made
of isomorphism classes of the boundary maps, you must make these maps fit to
each other in such a way there is \emph{equality} between appropriate kernels
and images. So that when you change the cycles representing the homology
classes during the computation of the component \(d_3\) of the resolution for
example, then the computation of \(d_2\) could also be modified.

\subsubsection{Second Aramova-Herzog example.}

On one hand it is significantly simpler than the first one: the concerned
module is a \(\fk\)-vector space of finite dimension 3, so that any computation
is elementary. On the other hand it is a little harder: the interesting
differential to be constructed is quadratic. Note in particular it was not
obvious in the previous example to obtain the effective homology: the concerned
module was a \(\fk\)-vector space of infinite dimension, but the standard
methods of effective homology know how to overcome such a problem; in fact they
were invented exactly to \emph{overcome} such a problem, see~\cite{RBSR6}.

The underlying ground ring now is \(\fR = \fk[x_1, x_2]\) and we consider the
module \(M = \mbox{\(<\! x_1, x_2 \!>\)} / \mbox{\(<\! x_1^2, x_2^2 \!>\)}\).
The module \(M\) is a \(\fk\)-vector space of dimension 3. The Koszul complex
is of dimension 3 in degrees 0 and 2, of dimension 6 in degree 1. The simplest
form of the effective homology is well described by this figure.
\[
\xymatrix@R0pt@C10pt{
 & \ar@{-}[rrrrrr] \ar@{-}[dddddddddddd] && \ar@{-}[dddddddddddd] &&
 \ar@{-}[dddddddddddd] && \ar@{-}[dddddddddddd]
\\
  && \Ksz_0(M) = \fk^3 && \Ksz_1(M) = \fk^6 && \Ksz_2(M) = \fk^3 &
\\
 &\ar@{-}[rrrrrr]&&&&&&
\\
 R_1 &&&& {x_1\,dx_2} \ar@{<->}[dddll] && {- x_1 \, dx_1.dx_2} \ar@{<->}[dddll]
\\
 && && && {x_2 \, dx_1.dx_2} \ar@{<->}[dddll] &
\\
 &\ar@{-}[rrrrrr]&&&&&&
\\
 R_2 && {x_1 x_2} && {x_1 x_2 \, dx_1} && &
\\
 && && {x_1 x_2 \, dx_2} && &
\\
 &\ar@{-}[rrrrrr]&&&&&&
\\
 R_3&& x_1 && {x_2 \, dx_1 - x_1 \, dx_2} && {x_1 x_2 \, dx_1.dx_2} &
\\
 && x_2 && {x_1 \, dx_1} && &
\\
 && && {x_2 \, dx_2} && &
\\
 &\ar@{-}[rrrrrr]&&&&&&
}
\]

Each column corresponds to a component of the Koszul complex and the (almost)
canonical basis  is shared in boundary preimages, cycles homologous to zero,
and homology classes, each homology class being represented by a cycle not at
all homologous to zero. The effective homology:
\[\rho = \raisebox{4pt}{\framebox{\xymatrix@1{
 \scriptstyle h\ \textstyle \autoarrow{0.8}\ \textrm{Ksz}(M) \ar@<-10\ul>[r]_-f & H
 \ar@<-10\ul>[l]_-g }}}
 \]
is read on the figure as follows. The map \(g\) consists in representing the
homology classes by the cycles listed on the bottom row \(R_3\). The map~\(f\)
is the inverse projection which forgets the  basis vectors of the rows \(R_1\)
and \(R_2\). The differentials and the homotopy operator \(h\) are
simultaneously represented by bidirectional arrows. The chosen supplementary of
the homology groups -- in fact of the representing cycles -- are shared in two
components (\(R_1\) and \(R_2\)) isomorphic through the differential in the
decreasing direction, through the homotopy operator in the increasing
direction. This diagram expresses in a very detailed way the Betti numers are
\((2, 3, 1)\).

The chain-complex \(H\) is \([0 \leftarrow \fk^2 \leftarrow \fk^3 \leftarrow
\fk \leftarrow 0]\) with null differentials. We have to install the right
differential on \(H \otimes \fR\). With the same notations as in the previous
section, the differential \(d_2\) of the minimal resolution is obtained by a
unique non-null term of the series \((\Sigma)\) following the path:
\[
\begin{array}{rcl}
 && \fh_{2,1}
\\
 (g_2 \otimes \id{\fR}) & \mapsto & x_1 x_2 \, dx_1.dx_2
\\
 \partial' & \mapsto & x_1 x_2 \, dx_2 \otimes x_1 - x_1 x_2 \, dx_1 \otimes
 x_2
\\
 -(h_1 \otimes \id{\fR}) & \mapsto & - x_2 \, dx_1.dx_2 \otimes x_1 - x_1 \,
 dx_1.dx_2 \otimes x_2
\\
 \partial' & \mapsto & - x_2 \, dx_2 \otimes x_1^2 + (x_2 \, dx_1 - x_1 \,
 dx_2) \otimes x_1 x_2 + x_1 \, dx_1 \otimes x_2^2
\\
 (f_1 \otimes \id{\fR}) & \mapsto & - \fh_{1,3} \otimes x_1^2 - \fh_{1,1} \otimes
 x_1 x_2 + \fh_{1,2} \otimes x_2^2,
\end{array}
\]
that is, the same result as in~\cite{ARHR}, except innocent sign changes and
permutations. All the other terms produced by the series \((\Sigma)\) are null.

The ``path'' described above makes also obvious the nilpotency argument which
guarantees the convergence of the series \((\Sigma)\): in \(M \otimes \wedge V
\otimes \fR\), the central term \(\wedge V\) ``inhales'' the monomials from the
left hand factor \(M\) and partly ``exhales'' them to the right hand side after
some processing, giving back also something on the left hand side but with a
strictly inferior degree. After a finite number of steps, certainly nothing
anymore on the left hand side. This is particularly clear in the homogeneous
case, a little more difficult but interesting in the general case: the Groebner
monomial orders again play an important role here.

You see in fact the nature of this example is \emph{essentially} the same as
for our initial ``minimal non-trivial'' example.

\subsubsection{The favourite Kreuzer-Robbiano example.}

Martin Kreuzer and Lorenzo Robbiano use a little more complicated toy example
in their book~\cite[Chapter 4]{KRRB}, in fact close to the first Aramova-Herzog
example. Again the ring \(\fR = \fk[x_1, x_2, x_3, x_4]\) but the ideal is
nomore monomial: \(I = \mbox{\(<\!\! x_2^3 - x_1^2 x_3,\)}\ x_1 x_3^2 - x_2^2
x_4, x_3^3 - x_2 x_4^2, x_2 x_3 - x_1 x4 \!\!>\). It is a Groebner basis for
DegRevLex, so that step 1 of the algorithm is void, but the ideal is nomore
monomial and step 3 is not. Keeping the leading terms, we consider the close
ideal \(I' = \ideal{x_2^3, x_1 x_3^2, x_3^3, x_2 x_3}\). It is a monomial ideal
and the effective homology of the Koszul complex \(\Ksz(\fR/I')\) is easily
computed; the Betti numbers are \((1, 4, 4, 1)\) and Kenzo gives for example as
a generator of the 3-homology the cycle \(- x_3^2 \, dx_1.dx_2.dx_3\). Applying
the homological perturbation lemma to take account of the difference between
\(I\) and \(I'\) gives the effective homology of \(\Ksz(\fR/I)\); the new Betti
numbers are certainly bounded by the previous ones, but in this simple case,
they are the same. The generator of the homology in dimension 3 is now \(-
x_3^2 \, dx_1.dx_2.dx_3 + x_2 x_4 \, dx_1.dx_2.dx_4 - x_1 x_3 \, dx_1.dx_3.dx_4
+ x_2^2 \, dx_2.dx_3.dx_4\). There remains to play the same game with the
components \(f\), \(g\) and \(h\) of the effective homology, and also with the
differential \(\partial'\) of the Aramova-Herzog bicomplex, exactly the same
game as before, nothing more, to obtain the minimal resolution:
\[
0 \longleftarrow \fR \stackrel{d_1}{\longleftarrow} \fR^4
\stackrel{d_2}{\longleftarrow} \fR^4 \stackrel{d_3}{\longleftarrow} \fR
\longleftarrow 0
\]
with the matrices:
\[
d_1 = \left[\begin{array}{c} x_1^2 x_3 - x_2^3, - x_1 x_3^2 + x_2^2 x_4, x_2
x_4^2 - x_3^3, - x_1 x_4 + x_2 x_3 \end{array}\right]
\]
\[
d_2 = \left[\begin{array}{cccc}
 0 & - x_3 & - x_4 & 0
 \\
 - x_3 & - x_1 & - x_2 & x_4
 \\
 x_1 & 0 & 0 & - x_2
 \\
 x_2 x_4 & - x_2^2 & - x_1 x_3 & - x_3^2
\end{array}\right]
\hspace{2cm} d_3 = \left[\begin{array}{c}
 - x_2 \\ - x_4 \\ x_3 \\ - x_1
\end{array}\right]
\]

\subsubsection*{Another toy example.}\label{25180}

Let us finally consider now the non-homogeneous ideal:
\[
I =\ \ideal{t^5 - x, t^3 y - x^2, t^2 y^2 - x z, t^3 z - y^2, t^2 x - y, t x^2
- z, x^3 - t y^2, y^3 - x^2 z, x y - t z }
\]

This ideal seems more complicated than the previous one, but in a sense in fact
it is not. This ideal is obtained by applying the DegRevLex Groebner process to
\(I =\ \ideal{x - t^5, y - t^7, z - t^{11} }\) and the simple arithmetic nature
of the toric generators allows us to expect a simple minimal resolution. But
the program ignores this expression of \(I\) and it is interesting to observe
the result of its study: the minimal resolution is in principle a machine to
analyze the \emph{deep} structure of an ideal or module. Macaulay2's
\boxtt{resolution} gives for \(\fR/I\) a resolution with Betti numbers \((1, 7,
11, 6, 1)\) which is not minimal\footnote{But the writer of this part of the
text is not at all a Macaulay2 expert; using the rich set of Macaulay2
procedures, it is certainly possible to compute the minimal resolution.}. On
the contrary, Singular's \boxtt{mres} computes the minimal resolution,
necessarily equivalent to ours; but to our knowledge, Singular does not give
any information about the connection between the homology of the Koszul complex
and this minimal resolution, in particular between the \emph{effective}
character of the homology of the Koszul complex and the \emph{effective}
character of the obtained resolution. No indication in~\cite{GRPF} about these
subjects.

The approximate monomial module \(\fR/I'\) has Betti numbers \((1, 9, 15, 8,
1)\). Applying the homological perturbation lemma between \(\Ksz(\fR/I')\) and
\(\Ksz(\fR/I)\) gives the effective homology of the last one. The Betti numbers
are, surprise, \((1, 3, 3, 1)\). For example a generator for the 3-homology is
\(- x^2 \, dt.dx.dy + t x \, dt.dx.dz - t^4 dt.dy.dz + dx.dy.dz\). The same
process as before using the Aramova-Herzog bicomplex now describes a possible
minimal resolution. The differentials can be:
\[
d_1 = \left[ - t^2 x + y, - t x^2 + z, - t^5 + x \right]
\]
\[
d_2 = \left[\begin{array}{ccc}
 0 & t^5 - x & tx^2 - z
 \\
 t^5 - x & 0 & - tx^2 + y
 \\
 - tx^2 + z & - t^2x + y & 0
\end{array}\right]
\]
\[
d_3 = \left[\begin{array}{c}
 - t^2x + y \\ tx^2 - z \\ - t^5 + x
\end{array}\right]
\]

With respect to the series \((\Sigma)\), each term of degree \(k\) in the
previous matrices comes from a term of the series with \(i = k-1\). Here all
the terms of the series are null for \(i \geq 5\): in fact the degree
corresponds to the number of applications of \(\partial'\).

\section{Simplicial sets.}

\subsection{Introduction.}

To illustrate in Section~\ref{70994} how the \emph{chain-complexes} can be
used, the notion of \emph{simplicial complex} was defined. The general
organization of traditional algebraic topology is roughly as explained in the
diagram:
\[
\textrm{Topology} \rightarrow \textrm{Combinatorial Topology} \rightarrow
\textrm{Chain-Complexes} \rightarrow \textrm{Homology Groups}
\]

\emph{Constructive} algebraic topology must improve this framework. On one
hand, \emph{locally effective} objects are systematically used to
\emph{implement}, theoretically or concretely, the infinite objects which are
quickly unavoidable. On the other hand, a systematic connection with
\emph{effective} objects must be maintained during the construction steps,
currently the only method allowing one to easily produce \emph{algorithms}
computing the traditional invariants: homology groups, homotopy groups,
Postnikov (pseudo-)invariants\ldots

The notion of simplicial complex is the most elementary method to settle a
connection between common ``general'' topology and homological algebra. The
``sensible'' spaces can be triangulated, at least up to homotopy, and instead
of using the notion of topological space, too ``abstract'', only the spaces
having the homotopy type of a CW-complex (see~\cite{LNWN}) are considered, and
all these spaces in turn have the homotopy type of a simplicial complex. So
that a lazy algebraic topologist can decide every space is a simplicial
complex.

But many common \emph{constructions} in topology are difficult to make explicit
in the framework of simplicial complexes. It soon became clear in the forties
the tricky and elegant notion of simplicial \emph{set} is much better. It is
the subject of this section. The reference~\cite{MAY} certainly remains the
basic reference in this subject; it is a book of Mathematics' Gold Age, when a
reasonable detail level was naturally required, and in this respect, this book
is perfect; in particular many explicit formulas, quite useful if you want to
\emph{constructively} work, can be found only in this book. A unique flaw: no
concrete examples; the present section must be understood just as a reading
help to Peter May's book, providing the ``obvious'' examples that are necessary
to understand the exact motivation of this subtle notion of simplicial set and
the related definitions; these examples are obvious, except for the beginner.
Combining both, you should be quickly able to work yourself with this wonderful
tool.

\subsection{The category $\Delta$.}

Some strongly structured  sets of indices are  necessary to define the notion
of \emph{simplicial object}; they are conveniently organized as the  category
$\Delta$. An object  of  $\Delta$  is a set $\und{m}$, namely the set of
integers $\und{m} = \{0,1,\ldots, m-1,m\}$; this set is canonically
\emph{ordered} with the usual order between integers.

A  $\Delta$-morphism  $\alpha:   \und{m} \rightarrow  \und{n}$   is an
\emph{increasing} map.   Equal  values are  permitted;  for  example a
$\Delta$-morphism  $\alpha:  \und{2}   \rightarrow  \und{3}$  could be defined
by $\alpha(0) = \alpha(1) = 1$ and $\alpha(2) = 3$. The set of
$\Delta$-morphisms from   \und{m}\   to   \und{n}  is   denoted     by
$\Delta(\und{m}, \und{n})$; the subset of injective (resp. surjective)
morphisms        is  denoted    by   $\dinj(\und{m},\und{n})$   (resp.
$\dsrj(\und{m},\und{n})$).

Some \emph{elementary}  morphisms are  important, namely  the simplest
non-surjective and  non-injective  morphisms. For  geometric   reasons
explained   later, the first ones  are  the \emph{face morphisms}, the second
ones are the \emph{degeneracy morphisms}.

\begin{dfn}
--- \emph{The \emph{face morphism} $\partial^m_i: \und{m-1} \rightarrow
\und{m}$ is defined for $m \geq 1$ and $0 \leq i \leq m$ by:
$$\begin{array}{rll}
\partial^m_i(j) =& j & \textrm{if $j < i$,} \\
\partial^m_i(j) =& j+1 & \textrm{if $j \geq i$.}
\end{array}$$
}
\end{dfn}

The face morphism $\partial^m_i$ is the unique injective morphism from
\und{m-1}\ to \und{m} such that the integer $i$ is not in the image. The face
morphisms generate the injective morphisms, in fact in a unique way if a growth
condition is required.

\begin{prp}
--- Any injective $\Delta$-morphism $\alpha \in \dinj(\und{m},
\und{n})$ has a unique expression: $$ \alpha = \partial^{n}_{i_n} \circ \ldots
\circ \partial^{m+1}_{i_{m+1}} $$ satisfying the relation $i_n > i_{n-1} >
\ldots > i_{m+1}$.
\end{prp}

\proof The index set $\{i_{m+1}, \ldots, i_n\}$ is exactly the difference set
$\und{n} -  \alpha(\und{m})$, that is, the set   of the integers where
surjectivity fails. \QED

Frequently the  upper index $m$ of  $\partial^m_i$  is omitted because clearly
deduced from  the context.  For example  the unique injective morphism $\alpha:
\und{2} \rightarrow \und{5}$ the  image of which is $\{0,2,4\}$  can   be
written   $\alpha    =  \partial_5  \partial_3
\partial_1$.

If two face morphisms  are composed in the   wrong order, they  can be
exchanged:   $\partial_i \circ    \partial_j =  \partial_{j+1}   \circ
\partial_i$ if  $j  \geq i$.   Iterating  this process allows  you  to
quickly  compute    for  example  $\partial_0   \partial_2 \partial_4
\partial_6 = \partial_9 \partial_6 \partial_3 \partial_0$.

\begin{dfn}
--- \emph{The \emph{degeneracy morphism} $\eta^m_i: \und{m+1} \rightarrow
\und{m}$ is defined for $m \geq 0$ and $0 \leq i \leq m$ by:
$$\begin{array}{rll} \eta^m_i(j) =& j & \textrm{if $j \leq i$,} \\ \eta^m_i(j)
=& j-1 & \textrm{if $j > i$.}
\end{array}$$
}
\end{dfn}

The degeneracy morphism  $\eta^m_i$ is the unique surjective  morphism from
\und{m+1}\   to \und{m}    such that  the  integer $i$    has two pre-images.
The  degeneracy    morphisms   generate   the  surjective morphisms, in fact in
a unique way if a growth condition is required.

\begin{prp}
--- Any surjective $\Delta$-morphism $\alpha \in \dsrj(\und{m},
\und{n})$ has a unique expression: $$ \alpha = \eta^{n}_{i_n} \circ \ldots
\circ \eta^{m-1}_{i_{m-1}} $$ satisfying the relation $i_n < i_{n+1} < \ldots <
i_{m-1}$.
\end{prp}

\proof The  index set  $\{i_n, \ldots,   i_{m-1}\}$ is  exactly  the set   of
integers   $j$  such that $\alpha(j)    =  \alpha(j+1)$, that is,  the integers
where injectivity fails. \QED

Frequently the  upper   index  $m$ of $\eta^m_i$   is  omitted because clearly
deduced from the context.  For  example the unique surjective morphism $\alpha:
\und{5} \rightarrow \und{2}$ such that $\alpha(0) = \alpha(1)$ and $\alpha(2) =
\alpha(3) = \alpha(4)$ can be expressed $\alpha = \eta_0 \eta_2 \eta_3$.

If two face  morphisms are composed in   the wrong order, they can  be
exchanged: $\eta_i   \circ \eta_j =  \eta_{j} \circ  \eta_{i+1}$ if $i \geq j$.
Iterating this process  allows  you to quickly  compute for example $\eta_3
\eta_3 \eta_2 \eta_2 = \eta_2 \eta_3 \eta_5 \eta_6$.

\begin{prp}
--- Any $\Delta$-morphism $\alpha$ can be $\Delta$-decomposed
in a unique way: $$ \alpha = \beta \circ \gamma $$ with $\beta$ injective and
$\gamma$ surjective.
\end{prp}

\proof The intermediate $\Delta$-object \und{k}\ necessarily satisfies $k+1 =
\textrm{Card}(\textbf{im}(\alpha))$. The growth condition then gives a unique
choice for $\beta$ and $\gamma$. \QED

\begin{crl}
--- Any $\Delta$-morphism $\alpha: \und{m} \rightarrow \und{n}$ has
a unique expression: $$ \alpha = \partial_{i_n} \circ \ldots \circ
\partial_{i_{k+1}}
         \circ \eta_{j_k} \circ \ldots \circ \eta_{j_{m-1}}
$$ satisfying the conditions $i_n > \ldots > i_{k+1}$ and $j_k < \ldots <
j_{m-1}$. \QED
\end{crl}

Finally  if face and  degeneracy  morphisms are composed  in the wrong order,
they can be exchanged:

$$
\begin{array}{rcll}
\eta_i \circ \partial_j &=& \textrm{id} & \mbox{if $j=i$ or $j=i+1$;}\\
       &=& \partial_{j-1} \circ \eta_{i} & \mbox{if $ j \geq i+2$;}\\
       &=& \partial_j \circ \eta_{i-1} & \mbox{if $j < i$.}
\end{array}
$$

All these commuting  relations  can be used   to convert an  arbitrary
composition of faces and degeneracies into the canonical expression:

$$ \alpha = \eta_9 \partial_6 \eta_3 \partial_7 \eta_9 \partial_8 \eta_6
\partial_2 \eta_4 \partial_9 =
\partial_7 \partial_6 \partial_2 \eta_2 \eta_4 \eta_6.
$$ This  relation means   the  image of $\alpha$   does   not contain the
integers  2,  6 and 7,  and  the  relations  $\alpha(2) =  \alpha(3)$,
$\alpha(4) = \alpha(5)$ and $\alpha(6) = \alpha(7)$ are satisfied.

\begin{crl}
---   A   \emph{contravariant}   functor   $X:   \Delta    \rightarrow
\textrm{\emph{CAT}}$ is   nothing but  a collection    $\{X_m\}_{m \in
\mathbb{N}}$  of objects    of  the target  category  \emph{CAT},  and
collections     of    \emph{CAT}-morphisms   $\{X(\partial^m_i):   X_m
\rightarrow   X_{m-1}\}_{m   \geq  1\,,\,0  \leq    i   \leq m}$   and
$\{X(\eta^m_i): X_m \rightarrow X_{m+1}\}_{m \geq  0\,,\,0 \leq i \leq m}$
satisfying the commuting relations: $$
\begin{array}{rcll}
X(\partial_i) \circ X(\partial_j) &=& X(\partial_j) \circ X(\partial_{i+1})
   & \mbox{\emph{if $i \geq j$,}}\\
X(\eta_i) \circ X(\eta_j) &=& X(\eta_{j+1}) \circ X(\eta_i)
   & \mbox{\emph{if $j \geq i$,}}\\
X(\partial_i) \circ X(\eta_j) &=& \textrm{\emph{id}} & \mbox{\emph{if
$i=j,j+1$,}}\\ X(\partial_i) \circ X(\eta_j) &=& X(\eta_{j-1}) \circ
X(\partial_i)
   & \mbox{\emph{if $j > i$,}}\\
X(\partial_i) \circ X(\eta_j) &=& X(\eta_j) \circ X(\partial_{i-1})
   & \mbox{\emph{if $i > j+1$.}}
\end{array}
$$
\end{crl}

In the five last relations, the  upper indices have been omitted. Such a
contravariant functor is a  \emph{simplicial object} in the category CAT.  If
$\alpha$ is    an  arbitrary $\Delta$-morphism,  it   is then sufficient to
express $\alpha$ as a composition of face and degeneracy morphisms; the image
$X(\alpha)$ is necessarily the composition of the images of  the corresponding
$X(\partial_i)$'s and $X(\eta_i)$'s; the above relations assure the definition
is coherent.

\subsection{Terminology and notations.}

\begin{dfn}
--- \emph{A \emph{simplicial set} is a simplicial object in the category
of sets.}
\end{dfn}

A    simplicial set    $X$  is given     by  a  collection  of    sets
$\{X(\und{m})\}_{m \in  \mathbb{N}}$    and   collections   of    maps
$\{X_\alpha\}$, the index $\alpha$ running the $\Delta$-morphisms; the usual
coherence properties must be  satisfied. As explained at the end of   the
previous   section,   it   is   sufficient  to define   the $X(\partial^m_i)$'s
and  the  $X(\eta^m_i)$'s with  the  corresponding commuting relations.

The set $X(\und{m})$ is usually denoted by $X_m$ and is called the set of
$m$-simplices of $X$; such  a simplex has the \emph{dimension} $m$. To be  a
little more precise,  these  simplices are  sometimes called \emph{abstract}
simplices,  to  avoid  possible confusions   with  the \emph{geometric}
simplices defined   a little  later.  An  (abstract) $m$-simplex is only
\emph{one} element of $X_m$.

If   $\alpha \in \Delta(\und{n},\und{m})$,  the corresponding morphism
$X(\alpha):   X_m \rightarrow  X_n$ is most    often simply denoted by
$\alpha^\ast: X_m \rightarrow X_n$ or  still more simply $\alpha: X_m
\rightarrow X_n$. In particular the faces and degeneracy operators are maps
$\partial_i:   X_m    \rightarrow X_{m-1}$   and   $\eta_i:  X_m \rightarrow
X_{m+1}$. If $\sigma$ is   an $m$-simplex, the  (abstract) simplex $\partial_i
\sigma$ is   its   $i$-th face, and the    simplex $\eta_i \sigma$ is its
$i$-th degeneracy; we will  see the last one is ``particularly'' abstract.

\subsection{The structure of simplex sets.}

\begin{dfn}\label{62937}
---  \emph{An  $m$-simplex $\sigma$  of   the  simplicial set  $X$  is
\emph{degenerate}  if there exists an integer  $n <  m$, an $n$-simplex $\tau
\in     X_n$   and     a    $\Delta$-morphism   $\alpha   \in
\Delta(\und{m},\und{n})$ such that $\sigma = \alpha(\tau)$. The set of
non-degenerate   simplices of  dimension  $m$ in    $X$  is denoted by
$X^{ND}_m$.}
\end{dfn}

Decomposing the morphism $\alpha  = \beta \circ \gamma$  with $\gamma$
surjective, we  see that  $\sigma =  \gamma (\beta  (\tau))$, with the
dimension of   $\beta(\tau)$ less or  equal   to $n$;  so  that in the
definition of  degenaracy, the connecting  $\Delta$-morphism  $\alpha$ can be
required to be surjective. The relation $\sigma = \alpha(\tau)$ with $\alpha$
surjective is    shortly   expressed  by saying   the $m$-simplex $\sigma$
\emph{comes from} the $n$-simplex $\tau$.

Eilenberg's lemma   explains  each degenerate   simplex   comes from a
canonical non-degenerate one.

\begin{lmm}
--- \textbf{\emph{(Eilenberg's lemma)}} If $X$ is a simplicial set and
$\sigma$  is an $m$-simplex   of  $X$, there   exists a unique  triple
$T_\sigma  = (n,  \tau, \alpha)$ satisfying  the following conditions:

{\setlength{\parindent}{40pt} 1.  The   first component   $n$  is a  natural
number  $n   \leq m$; {\setlength{\parskip}{0pt}

2.    The second component   $\tau$  is  a non-degenerate  $n$-simplex
\mbox{$\tau \in X\ND_n$};

3. The third component $\alpha$ is a $\Delta$-morphism $\tau \in \dsrj(\und{m},
\und{n})$;

4. The relation $\sigma = \alpha(\tau)$ is satisfied. }}
\end{lmm}

\begin{dfn}
--- \emph{This triple $T_\sigma$ is called the \emph{Eilenberg triple}
of $\sigma$.}
\end{dfn}

\proof Let $\mathcal{T}$ be the set of  triples $T =  (n, \tau, \alpha)$ such
that     $n  \leq   m$,       $\tau  \in    X_n$  and   $\alpha \in
\Delta(\und{m},\und{n})$  satisfy  $\sigma  =  \alpha(\tau)$.  The set
$\mathcal{T}$ certainly contains the triple $(m, \sigma, \textrm{id})$ and
therefore is non empty. Let  $(n_0, \tau_0, \alpha_0)$ be an element of
$\mathcal{T}$ where the  first component, the  integer $n_0$,  is minimal. We
claim $(n_0, \tau_0, \alpha_0)$ is the Eilenberg triple.

Certainly $n_0 \leq  m$. The $n_0$-simplex $\tau_0$ is non-degenerate;
otherwise   $\tau_0  = \beta(\tau_1)$  with  the   dimension  $n_1$ of $\tau_1$
less than $n_0$,  but   then $(n_1, \tau_1,   \beta\alpha_0)$ would be a triple
with $n_1 < n_0$.  Finally $\alpha_0$ is surjective, otherwise $\alpha_0 =
\beta \gamma$ with $\gamma \in \dsrj(m,n_1)$ and $n_1 < n_0$; but again the
triple $(n_1, \beta(\tau_0), \gamma)$ would be a  triple denying the required
property  of $n_0$. The existence of an Eilenberg triple is proved and
uniqueness remains to be proved.

Let $(n_1,   \tau_1, \alpha_1)$ be    another Eilenberg  triple.   The
morphisms  $\alpha_0$  and  $\alpha_1$ are   surjective and respective sections
$\beta_0 \in   \dinj(\und{n_0},\und{m})$  and $\beta_1   \in
\dinj(\und{n_1},\und{m})$   can be  constructed:   $\alpha_0 \beta_0 =
\textrm{id}$  and $\alpha_1 \beta_1  =   \textrm{id}$. Then $\tau_0  =
(\alpha_0 \beta_0)      (\tau_0)   =   \beta_0(\alpha_0(\tau_0))     =
\beta_0(\sigma)  =   \beta_0 (\alpha_1  (\tau_1)) = (\alpha_1 \beta_0)
(\tau_1)$; but $\tau_0$  is     non-degenerate,   so that  $n_1      =
\mathrm{dim}(\tau_1) \geq n_0  = \mathrm{dim}(\tau_0)$;  the analogous relation
holds when $\tau_0$ and $\tau_1$  are exchanged, so that $n_1 \leq n_0$ and the
equality $n_0 = n_1$ is proved.

The  relation $\tau_0 =   \beta_0  (\alpha_1 (\tau_1))$  with $\tau_0$
non-degenerate   implies  $\alpha_1 \beta_0  = \textrm{id}$, otherwise
$\alpha_1 \beta_0     =      \gamma  \delta$  with    $\delta      \in
\dsrj(\und{n_1},\und{n_2})$  and $n_2 <  n_1 =  n_0$, but this implies $\tau_0$
comes   from  $\gamma(\tau_1)$    of  dimension   $n_2$ again contradicting the
non-degeneracy  property    of $\tau_0$; therefore $\alpha_1 \beta_0 =
\textrm{id}$ but this equality implies $\tau_0 = \tau_1$.

If $\alpha_0  \neq  \alpha_1$,   let $i$  be   an integer   such  that
$\alpha_0(i) = j \neq \alpha_1(i)$;  then the section $\beta_0$ can be chosen
with $\beta_0(j) = i$; but  this implies $(\alpha_1 \beta_0)(j) \neq j$, so
that the relation $\alpha_1  \beta_0 = \textrm{id}$ would not hold.  The last
required equality  $\alpha_0  = \alpha_1$  is also proved. \QED

Each  simplex comes  from    a  unique non-degenerate   simplex,   and
conversely,  for  any non-degenerate $m$-simplex  $\sigma \in X\ND_m$, the
collection \mbox{$\{\alpha(\sigma)\,;\,\alpha \in \dsrj  (\und{n},
\und{m})\,;$}\break\mbox{$n \geq m\}$} is a perfect description of all
simplices  coming from $\sigma$, that  is, of all degenerate simplices
\emph{above}   $\sigma$. This  is   also  expressed in   the following formula,
describing the structure of the simplex set of any simplicial set $X$:

$$ \coprod_{m \in \mathbb{N}} X_m = \coprod_{m \in \mathbb{N}}\ \
\coprod_{\sigma \in X\ND_m}\ \ \coprod_{n \geq m}
\dsrj(\underline{n},\underline{m})(\sigma). $$

\subsection{Examples.}

\subsubsection{Discrete simplicial sets.}\label{42832}

\begin{dfn}
--- \emph{A simplicial set $X$ is \emph{discrete} if $X_m = X_0$ for
every $m \geq 1$, and if for every $\alpha \in \Delta(\und{m}, \und{n})$, the
induced map $\alpha^\ast: X_n \rightarrow X_m$ is the identity.}
\end{dfn}

The reason of this definition is that the \emph{realization} (see
Section~\ref{40119}) of such a simplicial set is the discrete point set $X_0$;
the Eilenberg triple of any simplex $\sigma \in X_m = X_0$ is $(0, \sigma,
\alpha)$ where the map $\alpha$ is the unique element of $\Delta(\und{m},
\und{0})$.

\subsubsection{The simplicial complexes.}

A \emph{simplicial  complex} $K =  (V,S)$ is a  pair where $V$  is the
\emph{vertex set}  (an arbitrary set,  finite or not), and  $S \subset
\mathcal{P}_F(V)$  is a set of finite  sets of vertices satisfying the
properties:

{\setlength{\parindent}{40pt} 1.  For any $v \in V$,  the one element   subset
$\{v\}$ of  $V$ is an element of $S$;{\setlength{\parskip}{0pt}

2. For any $\tau \subset \sigma \in S$, then $\tau \in S$. }}

The \emph{simplex} $\sigma  \in S$ \emph{spans} its  elements. If $S =
\mathcal{P}_F(V)$, then $K$ is  the \emph{simplex} freely generated by $V$, or
more simply the simplex spanned by $V$.

The terminology is a little incoherent because a simplicial \emph{set} is an
object more  sophisticated than a simplicial \emph{complex}, but this
terminology is  so well established that  it is probably too late to modify it.

The simplicial complex $K = (V,S)$ is \emph{ordered} if the vertex set $V$  is
provided with  a  \emph{total} order\footnote{Other situations where  the order
is not  total   are  also interesting   but will  be considered later.}. Then a
simplicial   set again  denoted by $K$  is canonically associated;  the simplex
set $K_m$   is  the   set  of \emph{increasing} maps $\sigma: \und{m}
\rightarrow K$ such that  the image of $\und{m}$ is an element of $S$; note
that such a map $\sigma$ is   not necessarily injective.  If   $\alpha$ is a
$\Delta$-morphism $\alpha \in \Delta(\und{n},  \und{m})$ and $\sigma$  is an
$m$-simplex $\sigma \in  K_m$,  then  $\alpha(\sigma)$   is naturally defined
as $\alpha(\sigma) = \sigma \circ \alpha$. A simplex  $\sigma \in K_m$ is
non-degenerate if  and  only if  $\sigma   \in \dinj(\und{m},  V)$; if $\sigma
\in K_m = \Delta(\und{m}, V)$, the Eilenberg triple $(n, \tau, \alpha)$
satisfies  $\sigma   =  \tau \circ  \alpha$    with $\alpha$ surjective and
$\tau$ injective.

This in particular  works  for  $K=(\und{d},\mathcal{P}(\und{d}))$ the simplex
freely generated by  $\und{d}$   provided with the  canonical vertex  order. We
obtain in  this  way  the   canonical structure  of simplicial set for the
\emph{standard $d$-simplex} $\Delta^d$.

\subsubsection{The spheres.}\label{96741}

Let $d$ be a natural number. The simplest simplicial version $S = S^d$ of  the
$d$-sphere is defined   as  follows: the set of  $m$-simplices $S_m$ is   $S_m
= \{\ast_m\}  \coprod \dsrj(\und{m},  \und{d})$;  if $\alpha \in
\Delta(\und{n},  \und{m})$ and $\sigma$ is  an $m$-simplex $\sigma   \in S_m$,
then  $\alpha(\sigma)$   depends  on the   nature of~$\sigma$:

\begin{enumerate}
\item If $\sigma = \ast_m$, then $\alpha(\sigma) = \ast_n$;
\item Otherwise $\sigma \in \dsrj(\und{m}, \und{d})$ and if $\sigma
\circ \alpha$ is surjective, then $\alpha(\sigma) = \sigma \circ \alpha$, else
$\alpha(\sigma) = \ast_n$ (the emergency solution when the natural solution
does not work).
\end{enumerate}

This is nothing but the canonical quotient $S^d  = \Delta^d / \partial
\Delta^d$,  at   least if  $d  > 0$;   more  generally  the  notion of
simplicial subset is naturally defined and a quotient then appears. In the case
of the construction of $S^d  = \Delta^d / \partial \Delta^d$, the subcomplex
$\partial \Delta^d$  is  made of the  simplices $\alpha \in \Delta(\und{m},
\und{d})$ that are not surjective.

The   Eilenberg triple of  $\ast_m$  is  $(0, \ast_0, \alpha)$ where $\alpha$
is the  unique  element of   $\Delta(\und{m}, \und{0})$. The Eilenberg triple
of $\sigma \in  \dsrj(\und{m}, \und{d}) \subset S_m$ is  $(d,  \textrm{id},
\sigma)$. There    are only two  non-degenerate simplices, namely   $\ast_0 \in
S_0$  and $\textrm{id}(\und{d}) \in S_d$, even if $d=0$.

\subsubsection{Classifying spaces of discrete groups.}\label{33788}

Let  $G$  be a  (discrete) group.   Then a simplicial   version of its
classifying space  $BG$ can be  given. The set of $m$-simplices $BG_m$ is the
set of ``$m$-bars'' $\sigma = [g_1 | \ldots | g_m]$ where every $g_i$ is an
element of $G$. It is  simpler in this situation to define the structure
morphisms only for the face and degeneracy operators:

{\setlength{\parindent}{40pt} 1.  $\partial_0[g_1   |   \ldots    |  g_m]   =
[g_2    |  \ldots   | g_m]$;{\setlength{\parskip}{0pt}

2. $\partial_m[g_1 | \ldots | g_m] = [g_1 | \ldots | g_{m-1}]$;

3. $\partial_i[g_1 | \ldots | g_m] = [\ldots  | g_{i-1} | g_{i}g_{i+1}
                     | g_{i+2} | \ldots]$ if $0 < i < m$;

4.  $\eta_i[g_1 |  \ldots |   g_m]  = [\ldots | g_i  | e_G |g_{i+1}  |
\ldots]$, where $e_G$ is the unit element of $G$. }}

In particular $BG_0 = \{[\,]\}$ has only one element.

The $m$-simplex $[g_1 | \ldots | g_m]$ is degenerate if and only if one of the
$G$-components is the unit element.

The various commuting relations must  be verified; this works but does not give
obvious indications on the very  nature of this construction; in  fact there is
a  more conceptual  description.  Let us define the simplicial set $EG$ by
$EG_m = \textrm{SET}(\und{m},G)$, that is, the maps  of   \und{m}\ to $G$
without  taking account  of   the ordered structure  of \und{m} (the group  $G$
is not  ordered); if $\alpha \in \Delta(\und{n}, \und{m})$ there is  a
canonical way to define $\alpha: EG_m \rightarrow EG_n$; it would be more or
less coherent to write $EG = G^\Delta$.

There is a canonical left action of the group $G$ on $EG$, and $BG$ is the
natural  quotient of $EG$  by  this action.  A simplex $\sigma \in EG_m$  is
nothing but a $(m+1)$-tuple   $(g_0, \ldots, g_m)$  and the action  of  $g$
gives the  simplex   $(gg_0,  \ldots, gg_m)$.  If two simplices are
$G$-equivalent, the  products $g_{i-1}^{-1}g_i$ are  the same; the  quotient
$BG$-simplex $[g_1,  \ldots,  g_m]$ denotes  the equivalence class of   all the
$EG$-simplices $(g,    gg_1, gg_1g_2, \ldots)$, which  can be imagined   as a
simplex where the  \emph{edge} between the vertices $i-1$ and $i$ ($i > 0$) is
labeled by $g_i$ to be considered as a (right) operator  between the adjacent
vertices.  Then the boundary and degeneracy operators are clearly explained and
it is even not necessary to    prove the commuting  relations, they   can be
deduced of the coherent structure of $EG$.

\subsubsection{The Eilenberg-MacLane spaces.}\label{70076}

The   previous example constructs  an  \emph{Eilenberg-MacLane} space, that is,
a space  with   only   one  non-zero homotopy  group.   The \emph{realization}
process (see later)  applied to the simplicial set $BG$ produces a   model  for
$K(G,1)$: all the homotopy groups are null except $\pi_1$ canonically
isomorphic to $G$.    The construction can    be generalized to  construct
$K(\pi,d)$, $d>1$, when   $\pi$ is an \emph{abelian} group. This    requires
the simplicial  definition of homology groups, explained in another lecture
series. See also \cite[Chapter V]{MAY} where these questions are carefuly
detailed.

Let $\pi$ be a   fixed abelian group, and   $d$ a natural  number. The
simplicial set  $E(\pi,   d)$ is  defined    as follows. The  set   of
$m$-simplices $E(\pi,   d)_m$,  shortly denoted  by  $E_m$,  is $E_m =
C^d(\Delta^m, \pi)$,  the  group of  \emph{normalized} $d$-cochains on the
standard $m$-simplex with values in $\pi$. Such a cochain $\sigma$ is nothing
but a map  $\sigma: \Delta^m_d \rightarrow \pi$, defined on the (degenerate or
not)    $d$-simplices of $\Delta^m$, null for   the degenerate simplices.   If
$\alpha$ is   a $\Delta$-morphism $\alpha: \und{n}  \rightarrow  \und{m}$, this
map defines   a simplicial  map $\alpha_\ast: \Delta^n \rightarrow \Delta^m$
which in turns defines a pullback map     $\alpha^\ast: C^d(\Delta^m,  \pi)
\rightarrow C^d(\Delta^n, \pi)$ between $m$-simplices and $n$-simplices of
$E_m$.

The simplicial  set $E(\pi,d)$   so defined contains    the simplicial subset
$K(\pi, d)$, made  only of the \emph{cocycles},  those cochains the coboundary
of     which ($d:    C^d(\Delta^m, \pi)  \rightarrow C^{d+1}(\Delta^m,   \pi$)
is  null.    In   fact $E(\pi,  d)$  is   a \emph{simplicial group}, that is, a
simplicial object in the category of groups, and  $K(\pi, d)$  is  a simplicial
subgroup.  The  quotient simplicial group $E(\pi,d) / K(\pi, d)$  is
canonically isomorphic to $K(\pi, d+1)$  and this  structure   defines the
Eilenberg-MacLane fibration:

$$ K(\pi, d) \hookrightarrow E(\pi, d) \rightarrow K(\pi, d+1) $$

See later   the section  about \emph{simplicial fibrations}    for some
details.

\subsubsection{Simplicial loop spaces.}\label{87671}

Let $X$ be a simplicial set. We can construct a new simplicial set $DT(X)$ (the
acronym $DT$ meaning Dold-Thom) from X, where $DT(X)_m$ is the free
\(\bZ\)-module generated by the $m$-simplices $X_m$; the operators of $DT(X)$
are also ``generated" by the operators of $X$. This is a simplicial version of
the Dold-Thom construction, producing a new simplicial set $DT(X)$, the
homotopy groups of which being the homology groups of the initial $X$. The
simplicial set $DT(X)$ is also of \emph{simplicial group}; its simplex
\emph{sets} are nothing but the chain groups at the origin of the simplicial
homology of $X$, but in $DT(X)$, each simplicial ``chain" of $X$ is \emph{one}
(abstract) simplex. See \cite[Section 22]{MAY}.

The same construction can be undertaken, but instead of using the abelian group
generated by the simplex sets $X_m$, we could consider the free
\emph{non-commutative} group generated by $X_m$. This also works, but then the
obtained space is a simplicial model for the \emph{James construction} of
$\Omega \Sigma X$, the loop space of the (reduced) suspension of $X$. See
\cite{CRML} for the James construction in general and~\cite{DNCT1} for the
simplicial case.

It is even possible to construct the ``pure" loop space $\Omega X$, without any
suspension. This is due to Daniel Kan~\cite{KAN} and works as follows. It is
necessary to assume $X$ is reduced, that is with only one vertex: the
cardinality of $X_0$ is 1. Let $X^\ast_m$ the set of all $m$-simplices, except
those that are 0-degenerate: $X^\ast_m = X_m - \eta_0(X_{m-1})$; this makes
sense for $m \geq 1$. Then let $GX_m$ be the free \emph{non-commutative} group
generated by $X^\ast_{m+1}$; to avoid possible confusions, if $\sigma \in
X^\ast_{m+1}$, let us denote by $\tau(\sigma)$ the corresponding
\emph{generator} of $GX_m$. The simplicial object $GX$ to be defined is a
simplicial \emph{group}, so that it is sufficient to define face and degeneracy
operators for the generators:

$$
\begin{array}{rcll}
 \partial_i \tau(\sigma) &=& \tau(\partial_{i+1} \sigma), & \mbox{if $1 \leq i \leq m$;}\\
 \partial_0 \tau(\sigma) &=& \tau(\partial_1 \sigma)\tau(\partial_0 \sigma)^{-1}; &\\
 \eta_i \tau(\sigma) &=& \tau(\eta_{i+1} \sigma), &\mbox{if $0 \leq i \leq m$}.
\end{array}
$$

These definitions are coherent, and the simplicial set $GX$ so obtained is a
simplicial version of the loop space construction. See \cite[Chapter VI]{MAY}
for details and related questions, mainly the \emph{twisted Eilenberg-Zilber
Theorem}, at the origin of the general solution described in \cite{SRGR4,RBSR6}
for the computability problem in algebraic topology.

\subsubsection{The singular simplicial set.}\label{49261}

Let $X$  be  an arbitrary  topological space. Then  the \emph{singular
simplicial set} associated with $X$ is constructed  as follows. The set of
$m$-simplices $SX_m$  is  made  of  the continuous maps    $\sigma: \Delta^m
\rightarrow X$;  \emph{one} (abstract)  simplex is \emph{one} continuous map
but no topology is  installed on $SX_m$; in particular when   $SX$ will be
\emph{realized}  in   the following section,  the \emph{discrete} topology must
be used.  The   source of the abstract $m$-simplex $\sigma$   is the geometric
$m$-simplex $\Delta^m \subset \bR^{\underline{m}}$ provided with the
traditional topology. If $\alpha \in \Delta({\und{n}, \und{m}})$ is  a
$\Delta$-morphism, this $\alpha$ defines a natural continuous  map
$\alpha_\ast:  \Delta^n \rightarrow \Delta^m$ between geometric simplices, and
this allows us to naturally define $\alpha^\ast(\sigma) = \sigma \circ
\alpha_\ast$. An enormous simplicial set is so defined if $X$ is an arbitrary
topological space; it is at the origin of the \emph{singular homology} theory.

\subsection{Realization.}\label{40119}

If $K = (V,S)$ is a simplicial complex, the realization $|K|$ is a subset of
$\bR^{(V)}$, the \(\bR\)-vector space generated by the vertices $v \in V$; a
point $x \in \bR^{(V)}$ is a function $x: V \rightarrow \bR$ almost everywhere
null, that is, the set of $v$'s where $x$ is non-null is finite. Such a
function can also be denoted by $x = \{x_v\}_{v \in V}$, the set of indexed
values, or also the linear notations $x = \sum x_v.e_v$ or $x = \sum x_v.v$ can
be used. Then $|K|$ is the set of $x$'s in $\bR^{(V)}$ satisfying the following
conditions:

{\setlength{\parindent}{40pt} 1. For every $v \in V$, the inequality $0 \leq
x_v \leq 1$ holds;{\setlength{\parskip}{0pt}

2. The relation $\sum_{v \in V} x_v = 1$ is satisfied;

3. The set $\{v \in V \st x_v \neq 0\}$ is a simplex $\sigma \in S$. }}

The right topology to install on $|K|$ is induced by all the finite dimensional
spaces $\bR^\sigma$ for $\sigma \in S$. In this way the realization $|K|$ is a
CW-complex. In particular, if $\Delta^m$ is the simplex freely generated by
$\und{m}$, the realization is the standard geometric $m$-simplex again denoted
by $\Delta^m$, provided with its ordinary topology. In general the topology of
$|K|$ is induced by its (geometric) simplices.

If $\alpha: \und{m} \rightarrow \und{n}$ is a $\Delta$-morphism, then $\alpha$
defines a covariant induced map $\alpha_\ast: \Delta^m \rightarrow \Delta^n$
(between the ``simplicial" simplices or the geometric realizations, as you
like) and for any simplicial set $X$ a contravariant  induced map $\alpha^\ast:
X_n \rightarrow X_m$. From now on, unless otherwise stated, $\Delta^m$ denotes
the \emph{geometric} standard simplex, that is, the convex hull of the
canonical basis of $\bR^{\underline{m}}$.

If $X$ is a simplicial set, the (\emph{expensive}) realization $|X|$ of $X$ is:

$$ |X| = \coprod_{m \in \bN} X_m \times \Delta^m \ / \ \approx. $$

Each component of the coproduct is the product of the discrete set of
$m$-simplices by the geometric $m$-simplex; in other words, each abstract
simplex $\sigma$ in $X_m$ gives birth to a geometric simplex $\{\sigma\} \times
\Delta^m$, and they are attached to each other following the instructions of
the equivalence relation $\approx$, to be defined. Let $\alpha \in
\Delta(\und{m}, \und{n})$ be some $\Delta$-morphism, and let $\sigma$ be an
$n$-simplex $\sigma \in X_n$ and $t \in \Delta^m \subset \bR^{\underline{m}}$.
Then the pairs $(\alpha^\ast \sigma, t)$ and $(\sigma,  \alpha_\ast t)$ are
declared equivalent.

It is not obvious to understand what is the topological space so obtained. A
description a little more explicit but also a little more complicated explains
more satisfactorily what should be understood.

The \emph{cheap} realization $\|X\|$ of the simplicial set $X$ is:

$$ \|X\| = \coprod_{m \in \bN} X^{ND}_m \times \Delta^m \ / \ \approx $$

\noindent where the equivalence relation $\approx$ is defined as follows. Let
$\sigma$ be a non-degenerate $m$-simplex and $i$ an integer $0 \leq i \leq m$;
let also $t \in \Delta^{m-1}$; the abstract $(m-1)$-simplex $\partial_i^\ast
\sigma$ has a well defined Eilenberg triple $(n, \tau, \alpha)$; then we decide
to declare equivalent the pairs $(\sigma, \partial_{i\ast} (t)) \approx (\tau,
\alpha_\ast(t))$.

For example let $S = S^d$ be the claimed simplicial version of the $d$-sphere
described in Section \ref{96741}. This simplical set $S$ has only two
non-degenerate simplices, one in dimension 0, the other one in dimension $d$.
The cheap realization needs a point $\Delta^0$ and a geometric $d$-simplex
$\Delta^d$ corresponding to the abstract simplex $\textrm{id} \in
\Delta(\und{d}, \und{d})$; then if $t \in \Delta^{d-1}$ and $0 \leq i \leq d$,
the equivalence relation asks for the Eilenberg triple of
$\partial_i(\textrm{id}) = \ast_{d-1}$ which is $(0, \ast_0, \eta)$, the map
$\eta$ being the unique element of $\Delta(\und{d-1}, \und{0})$. Finally the
initial pair $(\textrm{id}, \partial_{i\ast}t)$ in the realization process must
be identified with the pair $(\ast_0, \Delta^0)$; in other words $\|S\| =
\Delta^d / \partial\Delta^d$, homeomorphic to the unit $d$-ball with the
boundary collapsed to one point.

\begin{prp}
--- Both realizations, the expensive one and the cheap one, of
a simplicial set $X$ are canonically homeomorphic.
\end{prp}

\proof The homeomorphism $f: |X| \rightarrow \|X\|$ to be constructed maps the
equivalence class of the pair $(\sigma,t) \in X_m \times \Delta^m$ to the
(equivalence class of the) pair $(\tau, \alpha_\ast(t)) \in X_n \times
\Delta^n$ if the Eilenberg triple of $\sigma$ is $(n, \tau, \alpha)$. The
inverse homeomorphism $g$ is induced by the canonical inclusion $\coprod
X^{ND}_m \times \Delta^m \hookrightarrow \coprod X_m \times \Delta^m$. These
maps must be proved coherent with the defining equivalence relations and
inverse to each other; their continuity is an obvious consequence of the
definition simplex by simplex.

If $\alpha = \beta\gamma$ is a $\Delta$-morphism expressed as the composition
of two other $\Delta$-morphisms, then an equivalence $(\sigma, \beta_\ast
\gamma_\ast t) \approx (\gamma^\ast \beta^\ast \sigma, t)$ can be considered as
a consequence of the relations $(\sigma, \beta_\ast \gamma_\ast t) \approx
(\beta^\ast \sigma, \gamma_\ast t)$ and $(\beta^\ast \sigma, \gamma_\ast t)
\approx (\gamma^\ast \beta^\ast \sigma, t)$, so that it is sufficient to prove
the coherence of the definition of $f$ with respect to the \emph{elementary}
$\Delta$-operators, that is, the face and degeneracy operators.

Let us assume the Eilenberg triple of $\sigma \in X_m$ is $(n, \tau, \alpha)$,
so that $f(\sigma, t) = (\tau, \alpha_\ast t)$. We must in particular prove
that $f(\eta_i^\ast \sigma, t)$ and $f(\sigma, \eta_{i\ast} t)$ are coherently
defined. The second image is the equivalence class of $(\tau, \alpha_\ast
\eta_{i\ast} t)$; the Eilenberg triple of $\eta_i^\ast \sigma$ is $(n, \tau,
\alpha\eta_i)$ so that the first image is the equivalence class of $(\tau,
(\alpha \eta_i)_\ast t)$ and both image representants are even equal.

Let us do now the analogous work with the face operator $\partial_i$ instead of
the degeneracy operator $\eta_i$. Two cases must be considered. If ever the
composition $\alpha \partial_i \in \Delta(\und{m-1}, \und{n})$ is surjective,
the proof is the same. The interesting case happens if $\alpha \partial_i$ is
not surjective; but its image then forgets exactly one element $j$ ($0 \leq j
\leq n$) and there exists a unique surjection $\beta \in \Delta(\und{m-1},
\und{n-1})$ such that $\alpha \partial_i = \partial_j \beta$. The abstract
simplex $\partial_j^\ast \tau$ gives an Eilenberg triple $(n', \tau', \alpha')$
and the unique possible Eilenberg triple for $\partial_i^\ast \sigma$ is $(n',
\tau', \beta\alpha')$. Then, on one hand,  the $f$-image of $(\sigma,
\partial_{i\ast}t)$ is $(\tau, \alpha_\ast \partial_{i\ast} t) = (\tau,
\partial_{j\ast}\beta_\ast t)$; on the other hand the $f$-image of
$(\partial^\ast_i \sigma, t)$ is $(\tau', \alpha_\ast \beta_\ast t)$; but
according to the definition of the equivalence relation $\approx$ for $\|X\|$,
both $f$-images are equivalent. The coherence of $f$ is proved.

Let $\sigma \in X^{ND}_m$, $0 \leq i \leq m$, $t \in \Delta^{m-1}$ and $(n,
\tau, \alpha)$ (the Eilenberg triple of $\partial_i^\ast \sigma$) be the
ingredients in the definition of the equivalence relation for $\|X\|$; the
pairs $(\sigma, \partial_{i\ast} t)$ and $(\tau, \alpha_\ast t)$ are declared
equivalent in $\|X\|$; the map $g$ is induced by the canonical inclusion of
coproducts, so that we must prove the same pairs are also equivalent in $|X|$.
But this is a transitive consequence of $(\sigma, \partial_{i\ast} t) \approx
(\partial^\ast_i \sigma, t) = (\alpha^\ast \tau, t) \approx (\tau, \alpha_\ast
t)$. We see here we had only described the binary relations \emph{generating}
the equivalence relation $\approx$; the defining relation is not necessarily
stable under transitivity. The coherence of $g$ is proved.

The relation $fg = \textrm{id}$ is obvious. The other relation $gf =
\textrm{id}$ is a consequence of the equivalence in $|X|$ of $(\sigma, t)
\approx (\tau, \alpha_\ast t)$ if the Eilenberg triple of $\sigma$ is $(n,
\tau, \alpha)$. \QED

\subsubsection{Examples.}\label{84409}

Let us consider the construction of the classifying space of the group $G =
\bZ_2$ described in Section~\ref{33788}. The universal ``total space" EG has
for every $m \in \bN$ exactly two non-degenerate $m$-simplices
$(0,1,0,1,\ldots)$ and $(1,0,1,0,\ldots)$. The only non degenerate faces are
the $0$-face and the $m$-face. For example the faces of $(0,1,0,1)$ are
$(1,0,1) \in EG^{ND}_2$, $(0,0,1) = \eta_0 (0,1)$, $(0,1,1) = \eta_1 (0,1)$ and
$(0,1,0) \in EG^{ND}_2$. Each non-degenerate $m$-simplex is attached to the
$(m-1)$-skeleton of $EG$ like each hemisphere of $S^m$ is attached to the
equator $S^{m-1}$ and $EG$ is nothing but the infinite sphere $S^\infty$. The
details are not so easy; the key point consists in proving the geometric
$m$-simplex corresponding for example to $\sigma = (0,1,0,1,\ldots)$ with a few
identification relations on the boundary, following the instructions read from
the various iterated faces of $\sigma$, is again homeomorphic to the $m$-ball,
its boundary to the $(m-1)$-sphere; the simplest case is $\Delta^2 / \partial_1
\Delta^2 \cong D^2$, for $\partial_1 \Delta^2 = \Delta^1$ is contractible, and
this can be extended to the higher dimensions.

The classifying space $BG$ is the quotient space of $EG$ by the canonical
action of $\bZ_2$, that is, the quotient space of $S^\infty$ by the
corresponding action; so that $BG$ is homeomorphic to the infinite real
projective space $P^\infty \bR$; the $m$-skeleton (throw away all the
non-degenerate simplices of dimension $> m$ and also \emph{their} degeneracies)
is a combinatorial description of $P^m \bR$. If $\sigma_m = [1|1|\ldots|1|1]$
denotes the unique non-degenerate simplex of $BG$; then $\partial_0 \sigma_m =
\sigma_{m-1}$, $\partial_1 \sigma_m = \eta_0 \sigma_{m-2}$, \ldots,
$\partial_{m-1} \sigma_m = \eta_{m-2} \sigma_{m-2}$ and $\partial_m \sigma_m =
\sigma_{m-1}$.

Let us also consider the case of the singular simplicial set of a topological
space~$X$ (see Section \ref{49261}). There is a canonical continuous map $f:
|SX| \rightarrow X$ defined as follows; if $(\sigma, t)$ represents an element
of $|SX|$, this means the (abstract) simplex~$\sigma$ is a continuous map
$\sigma: \Delta^m \rightarrow X$, but $t$ is an element of the geometric
simplex $\Delta^m$, so that it is tempting to define $f(\sigma, t) =
\sigma(t)$; it is easy to verify this definition is coherent with the
equivalence relation defining $|SX|$. This map is always a weak homotopy
equivalence, and is an ordinary homotopy equivalence if and only if $X$ has the
homotopy type of a CW-complex.

\subsubsection{Simplicial maps.}

A natural notion of \emph{simplicial map} $f: X \rightarrow Y$ between
simplicial sets  can be defined. The map  $f$ must be a system $\{f_m: X_m
\rightarrow  Y_m\}_{m  \in  \mathbb{N}}$ satisfying the  commuting relations
$\alpha^\ast_X \circ  f_m  = f_n \circ \alpha^\ast_Y$   if $\alpha$ is a
$\Delta$-morphism $\alpha \in \Delta(\und{m}, \und{n})$. If $f: X \rightarrow
Y$ is such  a simplicial map, a realization $|f|: |X| \rightarrow |Y|$, a
continuous map, is canonically defined.

\subsection{Associated chain-complexes.}\label{18426}

In the same way simplicial \emph{complexes} produce chain-complexes, see
Section~\ref{43994}, simplicial sets also produce chain-complexes.

\begin{dfn}---
\emph{Let \(X\) be a simplicial set. The \emph{chain-complex \(C_\ast(X) =
C_\ast(X ; \fR)\) associated with \(X\)} is defined as follows:
\begin{itemize}\sitemsep
\item
\(C_m(X)\) is the free \(\fR\)-module generated by \(X_m\), the set of
\(m\)-simplices of \(X\);
\item
The differential \(d: C_{m}(X) \rightarrow C_{m-1}(X)\) is the \(\fR\)-linear
morphism defined by \(d(\sigma) = \sum_{i = 0}^m (-1)^i \partial_i(\sigma)\) if
\(\sigma \in X_m\).
\end{itemize}
}
\end{dfn}

In algebraic topology, most often some ground ring \(\fR\) is underlying,
frequently \(\fR = \bZ, \bQ\ \textrm{or}\ \bZ_p\) for a prime \(p\). The
chain-complex so defined is the \emph{standard} chain-complex, not taking
account of possible shorthands due to \emph{degenerate} simplices. From a
theoretical point of view, this chain-complex is frequently more convenient,
because the technicalities about the nature, degenerate or not, of every
simplex are not necessary. On the contrary, for concrete calculations,
typically for finite simplicial sets, the \emph{normalized} associated
chain-complex can be more convenient. The right statement of the
Eilenberg-Zilber Theorem, see Section~\ref{66038}, also requires normalized
chain-complexes, and it is so important this will become the \emph{default}
option.

\begin{dfn}---
\emph{Let \(X\) be a simplicial set. The \emph{normalized} chain-complex
\(C^N_\ast(X) = C^N_\ast(X ; \fR)\) associated with \(X\) is defined as
follows:
\begin{itemize}\sitemsep
\item
\(C^N_m(X)\) is the free \(\fR\)-module generated by \(X\ND_m\), the set of
\emph{non-degenerate} \(m\)-simplices of \(X\);
\item
The differential \(d: C^N_{m}(X) \rightarrow C^N_{m-1}(X)\) is the
\(\fR\)-linear morphism defined by \(d(\sigma) = \sum_{i = 0}^m (-1)^i
\partial_i(\sigma)\) if \(\sigma \in X_m\) where every possible occurence of a
degenerate simplex in the alternate sum is cancelled.
\end{itemize}
}
\end{dfn}

See page~\pageref{24484} where the minimal triangulation of the real projective
plane, called \boxtt{short-P2R} was used to compute the homology of this
projective plane. The happy event is that for every simplicial set \(X\), both
chain-complexes \(C_\ast(X)\) and \(C^N_\ast(X)\) have the same homology.

\begin{thr} \emph{\textbf{(Normalization Theorem)}}
The graded submodule \(C^D_\ast(X)\) generated by the degenerate simplices is a
sub\emph{complex} of \(C_\ast(X)\): the boundary of a degenerate simplex is a
combination of degenerate simplices; this chain-complex is acyclic. The
canonical isomorphism \(C^N_\ast(X) \cong C_\ast(X) / C^D_\ast(X)\) induces a
canonical isomorphism \(H_\ast(C_\ast(X)) \cong H_\ast(C^N_\ast(X))\).
\end{thr}

The \emph{right} definition of \(C^N_\ast(X)\) in fact is \(C^N_\ast(X) :=
C_\ast(X) / C^D_\ast(X)\). The inductive proof~\cite[VIII.6]{MCLN2} can easily
be arranged to prove:

\begin{thr}\label{06207}
A general algorithm computes:
\[
X \mapsto [\rho_X: C_\ast(X) \rrdc C^N_\ast(X)]
\]
where:
\begin{enumerate}
\item
\(X\) is a simplicial set;
\item
\(\rho_X\) is a chain-complex reduction.
\end{enumerate}
\end{thr}

If a simplex \(\sigma\) is an \(m\)-simplex, the induction can be chosen going
from 0 to~\(m\) or symmetrically from \(m\) to 0. So that there are \emph{two}
such canonical general algorithms. See also~\cite{RBSR4} for a
\emph{categorical} programming of this algorithm.

\subsection{Products of simplicial sets.}

\begin{dfn}\label{43867}
--- \emph{If $X$ and $Y$ are two simplicial sets, the \emph{simplicial
product} $Z  = X \times Y$ is  defined by $Z_m  = X_m  \times Y_m$ for every
natural number $m$, and  $\alpha^\ast_Z = \alpha^\ast_X \times \alpha^\ast_Y$
if $\alpha$ is a $\Delta$-morphism.}
\end{dfn}

The definition of the  product  of two  simplicial sets  is  perfectly trivial
and is  however  at the  origin  of several landmark problems  in algebraic
topology, for   example the deep  structure of   the twisted Eilenberg-Zilber
Theorem,   still quite   mysterious,   and also the enormous field around the
Steenrod algebras.

Every  simplex  of  the product $Z  =   X \times  Y$  is a \emph{pair}
$(\sigma, \tau)$  made of one simplex in  $X$ and one simplex  in $Y$; both
simplices must have the  \emph{same dimension}. It is tempting at this point,
because of the  ``product'' ambience, to denote by $\sigma \times \tau$ such a
simplex in  the product but  \emph{this would be a terrible error!} This is not
at all the right  point of view; the pair $(\sigma, \tau) \in Z_m$ is the
unique simplex in $Z$ whose respective \emph{projections} in $X$ and $Y$ are
$\sigma$ and $\tau$ and this is the   reason why the pair  notation  $(\sigma,
\tau)$  is the only one which is  possible.  For   example the diagonal  of   a
square  is  a 1-simplex, the   unique simplex  the projections  of   which  are
both factors of the square; on the contrary, the ``product'' of the factors is
simply the square, which does not have the dimension 1 and which is even not a
simplex.

\begin{thr}
---  If $X$ and $Y$ are  two simplicial sets and $Z  =  X \times Y$ is
their simplicial product,  then there exists a canonical homeomorphism between
$|Z|$ and $|X| \times |Y|$, the last product being the product of CW-complexes
(or also of $k$-spaces~\cite{STNR2}).
\end{thr}

If you  consider the  product  $|X| \times  |Y|$   as the product   of
topological spaces, the same accident as for CW-complexes (see~\cite{LNWN}) can
happen.

\proof There are natural simplicial projections   $X \times Y \rightarrow  X$
and $Y$  which define a canonical  continuous map $\phi:  |X \times Y|
\rightarrow |X| \times |Y|$. The interesting question is to define its inverse
$\psi: |X| \times |Y| \rightarrow |X \times Y|$.

First of  all, let us  detail  the case of  $X  = \Delta^2$ and  $Y  =
\Delta^1$ where the essential points are visible. The first factor $X$ has
dimension 2,  and the second one $Y$  has dimension 1 so  that the product $Z$
shoud have dimension 3. What about the 3-simplices of $Z$? There are 3  such
\emph{non-degenerate} 3-simplices, namely  $\rho_0 = (\eta_0 \sigma, \eta_2
\eta_1 \tau)$, $\rho_1 = (\eta_1 \sigma, \eta_2 \eta_0 \tau)$ and  $\rho_2 =
(\eta_2 \sigma,  \eta_1 \eta_0 \tau)$, if $\sigma$ (resp. $\tau$)  is the
unique non-degenerate 2-simplex (resp. 1-simplex) of $\Delta^2$ (resp.
$\Delta^1$). This is nothing  but the decomposition   of   a prism $\Delta^2
\times \Delta^1$ in  three tetrahedrons.

Note no non-degenerate 3-simplex is present in $X$ and $Y$ and however some
3-simplices must be produced for  $Z$; this is possible thanks to the
\emph{degenerate} simplices  of $X$ and  $Y$ where they are  again playing a
quite tricky role in our  workspace; in particular a pair of \emph{degenerate}
simplices     in  the    factors  can    produce   a \emph{non-degenerate}
simplex in the product!  This happens when there is no common degeneracy in the
factors.

For example  the tetrahedron $\rho_0  = (\eta_0  \sigma, \eta_2 \eta_1 \tau)$
inside $Z$ is  \emph{the} unique 3-simplex the first projection of  which is
$\eta_0 \sigma$, and  the second   projection is $\eta_2 \eta_1 \tau$;  the
first projection is a  tetrahedron collapsed on the triangle $\sigma$,
identifying two points  when the sum of barycentric coordinates of  index 0 and
1  (the indices where injectivity fails in $\eta_0$) are equal; the  second
projection is a tetrahedron collapsed on an  interval, identifying  two  points
when the sum  of barycentric coordinates of index 1, 2 and 3 are equal.

Let us take a point of  coordinates $r = (r_0,  r_1, r_2, r_3)$ in the simplex
$\rho_0$. Its first projection is the point  of $X = \Delta^2$ of barycentric
coordinates  $s = (s_0  = r_0 + r_1,  s_1 = r_2, s_2  = r_3)$;  in the same way
its  second projection is  the  point of $Y = \Delta^1$ of barycentric
coordinates $t = (t_0 = r_0,  t_1 = r_1 + r_2 + r_3)$. So that:

$$ \phi(\rho_0, (r_0,r_1,r_2,r_3)) =
   ((\sigma, (r_0+r_1, r_2, r_3)), (\tau,(r_0, r_1+r_2+r_3)))
$$

In the same way:

\begin{eqnarray*}
\phi(\rho_1, (r_0,r_1,r_2,r_3)) &=&
   ((\sigma, (r_0, r_1+r_2, r_3)), (\tau, (r_0+r_1, r_2+r_3)))\\
\phi(\rho_2, (r_0,r_1,r_2,r_3)) &=&
   ((\sigma, (r_0, r_1, r_2+r_3)), (\tau, (r_0+r_1+r_2, r_3)))
\end{eqnarray*}

The challenge  then   consists  in deciding   for an  arbitrary  point
$((\sigma, (s_0, s_1,  s_2)), (\tau, (t_0, t_1)))  \in |X| \times |Y|$ what
simplex  $\rho_i$ it comes from and   what a good $\phi$-preimage $(\rho_i, r)$
could be. You obtain the  solution in comparing the sums $u_0  = s_0$, $u_1 =
s_0 + s_1$, $u_2  = t_0$ ;  the sums  $s_0 + s_1 +s_2$ and $t_0 + t_1$ are
necesarily equal to 1  and do not play any role. You see in the three cases,
the values of $u_i$'s are:

$$
\begin{array}{rl}
((\eta_0 \sigma, \eta_2 \eta_1 \tau), r) \Rightarrow &
   u_0 = r_0+r_1, u_1 = r_0+r_1+r_2, u_2 = r_0, \\
((\eta_1 \sigma, \eta_2 \eta_0 \tau), r) \Rightarrow &
   u_0 = r_0, u_1 = r_0+r_1+r_2, u_2 = r_0+r_1, \\
((\eta_2 \sigma, \eta_1 \eta_0 \tau), r) \Rightarrow &
   u_0 = r_0, u_1 = r_0+r_1, u_2 = r_0 +r_1+r_2,
\end{array}
$$

\noindent so  that you can  guess the degeneracy operators to  be applied to
the factors  $\sigma$ and   $\tau$ from  the  order of  the $u_i$'s;  more
precisely, sorting the  $u_i$'s puts the $u_2$  value in position 0, 1 or 2,
and this gives  the index for the degeneracy   to be applied  to $\sigma$; in
the same way the $u_0$ and $u_1$ values must be installed in positions ``1 and
2'',  or ``0  and 2'', or  ``0  and 1'' and  this gives the  two indices (in
reverse order)  for the degeneracies  to be applied to $\tau$. It's  a question
of \emph{shuffle}. Furthermore you can find the components $r_i$  from the
differences between successive $u_i$'s. Now we can construct the map $\psi$:

$$
\begin{array}{rcll}
\phi((\sigma,s)(\tau,t)) &=& (\rho_0, (u_2, u_0-u_2, u_1-u_0, 1-u_1))
   & \textrm{if $u_2 \leq u_0 \leq u_1$}, \\
                         &=& (\rho_1, (u_0, u_2-u_0, u_1-u_2, 1-u_1))
   & \textrm{if $u_0 \leq u_2 \leq u_1$}, \\
                         &=& (\rho_2, (u_0, u_1-u_0, u_2-u_1, 1-u_2))
   & \textrm{if $u_0 \leq u_1 \leq u_2$}.
\end{array}
$$

There  seems an ambiguity  occurs  when there  is an equality  between $u_2$
and  $u_0$  or $u_1$, but   it  is easy   to  see both possible preimages are
in fact the same in $|Z|$.

Now this can be extended to the general case, according to the following
recipe. Let $\sigma \in X_m$ and $\tau \in Y_n$ be two simplices, $s \in
\Delta^m$ and $t \in \Delta^n$ two geometric points. We must define
$\psi((\sigma, s), (\tau, t)) \in |Z| = |X \times Y|$. We set $u_0 = s_0$, $u_1
= s_0+s_1$, \ldots, $u_{m-1} = s_0+\ldots+s_{m-1}$, $u_m = t_0$, $u_{m+1} =
t_0+t_1$, \ldots, $u_{m+n-1} = t_0+\ldots+t_{n-1}$. Then we sort the $u_i$'s
according to the increasing order to obtain a sorted list $(v_0 \leq \ldots
\leq v_{m+n-1})$. In particular $u_m = v_{i_0}, \ldots, u_{m+n-1} =
v_{i_{n-1}}$ with $i_0 < \ldots < i_{n-1}$, and $u_0 = v_{j_0}, \ldots, u_{m-1}
= v_{j_{m-1}}$ with $j_0 < \ldots < j_{m-1}$. Then:

$$
\begin{array}{l}
\psi((\sigma, s), (\tau, t)) = \\ ((\eta_{i_{n-1}} \ldots \eta_{i_0} \sigma,
 \eta_{j_{m-1}} \ldots \eta_{j_0} \tau),
 (v_0, v_1-v_0, \dots, v_{m+n-1}-v_{m+n-2}, 1-v_{m+n-1})).
\end{array}
$$

Now  it is easy  to prove  $\psi \circ \phi  =  \textrm{id}_{|Z|}$ and $\phi
\circ \psi = \textrm{id}_{|X|  \times |Y|}$, following the proof structure
clearly visible  in  the case of  $X  = \Delta^2$ and $Y  = \Delta^1$.

It  is  also necessary   to prove  the  maps  $\phi$   and $\psi$  are
continuous. But $\phi$  is   the product of  the  realization   of two
simplicial maps and is therefore continuous. The map $\psi$ is defined in a
coherent way for each \emph{cell} $\sigma \times \tau$ (this time it is really
the   \emph{product} $|\sigma| \times |\tau| \subset  |X| \times  |Y|$) and is
clearly  continuous on each  cell; because of the definition of the
CW-topology, the map $\psi$ is continuous. \QED

If three simplicial sets $X$, $Y$ and $Z$ are given, there is only one natural
map $|X \times Y \times Z| \rightarrow |X| \times |Y| \times |Z|$ so that
``both" inverses you construct by applying twice the previous construction of
$\psi$, the first one going through $|X \times Y| \times |Z|$, the second one
through $|X| \times |Y \times Z|$ are necessarily the same: the
$\psi$-construction is \emph{associative}, which is interesting to prove
directly; it is essentially the associativity of the Eilenberg-MacLane
formula~\cite[Theorem~5.2]{ELMC1}.

\subsubsection{The case of simplicial groups.}

Let $G$ be  a \emph{simplicial group}. The object  $G$ is a simplicial object
in the group category; in other words each simplex set $G_m$ is provided with a
group structure  and     the  $\Delta$-operators $\alpha^\ast: G_m \rightarrow
G_n$ are group homomorphisms.

This gives  in particular  a  continuous canonical map  $|G  \times G|
\rightarrow  |G|$;  then identifying $|G  \times   G|$ and $|G| \times |G|$, we
obtain a ``continuous''  group structure for $|G|$; the word \emph{continuous}
is put between quotes, because this does not work in general  in the
topological  sense:  this works  always only in the category of ``CW-groups''
where the group structure  is  a map  $|G| \times |G|  \rightarrow |G|$, the
source  of which  being evaluated in the CW-category;  because of  this
definition of   product, it is then true that $|G| \times |G| = |G \times G|$.
The composition rule so defined on $|G|$ satisfies the group axioms; in
particular the associativity property comes from the considerations about the
associativity of the $\psi$-construction in the previous section.

\subsection{Kan extension condition.}\label{40872}

Let us consider the standard simplicial model $S^1$ of the circle, with one
vertex and one non-degenerate 1-simplex $\sigma$. This unique 1-simplex clearly
represents a generator of $\pi_1(S^1)$, but its double cannot be so
represented. This has many disadvantages and correcting this defect was
elegantly solved by Kan.

\begin{dfn}
--- \emph{A \emph{Kan $(m,i)$-hat} (Kan hat in short)
in a simplicial set $X$ is a $(m+1)$-tuple $(\sigma_0, \ldots, \sigma_{i-1},
\sigma_{i+1}, \ldots, \sigma_{m+1})$ satisfying the relations $\partial_j
\sigma_k = \partial_{k-1} \sigma_{j}$ if $j < k$, $j \neq i \neq k$.}
\end{dfn}

For example the pair $(\partial_0 \textbf{id}, \partial_1 \textbf{id},
\partial_2 \textbf{id},)$ is a Kan $(3,3)$-hat in the standard 3-simplex
$\Delta^3$ if $\textbf{id}$ is the unique non-degenerate 3-simplex. Also the
pair $(\sigma, \sigma)$ is a Kan $(2,1)$-hat of the above $S^1$.

\begin{dfn}
--- \emph{If $(\sigma_0, \ldots, \sigma_{i-1}, \sigma_{i+1},
\ldots, \sigma_{m+1})$ is a Kan $(m,i)$-hat in the simplicial set $X$, a
\emph{filling} of this hat is a simplex $\sigma \in X_{m+1}$ such that
$\partial_j \sigma = \sigma_j$ for $j \neq i$.}
\end{dfn}

The 3-simplex $\textbf{id}$ of $\Delta^3$ is a filling of the example Kan hat
in $\Delta^3$. The example Kan hat of $S^1$ has no filling. A Kan $(m,i)$-hat
is a system of $m$-simplices arranged like all the faces except the $i$-th one
of a \emph{hypothetical} $(m+1)$-simplex.

\begin{dfn}
--- \emph{A simplicial set $X$ satisfies the \emph{Kan extension
condition} if any Kan hat has a filling.}
\end{dfn}

The standard simplex $\Delta^d$ \emph{does not}\footnote{Claude Quitt\'e and
Zhao Gong Yun noticed the erroneous opposite statement in a previous version,
thanks!} satisfy the Kan condition: Consider the \((1,2)\)-hat of \(\Delta^2\)
\((\sigma_0 =
\partial_1(\text{id}), \sigma_1 =
\partial_0(\text{id}))\) ; the only necessary incidence relation \(\partial_0 \sigma_0 =
\partial_0 \sigma_1\) is satisfied, but the wrong order in the indices prevents from filling this hat.
Exercise: prove the same for \(\Delta_1\). Most elementary simplicial sets do
not satisfy the Kan condition.

The  simplicial sets    satisfying the Kan   extension condition  have numerous
interesting properties; for example their homotopy groups can be
combinatorially defined    \cite[Chapter  1]{MAY},   a  canonical
\emph{minimal}  version is  included,  also  satisfying the  extension
condition \cite[Section   9]{MAY},  a  simple   decomposition  process produces
a   Postnikov tower \cite[Section~8]{MAY}.

The simplicial groups are important in this respect: in fact a simplicial group
always satisfies the Kan extension condition \cite[Theorem 17.1]{MAY}. For
example the simplicial description of $P^\infty \bR$ (see Section~\ref{84409})
is a simplicial group and therefore satisfies the Kan condition, which is not
so obvious; it is even minimal. The singular complex $SX$ of a topological
space $X$ also satisfies the Kan condition but in general is not minimal. These
simplicial sets satisfying the Kan condition are so interesting that it is
often useful to know how to \emph{complete} an arbitrary given simplicial set
$X$ and produce a new simplicial set $X'$ with the same homotopy type
satisfying the Kan condition. The Kan-completed $X'$ can be constructed as
follows.

Let us define first an elementary completion $\chi(X)$ for $X$. For each Kan
$(m,i)$-hat of $X$, we decide to add the hypothetical $(m+1)$-simplex (even if
a ``solution" preexists), and the ``missing" $i$-th face; such a completion
operation does not change the homotopy type of $X$. Doing this completion
construction for every Kan hat of $X$, we obtain the first completion
$\chi(X)$. Then we can define $X_0 = X$, $X_{i+1} = \chi(X_i)$ and $X' =
\textrm{lim}_\rightarrow X_i$ is the desired Kan completion. You can also run
this process in considering only the failing hats.

\subsection{Simplicial fibrations.}\label{92270}

A \emph{fibration} is a map $p: E \rightarrow B$ between a \emph{total space}
$E$  and a \emph{base  space}  $B$ satisfying a few properties describing more
or less  the  total space    $E$ as a  \emph{twisted product}    $F \times_\tau
B$.  In  the  simplicial context,  several definitions   are possible.    The
notion   of \emph{Kan}   fibration corresponds to  a   situation where  a
simplicial homotopy   lifting property is satisfied; to state this property,
the elementary datum is a Kan hat in the total space and a  given filling of
its projection in the base space;  the  Kan fibration  property  is satisfied
if  it is possible to fill the Kan hat in the total space in a coherent way
with respect to the given filling in   the  base space.   This notion is  the
simplicial version of     the  notion of \emph{Serre    fibration},  a
projection where  the homotopy lifting property   is satisfied for the maps
defined on polyhedra.  The   reference  \cite{MAY} contains   a detailed study
of  the  basic  facts  around Kan  fibrations,  see \cite[Chapters I and
II]{MAY}.

We will examine with a little more details the notion of \emph{twisted
cartesian product},  corresponding to the  topological notion of fibre bundle.
It is a  key notion in topology,  and the simplicial framework is  particularly
favourable for several    reasons. In particular the Serre spectral sequence
becomes well structured  in this environment, allowing us to extend   it  up to
a \emph{constructive} version, one of the main subjects  of another lecture
series of  this Summer School. We give here the basic necessary definitions for
the notion of twisted cartesian product.

A reasonably general situation consists  in considering the case where a
simplicial group $G$ acts on the  fibre space, a simplicial set $F$, the fibre
space. As   usual  this means  a    map $\phi: F \times G \rightarrow F$ is
given; source  and target are simplicial  sets, the first one being  the
product  of $F$ by    the  simplicial set   $G$, underlying the  simplicial
group; the map $\phi$  is a simplicial map; furthermore each  component
$\phi_m: (F \times  G)_m =  F_m \times G_m \rightarrow F_m$ must satisfy the
traditional  properties of the right actions of a group on a set. We  will use
the shorter notation $f.\,g$ instead of  $\phi(f,  g)$. Let  also   $B$ be  our
base  space,  some simplicial set.

\begin{dfn}\label{51596}
--- \emph{A  \emph{twisting operator} $\tau: B    \rightarrow G$ is  a
family  of maps  $\{\tau_m:   B_m \rightarrow  G_{m-1}\}_{m \geq   1}$
satisfying the following properties.
\begin{enumerate}\sitemsep
\item
$\partial_0 \tau(b) = \tau(\partial_1 b) \tau(\partial_0 b)^{-1}$;
\item
$\partial_i \tau(b) = \tau(\partial_{i+1}(b))$ if $i \geq 1$;
\item
$\eta_i \tau(b) = \tau(\eta_{i+1}b)$;
\item
$\tau(\eta_0 b) = e_m$  if $b \in B_m$,  the unit element of $G_m$ being $e_m$.
\end{enumerate}}
\end{dfn}

In particular it is not  required  $\tau$ is a \emph{simplicial  map}, and in
fact,  because  of  the  degree  -1 between source   and target dimensions,
this does not make sense.

\begin{dfn}\label{69699}
--- \emph{If a twisting operator $\tau:  B \rightarrow G$ is given, the
corresponding \emph{twisted cartesian product}  $E = F \times_\tau B$ is  the
simplicial  set defined  as follows.  Its set of $m$-simplices $E_m$ is the
same  as for the  non-twisted product  $E_m = F_m  \times B_m$; the face and
degeneracy operators are  also the same as for the non-twisted  product   with
only one   exception:  $\partial_0(f, b) = (\partial_0 f . \, \tau(b),
\partial_0 b)$.}
\end{dfn}

  The twisting operator $\tau$, the  unique ingredient at the origin of
a  difference between  the non-twisted product  and the $\tau$-twisted one,
acts in the following way:  the twisted product is constructed in a recursive
way  with respect to the base  dimension. Let $B^{(k)}$ be the $k$-skeleton of
$B$ and let us suppose  $F \times_\tau B^{(k)}$ is already constructed. Let
$\sigma$ be a $(k+1)$-simplex of $B$; we must describe how  the product $F
\times \sigma$   is to be  attached to $F \times B^{(k)}$; what  is above the
faces $\partial_i  \sigma$ for $i \geq  1$ is naturally  attached;  but  what
is  above the  $0$-face is shifted by the  translation defined by  the
operation of $\tau(b)$. It is  not   obvious such an  attachment  is coherent,
but  the various formulas of Definition~\ref{69699} are exactly the relations
which must be satisfied by  $\tau$ for consistency.  It was not obvious,
starting from  scratch,  to  guess  this    is  a  good framework for working
simplicially about  fibrations; this  was  invented  (discovered?) by Daniel
Kan \cite{KAN};   the previous work  by Eilenberg   and MacLane
\cite{ELMC1,ELMC2} in the  particular case of the fibrations  relating the
elements   of  the  Eilenberg-MacLane    spectra   was probably determining.

\subsubsection{The simplest example.}

Let us  describe in this  way the exponential fibration $\textbf{exp}: \bR
\rightarrow S^1: t \rightarrow e^{2 \pi it}$. We take for $S^1$ the model with
one  vertex  $\ast_0$  and    one non-degenerate  edge $\textrm{id(\und{1})} =
\sigma$ (see  Section~\ref{96741}). For $\bR$, we choose $\bR_0 = \bZ$ and
$\bR_1^{ND}  = \bZ$, that  is one vertex $k_0$ and one non-degenerate edge
$k_1$ for each  integer $k \in  \bZ$; the faces are defined by $\partial_i(k_1)
= (k+i)_0$ ($i = 0  \ \textrm{or} \ 1$). The discrete (see Section~\ref{42832})
simplicial  group $\bZ$ acts on the fibre; for any  dimension $d$, the group of
\(d\)-simplices is $\bZ$ with the natural structure, and $k_i \, . \, g =
(k+g)_i$ for $i =  0 \ \textrm{or}  \ 1$. It is then clear that the  right
twisting operator for the exponential  fibration is $\tau(g) = 1$ for $g \in
\bR_1^{ND}$.

\subsubsection{Fibrations between $K(\pi, n)$'s.}

Let us recall (see Section \ref{70076}) $E(\pi,d)$ is the simplicial set
defined by $E(\pi,d)_m = C^d(\Delta^m, \pi)$ (only \emph{normalized} cochains)
and $K(\pi,n)$ is the simplicial subset made of the \emph{cocycles}. The maps
between simplex sets to be associated with $\Delta$-morphisms are naturally
defined. A simplicial projection $p: E(\pi,d) \rightarrow K(\pi, d+1)$
associating to an $m$-cochain $c$ its coboundary $\delta c$, necessarily a
cocycle, is also defined. The simplicial set $\Delta^m$ is contractible, its
cochain-complex is acyclic and the kernel of $p$, the potential \emph{fibre
space}, is therefore the simplicial set $K(\pi,d)$. The base space is clearly
the quotient of the total space by the fibre space (\emph{principal}
fibration), and a systematic examination of such a situation (see \cite[Section
18]{MAY}) shows $E(\pi,d)$ is necessarily a twisted cartesian product of the
base space $K(\pi, d+1)$ by the fibre space $K(\pi,d)$.

It is not so easy to guess a corresponding twisting operator. A solution is
described as follows; let $z \in Z^{d+1}(\Delta^m, \pi)$ a base $m$-simplex;
the result $\tau(z) \in Z^{d}(\Delta^{m-1}, \pi)$ must be a $d$-cocyle of
$\Delta^{m-1}$, that is a function defined on every $(d+1)$-tuple $(i_0,
\ldots, i_{d})$, with values in $\pi$, and satisfying the cocycle condition;
the solution $\tau(z) (i_0, \ldots, i_{d}) = z(0, i_0+1, \ldots, i_{d}+1) -
z(1, i_0+1, \ldots, i_{d}+1)$ works, but seems a little mysterious. The good
point of view consists in considering the notion of \emph{pseudo-section} for
the studied fibration; an actual section cannot exist if the fibration is not
trivial, but a pseudo-section is essentially as good as possible; see the
definition of pseudo-section in \cite[Section 18]{MAY}. When a pseudo-section
is found, a simple process produces a twisting operator; in our example, the
twisting operator comes from the pseudo-section $\rho(z)(i_0, \ldots, i_d) =
z(0, i_0+1, \ldots, i_d+1)$, quite natural.

The fibrations between Eilenberg-MacLane spaces are a particular case of
universal fibrations associated with simplicial groups. See \cite[Section
21]{MAY}.

\subsubsection{Simplicial loop spaces.}

A simplicial set $X$ is \emph{reduced} if its 0-simplex set $X_0$ has only one
element. We have given in Section~\ref{87671} the Kan combinatorial version
$GX$ of the loop space of $X$. This loop space is the fibre space of a
\emph{co-universal} fibration:

$$ GX \hookrightarrow GX \times_\tau X \rightarrow X. $$

Only the twisting operator $\tau$ remains to be defined. The definition is
simply\ldots\ $\tau(\sigma) := \tau(\sigma)$ for both possible meanings of
$\tau(\sigma)$; the first one is the value of the twisting operator to be
defined for some simplex $\sigma \in X_{m+1}$ and the second one is the
generator of $GX_m$ corresponding to $\sigma \in X_{m+1}$, the unit element of
$GX_m$ if ever $\sigma$ is 0-degenerate (see Section~\ref{87671}). The
definition of the face operators for $GX$ are exactly those which are required
so that the twisting operator so defined is coherent.

It is again an example of \emph{principal fibration}, that is the fibre space
is equal to the structural group and the action $GX \times GX \rightarrow GX$
is equal to the group multiplication. This fibration is co-universal, with
respect to $X$; in fact, let $H \hookrightarrow H \times_{\tau'} X
\stackrel{p}{\rightarrow} X$ another \emph{principal} fibration above $X$ for
another twisting operator $\tau': X \rightarrow H$. Then the free group
structure of $GX$ gives you a unique group homomorphism $GX \rightarrow H$
inducing a canonical morphism between both fibrations.

If the simplicial space $X$ is 1-reduced (only one vertex, no non-degenerate
1-simplex), then an important result by John Adams \cite{ADMS} allows one to
compute the homology groups of $GX$ if the initial simplicial set $X$ is of
finite type; an intermediate ingredient, the \emph{Cobar construction}, is the
key point. One of the main problems in Algebraic Topology consists in solving
the analogous problem for the iterated loop spaces $G^nX$ when $X$ is
$n$-reduced; it is the problem of \emph{iterating the Cobar construction}; one
of the lecture series of this Summer School is devoted to this subject.

\section{Serre spectral sequence.}\label{66038}

\subsection{Introduction.}

We begin now the part of this text devoted to Algebraic Topology. The general
idea is that Topology is difficult; on the contrary, Algebra is easy, an
appreciation certainly shared by Deligne, Faltings, Wiles, Lafforgue\ldots\ Let
us be serious; as explained in the introduction of this text, the matter is not
at all to switch from Topology to Algebra, the actual subject is to make
topology \emph{constructive}, in particular the natural problem of
\emph{classification}. Because common algebra has a naturally constructive
framework, it is understood switching from topology to algebra could be useful.
The goal of \emph{constructive} algebraic topology consists in organizing the
translation process in such a way that common constructive algebra actually
allows you to constructively work in topology.

The Eilenberg-Zilber Theorem is unavoidable in Algebraic Topology, it allows to
compute \(H_\ast(X \times Y)\) when \(H_\ast(X)\) and \(H_\ast(Y)\) are known.
In a sense it is the last case where ``ordinary'' algebraic topology succeeds:
\emph{ordinary} homology groups of the ingredients \(X\) and \(Y\) are
sufficient to determine the homology groups of the product \(X \times Y\). The
next natural case concerns the Serre spectral sequence: if the product \(X
\times_\tau Y\) is \emph{twisted}, then the ordinary homology groups of \(X\)
and \(Y\) in general are not sufficient to design an algorithm determining
\(H_\ast(X \times_\tau Y)\).

We present in this section both Eilenberg-Zilber Theorems, the original one,
non-twisted, and also the twisted one, in fact due to Edgar
Brown~\cite{BRWNE2}, put under its modern form by Shih Weishu~\cite{SHIH} and
Ronnie Brown~\cite{BRWNR1}. The \emph{effective} Serre spectral sequence is
then an obvious consequence of the \emph{twisted} Eilenberg-Zilber Theorem.

\begin{uos}---
In the part of this text devoted to Algebraic Topology, if \(X\) is a
simplicial set, \(C_\ast(X)\) denotes the \emph{normalized} chain-complex
\(C^N_\ast(X)\) canonically associated with \(X\); that is, \(C_\ast(X)\)
denotes which should be denoted by \(C^N_\ast(X) := C_\ast(X) / C^D_\ast(X)\).
\end{uos}

Because of Theorem~\ref{06207}, this choice has no incidence upon the
theoretical nature of the results. For concrete calculations, one or other
choice can significantly change computing time and/or space.

\subsection{The Eilenberg-Zilber Theorem.}

If $X$ and $Y$ are two simplicial sets, the \emph{cartesian product} $X \times
Y$  is  naturally defined by  $(X \times  Y)_n = X_n  \times Y_n$, and the face
and degeneracy operators are  the  products of the corresponding operators of
each  factor simplicial set; see Definition~\ref{43867}. If $\sigma \in X_n$
and $\tau \in Y_n$ are two $n$-simplices,  the    notation $(\sigma,\tau)$ must
be preferred  to the tempting notation $\sigma \times  \tau$: the pair notation
$(\sigma,\tau)$ has the advantage to clearly  mean  this is  the $n$-simplex
whose  first (resp. second) \emph{projection} is  $\sigma$ (resp. $\tau$). The
``product'' $\sigma \times  \tau$, even  if both  simplices have not  the same
dimension, should normally denote  the element of  $C_\ast(X \times Y)$ which
is the Eilenberg-MacLane image  of the element  $\sigma \otimes  \tau \in
C_\ast  X \otimes C_\ast Y$,  that is, the geometrical decomposition in
simplices of the geometrical product of $\sigma$ and $\tau$.

\begin{thr}\label{20016}\textbf{\emph{(Eilenberg-Zilber Theorem)}}---
A general algorithm computes:
\[
(X, Y) \mapsto [\rho_{X,Y}: C_\ast(X \times Y) \rrdc C_\ast(X) \otimes
C_\ast(Y)].
\]where:
\senumerate
\item
\(X\) and \(Y\) are simplicial sets;
\item
\(\rho_{X,Y}\) is a \emph{reduction} from the chain-complex of the product
\(C_\ast(X \times Y)\) over the tensor product of chain-complexes \(C_\ast(X)
\otimes C_\ast(Y)\).
\end{enumerate}
\end{thr}

Let us recall this theorem requires considering \emph{normalized}
chain-complexes. It is frequently presented as a consequence of the theorem of
acyclic models \cite{SPNR}, which is not very explicit; however this method can
be made effective~\cite{RBSR4}. It is simpler to use the effective formulas for
the Eilenberg-Zilber reduction $\rho_{X,Y} = (f,g,h)$ known as the
Alexander-Whitney ($f$), Eilenberg-MacLane~($g$) and Shih~($h$) operators. They
come from the recursive  definition of these operators (see \cite{ELMC1} and
\cite{ELMC2}, or \cite{SHIH}). It is in the papers~\cite{ELMC1,ELMC2} that
(homological) reductions\footnote{They were called \emph{contractions}, but it
was a serious terminological imprecision: reduction is reserved for
simplification in an \emph{algebraic} framework, and contraction in a
\emph{topological} framework. And it is essential to understand that a
chain-complex associated with a topological object in general loses the
topological nature of the object.} between chain-complexes appeared for the
first time. Only the last requirement \(h^2 = 0\) was missing.


The     Eilenberg-MacLane  and Shih  operators     have  an  essential
``exponential" nature. It is not  a question of method of computation, it is a
question of very nature:  the number of \emph{different terms} produced by the
Eilenberg-MacLane operator working on a tensor product of  bi-degree $(p,q)$ is
the  binomial coefficient~$p+q \choose p$. So that any   algorithm going
through such  a  formula is  necessarily of exponential    complexity.
Furthermore this  formula    is   unique \cite{PROT1},  and the difficulty met
here is  therefore quite essential. In a sense, ``classical'' algebraic
topology, typically the work around Steenrod operations,  consists in avoiding
the definitively exponential complexity of the Eilenberg-MacLane formula in
order to be able to reach  higher dimensions; this  text on the contrary
focuses on \emph{arbitrary} spaces in low dimensions (something like $< 12$)
where much interesting work is also to be done. A consequence of these
considerations is that our computing  methods will  certainly not lead to high
sphere homotopy groups; we are  processing the \emph{orthogonal}  problem: we
are not concerned by  high dimensional invariants  of \emph{known} objects,  we
are only interested by the first invariants of \emph{random} objects.

Interpreting the Eilenberg-Zilber Theorem in the framework of objects with
effective homology requires \emph{composition} of equivalences.

\begin{prp}\label{13837}---
A general algorithm computes:
\[
[\varepsilon: A_\ast \lrdc B_\ast \rrdc C_\ast, \varepsilon': C_\ast \lrdc
D_\ast \rrdc E_\ast] \mapsto \varepsilon'': [A_\ast \lrdc F_\ast \rrdc E_\ast]
\]
where:
\begin{enumerate}\sitemsep
 \item
 \(\varepsilon\) and \(\varepsilon'\) are two given equivalences between
 chain-complexes, the ``target'' of~\(\varepsilon\) being the ``source'' of
 \(\varepsilon'\).
 \item
 \(\varepsilon''\) is an equivalence between the extreme chain-complexes which
 must be considered as the composition \(\varepsilon'' = \varepsilon \circ
 \varepsilon'\).
 \end{enumerate}
\end{prp}

\proof Instead of a complex direct proof, a small collection of quite
elementary lemmas gives the answer.

\begin{lmm}---
The cone of an identity chain map \(\emph{\(\id{}\)}: C_\ast \leftarrow
C_\ast\) is acyclic; more precisely a simple algorithm constructs a reduction
\emph{\(\Cone(\id{}) \rrdc 0\)}.
\end{lmm}

\proof Apply Lemma~\ref{47522} to the short exact sequence:
\[
0 \leftarrow 0 \leftarrow C_\ast \stackrel{\textrm{\scriptsize id}}{\leftarrow}
C_\ast \leftarrow 0.
\]\QED

\begin{lmm}---
Let \(\rho = (f,g,h): D_\ast \rrdc C_\ast\) be a reduction. Then \(\Cone(f)\)
is acyclic; more precisely, an algorithm constructs a reduction \(\Cone(f)
\rrdc 0\).
\end{lmm}

\proof Applying the Cone Reduction Theorem~\ref{47942} to \(\Cone(f)\), using
the given reduction \(\rho\) for the source \(D_\ast\) over \(C_\ast\) and the
trivial identity reduction \(C_\ast \rrdc C_\ast\) for the target \(C_\ast\)
produces a reduction \(\Cone(f) \rrdc \Cone(\id{C_\ast})\). Composing
(Proposition~\ref{28461}) this reduction with the reduction
\(\Cone(\id{C_\ast}) \rrdc 0\) of the previous lemma gives the result. \QED

\begin{dfn}---
\emph{If \(f: B_\ast \rightarrow C_\ast\) and \(f': C_\ast \leftarrow D_\ast\)
are two chain-complex morphisms, the \emph{bicone} \(\textrm{BiCone}(f,f')\) is
constructed from \(\Cone(f)\) and \(\Cone(f')\) by identification of both
target chain-complexes \(C_\ast\).}
\end{dfn}

It is an \emph{amalgamated} sum of both cones \emph{along} the common component
\(C_\ast\).

\begin{lmm}---
Let \(\rho = (f,g,h): B_\ast \rrdc C_\ast\) and \(\rho' = (f',g',h'): D_\ast
\lrdc C_\ast\) be two reductions. An algorithm constructs a reduction
\(\textrm{\emph{BiCone}}(f,f')^{[-1]} \rrdc B_\ast\) and another one
\(\textrm{\emph{BiCone}}(f,f')^{[-1]} \rrdc D_\ast\).
\end{lmm}

\proof The bicone \(\textrm{BiCone}(f,f')\) can be interpreted as \(\Cone(f:
B_\ast \rightarrow \Cone(f'))\), calling again \(f\) the chain-complex morphism
with the same source as \(f\) and going to \(C_\ast\) which is also a
sub-\emph{chain-complex} of \(\Cone(f')\). This allows us to apply again the
Cone Reduction Theorem to the trivial identity reduction over \(B_\ast\) and
the reduction to 0 of \(\Cone(f')\). The desuspension process for the bicone is
necessary, because in a cone, the source is suspended.\QED

\noindent\textsc{Proof of Proposition~\ref{13837}.} It is a consequence of the
next diagram and composition of reductions:
\[
A_\ast \lrdc B_\ast \lrdc \textrm{BiCone}(f,f')^{[-1]} \rrdc D_\ast \rrdc
E_\ast.
\]\QED

\begin{crl}---
A general algorithm computes:
\[
(X_{EH}, Y_{EH}) \mapsto (X \times Y)_{EH}
\]
where:
 \senumerate
 \item
 \(X_{EH}\) and \(Y_{EH}\) are simplicial sets \emph{with effective homology};
 \item
 \((X \times Y)_{EH}\) is a version \emph{with effective homology} of the
 product \(X \times Y\).
\end{enumerate}
\end{crl}

\proof Let  $(X,C_\ast(X),EC_X,\varepsilon_X)$   and    $(Y,C_\ast(Y),EC_Y,
\varepsilon_Y)$ be two simplicial sets with effective homology. Eilenberg and
Zilber  give  an  equivalence $\varepsilon_1:   \mbox{\(C_\ast(X \times Y)\)}
\stackrel{=}{\lrdc} C_\ast(X \times Y) \Rightarrow   C_\ast(X)  \otimes
C_\ast(Y)$ (the left reduction  is trivial); Proposition~\ref{49390} gives also
an equivalence $\varepsilon_2$ between $C_\ast(X) \otimes C_\ast(Y)$ and $EC_X
\otimes EC_Y$. Composing these homotopy equivalences (Proposition \ref{13837}),
we obtain the wished homotopy equivalence between $C_\ast(X \times  Y)$  and
the \emph{effective} chain-complex $EC_X \otimes EC_Y$ \QED

The  K\"unneth   Theorem is  not used;  it allows you  to \emph{guess} the
homology groups  of $EC_X \otimes  EC_Y$ if you  know the homology   groups of
factors,  but we  are not  concerned  by this question: the chain-complexes
$EC_X$ and $EC_Y$ are effective, so that $EC_X \otimes EC_Y$ is also effective,
and this  is sufficient. We are on the contrary  essentially interested by an
\emph{explicit} homology equivalence between $C_\ast(X \times Y)$ and $EC_X
\otimes EC_Y$, and the explicit definition of  the Eilenberg-Zilber reduction
is  the key point.

Let us finish this presentation of the Eilenberg-Zilber Theorem by a typical
application. It is elementary to compute the homology of the real projective
plane \(P^2\bR\); this was done by means of Kenzo at page~\pageref{49275}, but
once the simplicial set technique is known, pen and paper are enough. The
minimal\footnote{Minimal by the number of simplices, but the Kan condition, see
Section~\ref{40872}, is not satisfied, so that this minimal description of
\(P^2\bR\) is not minimal in the sense of Kan~\cite[\S 9]{MAY}.} simplicial
descrition has three non-degenerate simplices: one vertex, the base point, one
edge, the equivalence class of the equator and one triangle. The normalized
chain-complex is:
\[
C_\ast(P^2\bR) = [\cdots 0 \leftarrow \bZ \stackrel{0}{\leftarrow} \bZ
\stackrel{\times 2}{\leftarrow} \bZ \leftarrow 0 \cdots]
\]
with \(\times 2\) between degrees 2 and 1.

If you work in the style of traditional algebraic topology, you deduce the
homology groups \(H_\ast(P^2\bR) = (\bZ, \bZ_2)\) in degrees 0 and 1, the
others being null.

Now your client orders a \emph{construction} \(X := P^2\bR \times P^2\bR\) and
asks for \(H_\ast(X) = ??\). In traditional style, you will try to deduce the
homology groups of \(X\) from those of \(P^2\bR\); the answer is the K\"unneth
formula~\cite[Section 5.3]{SPNR}:
\[
 H_n(X \times Y) = \left(\oplus_{p=0}^n H_p(X) \otimes H_{n-p}(Y)\right) \oplus
 \left(\oplus_{p=0}^{n-1} \Tor^\bZ_1(H_p(X), H_{n-1-p}(Y))\right).
\]
The bad student forgets the torsion terms, which require some lucidity.
Furthermore the sum decomposition is not canonical. The result is \(H_\ast(X) =
(\bZ, \bZ_2 \oplus \bZ_2, \bZ_2, \bZ_2)\), the last \(\bZ_2\) being the only
contribution of torsion terms.

In the spirit of effective homology, you observe the normalized chain-complex
\(C_\ast(P^2\bR)\) is already effective, so that a trivial effective homology
is enough. If ever you are interested by the standard homology groups, you can
ask a machine, but here it is obvious, you obtain the right homology groups of
\(P^2\bR\). Now what about \(H_\ast(X)\)? First you \emph{simplicially}
construct \(X = P^2\bR \times P^2\bR\); it is not so simple, but a machine does
it automatically; the numbers of non-degenerate simplices are \((1, 3, 9, 12,
6)\); the corresponding normalized chain-complex \(C_\ast(X)\) is relatively
complex, but Eilenberg and Zilber explain to us there is a reduction
\(C_\ast(X) \rrdc C_\ast(P^2\bR) \otimes C_\ast(P^2\bR)\) and the last
chain-complex is elementarily computed; presented as a bicomplex, it is:
\[
\xymatrix{
 \bZ \ar[d]_{\times 2} & \bZ \ar[l]_0 \ar[d]|{\times (-2)} & \bZ \ar[l]_{\times
 2} \ar[d]_{\times 2}
 \\
 \bZ \ar[d]_0 & \bZ \ar[l]_0 \ar[d]_0 & \bZ \ar[l]_{\times 2} \ar[d]_0
 \\
 \bZ & \bZ \ar[l]_0 & \bZ \ar[l]_{\times 2}
}
\]
giving the expected homology groups. This is nothing but the standard
calculation giving \(\Tor^\bZ_1(\bZ_2, \bZ_2) = \bZ_2\), so that you can wonder
why such a presentation? The crucial point is the following: your client will
probably tomorrow undertake a new construction, more or less complicated, using
the \emph{geometry} of \(X = (P^2\bR)^2\); if the construction is a little
complicated, then it is not possible to \emph{describe} it at the homological
level and you cannot continue. We will soon see critical examples. In effective
homology, the simplicial model of \(X\) remains present in your environment
with the Eilenberg-Zilber connection with the chain-complex displayed above.
Whatever construction is imagined by your client, \emph{you will be ready to
construct an algorithm} providing the effective homology of the resulting
object.

\subsection{The twisted Eilenberg-Zilber Theorem.}

Let \(F \hookrightarrow [E = F \times_\tau B] \rightarrow B\) be a twisted
product, that is, a simplicial fibration defined by a \emph{base space} \(B\),
some simplicial set, a \emph{fibre space} \(F\), another simplicial set, and
some twisting operator \(\tau: B \rightarrow G\), the target \(G\) being some
simplicial group acting over the fibre space \(F\). See Section~\ref{92270} for
details and examples. If \(\tau\) is trivial, the Eilenberg-Zilber Theorem
gives a reduction \(C_\ast(E) \rrdc\ \mbox{\(C_\ast(F) \otimes C_\ast(B)\)}\).
The so-called ``twisted'' Eilenberg-Zilber Theorem constructs an analogous
reduction \(C_\ast(E) \rrdc C_\ast(F) \otimes_t C_\ast(E)\), the index of
\(\otimes_t\) meaning the \emph{differential} of the usual chain-complex tensor
product \(C_\ast(F) \otimes C_\ast(B)\) being (deeply) modified.

\begin{thr}\label{80728} \textbf{\emph{(Twisted Eilenberg-Zilber Theorem)}}---
An algorithm computes:
\[
 \Phi \mapsto \rho
\]
where:
 \senumerate
\item
\(\Phi\) is a simplicial fibration \(\Phi = \{F \hookrightarrow [E = F
\times_\tau B] \rightarrow B\}\).
\item
\(\rho: C_\ast(E) \rrdc C_\ast(F) \otimes_t C_\ast(B)\) is a reduction of the
(normalized) chain-complex of the total space of the fibration over a
chain-complex \(C_\ast(F) \otimes_t C_\ast(B)\); the underlying \emph{graded
module} of the latter is the same as for \(C_\ast(F) \otimes C_\ast(B)\), but
the differential is modified to take account of the twisting operator \(\tau\).
\end{enumerate}
\end{thr}

\proof The ordinary  (non-twisted) Eilenber-Zilber Theorem gives  a reduction
between  the non-twisted cartesian   and tensor products,  the twisting
operator being null.  But we must take account  of the twisting operator
$\tau$; this twisting operator does not change the  underlying top  graded
module, only the differential  is modified: the  0-face  operator  is twisted,
see Definition~\ref{69699}. The basic perturbation lemma may be applied if the
nilpotency condition is satisfied.

If $(f,b)$ is a simplex of $E$, the component $b$ has a unique form $b =  \eta
b'$    where $b'$  is   non-degenerate   and   $\eta$   is a multi-degeneracy
operator; if $b$  is  non-degenerate then $b'=b$ and $\eta$  is the identity,
no   degeneracy at all. Following Serre,  the filtration degree   of $(f,b)$ is
the dimension   of $b'$,  the ``base dimension''. The Shih   homotopy operator
of  Eilenberg-Zilber    is \emph{natural}, and when it works on $(f,b)$ it is
equal  to the one which is defined on $F  \times b'$, just above  the simplex
$b'$; therefore the Shih operator does not increase the filtration degree.

On the contrary the  perturbation $\hdl(f,b) = (\partial_0  f \ldotp \tau(b),
\partial_0 b)  -  (\partial_0    f,  \partial_0 b)$   has   a filtration degree
smaller than the filtration   degree of $(f,b)$. If $b$ is non-degenerate, it
is obvious. If  $b$ is degenerate and if the $\eta$ in  the expression $b =
\eta  b'$ does not  contain a $\eta_0$, then  $\partial_0  \eta b'  = \eta'
\partial_0 b'$,  because  of the commuting  relation $\partial_0 \eta_i =
\eta_{i-1} \partial_0$ if $i
>0$; the filtration degree of $(f', \partial_0 b)$  is again less than
the one  of $(f,b)$. Finally,  if the multi-degeneracy operator $\eta$ contains
a $\eta_0$, then $\tau(b)$ is trivial, see Definition~\ref{51596}, and the
perturbation is null. The nilpotency hypothesis is satisfied.

The basic perturbation lemma is then applied and produces the wished reduction.
\QED

The following technical proposition is the key point allowing one to use the
twisted Eilenberg-Zilber Theorem to obtain a version with effective homology of
the Serre spectral sequence.

\begin{prp}\label{89844}---
Let $\Phi=(B,F,G,\tau,E) = \{F \hookrightarrow [E = (F \times_\tau B)]
\rightarrow B\}$ be a simplicial fibration. Let $\rho: C_\ast(F \times B)
\Rightarrow C_\ast(F) \otimes C_\ast(B)$ (resp. $\rho': C_\ast(F \times_\tau B)
\Rightarrow C_\ast(F) \otimes_t C_\ast(B)$) be the non-twisted (resp. twisted)
reduction given by the Eilenberg-Zilber (resp. twisted Eilenberg-Zilber)
Theorem. Let $d$ (resp. $d'$) be the differential of $C_\ast(F) \otimes
C_\ast(B)$ (resp. $C_\ast(F) \otimes_t C_\ast(B)$) and let $\delta = d'-d$ be
the bottom differential perturbation computed by the twisted Eilenberg-Zilber
Theorem. Then, \emph{if $B$ is $1$-reduced}, the bottom perturbation $\delta$
decreases the filtration degree at least by $2$.
\end{prp}

A simplicial set \(B\) is \emph{1-reduced} if it has only one vertex and no
non-degenerate 1-simplex, therefore, only one 1-simplex, the unique degeneracy
of the unique vertex.

The conclusion of the proposition is to be understood as follows: if $b$ (resp.
$f$) is a $p$-simplex (resp. $q$-simplex) of $B$ (resp. $F$), then:

$$ \delta (f \otimes b) = \sum_{r=2}^p \delta_r(f \otimes b) $$

\noindent where $\delta_r(f \otimes b) \in C_{q+r-1}(F) \otimes C_{p-r}(B)$.
Note it is not possible to coherently choose one of both possible notations $(f
\otimes b)$ and $(f \otimes_t b)$: in fact $\delta = d'-d$ and $d$ (resp.
\(d'\)) is to be applied to $(f \otimes b)$ (resp. $(f \otimes_t b)$).

\proof Let $\rho=(AW,EML,SH)$ the ordinary Eilenberg-Zilber reduction between
$C_\ast(F \times B)$ and $C_\ast(F) \otimes C_\ast(B)$. If $\hdl = \hd\,'-\hd$
is the top perturbation, the explicit formula for the bottom perturbation in
the proof of Theorem, cf. page~\pageref{35472}, gives:

$$ \delta(f \otimes b) = (AW \circ (\sum_{i=0}^\infty (-1)^i (\hdl \circ SH)^i)
                      \circ \hdl \circ EML) (f \otimes b).
$$

We have observed in the previous proof the top perturbation $\hdl$ decreases
the filtration degree at least by 1; furthermore, the Shih operator does not
increase this filtration degree; therefore, the components with $i \geq 1$ in
the expression just above satisfy the wished condition. The main work concerns
only the $i=0$ component.

The Eilenberg-MacLane operator working on $f \otimes b$ ($f$ a non-degenerate
$q$-simplex of $F$, $b$ a non-degenerate $p$-simplex of B) produces a set of
terms, shuffles of the form $\pm(\eta f, \eta' b)$ for some multi-degeneracy
operators $\eta$ and $\eta'$. If $\eta'$ contains a $\eta_0$, then the
corresponding twist is trivial and there is no perturbation. We can organize
the other terms as follows: $\pm(\eta f, \eta' \eta''b)$ where $\eta$ contains
a $\eta_0$, $\eta''$ is a composition of consecutive degeneracies $\eta'' =
\eta_k \eta_{k-1} \ldots \eta_2 \eta_1 = \eta_1^k$, and $\eta'$ is another
composition $\eta' = \eta_{i_\ell} \ldots \eta_{i_1}$ with $i_1 \geq k+2$ and
$k+\ell=q$; the integer $k+1$ is the first missing index in the degeneracies of
the second component. We have then the expression:
\[
(\hdl \circ EML) (f \otimes b) = \sum \pm [(\partial_0 \eta f \ldotp
\tau(\eta'\eta''b), \eta'_{-1}\eta''_{-1}\partial_0b) -(\partial_0 \eta f,
\eta'_{-1}\eta''_{-1}\partial_0b)].
\]
In the expression above, a term $\eta'_{-1}$ denotes the multi-degeneracy
operator $\eta'$ where all the indices have been replaced by the same minus
one; in particular $\eta''_{-1} = \eta_{k-1} \ldots \eta_0$. There remains to
apply the Alexander-Whitney operator:

$$ AW(f',b') = \sum_{j=0}^{p+q-1} \partial_{j+1}^{p+q-1-j}f' \otimes
                               \partial_0^j b'.
$$

If $j > k$, then there are at least two operators $\partial_0$ which remain
alive in the right component; this comes from the relation $\partial_0^j
\eta_{k-1} \ldots \eta_0 = \partial_0^{j-k}$. In such a case, the term becomes
something like $\pm(\ldots,\eta'''\partial_0^m b)$ with $m \geq 2$, and the
result is obtained.

If $j \leq k$, the twisting modifier $\tau(\eta'\eta''b)$ becomes by
Alexander-Whitney $\tau(\partial_{j+2}^{p+q-1-j}\eta'\eta''b)$, because the
face index is increased by one when entered inside the $\tau$ argument. On one
hand the inequality $p+q-1-j \geq p+q-1-k = p-1+\ell$ is satisfied; on the
other hand all the indices $i_\ell,\ldots,i_1$ are $> k+1 \geq j+1$, so that
the following relation is satisfied:

$$
\partial_{j+2}^{p+q-1-j} \eta' \eta'' = \partial_{j+2}^{p-1+k-j} \eta''.
$$

\noindent But we have also the relation:

$$
\partial_{j+2}^{p-1+k-j} \eta_k \ldots \eta_1 b=
\eta_j \ldots \eta_1 \partial_2^{p-1} b; $$

\noindent finally, the $p$-simplex $b$ gives a 1-simplex $\partial_2^{p-1} b$,
dimension 1, necessarily the \(\eta_0\)-degeneracy of the base-point, for the
base space $B$ is 1-reduced; the corresponding twist is trivial and the
associated bottom perturbation is null. \QED

The previous demonstration is a little technical but more elementary than the
original ones \cite{BRWNE2,SHIH} (see also~\cite{GGNH1}), where the interesting
notion of \emph{twisting cochain} is required and used to make more conceptual
the result. The present demonstration is sufficient to give a
\emph{certificate} for the corresponding computer program.

\subsection{The \emph{effective} version of the Serre spectral sequence.}

Let $F \hookrightarrow [E = F \times_\tau B] \rightarrow B$ be a simplicial
fibration. The Serre spectral sequence gives a set of relations between the
homology groups of $F$, $E$ and $B$. In some particular cases, this spectral
sequence gives a method allowing you to deduce the homology  groups of one of
the components, $E$  for example, when  the homology  groups  of the others
($B$ and $F$) are given. An example of this sort has been given in
Section~\ref{11076} where \(H_\ast(B)\) was deduced from \(H_\ast(F)\) and
\(H_\ast(E)\), known.

But in the general case, the Serre spectral sequence \emph{is not} an
algorithm; see for example \cite[pp 6 and 28]{MCCL} for a serious warning about
this question, unfortunately not formalized: a \emph{computational} environment
is required there to obtain a \emph{mathematical} statement of the obstacle.
Section~\ref{62859} was devoted to the first historical example where the
spectral sequence method \emph{failed} to compute a sphere homotopy group.

We show here the \emph{effective} homology methods give very easily a
\emph{constructive} version of the Serre spectral sequence. For example the
Kenzo program ``stupidly'' computes in one minute \(\pi_6 S^3 = \bZ_{12}\).

\begin{thr}---
An algorithm computes:
\[
(F_{EH}, B_{EH}, \tau) \mapsto E_{EH}
\]
where:
 \senumerate
 \item
 \(F_{EH} = (F, C_\ast(F), EC_\ast^F, \varepsilon_F)\) is a version \emph{with
 effective homology} of the fibre space \(F\);
 \item
 \(B_{EH} = (B, C_\ast(B), EC_\ast^B, \varepsilon_B)\) is a version \emph{with
 effective homology} of the base space \(B\); we assume the base space \(B\) is
 \emph{1-reduced}: only one vertex, no non-degenerate 1-simplex;
 \item
 \(\tau: B \rightarrow G\) is a twisting operator with values in a simplicial
 group \(G\) acting over the fibre space \(F\), defining the \emph{twisted
 product} \(E = F \times_\tau B\);
 \item
 \(E_{EH} = (E, C_\ast(E), EC_\ast^E, \varepsilon_E)\) is a version \emph{with
 effective homology} of the total space \(E\).
 \end{enumerate}
\end{thr}

In  other  words, you can  \emph{compute}  the  homology groups of the total
space $E$, no  mysterious unreachable differential, no extension problem at
abutment; see \cite[pp 6  and 28]{MCCL}. More important, if the total space $E$
is one  of  the elements of a  new ``reasonable'' construction,  the  object
$E_{EH}$ can   again  be used to  obtain  a version with effective homology of
the  new constructed object, and so on.

\proof We must construct the equivalence \(\varepsilon_E: C_\ast(E) \eqvl
EC_\ast^E\). It is obtained as the composition of two equivalences,
\(\varepsilon_E := \varepsilon' \circ \varepsilon''\), see
Proposition~\ref{13837} for the construction of such a composition.

The first equivalence \(\varepsilon'\) is produced by the twisted
Eilenberg-Zilber Theorem:
\[
\varepsilon' = \{C_\ast(F \times_\tau B) \stackrel{=}{\lrdc} C_\ast(F
\times_\tau B)
  \rrdc C_\ast(F) \otimes_t C_\ast(B)\}
\]
where the left reduction is trivial. When $\varepsilon'$ and in particular the
\emph{twisted} tensor product $C_\ast(F) \otimes_t C_\ast(B)$ are constructed,
then we can construct the second necessary homotopy
equivalence~$\varepsilon''$, by applying the basic perturbation lemma to the
difference between \mbox{$C_\ast(F) \otimes_t C_\ast(B)$} and $C_\ast(F)
\otimes C_\ast(B)$. Two equivalences are available:
\[
\begin{array}{rcl}
\varepsilon_F &=& \{C_\ast(F) \lrdc \hC^F_\ast \rrdc EC^F_\ast\} \\
\varepsilon_B &=& \{C_\ast(B) \lrdc \hC^B_\ast \rrdc EC^B_\ast\}
\end{array}
\]
\noindent and we can construct their (non-twisted) tensor product (Proposition
\ref{49390}):
\[
 \varepsilon_{FB} = \{C_\ast(F) \otimes C_\ast(B) \lrdc \hC^F_\ast \otimes
 \hC^B_\ast \rrdc EC^F_\ast \otimes EC^B_\ast.
\}
\]

A \emph{filtration degree} is defined on the three tensor products according to
the degree with respect the second factor $C_\ast(B)$, $\hC_B$ or $EC_B$. Let
us introduce on the bottom left-hand chain-complex of this homotopy equivalence
the necessary perturbation to obtain the twisted tensor product $C_\ast(F)
\otimes_t C_\ast(B)$; the base space~$B$ is 1-reduced and according to
Proposition \ref{89844}, this perturbation decreases the filtration degree at
least by 2.

The left reduction of $\varepsilon_{FB}$ describes the left hand chain-complex
\(C_\ast(F) \otimes C_\ast(B)\) as a subcomplex of the top chain-complex $\hC_F
\otimes \hC_B$, and we can apply the easy perturbation lemma to the left
reduction; the perturbation can be so transferred to the top chain-complex
\(\hC^F_\ast \otimes \hC^B_\ast\), obtaining the same graded module with
another differential \(\hC^F_\ast \otimes_{t'} \hC^B_\ast\) with the same
property (Proposition~\ref{89844}) about the filtration degree for the
difference between the new and the old differential: this perturbation is
nothing but a \emph{copy} of the starting perturbation on a subcomplex of
\(\hC^F_\ast \otimes \hC^B_\ast\). The perturbation over \(\hC^F_\ast \otimes
\hC^B_\ast\) decreases the filtration degree at least by 2; the homotopical
component of the right reduction of $\varepsilon''$ increases the filtration
degree at most by one; the nilpotency hypothesis is satisfied. The basic
perturbation lemma can therefore be applied to the right reduction and the
perturbation obtained for the top chain-complex and the equivalence
\(\varepsilon''\) is obtained.

Both components \(EC^F_\ast\) and \(EC^B_\ast\) are \emph{effective}
chain-complexes; their tensor product, whatever is the differential, is
effective too. We have obtained a version \emph{with effective homology} of the
total space \(E\). \QED.

\section{The Eilenberg-Moore spectral sequence.}

\subsection{Introduction.}

Let \(F \hookrightarrow [E = F \times_\tau B] \rightarrow B\) be a fibration;
the \emph{total space} \(E\) is a twisted product of the \emph{base space}
\(B\) by the \emph{fibre space} \(F\), the \emph{twist} being defined by an
appropriate \emph{twisting operator} \(\tau: B \rightarrow G\), see
Definition~\ref{51596} which in particular explains the role of the
\emph{structural group} \(G\). As usual in constructive topology, we are
working inside the \emph{simplicial} framework.

The Serre spectral sequence, or more exactly the \emph{effective homology}
version of the Serre spectral sequence, see the previous section, allows us to
compute the effective homology of the total space when the effective homologies
of the base space and the fibre space are given; it is essentially a
\emph{product} operator. This is valid only if the base space is simply
connected, more exactly in our simplicial framework, if the base space is
1-reduced.

The Eilenberg-Moore spectral sequence corresponds to a \emph{division}. Because
the notion of twisted product is not symmetric with respect to both factors, in
fact \emph{two} Eilenberg-Moore spectral sequences are to be defined, but they
are similar. We will explain the \emph{Cotor} spectral sequence, expressing the
homology of the fibre space \(F\) as a ``Cotor'' operation between the
homologies of the base space and the total space. The symmetric Tor spectral
sequence describes the homology of the base space as a ``Tor'' involving the
homologies of the total space, the fibre space \emph{and the structural group}.
We give here a reasonable level of details for the Cotor spectral sequence and
will briefly explain how the symmetric result for the Tor spectral sequence is
obtained.

\subsection{Coalgebra and comodule structures.}

The notions of algebra and module are common. The Cotor spectral sequence needs
the symmetric notions of coalgebra and comodule.

\begin{dfn}---
\emph{A \emph{differential coalgebra} is a chain-complex \(C_\ast\) provided
with a coproduct \(\Delta: C_\ast \rightarrow C_\ast \otimes C_\ast\) and a
counit \(\eta: C_0 \rightarrow \fR\), satisfying the rules that are required
for differential algebras, with the ``arrows reversed''.}
\end{dfn}

See for example~\cite[VI.9]{MCLN2}. In particular the tensor product \(C_\ast
\otimes C_\ast\) is itself a chain-complex (Definition~\ref{05496}) and the
coproduct \(\Delta\) must be \emph{compatible with the differentials} of
\(C_\ast\) and \(C_\ast \otimes C_\ast\); in other words the coproduct is a
\emph{chain-complex morphism}. The coproduct is \emph{homogeneous}: the (total)
degree of a coproduct is equal to the degree of the argument: \(|\Delta(x)| =
|x|\). The coproduct is coassociative: \((\id{} \otimes \Delta) \circ \Delta =
(\Delta \otimes \id{}) \circ \Delta\). The counit satisfies \((\eta \otimes
\id{}) \circ \Delta = (\id{} \otimes \eta) \circ \Delta = \id{}\); in these
equalities, you have to identify \(C_\ast \otimes \fR = \fR \otimes C_\ast =
C_\ast\). All the tensor products, unless otherwise stated, are
\(\otimes_\fR\).

\begin{dfn}---
\emph{Let \(X\) be a simplicial set. The \emph{canonical coalgebra structure}
of \(C_\ast(X)\) is defined by the coproduct \(\Delta(\sigma) = \sum_{i=0}^n
\partial_{i+1} \cdots \partial_n \sigma \otimes \partial_0^i \sigma\) and the
counit~\(\eta\) is defined by \(\eta(\sigma) = 1_\fR\) if \(\sigma \in X_0\),
\(\eta(\sigma) = 0\) otherwise.
 }
\end{dfn}

The formulas that are given for individual simplices must as usual be linearly
extended to combinations of simplices. The coproduct is easily understood in
the simplicial complex case:
\[
\Delta(0123) = 0 \otimes 0123 + 01 \otimes 123 + 012 \otimes 23 + 0123 \otimes
3
\]
where for example \(123\) denotes the simplex spanned by the vertices 1, 2 and
3. And there is a unique way to extend this game rule to the general case of
simplicial sets. This coproduct is known as the \emph{Alexander-Whitney}
coproduct. The counit consists in deciding the image of a 0-simplex is the unit
of the ground ring, and is null in higher dimensions. Usual algebraic topology
uses this \emph{coproduct} to install by duality a \emph{product} on the
cohomology, providing to this cohomology an \emph{algebra structure}.

\begin{dfn}---
\emph{If \(C_\ast\) is a differential coalgebra, a \emph{differential} (right)
\emph{comodule} is a chain-complex \(M_\ast\) provided with an external
coproduct \(\Delta_M: M \rightarrow M \otimes C_\ast\), satisfying the rules
that are required for differential algebras, with the ``arrows reversed''.
 }
\end{dfn}

The external coassociativity rule becomes \((\Delta_M \otimes \id{C}) \circ
\Delta_M = (\id{M} \otimes \Delta_C) \circ \Delta_M\). The external counit rule
is \((\id{M} \otimes \eta) \circ \Delta_M = \id{M}\) where an identification
\(M \otimes \fR = M\) is necessary.

For example, if \(f: X \rightarrow Y\) is a simplicial morphism between two
simplicial sets, there is a canonical way to provide \(C_\ast(X)\) with a
\(C_\ast(Y)\)-comodule structure. Decide \(\Delta_{X,Y} = (\id{C_\ast{X}}
\otimes f) \circ \Delta_X\): it is a process of coextension of scalars, the
coproduct \(\Delta_X\) of the \emph{coalgebra} \(C_\ast(X)\) is ``extended'' to
a comodule coproduct~\(\Delta_{X,Y}\).

\subsection{The Cobar construction.}

If you intend to make \emph{divisions}, a good idea could consist in firstly
studying \emph{inverses}: a division is most often nothing but a multiplication
by an inverse. In topology, when you want to consider in a fibration \(F
\hookrightarrow E \rightarrow B\) the fibre space \(F\) as a (twisted) quotient
of \(E\) by \(B\), it is natural to look for some ``inverse'' of \(B\).

\begin{dfn}---
\emph{Let \(B\) be a \emph{pointed} topological space. The \emph{path space}
\(P B\) of \(B\) is the space of all the continuous maps \(PB =
\mathcal{C}((I,0), (B, \ast))\). The loop space~\(\Omega B\) of~\(B\) is the
space of all the continuous maps \(\Omega B = \mathcal{C}((I,0,1), (B, \ast,
\ast))\). A canonical fibration \(\Omega B \hookrightarrow P B \rightarrow B\)
is defined.
 }
\end{dfn}

The space \(B\) is \emph{pointed}, that is, \(B\) is a shorthand for \(B = (B,
\ast)\) where \(\ast\) is some distinguished point of \(B\), its \emph{base
point}. In the case of a path \(\gamma \in P B\), the image of \(0 \in I =
[0,1]\) must be the base point \(\ast\). The same for 0 and 1 in the case of a
\emph{loop}. A loop is a path, but in general a path is not a loop. The
canonical projection \(\pr: PB \rightarrow B\) is defined by \(\pr(\gamma) =
\gamma(1)\). The fibre \emph{above the base point} is the loop space \(\Omega
B\); it is a Hurewicz fibration: the canonical projection \(\pr: PB \rightarrow
B\) satisfies the \emph{homotopy lifting} property, see~\cite[2.2]{SPNR}. The
total space \(PB\) is contractible: every path can be retracted along itself to
the trivial path \(\gamma_0 \equiv \ast\). So that the total space has the
homotopy type of a point, the unit in the world of topology. And the loop space
\(\Omega B\) is therefore a sort of \emph{inverse} \(\Omega B =\)
``\(B^{-1}\)''.

These constructions were invented by Jean-Pierre Serre when he designed
appropriate tools allowing him to ``compute'' homotopy groups of spheres. But
we cannot work with \emph{general} topology on a computer, and the analogous
process in \emph{combinatorial} topology was discovered (or invented?) by
Daniel Kan and was sketched in Section~\ref{87671}. We summarize the
corresponding result in this theorem.

\begin{thr}\label{13368}---
A functor \(\Omega\) can be defined on the category of reduced simplicial sets.
If \(B\) is such a simplicial set, then a canonical twisting operator \(\tau\)
is defined \(\tau: B \rightarrow \Omega B\) defining a twisting product \(PB =
\Omega B \times_\tau B\) which is contractible.
\end{thr}

It is the simplicial version of the Hurewicz fibration \(\Omega B
\hookrightarrow PB \rightarrow B\). The chapter VI of~\cite{MAY} gives all the
possible details about this question. In this way we have a simple process to
construct the ``inverse'' of a \emph{base} space.

The next step must go from combinatorial topology to algebra. It happens in a
sense a differential coalgebra is the translation in algebra of a topological
space\footnote{A deeper study shows this is not exact; some essential
information in general is \emph{lost} in this translation process. If you want
to keep the \emph{whole} homotopy type of the space \(B\), you must endow the
chain complex \(C_\ast(B)\) not only with the coalgebra structure, but with
some \(E_\infty\)-coalgebra structure, \(E_\infty\) being an appropriate
algebraic \emph{operad}, which can be understood as the completion of the
Steenrod operations. See~\cite{BRFR,MNDL2} for this essential point.}.

\begin{uos}\label{62309}---
Our differential coalgebras \(C_\ast\) are assumed from now on
\mbox{\emph{1-reduced}}. This means the 0-component \(C_0\) is isomorphic to
the ground ring \(\fR\) by the coaugmentation \(\eta\) of the coalgebra, and
the 1-component \(C_1\) is null.
\end{uos}

Many definitions and results given here can be extended to significantly more
general situations, but our main result is concerned by the 1-reduced case, and
stating now this restriction makes easier the exposition.

\begin{dfn}---
\emph{Let \(C\) be a differential coalgebra, and \(M\) (resp. \(N\)) be a right
(resp. left) \(C\)-comodule. The \emph{Cobar construction} \(\Cobar^C(M, N)\)
is a \emph{bicomplex} defined as follows:
\[
\Cobar^C(M,N) = \oplus_{p=0}^\infty (M \otimes \oC^{\otimes p} \otimes
N)^{[-p]}
\]
where \(\oC\) is the \emph{coaugmentation ideal} of \(C\); the differential
structure of \(\Cobar^C(M,N)\) comes from two differentials, the vertical
differential \(d_v\) is deduced from the component differentials, and the
horizontal differential \(d_h\) is deduced from the various coproducts.
 }
\end{dfn}

It is a first quadrant bicomplex, the \emph{horizontal} degree is \(p\) and the
\emph{vertical} degree is deduced from the grading of \((M \otimes \oC^p
\otimes N)\), each factor being graded. The total degree is the
\emph{difference} between vertical and horizontal degree, this is necessary
because of the nature of the coproduct which \emph{unfolds} an element into a
sum of tensor products. This point is recalled by the exponent \([-p]\) in the
initial formula; you can consider this difference as a \emph{desuspension}
process. A purist usually prefers install this bicomplex in the second
quadrant, but the notation becomes a little heavier.

The reader notes we allow us not to indicate the grading property by the usual
\mbox{\(\ast\)-index}, in order to trim the notation when it is possible. Let
us detail a little more this definition of the Cobar construction. The
coaugmentation \(\eta: C \rightarrow \fR\) has a kernel \(\oC\); because of the
restriction~\ref{62309}, the coaugmentation ideal \(\oC\) is nothing but \(C\)
with the 0-component cancelled. The grading of \(\oC\) therefore begins in
degree 2. The differential of \(\oC\) in degree 2 is null as in \(C\) itself,
because of the absence of a 1-component.

The components \(M\), \(\oC\) and \(N\) in the formula defining the Cobar are
chain-complexes, so that their tensor products are chain-complexes too, see
Definition~\ref{05496}; so is defined the \emph{vertical} differential of the
chain complex, signs being deduced from the Koszul rule, except the role of
\((-1)^n\) explained later:
\[
\begin{array}{rcl}
 d_v(a \otimes c_1 \otimes \cdots \otimes c_n \otimes b) &=&
 \phantom{+}\ (-1)^n da \otimes c_1 \otimes \cdots \otimes c_n \otimes b
 \\
 && +\ (-1)^{n + |a|} a \otimes dc_1 \otimes \cdots \otimes c_n \otimes b
 \\
 && +\ \cdots\ \cdots
 \\
 && +\ (-1)^{n + |ac_1\cdots c_{n-1}|} a \otimes c_1 \cdots \otimes dc_n \otimes b
 \\
 && +\ (-1)^{n+|ac_1\cdots c_{n}|} a \otimes c_1 \cdots \otimes c_n \otimes db
\end{array}
\]

The coalgebra \(C\) and the comodules \(M\) and \(N\) are provided with
coproducts. The ideal \(\oC\) inherits a pseudo-coproduct again denoted by
\(\Delta: \oC \rightarrow \oC \otimes \oC\) by cancelling in the original
coproduct the factors of degree 0 in the result. For example in the case of the
standard \(s\)-simplex \(\Delta^s\) with the 1-skeletton collapsed on the base
point to satisfy the 1-reduced requirement, we would have: \( \Delta(0123) =
0\) because the 0- and 1-simplices do not exist anymore in
\(\oC_\ast^N(\Delta^n)\); on the contrary \(\Delta(01234) = 012 \otimes 234\).
The same process for~\(M\) and \(N\) gives pseudo-coproducts \(\Delta: M
\rightarrow M \otimes \oC\) and \(\Delta: N \rightarrow \oC \otimes N\). Then
the horizontal differential is defined, if \(m \in M\), \(c_i \in \oC\) and \(b
\in N\), by the formula:
\[
\begin{array}{rcl}
 d_h(a \otimes c_1 \otimes \cdots \otimes c_n \otimes b) &=&
 \phantom{+\ } \Delta(a) \otimes c_1 \otimes \cdots \otimes c_n \otimes b
 \\
 && -\ a \otimes \Delta(c_1) \otimes \cdots \otimes c_n \otimes b
 \\
 && \pm\ \cdots\ \ \cdots
 \\
 && +\ (-1)^{n} a \otimes c_1 \otimes \cdots \otimes \Delta(c_n) \otimes b
 \\
 && +\ (-1)^{n+1} a \otimes c_1 \otimes \cdots \otimes c_n \otimes \Delta(b).
\end{array}
\]

The Cobar carries in particular a \emph{cosimplicial} structure. The
\(p\)-simplices are the elements of \(M \otimes \oC^{\otimes p} \otimes N\) and
the coface operator \(\partial_i: M \otimes \oC^{\otimes p} \otimes N
\rightarrow M \otimes \oC^{\otimes (p+1)} \otimes N\) is defined by applying
the coproduct to the \(i\)-th copy of \(\oC\), or to \(M\) (resp. \(N\)) if \(i
= 0\) (resp. \(i = p + 1\)). As for the simplicial structures, a cosimplicial
structure is defined by a \emph{covariant} functor to the considered category,
here the \(\fR\)-modules. The compatibility of the coproduct with the
differentials guarantees every \(\partial_i\) is also compatible with the
differentials. As in the simplicial case, the alternate sum of the coface
operators is a differential\footnote{More precisely, the cosimplicial structure
should be installed on \(\oplus_{p=0}^\infty (M \otimes C^{\otimes p} \otimes
N)^{[-p]}\); our horizontal differential is the differential obtained for the
\emph{normalized} chain-complex, the normalization consisting in this case in
replacing every occurence of \(C\) by \(\oC\), see~\cite[X.2.2]{MCLN2}.}, our
horizontal differential~\(d_h\).

In this way every component of the horizontal differential \(d_h: M \otimes
\oC^{\otimes p} \otimes N \rightarrow M \otimes \oC^{\otimes (p+1)} \otimes N\)
is a chain-complex morphism; as usual to obtain a bicomplex, we must transform
the commutative squares into anticommutative squares, this is the role of the
factor \((-1)^n\) in the formula for the vertical differential, which can also
be considered as an effect of the desuspension process.

Many authors prefer to call \(d_v\) the \emph{tensorial} differential and
\(d_h\) the \emph{cosimplicial} differential. Our terminology, vertical and
horizontal refers to the bicomplex structure for our Cobar, which is very
important in effective homology.

\begin{thr}\label{93176}---
A general algorithm computes:
\[
(M_{EH}, C_{EH}, N_{EH}) \mapsto \Cobar^C(M,N)_{EH}
\]
where: \senumerate
\item
\(C_{EH}\) is a 1-reduced differential coalgebra with effective homology;
\item
\(M_{EH}\) (resp. \(N_{EH}\)) is a right (resp. left) \(C\)-comodule with
effective homolohy.
\item
The result \(\Cobar^C(M,N)_{EH}\) is a version \emph{with effective homology}
of the Cobar construction \(\Cobar^C(M, N)\).
\end{enumerate}
\end{thr}

\proof It is a simple application of the Bicomplex Reduction
Theorem~\ref{61564}. As usual, let us use the notation \(C \lrdc \hC \rrdc EC\)
for the given equivalence between \(C\) and some effective chain-complex
\(EC\), the same for \(M\) and \(N\). In a first step, we cancel the horizontal
differential of \(\Cobar^C(M, N)\), which is nothing but replacing the
\(C\)-coproduct by \(\Delta_0(x) = 1 \otimes x + x \otimes 1\), the unit 1
being defined by the coaugmentation, and for example the \(M\)-coproduct by
\(\Delta_0(x) = x \otimes 1\). We so obtain a simplified \(\Cobar^{C_0}(M_0,
N_0)\) which is a banal direct sum of tensor products and as a simple
consequence of Proposition~\ref{49390}, we obtain an equivalence:
\[
\Cobar^{C_0}(M_0, N_0) \lrdc \Cobar^{\hC_0}(\widehat{M}_0, \widehat{N}_0) \rrdc
\Cobar^{EC_0}(EM_0, EN_0)
\]
where the 0-index signals the coproduct is made or defined as trivial.

Now we reinstall into the initial Cobar the horizontal differential; this is a
perturbation. For the left hand member of the new equivalence to be
constructed, the so-called Easy Perturbation lemma~\ref{73716} must be firstly
applied. We obtain:
\[
\Cobar^{C}(M, N) \lrdc \widetilde{\Cobar}^{\raisebox{-5pt}{\scriptsize
\(\hC\)}}(\widehat{M}, \widehat{N}) \ \raisebox{-4pt}{\rule{2pt}{16pt}} \rrdc
\Cobar^{EC_0}(EM_0, EN_0)
\]
The tilde above the \(\widetilde{\Cobar}\) explains there is some similarity
between the new differential installed on the central object and a Cobar
structure, but it is not actually a Cobar construction\footnote{The notion of
\(A_\infty\)-structure is designed to handle such a situation~\cite{STSH}.}.
The vertical bar \raisebox{-4pt}{\rule{2pt}{16pt}} signals the right hand
reduction is no longer valid, because of the central perturbation.

\noindent
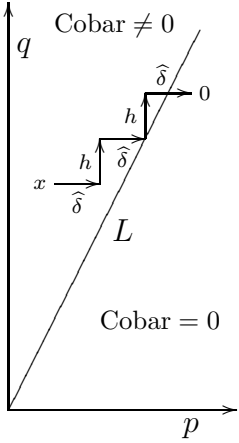
\begin{minipage}{0.22\textwidth}
\mbox{\begin{xy}<0.6cm,0cm>:<0cm,0.6cm>::
 (0,0)*{\hspace{0pt}} ;
 (4,0)*+!U{p} ; (0,8)*+!L{q} ;
 (1,5)*+!R{\scriptstyle x} ;
 (4,7)*+!L{\scriptstyle 0} ;
 (2,2)*!L{\mbox{\footnotesize\(\Cobar = 0\)}} ;
 (1,8.5)*!L{\mbox{\footnotesize\(\Cobar \neq 0\)}}
 \ar (0,0) ; (5,0)
 \ar (0,0) ; (0,9)
 \ar@{-}_{\textstyle L} (0,0) ; (4.2,8.4)
 \ar_{\hdl} (1,5) ; (2,5)
 \ar^h (2,5) ; (2,6)
 \ar_{\hdl} (2,6) ; (3,6)
 \ar^h (3,6) ; (3,7)
 \ar^<(0.35){\hdl} (3,7) ; (4,7)
 \end{xy}}
\end{minipage}\hfill
\begin{minipage}{0.60\textwidth}
\hspace*{20pt} For the right hand part of the equivalence, we must apply the
actual BPL and there remains to verify the nilpotency condition. It is here
that the 1-reduced property is required for our coalgebra~\(C\); because of
this property, the grading in column \(p\) begins at the ordinate \(2p\): the
bicomplex is null under a line \(L\) of slope 2. The path which is followed
when iterating the composition~\(h\hdl\) (notations of Theorem~\ref{07404}) is
a stairs of ``slope'' 1 which eventually goes under the line~\(L\), into an
area where the Cobar bicomplex is null: the nilpotency hypothesis is satisfied.
\end{minipage}

The basic perturbation lemma produces a reduction between
\(\widetilde{\Cobar}^{\raisebox{-5pt}{\scriptsize \(\hC\)}}(\widehat{M},
\widehat{N})\) and some pseudo-Cobar to be denoted as
\(\widetilde{\Cobar}^{\raisebox{-5pt}{\scriptsize \(EC\)}}(EM, EN)\). We
finally have an eauivalence:
\[\label{99700}
\Cobar^{C}(M, N) \lrdc \widetilde{\Cobar}^{\raisebox{-5pt}{\scriptsize
\(\hC\)}}(\widehat{M}, \widehat{N}) \rrdc
\widetilde{\Cobar}^{\raisebox{-5pt}{\scriptsize \(EC\)}}(EM, EN)
\]
Again because of the \emph{\(1\)-reduced} hypothesis for the coalgebra \(C\),
the right hand chain complex is \emph{effective}: the total degree is defined
as \(q-p\) and a homogeneous component of this chain complex is made of pieces
installed on a line of slope 1; but the intersection of this line with the
triangle ``Cobar \(\neq 0\)'' is finite, and the corresponding homogeneous
component therefore is \emph{effective}. \QED

It must be noted the perturbation lemma transforms \(\Cobar^{EC_0}(EM_0,
EN_0)\), a bicomplex, in fact without any horizontal differential, into
\(\widetilde{\Cobar}^{\raisebox{-5pt}{\scriptsize \(EC\)}}(EM, EN)\), a
multicomplex with arrows \(d^r_{p,q}\), maybe non-null for arbitrary values of
\(r\); see Definition~\ref{12080} and Theorem~\ref{61564}. The extra arrows so
defined are at the origin of the notion of \(A_\infty\)-coalgebra~\cite{STSH}.

A particular case of the Cobar contraction is crucial.

\begin{thr}\label{64957}---
The Cobar construction \(\Cobar^{C}(C, \fR)\) is a coresolution of \(\fR\).
\end{thr}

The prefix `co' in coresolution means it is an injective resolution \(\fR
\rightarrow \Cobar^{C}(C, \fR)\) instead of a projective resolution \(\fR
\leftarrow \Cobar^{C}(C, \fR)\).

We recall \(C\) is assumed reduced, so that \(C_0\) is isomorphic to the ground
ring~\(\fR\), which induces a left \(C\)-comodule structure \(\fR \rightarrow
C_0 \subset C = C \otimes \fR\).

\proof The contraction \(h\) of the complex \(\Cone(\fR \rightarrow \Cobar^C(C,
\fR))\) is defined by:
\[
h(c_0 \otimes c_1 \otimes \cdots \otimes c_n \otimes 1_\fR) = \eta(c_0) c_1
\otimes \ldots \otimes c_n \otimes 1_\fR
\]
The verification is a simple calculation. In particular \(\eta(c_0) \neq 0\)
only if \(c_0 \in C_0 \subset C\). \QED

\subsection{The \emph{effective} Eilenberg-Moore spectral sequence.}

In this section, a simplicial fibration \(F \hookrightarrow [E = F \times_\tau
B] \rightarrow B\) is given; in particular a group action \(G \times F
\rightarrow F\) for some simplicial group \(G\) and a twisting function \(\tau:
B \rightarrow G\) are present. In the Kenzo implementation, the twisting
function \(\tau\) \underline{\emph{is}} the fibration, more exactly, the
underlying principal fibration. Two equivalences \(C_\ast(E) \lrdc
\widehat{E}_\ast \rrdc EE_\ast\) and \(C_\ast(B) \lrdc \widehat{B}_\ast \rrdc
EB_\ast\) between the chain-complexes of \(E\) and \(B\) and \emph{effective}
chain-complexes \(EE_\ast\) and \(EB_\ast\) are given too, describing the total
space \(E\) and the base space \(B\) as simplicial sets \emph{with effective
homology}. Note no homological information is required for the structural group
\(G\) which can be any kind of \emph{locally effective} simplicial group. The
simplicial base space \(B\) \emph{is 1-reduced}: a unique vertex, the base
point, and no non-degenerate 1-simplex; the coalgebra \(C_\ast(B)\) is also
1-reduced and the coaugmentation ideal \(\overline{C_\ast(B)}\) begins only in
degree 2.

\begin{thr}\label{82517}---
A general algorithm computes:
\[
[F, G, B_{EH}, \tau, E_{EH}] \mapsto F_{EH}
\]
where \senumerate
\item
\(F\), \(G\), \(B\), \(\tau\) and \(E\) are as explained above;
\item
\(E_{EH}\) (resp. \(B_{EH}, F_{EH}\)) is a version with effective homology of
\(E\) (resp. \(B, F\)).
\end{enumerate}
\end{thr}

In other words, if the effective homology of the total space and of the base
space are known, an algorithm computes the effective homology of the fibre
space. Victor Gugenheim~\cite{GGNH2} computed an effective chain-complex, the
homology of which is guaranteed being the homology of the fibre space, but this
is not enough to iterate the process: an equivalence between the effective
chain-complex and the chain-complex of the fibre space is then necessary; it is
the key point to obtain a solution of the \emph{Adams' problem}.

The next proposition which has its own interest will be used.

\begin{prp}\label{00533}---
In the data of the Basic Perturbation Lemma~\ref{07404}, if the relation \(\hdl
g = 0\) holds, then the resulting perturbation \(\delta\) for the small
chain-complex is null: \(\delta = 0\). The same if \(f \hdl = 0\).
\end{prp}

It is a paradoxical result. Usually the BPL is used to construct a \emph{new}
interesting differential for the small graded module \(C_\ast\). It happens
sometimes we are interested by two \emph{different} reductions from the big
graded module \(\hC_\ast\) provided with two \emph{different} differentials
over the \emph{same} small chain-complex \(C_\ast\).

\proof Just glance at the formula \(\delta = f \hdl \phi g = f \psi \hdl g\) at
page~\pageref{80476}. \QED

\noindent\textsc{Proof of Theorem~\ref{82517}}. The Eilenberg-Zilber
Theorem~\ref{20016} produces a reduction \(C_\ast(F \times B) \rrdc C_\ast(F)
\otimes C_\ast(B)\) and the twisted Eilenberg-Zilber Theorem~\ref{80728}
another reduction \(C_\ast(F \times_\tau B) \rrdc C_\ast(F) \otimes_t
C_\ast(B)\). Theorem~\ref{64957} produces a reduction \(\fR \lrdc
\Cobar^{C_\ast(B)}(C_\ast(B), \fR)\); applying to the last reduction the
functor \(C_\ast(F) \otimes \!\! \functorobj\) gives again another reduction:
\[
(f,g,h): C_\ast(F) \lrdc \Cobar^{C_\ast(B)}(C_\ast(F) \otimes C_\ast(B), \fR).
\]
Note in the Cobar the \(C_\ast(B)\)-comodule structure of \(C_\ast(F) \otimes
C_\ast(B)\) is induced by the canonical projection \(C_\ast(F) \otimes
C_\ast(B) \rightarrow C_\ast(B)\), which projection requires the coaugmentation
\(C_\ast(F) \rightarrow \fR\). The same for the twisted tensor product below.

We intend to replace in the last reduction the ordinary tensor product
\(C_\ast(F) \otimes C_\ast(B)\) by the twisted one \(C_\ast(F) \otimes_t
C_\ast(B)\); this is a perturbation \(\hdl\) of the Cobar differential. Is the
condition \(\hdl g = 0\), which would allow us to apply
Proposition~\ref{00533}, satisfied? Let us call \(\ast_B\) the base point of
\(B\). The component \(g\) of the reduction maps \(C_\ast(F)\) onto the
\emph{sub-chain-complex} \(C_\ast(F) \otimes C_\ast(*_B)\) inside the 0-column
of the Cobar bicomplex; this a subcomplex not only inside the 0-column, but
also in the Cobar: the pseudo-coproduct \(C_\ast(\ast_B) \rightarrow
C_\ast(\ast_B) \otimes \overline{C_\ast(B)}\) is null, because of the
restriction to the coaugmentation ideal which cancels the 0-component. This
sub-chain-complex is left unchanged by the perturbation, for the \emph{base}
fibre of the total space is not modified by the twisting process. No
perturbation on the sub-chain-complex \(C_\ast(F) \otimes C_\ast(*_B)\) and the
condition \(\hdl g\) holds. Proposition~\ref{00533} produces a reduction:
\begin{equation}
(f',g,h'): C_\ast(F) \lrdc \Cobar^{C_\ast(B)}(C_\ast(F) \otimes_t C_\ast(B),
\fR). \tag{672}
\end{equation}
because the \(g\)-component is also unchanged.

There remains to apply Theorem~\ref{93176}. The base space \(B\) is given with
effective homology and the same for \(E\); in other words an equivalence:
\[
C_\ast(F \times_\tau B) \lrdc \widehat{E}_\ast \rrdc EE_\ast
\]
is given. Composing the left hand reduction with the twisted Eilenberg-Zilber
reduction \(C_\ast(F) \otimes_t C_\ast(B) \lrdc C_\ast(F \times_\tau B)\) gives
an equivalence:
\[
C_\ast(F) \otimes_t C_\ast(B) \lrdc \widehat{E}_\ast \rrdc EE_\ast,
\]
in other words the chain complex \(C_\ast(F) \otimes_t C_\ast(B)\) is with
effective homology. Theorem~\ref{93176} can be applied which produces an
equivalence:
\[
\Cobar^{C_\ast(B)}(C_\ast(F) \otimes_t C_\ast(B), \fR) \lrdc \rrdc
\widetilde{\Cobar}^{\raisebox{-5pt}{\scriptsize \(EB_\ast\)}}(EE_\ast, \fR).
\]
Finally composing the left hand reduction with the reduction (672) above gives
an equivalence:
\[
C_\ast(F) \eqvl \widetilde{\Cobar}^{\raisebox{-5pt}{\scriptsize
\(EB_\ast\)}}(EE_\ast, \fR)
\]
where the right hand chain-complex is effective. \QED

\subsection{Adams' problem.}\label{05508}

Our effective version of the Eilenberg-Moore spectral sequence gives a very
simple solution to Adams' problem. Frank Adams (\cite{ADMS}, see
also~\cite{ADHL}) designed an algorithm computing the homology groups of the
\emph{first} loop space \(\Omega X\) of a 1-reduced simplicial set; stated in
our framework, Adams'result is the following.

\begin{thr}\textbf{\emph{(Adams' Theorem)}} ---
Let \(X\) be a 1-reduced simplicial set. Then there exists a canonical
isomorphism between \(H_\ast(\Omega X ; \fR)\) and
\(H_\ast(\Cobar^{C_\ast(X)}(\fR, \fR))\).
\end{thr}

If \(X\) is a finite 1-reduced simplicial set, the Cobar is effective and the
homology groups are computable. Adams then asked for some analogous solution
for the \emph{iterated} loop space \(\Omega^n X\). Eighteen (\(!\)) years
later, Hans Baues~\cite{BAUS1} gave a solution for the second loop space
\(\Omega^2 X\); it depends on an ingenious possible geometrical model for the
second loop space; but again it is not possible to extend this model to the
third loop space \(\Omega^3(X)\)\ldots

The problem is in fact in the \emph{non-constructive} nature of Adams' solution
for the first loop space. Elementary homological algebra shows that for
reasonable ground rings \(\fR\) there \emph{exists} an equivalence
\(C_\ast(\Omega X) \eqvl \Cobar^{C_\ast(X)}(\fR, \fR)\), but the exact nature
of this equivalence is not studied.

Our \emph{effective} Eilenberg-Moore spectral sequence on the contrary will
\emph{constructively} prove the existence of this equivalence; and then the
iteration of the process is obvious, giving our solution to Adams' problem. So
simple that it is not difficult to implement it on a computer, leading to
programs computing homology groups of loop spaces otherwise so far unreachable.

We must make more precise Theorem~\ref{13368}.

\begin{thr}\label{85719}---
Let \(B\) be a 1-reduced \emph{locally effective} simplicial set. Then the path
space \(PB\) defined in Theorem~\ref{13368} has \emph{effective} homology.
\end{thr}

\proof See for example~\cite[Chapter VI]{MAY}. An \emph{explicit} contraction
is there \emph{constructed} for the chain-complex \(C_\ast(PB) = C_\ast(\Omega
B \times_\tau B)\) for the appropriate twist~\(\tau\) defining the path space.
It is easy to organize this contraction as a reduction \(C_\ast(\Omega B
\times_\tau B) \rrdc \fR\). \QED

\begin{crl}\label{65095}\textbf{\emph{(Effective Adams' Theorem)}}---
A general algorithm computes:
\[
B_{EH} \mapsto (\Omega B)_{EH}
\]
where \(B_{EH}\) is a 1-reduced simplicial set \emph{with effective homology}
(input) and \((\Omega B)_{EH}\) is a version \emph{with effective homology} of
the loop space (output).
\end{crl}

\proof Apply Theorem~\ref{82517} to the fibration: \(\Omega B \hookrightarrow
[E = \Omega B \times_\tau B] \rightarrow B\). The base space \(B\) is given
with its effective homology and the effective homology of the total space \(E =
PB = \Omega B \times_\tau B\) is computed by Theorem~\ref{85719}. \QED

\begin{crl}\thrtitle{Solution to Adams' problem} ---
A general algorithm computes:
\[
(n, X_{EH}) \mapsto (\Omega^n X)_{EH}
\]
where the input \(X_{EH}\) is an \(r\)-reduced simplicial set \emph{with
effective homology}, \(r \geq n\), and the output is a version \emph{with
effective homology} of the \(n\)-th loop space. In particular the
\emph{ordinary} homology of this iterated loop space is computable.
\end{crl}

The qualifier \emph{\(r\)-reduced} for \(X\) means in the simplicial structure
of \(X\) there is no non-degenerate simplex in dimension \(\leq r\) except the
base point in dimension~0.

\proof The simplicial model \(\Omega X\) for the \(r\)-reduced simplicial set
\(X\) is itself \((r-1)\)-reduced. It is sufficient to successively apply \(n\)
times Corollary~\ref{65095}. \QED

\subsection{Other Eilenberg-Moore spectral sequences.}

The reader can be puzzled by the non-symmetric presence of the ground ring
\(\fR\) in \(\Cobar^{C_\ast(B)}(E, \fR)\), the main ingredient in the
Eilenberg-Moore process. In fact our presentation is a particular case of a
more general situation.

\begin{dfn}---
\emph{Let \(F \hookrightarrow [E = F \times_\tau B] \rightarrow B\) be a
fibration and \(\beta: B' \rightarrow B\) be a simplicial map. These data
define an \emph{induced fibration} \(F \hookrightarrow E' \rightarrow B'\).
 }
\end{dfn}

The twisting function \(\tau\) is some ``degree'' -1 map \(\tau: B \rightarrow
F\), see Definition~\ref{51596}. The composition \(\tau' = \tau \beta\) also is
a twisting function, defining the induced fibration. Another point of view
consists in thinking of the total space \(E'\) as the cartesian product \(E' =
B' \times_B E\), where the set of \(n\)-simplices \(E'_n\) is \(E'_n =
\{(\sigma', \sigma) \in B'_n \times E_n \st f(\sigma') = \pr(\sigma)\}\) if
\(\pr\) is the projection \(\pr: E \rightarrow B\). Both definitions are
elementarily equivalent.

\begin{thr}\thrtitle{First effective Eilenberg-Moore spectral sequence} ---
A general algorithm computes:
\[
(B_{EH}, F, G, \tau, E_{EH}, B'_{EH}, \beta) \mapsto E'_{EH}
\]
where all the ingredients are as above, the \(EH\)-index meaning the
corresponding object is given (case of \(B\), \(E\) and \(B'\)) or produced
(case of \(E'\)) with \emph{effective homology}.
\end{thr}

Theorem~\ref{82517} is the particular case where \(\beta\) is the inclusion of
the base point in the base space \(\ast_B \hookrightarrow B\); the induced
fibration is then simply \(F \hookrightarrow F \rightarrow \ast_B\). The same
method constructs in the general case an equivalence:
\[
C_\ast(E') \eqvl \Cobar^{C_\ast(B)}(C_\ast(E), C_\ast(B')).
\]
Note in particular \(\beta\) defines, even if the map \(\beta: B' \rightarrow
B\) is not a fibration, a \(C_\ast(B)\)-comodule structure on \(C_\ast(B')\),
which makes coherent the definition of the Cobar. The proof is the same, you
just have to replace the right hand \(\fR\) in the various Cobars by
\(C_\ast(B')\).

What about the symmetric ``division''? If \(F \hookrightarrow E \rightarrow B\)
is a fibration, we could also be interested by something like \(B = E/F\) and
we would like to deduce the effective homology of the base space \(B_{EH}\)
from \(F_{EH}\) and \(E_{EH}\); possible? Yes, it is the second Eilenberg-Moore
spectral sequence. The general case works as follows. The main ingredients are
two simplicial sets \(E\) and \(E'\) and a simplicial group \(G\). A right
(resp. left) action is given \(\alpha: E \times G \rightarrow E\) (resp.
\(\alpha': G \times E' \rightarrow E'\)). This defines a cocartesian product
\(E \times_G E' := (E \times E') / \sim_G\) where the equivalence relation
\(\sim_G\) makes equivalent \((\alpha(\sigma, \gamma), \sigma') \sim_G (\sigma,
\alpha'(\gamma, \sigma'))\) when \(\sigma \in E_n\), \(\sigma' \in E'_n\) and
\(\gamma \in G\); think of the definition of a tensor product which, from a
categorical point of view, is analogous.

In the first Eilenberg-Moore spectral sequence, there must be a fibration
connecting the factor \(E\) of \(B' \times_B E\) with the base space \(B\). The
second spectral sequence depends on an analogous requirement: one action, for
example the first one \(\alpha: E \times G \rightarrow E\) must define a
principal fibration \(G \hookrightarrow E \rightarrow E/G\) where \(E/G\) is
nothing but \(E \times_G \{\ast\}\) for the trivial action \(G \times \{\ast\}
\rightarrow \{\ast\}\). If this condition is satisfied, an analogous
\emph{effective} spectral sequence is obtained.

\begin{thr}\label{50468}\thrtitle{Second \emph{effective} Eilenberg-Moore spectral sequence}
---
A general algorithm computes:
\[
(G_{EH}, E_{EH}, E'_{EH}, \alpha, \alpha') \mapsto (E \times_G E')_{EH}
\]
where the ingredients are as above, the \(EH\)-index meaning the corresponding
object is provided \emph{with effective homology}. The structural group \(G\)
is assumed 0-reduced.
\end{thr}

The proof is the same, the key intermediate ingredient being the Bar
construction \(\textrm{Bar}^{C_\ast(G)}(C_\ast(E), C_\ast(E'))\). In fact the
various multiplicative structures define a structure of differential algebra
over \(C_\ast(G)\), sometimes called the Pontrjagin structure, and
\(C_\ast(G)\)-module structures on \(C_\ast(E)\) and \(C_\ast(E')\). Note this
time the effective homology of the structural group is also \emph{required}: it
plays the role of the base space \(B\) in the symmetric situation.

It has been explained the loop space \(\Omega X\) can be considered as a
twisted inverse of the original space \(X\), for the appropriate twisted
product \(\Omega X \times_\tau X\) is contractible, has the homotopy type of a
point, and the point is the unit in the topological world. In the same way, the
classifying space construction~\cite[\S21]{MAY} allows one to construct the
\emph{universal fibration} \(G \hookrightarrow EG \rightarrow BG\) where the
total space \(EG\) is contractible and the base space is an ``inverse'' of G.
In particular if \(G\) is the Eilenberg-MacLane space \(K(\pi, 1)\), see
Section~\ref{70076}, then \(BG = K(\pi, 2)\) and more generally \(B^{n-1} G =
K(\pi, n)\). For sensible commutative groups \(\pi\), Theorem~\ref{50468} can
compute the effective homology of \(K(\pi,n)\). Coming back to the rough
explanations given in Section~\ref{62859} about the computation of homotopy
groups, it is easy to prove:

\begin{thr}---
A general algorithm computes:
\[
(n, X_{EH}) \rightarrow \pi_n X
\]
where \(X_{EH}\) is a 1-reduced simplicial set with effective homology and
\(\pi_n X\) is the \(n\)-th homotopy group of \(X\).
\end{thr}

This is a powerful generalization of Edgar Brown's Theorem~\cite{BRWNE1}: the
scope is much larger than in Edgar Brown's paper where the simplicial set is
assumed \emph{finite}, and the proof is more conceptual, so conceptual that the
machine implementation is not very difficult. See the Kenzo
documentation~\cite{DRSS}.

\section{The claimed Postnikov ``invariants''\protect\footnote{This section is a rough copy
of the paper~\cite{RBSR8}; which explains some redundancies and also some gaps
with respect to the previous sections; in spite of these defects, we think a
reader having reached this point of the text could be interested by this
section, making obvious a  surrealist error in the standard terminology of
Algebraic Topology.}}

\subsection{Introduction.}

\hfill\fbox{\parbox[t]{0.42\textwidth}{\footnotesize{\emph{
 \hspace*{0pt}\hfill As yet we are ignorant\\
 \hspace*{0pt}\hfill of an effective method of computing\\
 \hspace*{0pt}\hfill the cohomology of a Postnikov complex\\
 \hspace*{0pt}\hfill from
 \(\pi_n\) and \(k^{n+1}\)~\emph{\cite{EDM}}.
}}}}\vspace{20pt}

When this paper is written, the so-called \emph{Postnikov invariants} (or
\(k\)-invariants) are roughly fifty years old~\cite{PSTN1}; they are a key
component of standard Algebraic Topology. This notion is so important that it
is a little amazing to observe some important \emph{gaps} are still present in
our working environment around this subject, still more amazing to note these
gaps are seldom considered. One of these ``gaps'' is unfortunately an
\emph{error}, widely spread, and easy to state: the terminology ``Postnikov
\emph{invariants}'' is incorrect: any sensible definition of the
\emph{invariant} notion leads to the following conclusion: the Postnikov
invariants are not\ldots\ invariants. This is true even in the simply connected
case and to make easier the understanding, we restrict our study to this case.

First, several interesting questions of \emph{computability} are arisen by the
very notion of Postnikov invariant. It is surprisingly difficult fo find
citations related to this computability problem, as though this problem was
unconsciously ``hidden'' (\(?\)) by the topologists. The only significant one
found by the authors is the EDM title quotation\footnote{Other possible
quotations are welcome.}. In fact there are \emph{two} distinct problems of
this sort.

On one hand, if a simply connected space is presented as a \emph{machine
object}, does there exist a general algorithm computing its Postnikov
invariants? The authors have designed a general framework for
\emph{constructive} Algebraic Topology, giving in particular such a general
algorithm~\cite{SRGR3,RBSR6}. In the text, this process is formalized as a
functor \(\bfSP: \SSEH \widetilde{\times} I \rightarrow \cP\) where \(\SSEH\)
is an appropriate category of \emph{computable} topological spaces, and \(\cP\)
is the \emph{Postnikov category}. We will explain later the nature of the
factor~\(I\), in fact the heart of our subject.

On the other hand, a \emph{converse} problem must be considered. When a
Postnikov tower is given, that is, a collection of homotopy groups and
\emph{relevant} Postnikov invariants, how to construct the corresponding
topological space? The computability problem stated in the title quotation is a
(small) part of this converse problem. Again, our notion of \emph{constructive}
Algebraic Topology entirely solves it. The resulting computer program
Kenzo~\cite{DRSS} allows us to give a simple concrete illustration. In fact it
will be explained it is not possible to properly \emph{state} this
problem\ldots\ without having a solution of it! Again a strange situation to
our knowledge not yet considered by the topologists. Our solution for the
converse problem will be formalized as a functor \(\bfPS: \cP \rightarrow
\SSEH\).

There is a lack of symmetry between the functors \(\bfSP: \SSEH
\widetilde{\times} I \rightarrow \cP\) and \(\bfPS: \cP \rightarrow \SSEH\).
Instead of our functor \(\bfSP: \SSEH \widetilde{\times} I \rightarrow \cP\), a
simpler functor \(\bfSP: \SSEH \rightarrow \cP\), without the mysterious factor
\(I\), is expected, but in the current state of the art, such a functor
\emph{is not} available. It is a consequence of the following \emph{open}
problem: let \(P_1, P_2 \in \cP\) be two Postnikov towers; does there exist an
algorithm deciding whether \(\bfPS(P_1)\) and \(\bfPS(P_2)\) have the same
homotopy type or not? The remaining uncertainty is measured by the factor
\(I\). And because of this uncertainty, the so-called Postnikov invariants are
not\ldots\ invariants: the context clearly says they should be invariants of
the \emph{homotopy type}, but such a claim is equivalent to a solution of the
above decision problem.

It is even possible this decision problem \emph{does not have} any solution; in
fact, our Postnikov decision problem can be translated into an arithmetical
decision problem, a subproblem of the general tenth Hilbert problem to which
Matiyasevich gave a negative answer~\cite{MTSV}. If our decision problem had in
turn a negative answer, it would be \emph{definitively} impossible to transform
the common Postnikov invariants into \emph{actual} invariants.

\subsection{The Postnikov category and the PS functor.}\label{mocoi}

Defining a functor \(\bfPS: \cP \rightarrow \SSEH\) in principle consists in
defining the \emph{source} category, here the \emph{Postnikov category}
\(\cP\), the \emph{target} category, the \emph{simplicial set category}
\(\SSEH\), and \emph{then}, finally, the functor \(\bfPS\) itself. It happens
this is not possible in this case~: the Postnikov category \(\cP\) and the
functor \(\bfPS\) are \emph{mutually recursive}. More precisely, an object \(P
\in \cP\) is a limit \(P = \lim P_n\), every \(P_n\) being also an element of
\(\cP\). Let \(\cP_n\), \(n \geq 1\),  be the Postnikov towers limited to
dimension~\(n\). The definition of \(\cP_{n+1}\) \emph{needs} the partial
functor \(\bfPS_n: \cP_n \rightarrow \SSEH\) where \(\bfPS_n = \bfPS|\cP_n\)
and this is why the definitions of \(\cP\) and \(\bfPS\) are mutually
recursive.

We work only with simply connected spaces, the homotopy (or \(\bZ\)-homology)
groups of which being of \emph{finite type}. It is essential, when striving to
define \emph{invariants}, to have exactly \emph{one} object for every
isomorphism class of groups of this sort, so that we adopt the following
definition. No \(p\)-adic objects in our environment, which allows us to denote
\(\bZ / d \bZ\) by \(\bZ_d\); in particular \(\bZ_0 = \bZ\).

\begin{dfn}\label{szglm} ---
\emph{A \emph{canonical group} (abelian, of finite type) is a product
\(\bZ_{d_1} \times \cdots \times \bZ_{d_k}\) where the non-negative integers
\(d_i\) satisfy the divisibility condition: \(d_i\) divides \(d_{i+1}\) for
\mbox{\(1 \leq i < k\)}.}
\end{dfn}

Every abelian group of finite type is isomorphic to \emph{exactly one}
canonical group, for example the group \(\bZ^2 \oplus \bZ_6 \oplus \bZ_{10}
\oplus \bZ_{15}\) is isomorphic to the unique canonical group \(\bZ_{30} \times
\bZ_{30} \times \bZ_0 \times \bZ_0\); but such an isomorphism is
\emph{not}\ldots\ canonical; for example, for the previous example, there
exists an infinite number of such isomorphisms, and we will see this is the key
point preventing us from qualifying the Postnikov invariants as invariants.

\begin{dfn} ---
\emph{The category \(\SSEH\) is the category of the \emph{simply connected
simplicial sets with effective homology} described in \cite{RBSR6}.}
\end{dfn}

The framework of the present paper does not allow us to give the relatively
complex definition of this category. Roughly speaking, an object of this
category is a \emph{machine object} coding a (possibly infinite) simply
connected simplicial set with \emph{known} homology groups; furthermore a
\emph{complete} knowledge of the homology is required: mainly every homology
class has a canonical representant cycle, an algorithm computes the homology
class of every cycle, and if two cycles \(c_0\) and \(c_1\) are homologous, an
algorithm computes a chain \(C\) with \(\partial C = c_1 - c_0\). For example
it is explained in~\cite{RBSR7} that \(X = \Omega(\Omega
(P^\infty(\bR)/P^3(\bR)) \cup_4 D^4) \cup_2 D^3\) is an object of \(\SSEH\) and
the Kenzo program does compute the first homology groups of it, in the detailed
form just briefly sketched. More generally every ``sensible'' simply connected
space with homology groups of finite type has the homotopy type of an object of
\(\SSEH\); this statement is precisely stated in~\cite{RBSR6}, the proof is not
hard, it is only a repeated application of the so-called \emph{homological
perturbation lemma}~\cite{BRWNR1} and the \emph{most detailed} proof is the
Kenzo computer program itself~\cite{DRSS}, a Common Lisp text of about 16,000
lines.

The definitions of the category \(\cP\) and the functor \(\bfPS\) are
\emph{mutually recursive} so that we need a starting point.

\begin{dfn} ---
\emph{The category \(\cP_1\) has a unique object, the void sequence \(()_{2
\leq n \leq 1}\), the trivial Postnikov tower, and the functor \(\bfPS_1\)
associates to this unique object the trivial element \mbox{\(\bast \in \SSEH\)}
with only a base point.}
\end{dfn}

The next definitions of the category \(\cP_n\) and the functor \(\bfPS_n\)
assume the category \(\cP_{n-1}\) \emph{and} the functor \(\bfPS_{n-1}:
\cP_{n-1} \rightarrow \SSEH\) are already available.

\begin{dfn} ---
\emph{An object \(P_n \in \cP_n\) is a sequence \(((\pi_m, k_m))_{2 \leq m \leq
n}\) where:
\begin{itemize}
\item
\(((\pi_m, k_m))_{2 \leq m \leq n-1}\) is an element \(P_{n-1} \in \cP_{n-1}\);
\item
The component \(\pi_n\) is a canonical group;
\item
The component \(k_n\) is a cohomology class \(k_n \in
H^{n+1}(\bfPS_{n-1}(P_{n-1}), \pi_n)\)~;
\end{itemize}
Let us denote \(X_{n-1} = \bfPS_{n-1}(P_{n-1})\). The cohomology class \(k_n\)
classifies a fibration~:
\[
 K(\pi_n, n) \hookrightarrow K(\pi_n, n) \times_{k_n} X_{n-1} \twoheadrightarrow
 X_{n-1} \stackrel{k_n}{\longrightarrow} K(\pi_n, n+1) = BK(\pi_n, n).
\]
\begin{itemize}
\item
Then the functor \(\bfPS_n\) associates to \(P_n = ((\pi_m, k_m))_{2 \leq m
\leq n} \in \cP_n\) a version \emph{with effective homology} \(X_n =
\bfPS_n(P_n)\) of the total space \(K(\pi_n, n) \times_{k_n} X_{n-1}\).
\end{itemize}}
\end{dfn}

In particular our version \emph{with effective homology} of the Serre spectral
sequence and our versions \emph{with effective homology} of the
Eilenberg-MacLane spaces \(K(\pi, n)\) allow us to construct a version also
with effective homology of the total space \(K(\pi_n, n) \times_{k_n}
X_{n-1}\), here denoted by \(X_n\). We will give a typical small Kenzo
demonstration at the end of this section.

A canonical forgetful functor \(\cP_n \rightarrow \cP_{n-1}\) is defined by
forgetting the last component of \(((\pi_m, k_m))_{2 \leq m \leq n}\), which
allows us to define \(\cP\) as the projective limit \(\displaystyle \cP =
\lim_\leftarrow \cP_n\). If \(X_{n-1}\) is a simplicial set, the
\((n-1)\)-skeletons of \(X_{n-1}\) and \(K(\pi_n, n) \times_{k_n} X_{n-1}\) are
the same (for the standard model of \(K(\pi_n, n)\)), so that if
\(\displaystyle P = \lim_\leftarrow P_n\), the limit \(\displaystyle \bfPS(P) =
\lim_\leftarrow \bfPS_n(P_n)\) is defined also as an object of \(\SSEH\). The
category \(\cP\) and the functor \(\bfPS: \cP \rightarrow \SSEH\) are now
\emph{properly} defined.

The homotopy groups \(\pi_m\)'s of a Postnikov tower \(((\pi_m, k_m))_{2 \leq
m}\) can be defined \emph{firstly} independently of the \(k_m\)'s, but \(k_n\)
can be \emph{properly defined} only when \(((\pi_m, k_m))_{2 \leq m < n}\) is
given \emph{and only if} the functor \(\bfPS_{n-1}\) is available in the
environment. In other words, if the problem of the title EDM quotation is not
solved, the very notion of Postnikov tower cannot be made \emph{effective}.

\subsection{Kenzo example.}\label{djmab}

Let us play the game consisting in constructing the beginning of a Postnikow
tower with a \(\pi_i = \bZ_2\) at each stage and the ``simplest''
\emph{non-trivial} Postnikov invariant. First \(P_1 = ()\) and \(X_1 =
\bfPS_1(P_1) = \bast\). As planned, we choose \(\pi_2 = \bZ_2\) and \(k_2 \in
H^3(X_1, \bZ_2) = 0\) is necessary null, no choice. So that we define \(P_2 =
((\bZ_2, 0))\) and \(X_2 = K(\bZ_2, 2)\). The Kenzo function \boxtt{k-z2} can
construct this space. We show a copy of the dialog between a Kenzo user and the
Lisp machine.

 \bmp
 \bmpi\verb|> (setf X2 (k-z2 2))|\empim
 \bmpi\verb|[K13 Abelian-Simplicial-Group]|\empi
 \emp

This dialog goes as follows. The Lisp \emph{prompt} is the bigger sign
`\texttt{>}'. The Lisp user enters a Lisp \emph{statement}, here
``\texttt{(setf X2 (k-z2 2))}''. The Maltese cross `\(\maltese\)' signals the
end of the statement to be executed, it is added here to help the reader, but
it is not visible on the user screen. When the Lisp statement is finished, Lisp
\emph{evaluates} it, the computation time can be a microsecond or a few days or
more, depending on the statement to be evaluated, and when the evaluation
\emph{terminates}, a Lisp object is \emph{returned}, most often it is the
``result'' of the computation. Here the \texttt{K13} object (the Kenzo object
\#13) is constructed and returned, it is an abelian simplicial group. A Lisp
statement ``\texttt{(setf some-symbol (some-function some-arguments))}'' orders
Lisp to make the function \texttt{some-function} work, using the arguments
\texttt{some-arguments}; this function creates some object which is
\emph{returned} (displayed) \emph{and} assigned to the symbol
\texttt{some-symbol}; in this way, the created object remains reachable through
the symbol locating it.

The \(\bZ\)-homology in dimensions 3 and 4 of \(X_2\) (the arguments 3 and 5
must be understood as defining \(3 \leq i < 5\)):

 \bmp
 \bmpi\verb|> (homology X2 3 5)|\empim
 \bmpi\verb|Homology in dimension 3 :|\empi
 \bmpi\verb|---done---|\empix
 \bmpi\verb|Homology in dimension 4 :|\empi
 \bmpi\verb|Component Z/4Z|\empi
 \bmpi\verb|---done---|\empi
 \emp

\noindent to be read \(H_3 = 0\) and \(H_4 = \bZ_4\). The universal coefficient
theorem implies \(H^4(X_2, \bZ_2) = \bZ_2\), there is only one non-trivial
possible \(k_3 \in H^4(X_2, \bZ_2)\) and the Kenzo function \boxtt{chml-clss}
(\underline{c}o\underline{h}omology \underline{cl}ass) constructs it.

 \bmp
 \bmpi\verb|> (setf k3 (chml-clss X2 4))|\empim
 \bmpi\verb|[K125 Cohomology-Class on K30 of degree 4]|\empi
 \emp

The attentive reader can be amazed to see this cohomology class defined on
\boxtt{K30} and not \(\boxtt{K13} = X_2\). The explanation is the following.
Let us consider the \emph{\underline{ef}fective
\underline{h}o\underline{m}ology} of \(X_2\):

 \bmp
 \bmpi\verb|> (efhm X2)|\empim
 \bmpi\verb|[K122 Equivalence K13 <= K112 => K30]|\empi
 \emp

This is a chain equivalence between the chain complex of the considered space
and some \emph{small} chain complex, here the chain complex \boxtt{K30}. In
fact it is a \emph{strong} chain equivalence, made of two \emph{reductions}
through the intermediate chain complex \boxtt{K112} (see~\cite{RBSR6} for
details). So that defining a cohomology class of \(X_2\) is equivalent to
defining such a class for \boxtt{K30}. A \emph{small} chain complex is a free
\(\bZ\)-chain complex of finite type in every dimension. The chain complex
\boxtt{K13} of the standard model of \(X_2 = K(\bZ_2, 2)\) is already of finite
type, but the complex \boxtt{K30} is much smaller. For example, in dimension 6,
\boxtt{K13} has \(27,\hspace{-3pt}449\) generators and \boxtt{K30} has only 5.

The \(k_3\) class allows us to define the fibration canonically associated:
\[
F_3 = \left\{K(\bZ_2, 3) \hookrightarrow K(\bZ_2, 3) \times_{k_3} X_2
\twoheadrightarrow X_2 \stackrel{k_3}{\longrightarrow} K(\bZ_2, 4)\right\}
\]

We have now the Postnikov tower \(P_3 = ((\bZ_2, 0), (\bZ_2, k_3))\) with \(X_3
= \bfPS(P_3) = K(\bZ_2, 3) \times_{k_3} X_2\). The Kenzo program can construct
our fibration \(F_3\) and its total space \(X_3\).

 \bmp
 \bmpi\verb|> (setf F3 (z2-whitehead X2 k3))|\empim
 \bmpi\verb|[K140 Fibration K13 -> K126]|\empix
 \bmpi\verb|> (setf X3 (fibration-total F3))|\empim
 \bmpi\verb|[K146 Kan-Simplicial-Set]|\empi
 \emp

The fibration is modelled as a \emph{twisting operator} \(\tau_3: X_2
\rightarrow K(\bZ_2, 3)\) which is nothing but an avatar of \(k_3\), and we can
verify the target of \(\tau_3\) is really \(K(\bZ_2, 3)\).

 \bmp
 \bmpi\verb|> (k-z2 3)|\empim
 \bmpi\verb|[K126 Abelian-Simplicial-Group]|\empi
 \emp

We continue to the next stage of our Postnikov tower. We ``choose'' again
\(\pi_4 = \bZ_2\), but what about the next Postnikov invariant \(k_4\)? We must
choose some \(k_4 \in H^5(X_3, \bZ_2)\), so that we are in front of the problem
stated in the framed EDM title quotation. Fortunately, the Kenzo program
\emph{knows} how to compute the \emph{necessary} \(H^5\), the Kenzo program
knows a (simple) solution for the EDM problem. In fact it knows the
\emph{effective} homology of the fibre space \(K(\bZ_2, 3)\):

 \bmp
 \bmpi\verb|> (efhm (k-z2 3))|\empim
 \bmpi\verb|[K268 Equivalence K126 <= K258 => K254]|\empi
 \emp

\noindent In the same way, it knows the \emph{effective} homology of \(X_2 =
K(\bZ_2, 2)\), and the implicitly used \emph{effective homology version} of the
Serre spectral sequence, available in Kenzo, determines the effective homology
of the twisted product \(X_3\):

 \bmp
 \bmpi\verb|> (efhm X3)|\empim
 \bmpi\verb|[K358 Equivalence K146 <= K348 => K344]|\empi
 \emp

\noindent The chain-complex \boxtt{K344} is of finite type, its homology groups
are computable, and in this way Kenzo can compute the \(\bZ\)-homology groups
of \(X_3\).

 \bmp
 \bmpi\verb|> (homology X3 2 6)|\empim
 \bmpi\verb|Homology in dimension 2 :|\empi
 \bmpi\verb|Component Z/2Z|\empi
 \bmpi\verb|---done---|\empix
 \bmpi\verb|Homology in dimension 3 :|\empi
 \bmpi\verb|---done---|\empix
 \bmpi\verb|Homology in dimension 4 :|\empi
 \bmpi\verb|Component Z/2Z|\empi
 \bmpi\verb|---done---|\empix
 \bmpi\verb|Homology in dimension 5 :|\empi
 \bmpi\verb|Component Z/4Z|\empi
 \bmpi\verb|---done---|\empi
 \emp

Finally the universal coefficient theorem implies \(H^5(X^3, \bZ_2) = \bZ_2
\oplus \bZ_2\), and there are exactly four ways to add a new stage at our
Postnikov tower with \(\pi_4 = \bZ_2\). Four possible Postnikov
\emph{invariants} \(k_4\). In this simple case, rather misleading, it is true
such a \(k_4\) is an \emph{invariant} of the homotopy type of the resulting
space, but in the general case, we will see the situation is much more
complicated; this will be explained Section~\ref{tbpxw}

The \boxtt{chml-clss} Kenzo function constructs in such a case the
cohomology-class ``dual'' to the generator of \(H_5(X_3, \bZ) = \bZ_4\).

 \bmp
 \bmpi\verb|> (setf k4 (chml-clss X3 5))|\empim
 \bmpi\verb|[K359 Cohomology-Class on K344 of degree 5]|\empi
 \emp

\noindent and the process can be iterated as before, giving the fibration
\(F_4\) associated to \(k_4\), and the total space \(X_4 = \bfPS_4(P_4) =
K(\bZ_2, 4) \times_{k_4} X_3\) with \(P_4 = ((\bZ_2, 0), (\bZ_2, k_3), (\bZ_2,
k_4))\).

Constructing the next stage of the Postnikov tower \emph{needs} the knowledge
of \(H^6(X_4, \bZ_2)\), again a particular case of the EDM problem, and Kenzo
computes in a few seconds \(H^6(X_4, \bZ_2) = \bZ_2^4\)~: 16 different choices
for the next Postnikov invariant \(k_5\); again Kenzo knows how to directly
construct the ``simplest'' non-trivial invariant \(k_5\), in a sense which
cannot be detailed here\footnote{Depending on the Smith reduction of the
boundary matrices of the small chain complex which is the main component of the
effective homology of \(X_4\).}; the other cohomology classes could be
constructed and used as well, but the computations would be more complicated.
Then \(F_5\) and \(X_5\) are constructed, but this time a few hours of
computation are necessary to obtain \(H^7(X_5, \bZ_2) = \bZ_2^5\)~: there are
32 different choices for the next invariant \(k_6\) and again, in this
``simple'' case, such a \(k_6\) actually is an invariant of the homotopy type
of the resulting space, see Section~\ref{tbpxw}.

  And so on.

\subsection{Morphisms between Postnikov towers.}

\subsubsection{The definition.}

We have presented the Postnikov towers as being the objects of the Postnikov
category \(\cP\), so that we must also describe the \(\cP\)-morphisms. The
standard considerations around homotopy groups and Kan minimal models, see for
example~\cite{MAY}, lead to the following definition.

\begin{dfn}\label{rgynw} ---
\emph{Let \(P = ((\pi_n, k_n))_{n \geq 2}\) and \(P' = ((\pi'_n, k'_n))_{n \geq
2}\) be two Postnikov towers. A \emph{morphism} \(f: P \rightarrow P'\) is a
collection of group morphisms \(f = (f_n: \pi_n \rightarrow \pi'_n)_{n \geq
2}\) satisfying the following \emph{recursive} coherence property for every
\(n\). The sub-collection \((f_i)_{2 \leq i \leq n-1}\), if coherent, defines a
continuous map \(\phi_{n-1}: X_{n-1} (= \bfPS(P_{n-1})) \rightarrow X'_{n-1} (=
\bfPS(P'_{n-1})) \) between the \((n-1)\)-th stages of the respective Postnikov
towers. So that two canonical maps are defined:
\begin{itemize}
\item
The map \(\phi_{n-1}\) induces in a contravariant way a map
\(\phi^{\ast}_{n-1}: H^{n+1}(X'_{n-1}, \pi'_n) \rightarrow H^{n+1}(X_{n-1},
\pi'_n)\) between the cohomology groups;
\item
The map \(f_n\) induces in a covariant way a map \({f_n}_\ast: H^{n+1}(X_{n-1},
\pi_n) \rightarrow H^{n+1}(X_{n-1}, \pi'_n)\).
\end{itemize}
Then the equality \(\phi^{\ast}_{n-1}(k'_n) = {f_n}_\ast(k_n)\) is required.}

\emph{If so, a  continuous map \(\phi_n: X_n \rightarrow X'_n\) is defined,
which allows one to continue the recursive process. The projective limit
\(\displaystyle \phi = \lim_\leftarrow \phi_n\) then is a continuous map
\(\phi: X = \bfPS(P) \rightarrow X' = \bfPS(P')\).}\hfill\(\Box\)
\end{dfn}

\subsubsection{First example.}\label{izeiy}

This definition implies some \emph{isomorphisms} between \emph{different}
Postnikov towers can exist. Let us examine when a collection \(f = (f_n: \pi_n
\rightarrow \pi'_n)_{n \geq 2} : ((\pi_n, k_n))_{n \geq 2} \rightarrow
((\pi'_n, k'_n))_{n \geq 2}\) is an isomorphism. On one hand the coherence
condition stated above must be satisfied, on the other hand every \(f_n\) must
be a group isomorphism; if this is the case the obvious inverse \(g =
(f_n^{-1})_{n \geq 2}\) also satisfies the coherence condition and actually is
an inverse of \(f\).

The simplest example where a non-trivial isomorphism happens is the following.
Let us consider the small Postnikov tower \(P = ((\bZ, 0), (\bZ, k_3))\) where
\(k_3 \in H^4(K(\bZ,2))\) is \(k_3 = c_1^2\), the square of the canonical
generator \(c_1 \in H^2(K(\bZ, 2), \bZ)\), the first universal Chern class. The
corresponding space \(X = \bfPS(P)\) is the total space of a well defined
fibration:
\[
K(\bZ,3) \hookrightarrow X \twoheadrightarrow K(\bZ, 2)
\stackrel{c_1^2}{\longrightarrow} K(\bZ, 4)
\]
The same construction is valid replacing \(k_3\) by \(k'_3 = -k_3\); the
Postnikov tower \(P' = ((\bZ, 0), (\bZ, k'_3))\) produces a \emph{different}
fibration:
\[
K(\bZ,3) \hookrightarrow X' \twoheadrightarrow K(\bZ, 2)
\stackrel{-c_1^2}{\longrightarrow} K(\bZ, 4)
\]
It is important to understand the fibrations not only are different but they
are even non-isomorphic: their classifying maps are not \emph{homotopic}. Yet
the spaces \(X = \bfPS(P)\) and \(X' = \bfSP(P')\) are the same, that is, they
have the same homotopy type; the following diagram is induced by the group
morphism \(\varepsilon_4: K(\bZ, 4) \stackrel{K(-1, 4)}{\longrightarrow} K(\bZ,
4)\) associated to the symmetry \(-1: n \mapsto -n\) in \(\bZ\), and the same
for \(\varepsilon_3\).
\[
\begin{CD}
 K(\bZ,3) @>>> X     @>>> K(\bZ, 2) @>c_1^2>> K(\bZ, 4) \\ @V\varepsilon_3
 V\cong V         @V\varepsilon_3 \widetilde{\times} = V\cong V   @V=VV @V\varepsilon_4 V\cong V \\
 K(\bZ,3) @>>> X'    @>>> K(\bZ, 2) @>-c_1^2>> K(\bZ, 4) \\
\end{CD}
\]
The \(\cong\) sign between \(X\) and \(X'\) is particularly misleading. It is
correct from the topological point of view: both spaces \(X\) and \(X'\)
actually are homeomorphic and \(\varepsilon_3 \widetilde{\times}\!\! =\) is
such a homeomorphism. The \(\cong\) sign is incorrect with respect to the
principal \(K(\bZ, 3)\)-structures: the actions of \(K(\bZ,3)\) on the fibres
of \(X\) and \(X'\) are \emph{not compatible}; the satisfied relation is only
\((\varepsilon_3 \widetilde{\times}\!\! =)(a \cdot x) = \varepsilon_3(a) \cdot
(\varepsilon_3 \widetilde{\times}\!\! =)(x)\) and the principal structures
would be compatible if \((\varepsilon_3 \widetilde{\times}\!\! =)(a \cdot x) =
a \cdot (\varepsilon_3 \widetilde{\times}\!\! =)(x)\) was satisfied, this is
why the classifying maps are opposite.

Maybe the same phenomenon for the Hopf fibration is easier to be understood.
Usually we take \(S^3\) as the unit sphere of \(\bC^2\) so that a
\emph{canonical} \(S^1\)-action is underlying and a \emph{canonical}
characteristic class on the quotient \(S^3 / S^1\) is deduced. But if you
reverse the \(S^1\)-action, why not, the space \(S^3\) is not modified, the
quotient \(S^3 / S^1\) is not modified either, but the characteristic class is
the opposite one. In other words, it is important not to forget the classifying
map characterizes the isomorphism class of a principal fibration, but not the
homotopy type of the total space!

\subsubsection{The key example.}

The next example of a Postnikov tower with two stages is still rather simple
but is sufficient to understand the essential failure of the claimed Postnikov
\emph{invariants}.

Let us consider the tower \(P(\ell,k) = ((\bZ^\ell, 0), (\bZ, k))\), the
parameter \(\ell\) being some positive integer, and \(k\), the unique
non-trivial Postnikov ``invariant'' being an element \(k \in H^4(K(\bZ^\ell,
2), \bZ)\). A canonical isomorphism \(K(\bZ^\ell, 2) \cong K(\bZ, 2)^\ell\) is
available. The cohomology ring of \(K(\bZ, 2) = P^\infty \bC\) is the
polynomial ring \(\bZ[X]\) where \(X = c_1\) is the first universal Chern
class, of degree 2, so that \(H^\ast(K(\bZ^\ell, 2), \bZ) = \bZ[X_i]\) with \(1
\le i \le \ell\), every generator \(X_i\) being of degree~2. Finally
\(H^4(K(\bZ^\ell, 2), \bZ) = \bZ[X_i]^{[2]}\), the exponent \([2]\) meaning we
must consider only the sub-module of the homogeneous polynomials of degree~2
with respect to the \(X_i\)'s. Every \(k \in \bZ[X_i]^{[2]}\) thus defines a
two stages Postnikov tower \(P(\ell,k) = ((\bZ^\ell, 0), (\bZ, k))\).

Two such \emph{different} Postnikov towers \(P(\ell, k)\) and \(P(\ell', k')\)
can be isomorphic. If so, the homotopy groups must me the same and \(\ell =
\ell'\) and it is enough to wonder whether \(P(\ell, k) \stackrel{???}{\cong}
P(\ell, k')\). A possible isomorphism \(f: P(\ell, k) \rightarrow P(\ell, k')\)
is made of \(f_2: \bZ^\ell \stackrel{\cong}{\rightarrow} \bZ^\ell\) and \(f_3:
\bZ \stackrel{\cong}{\rightarrow} \bZ\). The component \(f_3\) is a possible
simple sign change, as in the first example~\ref{izeiy}, but the component
\(f_2\) is a \(\bZ\)-linear equivalence acting on the variables \([X_i]_{1 \leq
i \leq \ell}\). The coherence condition given in Definition~\ref{rgynw} becomes
\({f_3}_\ast (k) = {f_2}^\ast (k')\): the \({f_3}_\ast\) allows one to make
equivalent two classes of opposite signs, and the \({f_2}^\ast\), much more
interesting, allows one to make equivalent two classes \(k, k' \in
\bZ[X_i]^{[2]}\) where \(k\) is obtained from \(k'\) by a \(\bZ\)-linear change
of variables. We have here identified \(f_2\) with \(\phi_2\), the induced
automorphism of \(K(\bZ^\ell, 2) = X_2\), the first stage of both Postnikov
towers, see Definition~\ref{rgynw}.

Algebraic Topology succeeds: the topological problem of homotopy equivalence
between \(\bfPS(P(\ell, k))\) and \(\bfPS(P(\ell, k'))\) is transformed into
the algebraic problem of the \(\bZ\)-linear equivalence, up to sign, between
the ``quadratic forms'' \(k\) and \(k'\).  And this provides a complete
solution, because this landmark problem firstly considered by Gauss has now a
complete solution, see for example~\cite{SERR4, WTSN, CSSL}.

\subsubsection{Higher dimensions.}\label{dpzvm}

But instead of working with the integer \(3 = 2*2-1\), we could consider
exactly the same problem with the Postnikov tower:
\[
\cP_{2d-1} \ni P(\ell, d, k) = ((\bZ^\ell, 0), (0, 0), \ldots, (0, 0), (\bZ,
k_{2d-1} = k))
\]
defined by integers \(\ell \geq 1\), \(d \geq 2\) and a cohomology class \(k
\in H^{2d}(K(\bZ^\ell, 2), \bZ) = \bZ[X_i]^{[d]}\). Instead of an equivalence
problem between homogeneous polynomials of degree~2, we meet the same problem
but with homogeneous polynomials of degree~\(d\). And when this paper is
written, this problem seems entirely\footnote{Jiri Matousek points out this
qualifier is not correct, thanks! This problem in fact is \emph{theoretically}
solved, cf.~\cite{grsg}: very general computability results for problems about
arithmetic groups in particular cover our problem. But as far as we know, these
results, of course important and interesting, did not yet lead to concrete
implementations; it is a nice challenge to attack this question. Sure such
implementations are today rather problematic, but the computer scientists have
already obtained so many concrete good results that it would be erroneous to
leave this challenge off research.} \emph{open} as soon as \(d \geq 3\). Now is
the right time to recall what the very notion of invariant is.

\subsection{Invariants.}

\subsubsection{Elementary cases.}

What is an \emph{invariant}? An invariant is a process \(\cI\) which associates
to every object~\(X\) of some type some other object \(\cI(X)\), the relevant
\emph{invariant}; in other words, an invariant is a function. This terminology
clearly says that \(\cI(X)\) \emph{does not change} (does not \emph{vary}) when
\(X\) is replaced by \(X'\), if \(X\) and \(X'\) are equivalent in some sense:
a possible relevant equivalence between \(X\) and \(X'\) should imply the
\emph{equality} -- not again some other equivalence -- between \(\cI(X)\) and
\(\cI(X')\).

For example one of the most popular invariants is the set of \emph{invariant
factors} of square matrices. The concerned equivalence relation is the
similarity. If \(K\) is a commutative field and \(A \in M_n(K)\) is an \((n
\times n)\)-matrix with coefficients in \(K\) representing some endomorphism of
\(K^n\), the invariant factors of \(A\) are a sequence of polynomials \(\phi(A)
= (\mu_1, \ldots, \mu_k)\) characterizing \emph{in this case} the similarity
class of the matrix \(A\): two matrices \(A\) and \(B\) are similar if and only
if \(\phi(A)\ \fbox{\(=\)}\ \phi(B)\). Another example is the minimal
polynomial \(\mu_1(A)\): if two matrices are similar, they have the \emph{same}
minimal polynomial. Idem for the characteristic polynomial which is the product
of the invariant factors, and so on. It is well known that for example the
characteristic polynomial does not characterize the similarity class, yet it is
an invariant: if two matrices are similar, they have the \emph{same}
characteristic polynomial. Sometimes the characteristic polynomial is
sufficient to disprove the similarity between two matrices, sometimes not. The
trivial invariant consists in deciding that \(\cI(A) = \bast\), some fixed
object,  for every matrix; not very interesting but it is undoubtedly an\ldots\
invariant. Symmetrically the tentative invariant \(\cI(A) = A\) is \emph{not}
an invariant, for there exist different (\(!\)) matrices\footnote{See
\mbox{\scriptsize\texttt{http://encyclopedia.thefreedictionary.com/invariant}}
for other typical examples. Another amusing bug of the standard terminology in
Algebraic Topology is the expression ``characteristic class'' in the classical
fibration theory: the usual characteristic classes are actual
invariants~(\(!\)) of the isomorphism class but, except in simple situations,
they do not characterize (\(!\)) this isomorphism class.}.

Algebraic Topology is in a sense an enormous collection of (algebraic)
invariants associated to topological spaces, invariants with respect to some
equivalence relation, frequently the homotopy equivalence. Typically a homotopy
group \(\pi_n\) is an invariant of this sort. Not frequently, with respect to
some appropriate equivalence relation, it is possible a \emph{complete}
invariant is available. For example the \(H_1\) is a complete invariant for the
homotopy type of a finite connected graph, the genus is a complete invariant
for the diffeomorphism type of a closed orientable real manifold of
dimension~2.

The last two examples, quite elementary, are interesting, because the difficult
logical problem underlying this matter is often forgotten and easily
illustrated in these cases. Let \(M_0\) and \(M_1\) be two closed orientable
2-manifolds that are diffeomorphic; if \(g\) denotes the genus, then
\(g(M_0)\,\fbox{\(=\)}\,g(M_1)\): the genus is an invariant; furthermore it is
a complete invariant, because conversely \(g(M_0)\,\fbox{\(=\)}\,g(M_1)\)
implies both manifolds are diffeomorphic. We have framed the `=' sign, because
the main problem in the continuation of the story is there.

Let us consider now the case of the finite graphs. In fact, it is \emph{false}
the~\(H_1\) functor is an invariant. If you take a triangle graph \(G_0 =
\triangle\) and a square graph \(G_1 = \Box\), same homotopy type, the careless
topologist thinks \mbox{\(H_1(G_0) = H_1(G_1) = \bZ\)} so that \(H_1\) looks
like an invariant of the homotopy type, but it is important to understand this
is deeply erroneous. With respect to any coherent formal definition of
mathematics, in fact \(H_1(\triangle) \neq H_1(\Box)\), these \(H_1\)-groups
are only \emph{isomorphic}. To obtain an actual invariant of the homotopy type,
you must consider the functor \(\mathbf{H}_1 = \textrm{IC} \circ H_1\), where
\(\textrm{IC}\) is the ``isomorphism class'' functor, always difficult to
properly define from a logical point of view, see for example~\cite{BRBK}. But
in the case of the \(H_1\)-group of a finite graph, it is a free \(\bZ\)-module
of finite type, it is particularly easy to determine whether two such groups
are isomorphic and every topologist \emph{implicitly} apply the IC functor
without generating any error.

Such a situation is so frequent that most topologists come to confuse both
notions of \emph{functor} and \emph{invariant}, and the case of the Postnikov
``invariants'' is rather amazing.

\subsection{The alleged Postnikov ``invariants''. }

\subsubsection{Terminology.}\label{mobnv}

We start with a sensible topological space, for example a finite simply
connected CW-complex \(E\). The textbooks explain how it is possible to define
or sometimes to ``compute'' the Postnikov invariants \((k_n(E))_{n \geq 3}\).
In our framework, the problem is the following:
\begin{prb}\label{rqcxf} ---
How to determine a Postnikov tower \(P = ((\pi_n, k_n))_{n \geq 2}\) such
that~\(E\) and \(\emph{\bfPS}(P)\) have the same homotopy type?
\end{prb}
This problem, thanks to the general \emph{Constructive} Algebraic Topology
framework of the authors, now has a positive \emph{and} constructive solution.
The aforementioned textbooks also describe ``solutions'', but which do
\emph{not} satisfy the constructive requirements which should yet be required
in this context. See also~\cite{SCHN} for another theoretical
\emph{constructive} -- and interesting -- solution, significantly more complex,
so that it has not yet led to concrete results, that is, to machine programs.

Most topologists think a positive solution for Problem~\ref{rqcxf} imply the
\(k_n\)'s of the result are ``invariants'' of the homotopy type of \(E\). This
is simply \emph{false}, for any reasonable understanding of the word
\emph{invariant}, and it is rather strange such an error remains present a so
long time in a so important field as basic Algebraic Topology. The \(k_n\)'s
could be called \emph{invariants} if  they solved the next problem.
\begin{prb} ---
\emph{Construct a functor \(\bfSP: \SSEH \rightarrow \cP\) satisfying the
following properties:
\begin{enumerate}
\item
Some original space \(E \in \SSEH\) and \(\bfPS \circ \bfSP (E)\) have the same
homotopy type;
\item
If \(E\) and \(E' \in \SSEH\) have the same homotopy type, then \(\bfSP(E)\
\fbox{\(=\)}\ \bfSP(E')\).
\end{enumerate}}
\end{prb}
The first point is a rephrasing of Problem~\ref{rqcxf}, and the second states
that if~\(E\) and~\(E'\) have the same homotopy type, then the images
\(\bfSP(E)\) and \(\bfSP(E')\) are \fbox{equal}, not only \emph{isomorphic}. In
other words the claimed ``invariant'' must not \emph{change} when the source
object changes in the same equivalence class; this is of course (\(?\)) the
very notion of invariant.

The non-constructive topologist easily solves the problem by replacing the
category \(\cP\) by the quotient \(\cP / \textrm{Iso}\), and then a correct
solution is obtained, but it is an artificial one. The category \(\SSEH /
H\textrm{-equiv}\) and the canonical projection \(\SSEH \rightarrow \SSEH /
H\textrm{-equiv}\) would be much simpler, but obviously without any interest.

The right interpretation of the \(k_n\)'s is the following: combined with the
standard homotopy groups \(\pi_n\), they are to be considered as
\emph{directions for use} allowing one to reconstruct a simple object with the
right homotopy type; another rephrasing of Problem~\ref{rqcxf}. But it can
happen two different objects \(E\) and \(E'\) with the \emph{same} homotopy
type produce \emph{different} ``directions for use'', so that these
``directions for use'' are not invariants of the homotopy type. In fact such an
accident is the most common situation, except for the topologists working only
with paper and pencil.

\subsubsection{The SP functor, first try.}\label{tbpxw}

Let us briefly describe the standard solution of Problem~\ref{rqcxf}, a
solution which can be easily made constructive thanks to~\cite{SRGR3, RBSR6,
SCHN}. Let \(E\) be some reasonable\footnote{That is, an \(\SSEH\)-object,
see~\cite{RBSR6}.} simply connected space. There are many ways to determine
the\footnote{In fact \emph{some} Postnikov tower\ldots} Postnikov tower \(P =
\bfSP(E)\) and one of them is illustrated here with the beginning of the
simplest case, the 2-sphere~\(S^2\). Hurewicz indicates \(\pi_2 = H_2 = \bZ\);
the invariant \(k_2\) is necessarily null. The next step invokes the Whitehead
fibration:
\[
 K(\bZ, 1) \hookrightarrow E^3 \twoheadrightarrow S^2
 \stackrel{c_1}{\longrightarrow} K(\bZ, 2).
\]
where \(c_1\) is the canonical cohomology class, in this case the first Chern
class of the complex structure of \(S^2\). The first stage of the Postnikov
tower is \(X_2 = K(\bZ, 2) = P^\infty \bC\) and the first stage of the
complementary Whitehead tower is the total space \(E^3 = S^3\): our fibration
is nothing but the Hopf fibration. Then \(\pi_3(S^2) = \pi_3(S^3) = H_3(S^3,
\bZ) = \bZ\), so that the next Postnikov invariant is some \(k_3 \in H^4(X_2,
\bZ) = H^4(K(\bZ, 2), \bZ) = \bZ\). How to determine this cohomology class?

In general we obtain a fibration:
\[
E^n \hookrightarrow E \rightarrow X_{n-1}
\]
where \(X_{n-1}\) is the \((n-1)\)-stage of the Postnikov tower containing the
homotopy groups \((\pi_i)_{2 \leq i \leq n-1}\), and \(E^n\) is the
complementary \(n\)-stage of the Whitehead tower~\cite[Proposition~8.2.5]{DDPR}
containing the homotopy groups \((\pi_i)_{i \geq n}\); in the Kan context
of~\cite[\S~8]{MAY}, \(E^n\) is the \(n\)-th Eilenberg subcomplex of \(E\). How
to deduce a cohomology class \(k_n \in H^{n+1}(X_{n-1}, \pi_n)\)? The
\((n-1)\)-connectivity of \(E^n\) produces a transgression morphism \(H^n(E^n,
\pi_n) \rightarrow H^{n+1}(X_{n-1}, \pi_n)\); the group \(H^n(E^n, \pi_n)\)
contains a fundamental Hurewicz class and the image of this class in
\(H^{n+1}(X_{n-1}, \pi_n)\) is the wished \(k_n\). In the particular case of
\(S^2\) this process leads to an isomorphism \(H^3(S^3, \bZ)
\stackrel{\cong}{\longrightarrow} H^4(K(\bZ, 2), \bZ)\) so that \(k_3\) is the
image of the fundamental cohomology class of \(S^3\), that is, the (\(?\))
generator \(c_1^2\) of \(H^4(K(\bZ, 2), \bZ)\). Sure?

As usual we have light-heartedly mixed intrinsic objects and isomorphism
classes of these objects. The \emph{correct} isomorphism to be considered for
our example is \(H^3(E^3, \pi_3\underline{(E^3)}) \cong
H^4(K(\pi_2\underline{(S^2)}, 2), \pi_3\underline{(E^3)})\) where \(E^3\) is
now the total space of the \emph{canonical} fibration
\(K(\pi_2\underline{(S^2)}, 1) \hookrightarrow E^3 \twoheadrightarrow S^2\);
this isomorphism actually is canonical. But no canonical ring structure for
\(\pi_3(E^3)\) so that speaking of \(c_1^2\) does not make sense. There is
actually a canonical element \(k_3 \in H^4(K(\pi_2\underline{(S^2)}, 2),
\pi_3\underline{(E^3)})\), but such an element deeply depends on \(S^2\) itself
and cannot be qualified as an \emph{invariant} of the \emph{homotopy type} of
\(S^2\). An actual invariant should be taken in the ``absolute'' (independent
of \(S^2\)) group \(H^4(K(\bZ, 2), \bZ)\), but such a choice depends on an
isomorphism \(\pi_3(E^3) \cong \bZ\); two such isomorphisms are possible so
that in this case the \(k_3\) is defined \emph{up to sign}: it is well known
the Hopf fibration and the ``opposite'' one produce the ``same'' total space.

This is the reason why in the definition of a Postnikov tower, see
Definition~\ref{szglm}, we have decided to have only \emph{one} group for each
isomorphism class; this is easy and can be done in a \emph{constructive} way.
The goal being to obtain \emph{invariants}, we had to design our Postnikov
towers as a catalogue of possible Postnikov towers, in such a way that there
are no redundant copies up to isomorphism in this collection; bearing this
point in mind, it was mandatory to have only one copy for every isomorphism
class of group. But this was not enough, for it is today impossible to take the
same precaution for the second components, the \(k_n\)'s, the so-called
Postnikov invariants.

For example if the concerned homotopy groups are finite, then the number of
possible \(k\)-invariants is finite, so that the related equivalence problem is
theoretically solved; this was already noted by Edgar Brown~\cite{BRWNE1},
which conversely implies~(\(!\)) he did not know how to solve the general case.
On the contrary, as soon as the homotopy groups have infinite automorphism
groups, there is no known way to transform the pseudo-invariants into actual
invariants.

We understand now the reason of the repetitive remark in Section~\ref{djmab}:
``In this particular case, the \(k_n\) actually is an invariant of the homotopy
type''; we decided to systematically choose \(\pi_n = \bZ_2\), but the
automorphism group of \(\bZ_2\) is \emph{trivial}; no non-trivial automorphism
of the constructed tower can exist and \emph{then} the \(k_n\)'s are actual
\emph{invariants}.

But if some user intends to use the Postnikov invariants to try to prove the
spaces \(E\) and \(E'\) have different homotopy types, the following accident
can happen. A calculation could respectively produce the Postnikov towers
\(((\bZ^\ell, 0), (0, 0),\ldots, (\bZ, k_{2d-1}))\) and \(((\bZ^\ell, 0),
\ldots, (\bZ, k'_{2d-1}))\) (see Section~\ref{dpzvm}). If fortunately
\(k_{2d-1} = k'_{2d-1}\) our user can be sure the homotopy types are the same
but if on the contrary \(k_{2d-1} \neq k'_{2d-1}\), then he has to decide
whether two homogeneous polynomials of degree \(d\) are linearly equivalent or
not and for \(d \geq 3\): no general solution is known. \emph{Maybe} they are
equivalent, \emph{maybe} not; because the alleged invariants may\ldots\ vary,
in general our user cannot conclude: the claimed invariants cannot play the
role ordinarily expected for invariants; qualifying them as invariants is
therefore a deep \emph{error}.

\subsubsection{The SP functor, second try.}

The right definition for the \(\bfSP\) functor is now clear. We must add to the
data some \emph{explicit} isomorphisms between the homotopy groups \(\pi_n(E)\)
of the considered space~\(E\) with the corresponding \emph{canonical} groups,
see Definition~\ref{szglm}.
\begin{dfn} ---
\emph{The product \(\SSEH \widetilde{\times} I\) is the set of pairs \((E,
\alpha)\) where:
\begin{enumerate}
\item
The \(E\) component is a simplicial set with effective homology \(E \in
\SSEH\);
\item
The \(\alpha\) component is a collection \((\alpha_n)_{n \geq 2}\) of
isomorphisms \(\alpha_n: \pi_n(E) \stackrel{\cong}{\longrightarrow} \pi_n\)
where \(\pi_n\) denotes the unique \emph{canonical} group isomorphic to
\(\pi_n\underline{(E)}\).
\end{enumerate}}
\end{dfn}

The previous discussions of this text can reasonably be considered as a
demonstration of the next theorem.

\begin{thr}\label{ooaag} ---
A functor \(\emph{\bfSP}: \SSEH \widetilde{\times} I \rightarrow \cP\) can be
defined.
\begin{enumerate}
\item
If \((E, \alpha) \in \SSEH \widetilde{\times} I\), then \(E\) and
\(\emph{\bfPS} \circ \emph{\bfSP}(E, \alpha)\) have the same homotopy type.
\item
If \(P \in \cP\) is a Postnikov tower, there exists a unique \(\alpha\) such
that \(\emph{\bfSP}(\emph{\bfPS}(P), \alpha) = P\).
\end{enumerate}
\end{thr}

So that it is tempting -- and correct -- to replace the \(\bfPS\) functor by
another one \(\bfPS: \cP \rightarrow \SSEH \widetilde{\times} I\) to obtain a
better symmetry. But the ordinary topologists work with elements in \(\SSEH\),
not in \(\SSEH \widetilde{\times} I\).

\subsection{The Postnikov invariants in the available literature.}

Most textbooks speaking of Postnikov invariants (or \(k\)-invariants) use the
\emph{invariant} terminology without justifying it, so that strictly speaking,
no mathematical error in this case. For example \cite[p. 279]{DDPR} defines the
Postnikov invariant through a transgression morphism\footnote{We used this
method in Section~\ref{tbpxw}} and explains ``The \(k^i\) precisely constitute
the stepwise obstructions\ldots''; the statement about this obstruction of
course is correct but it seems the terminology should therefore speak of
Postnikov \emph{obstructions}? Nothing is explained about the \emph{invariant}
nature of these obstructions.

Other books speak of these invariants as objects allowing to \emph{reconstruct}
the right homotopy type. For example, in \cite[p. 412]{HTCH}: ``The map \(k_n\)
is equivalent to a class in \(H^{n+2}(X_n; \pi_{n+1}(K))\) called the \(n\)-th
\emph{\(k\)-invariant} of \(X\). These classes specify how to construct \(X\)
inductively from Eilenberg-MacLane spaces''. To be compared with our
considerations about the interpretation in terms of ``directions for use'' at
the end of Section~\ref{mobnv}. Again, no indication in this book about the
justification of the \emph{invariant} terminology. The Section ``The Postnikov
Invariants'' of~\cite[V.3.B]{DDNN} can be analyzed along the same lines.

In~\cite[VI]{GRJD}, because of a sophisticated categorical environment, the
authors prefer to define the \emph{general} notion of \emph{Postnikov tower}
for a space \(X\), each one containing in particular its
\(k_n\)-invariants~\mbox{\cite[VI.5]{GRJD}}; finally
Theorem~\cite[VI.5.14]{GRJD} proves two such Postnikov towers for the same
\(X\) are \emph{weakly} equivalent. In other words one source object produces
in general a large infinite set of (different!) \(k_n\)-invariants, for every
relevant \(n\); yet some invariant theory is interesting when different objects
can produce the same invariants, not when an object produces different
invariants! In fact, as explained in our text, this cannot be currently
avoided, but why these authors do not make explicit the misleading status of
these claimed invariants?

The book~\cite{MAY} systematically uses the powerful notion due to Kan of
\emph{minimal simplicial Kan-model}, often allowing a user to work in a
``canonical'' way, allowing frequently the same user to easily detect a
non-unicity problem. In this way \mbox{\cite[p. 113]{MAY}} correctly signals
that the map \(B \rightarrow K(\pi, n+1)\) leading to a \(k_n\)-invariant is
defined up to a \(\pi\)-automorphism, which is not a serious drawback: the
decision problem about the possible equivalence of two \(k_n\)'s under such an
automorphism is easy when \(\pi\) is of finite type. But the author does not
mention the same problem with respect to the automorphisms of the base space
\(B\), the automorphisms leading to the corresponding open problem detailed
here Section~\ref{dpzvm}.

The same author in a more recent textbook~\cite{MAY2} again considers the same
question. He defines the notion of Postnikov system in Section~22.4; the
existence of \emph{some} Postnikov system is proved, the terminology
\(k\)-invariant is used one time, between quotes seeming imply this expression
is not really appropriate, but without any explanations.

Hans Baues \cite[p.33]{BAUS4} on the contrary correctly respects the necessary
symmetry between the source and the target of the classifying map; but the
author is aware of the underlying difficulty and it is interesting to observe
how he ``solves'' the raised problem:

{\newlength{\parquot}
 \setlength{\parquot}{\textwidth}
 \addtolength{\parquot}{-\parindent}
 \addtolength{\parquot}{-\parindent}
 \addtolength{\parquot}{-\parindent}
 \begin{center}
 \fbox{\parbox{\parquot}
 {Here \(k_n(Y)\) is actually an \emph{invariant} of the homotopy type of \(Y\) in
 the sense that a map \(f: Y \rightarrow Z\) satisfies:
 \[
 (P_{n-1}f)^\ast k_n(Z) = (\pi_n f)_\ast k_n(Y)
 \]
 in \(H^{n+1}(P_{n-1} X, \pi_n Y)\).}}
 \end{center}}

Clearly explained, the author says that the \emph{invariant} is variable, but
in a \emph{functorial} way. Baues' condition is essentially the \emph{coherence
condition} of our Definition~\ref{rgynw}. If the appropriate morphisms of the
category \(\SSEH \widetilde{\times} I\) were defined, the functorial property
of the map \(\bfSP\) (Theorem~\ref{ooaag}) would be exactly Baues' relation.
But it is not explained in Baues' paper why a functor may be qualified as an
invariant.

Probably the most lucid reference about our subject is~\cite{WHTH1}. Chapter~IX
is entirely devoted to Postnikov systems. We find p.~423:
\begin{quotation}
The term `invariant' is used somewhat loosely here. In fact \(k^{n+2}\) is a
cohomology class of a space \(X^n\), which has not been constructed in an
invariant way. This difficulty, however, is not serious, for, as we shall show
below, the construction of the space \(X^n\) can be made completely natural.
\end{quotation}
This text is essentially a rephrasing of Baues' explanation. Again the common
confusion between the notions of invariant and functor is observed. To make
``natural'' its invariants, George Whitehead uses enormous singular models, so
that the obtained \(k^{n+2}\) heavily depends on \(X\) itself and not only on
its homotopy type. In fact Section~\cite[IX.4]{WHTH1} shows Whitehead is in
fact also interested in being able to reconstruct the homotopy type of \(X\)
from the ``natural'' associated Postnikov tower, and this goal is obviously
reached, but this does not provide a general machinery allowing one to detect
\emph{different homotopy types} when the associated \emph{invariants} are
\emph{different}.

\vspace{20pt}

\noindent{\footnotesize{\begin{tabular}{r}
 \texttt{Julio.Rubio@unirioja.es}
 \\
 \texttt{Francis.Sergeraert@ujf-grenoble.fr}
\end{tabular}}}

\end{document}